\documentclass{article}
\usepackage{latexsym}
\usepackage[all]{xy}
\usepackage{amsmath}

\newtheorem{thm}{Th\'eor\`eme}[section]
\newtheorem{prop}[thm]{Proposition}
\newtheorem{lem}[thm]{Lemme}
\newtheorem{df}[thm]{D\'efinition}

\newtheorem{cor}[thm]{Corollaire}
\newtheorem{hyp}[thm]{Hypoth\`ese}

\begin{document}

$N^{\circ}$ d'ordre:

\vspace{10mm}

\begin{center}

\large \bf TH\`ESE \rm 
\vspace{8mm}

\normalsize Pr\'esent\'ee en vue de l'obtention du
\vspace{8mm}

\large TITRE DE DOCTEUR DE L'UNIVERSIT\'E PAUL SABATIER DE TOULOUSE 
\vspace{8mm}

\normalsize \textit{Sp\'ecialit\'e:}

Math\'ematiques pures 
\vspace{10mm}

\textit{par}

\large Bertrand TOEN 
\vspace{20mm}

\Large \bf $K$-TH\'EORIE ET COHOMOLOGIE DES 

CHAMPS ALG\'EBRIQUES:

TH\'EOR\`EMES DE RIEMANN-ROCH, 

$\mathcal{D}$-MODULES ET TH\'EOR\`EMES "GAGA".  \rm 

\vspace{15mm}

\normalsize 

\textit{Soutenue le $24$ Juin $1999$, devant le jury compos\'e de:} 
\vspace{5mm}

\end{center}

\begin{tabular}{lll}

\hspace{-20mm} L. BREEN & Professeur, Universit\'e de Paris $13$ & \\
\hspace{-20mm} L. MORET-BAILLY & Professeur, Universit\'e de Rennes $1$ & \\
\hspace{-20mm} C. SIMPSON & Directeur de Recherche C.N.R.S. (Toulouse) & Co-directeur de Th\`ese \\
\hspace{-20mm} J. TAPIA & Professeur, Universit\'e Paul Sabatier (Toulouse) & Directeur de Th\`ese \\
\hspace{-20mm} C. VOISIN & Directrice de Recherche C.N.R.S. (Paris $7$) & Rapporteur \\
\hspace{-20mm} C. WALTER & Professeur, Universit\'e de Nice-Sophia-Antipolis & Rapporteur

\end{tabular}

\vspace{20mm}

\begin{center}
Laboratoire Emile Picard, Universit\'e Paul Sabatier,

$118$, route de Narbonne, $31062$ Toulouse Cedex France.
\end{center}

\newpage

\large

\vspace*{60mm}
\hspace*{10mm} \textit{Je d\'edie cette th\`ese \`a E.}

\newpage

\tableofcontents

\newpage

\begin{center}
\underline{\bf Remerciements:}
\end{center}

Je tiens tout d'abord, \`a remercier tr\`es sinc\`erement Joseph Tapia et
Carlos Simpson, qui ont su encadrer ce travail de th\`ese avec beaucoup
d'humanit\'e. Je les remercie d'avoir accept\'e de partager avec moi un peu de
leur passion pour la g\'eom\'etrie alg\'ebrique : je leur dois un
enrichissement math\'ematique consid\'erable.

Ce fut un r\'eel plaisir de travailler sous leur direction.  \\

Je remercie tout particuli\`erement Charles Walter, qui est \`a l'origine
de la question du th\'eor\`eme de Riemann-Roch pour les champs
alg\'ebriques, et qui par l\`a m\^eme a motiv\'e cette th\`ese. \\

Je remercie Claire Voisin et Charles Walter pour m'avoir fait
l'honneur d'\^etre rapporteurs de cette th\`ese, ainsi que pour l'int\'er\^et
qu'ils ont port\'e \`a ce travail. \\

Je voudrais aussi remercier Andr\'e Hirschowitz pour s'\^etre int\'eress\'e \`a ce
travail, et m'avoir permis de m'exprimer sur ce
sujet lors du s\'eminaire "GAF 3",  ce qui fut pour moi tr\`es motivant.
\\

Un grand merci \`a Lawrence Breen et Laurent Moret-Bailly
pour avoir accept\'e de faire partie du jury. \\

Je voudrais exprimer ma gratitude \`a Marie Brouard et Yveline Panabi\`ere,
dont la gentillesse et l'efficacit\'e sont hors du commun. \\

Je remercie tous les math\'ematiciens de Toulouse avec qui j'ai eu
l'occasion de participer aux s\'eminaires, et gr\^ace \`a qui j'ai beaucoup
appris. Je pense en particulier \`a : C. Simpson, J. Tapia, M. Reversat,
P. Essydieux, B. Angl\`es et T. Fiedler.  \\

J'exprime aussi mes remerciements \`a l'institut Max Plank pour leur
accueil durant les mois d'Avril et Mai $1997$. Ce s\'ejour \`a \'et\'e aussi
productif et qu'agr\'eable. \\

Je remercie aussi tous les autres \'etudiants du "Frigo" et du "Goulag",
avec qui j'ai partag\'e de nombreuses discussions, d\'ebats, repas,
soir\'ees~\dots et qui ont cr\'e\'e une ambiance de travail plus qu'agr\'eable. \\

Enfin, je remercie Jacques, Nicole et Guillaume, pour leur soutien et
leur compr\'ehension durant toute la dur\'ee de cette th\`ese. \\

\newpage

\setcounter{section}{-1}

\begin{section}{Introduction}
\hspace{5mm}
Depuis les id\'ees novatrices d'Evariste Galois, les math\'ematiciens ont
appris qu'\'etant donn\'e un probl\`eme \`a r\'esoudre, il est important de
consid\'erer les solutions dans leur ensemble, ainsi que les relations
qui les lient entre elles. D'une certaine fa\c{c}on, la th\'eorie des
champs, imagin\'ee par Alexandre Grothendieck dans les ann\'ees $60$, est
une illustration de ce principe. En effet, dans les probl\`emes de
classification, l'utilisation des champs permet de prendre en compte
l'information concernant l'ensembles des solutions, cod\'e sous forme
"d'ensembles d'objets", mais aussi les relations entre ses
solutions, cod\'ees sous forme de "morphismes entre les objets".  Ainsi,
tout comme la th\'eorie de Galois nous dit que l'ensemble des solutions
d'une \'equation alg\'ebrique peut-\^etre consid\'er\'e comme un groupoide, dont les
morphismes sont d\'etermin\'es par l'action de son groupe de Galois, les
solutions aux probl\`emes de classifications apport\'ees par la th\'eorie
des champs se trouvent sous forme de groupoides (~plus pr\'ecis\'ement
de "cat\'egories fibr\'ees en groupoides"~). \\

L'utilisation avec succ\`es de la th\'eorie des champs aux probl\`emes de
modules (~i.e. de classification en g\'eom\'etrie alg\'ebrique~), ainsi que
l'apparition des solutions sous forme de "champs alg\'ebriques",
remonte \`a l'article fondateur de Pierre Deligne et David Mumford
\cite{dm}. Les deux auteurs y montrent qu'il est possible de
"g\'eom\'etriser" (~ou encore "d'alg\'ebriser"~) la notion de champs, et
de cr\'eer ainsi un cadre unique contenant, d'une part les sch\'emas,
mais aussi certains champs de nature g\'eom\'etrique (~les champs de modules par
exemple~) : celui des champs alg\'ebriques.

Tout comme un sch\'ema peut-\^etre vu comme une solution au probl\`eme de
modules du foncteur qu'il repr\'esente (~qui classifie ses "points"~),
un champ alg\'ebrique est un objet de nature g\'eom\'etrique, classifiant
non seulement des objets (~ses points~), mais aussi leurs
isomorphismes. Ainsi, de fa\c{c}on un peu naive, on peut se repr\'esenter un
champ alg\'ebrique comme un objet en groupoides dans la cat\'egorie des
sch\'emas. \\

L'introduction de tels objets apporte de nombreux avantages en th\'eorie
des modules. A titre d'exemple, examinons bri\`evement le cas du
champ des courbes de genre $g$ (~qui est l'exemple trait\'e dans
\cite{dm}~).

Consid\'erons le probl\`eme de modules suivant : chercher un objet de
nature g\'eom\'etrique $\cal M\mit_{g}$, qui repr\'esente le foncteur qui \`a
un sch\'ema $X$ associe le groupoide des familles de courbes propres lisses de
genre $g$ sur $X$.

Il est d\'emontr\'e dans \cite{dm} qu'il existe un champ alg\'ebrique $\cal
M\mit_{g}$, lisse sur $Spec \bf Z$, et une \'equivalence
naturelle entre $Hom(X,\cal M\mit_{g})$, et le groupoides des courbes
de genre $g$ sur $X$. De plus, il existe une compactification
$\cal M\mit_{g} \hookrightarrow \overline{\cal M\mit}_{g}$, o\`u
$\overline{\cal M\mit}_{g}$ est un champ irr\'eductible, propre et
lisse sur $Spec \bf Z$, classifiant les courbes stables de genre $g$.

Ainsi, par d\'efinition m\^eme, il existe une courbe stable universelle
$\pi~:~\overline{\cal C\mit}_{g}~\longrightarrow~\overline{\cal M\mit}_{g}$,
correspondant au morphisme identit\'e de $\overline{\cal M\mit}_{g}$. De
plus, si
$\omega$ est le fibr\'e en droites canonique relatif
de $\overline{\cal C\mit}_{g}$ sur $\overline{\cal M\mit}_{g}$,
on peut alors d\'efinir des fibr\'es vectoriels $\cal V\mit_{m}:=
\pi_{*}(\omega^{\otimes m})$ sur $\overline{\cal M}_{g}$.

Les avantages d'une telle situation en comparaison avec celle dont on
dispose \`a l'aide de l'espace de modules grossier $\overline{M}_{g}$ (~
qui est le sch\'ema dont les points g\'eom\'etriques classifient les courbes
stables de genre $g$~) sont nombreux.

\begin{itemize}
\item
Tout d'abord le sch\'ema $\overline{M}_{g}$ n'est qu'une solution
partielle au probl\`eme de modules initial consistant \`a classifier
toutes les familles de courbes de genre $g$. Pour cette raison, la
courbe universelle de genre $g$ sur $\overline{M}_{g}$ ne peut pas
exister (~il manque une coh\'erence entre les points g\'eom\'etriques de
$\overline{M}_{g}$, du \`a l'existence de courbes avec des groupes
d'automorphismes non triviaux~).
\item
Bien que $\overline{\cal M}_{g}$ soit un champ lisse sur $Spec \bf Z$,
le sch\'ema $\overline{M}_{g}$ ne l'est pas. Par exemple, cela pose
quelques probl\`emes si l'on souhaite utiliser des outils analytiques tels
que la th\'eorie de Hodge sur $\overline{M}_{g}\otimes~\bf~C$, ou encore
la th\'eorie des intersections (~au moins \`a coefficients entiers~).
\item
Les fibr\'es vectoriels $\cal V\mit_{m}$ n'existent pas sur
$\overline{M}_{g}$ (~pour la m\^eme raison que la courbe universelle
n'existe pas~). Ainsi, l'\'etude de ces fibr\'es ne peut se faire que si l'on
accepte
de travailler sur le champ $\overline{\cal M}_{g}$.
\end{itemize}

Les trois exemples pr\'ec\'edents, qui ne sont que des exemples parmi
tant d'autres, montrent qu'il est capital de g\'en\'eraliser les outils de
la g\'eom\'etrie alg\'ebrique (~th\'eorie de Hodge, th\'eorie des
intersections, \\
$K$-th\'eorie et classes caract\'eristiques ...~) au cadre plus g\'en\'eral des champs
alg\'ebriques, si l'on veut profiter de la richesse de tels objets. \\

Dans cette th\`ese, nous nous sommes int\'eress\'es \`a essentiellement deux
probl\`emes
\begin{itemize}
\item L'\'etude de la $K$-th\'eorie alg\'ebrique des champs alg\'ebriques, et
son application \`a des formules de Riemann-Roch ainsi que des formules
d'indices de $\cal D$-modules.
\item Les relations entre champs alg\'ebriques et champs analytiques :
th\'eor\`emes GAGA et crit\`eres d'alg\'ebrisations.
\end{itemize}

Expliquons en quelques mots les raisons de ces deux choix. \\

Comme nous l'avons aper\c{c}u dans l'exemple pr\'ec\'edent, il existe
naturellement de nombreux fibr\'es vectoriels sur les champs de modules.
Il se trouve que l'\'etude de ces fibr\'es pr\'esente des int\'er\^ets g\'eom\'etriques
et arithm\'etiques.

Par exemple, dans l'exemple de $\overline{\cal M}_{g}$, les sections des fibr\'es
$\cal V\mit_{m}$ s'identifient aux "formes modulaires de poids $m$". De
m\^eme, il existe des fibr\'es vectoriels $\cal V\mit_{m}$ sur le champ (~
compactifi\'e~) $\overline{\cal A}_{1}$
des courbes elliptiques, dont les sections sont les formes modulaires
de poids $m$ sur le demi-plan de Poincar\'e.
On conna\^it depuis longtemps l'int\'er\^et de tels objets en arithm\'etique.

D'autre part, la formule bien connu pour des vari\'et\'es,
\mbox{$C_{1}(\cal V\mit_{m})=Div(f)$}, pour une section m\'eromorphe $f$ de $\cal
V\mit_{m}$ sur $\overline{\cal A}_{1}$, s'interpr\`ete comme la formule
classique (~\cite[$VII$, $3.1$]{se}~)
$$\frac{m}{24}=\frac{1}{4}v_{i}(f) + \frac{1}{6}v_{\rho}(f) +
\frac{1}{2}.\sum_{x \neq i,\rho}v_{x}(f)$$
Ceci permet l\'egitimement de penser que le formalisme des classes
caract\'eristiques a un int\'er\^et dans ce contexte.

Dans le m\^eme genre d'id\'ees, les formules de Riemann-Roch appliqu\'ee aux
fibr\'es $\cal V\mit_{m}$, sont aussi un outil pour \'etudier les dimensions
d'espaces
de formes modulaires. \\

Un autre champ d'application \'eventuel du formalisme des classes
caract\'eristiques
et des formules de type Riemann-Roch est celui de la g\'eom\'etrie
\'enum\'erative de M. Kontsevich et Y. Manin \cite[$2$]{ko} (~d'autant plus que
la transformation de Riemann-Roch est souhait\'ee pour certaines
constructions \cite[$5.4$]{bf}~). \\

L'\'etude des $\cal D$-modules sur les champs alg\'ebriques pr\'esente, elle aussi,
plusieurs int\'er\^ets. Le premier qui vient \`a l'esprit est l'\'etude
topologique des champs alg\'ebriques. Par exemple, les th\'eor\`emes d'indices
pour les $\cal D$-modules holon\^omes permettent de d\'emontrer des
formules de Gauss-Bonnet pour des caract\'eristiques d'Euler pond\'er\'ees (~
analogues \`a \cite{mac}~).

Par ailleurs, il est montr\'e dans \cite{jos}, que le
th\'eor\`eme de Riemann-Roch pour les $\cal D$-modules \'equivariants est
utile en th\'eorie des repr\'esentations. Dans la m\^eme veine d'id\'ee,
la correspondance de Langlands g\'eom\'etrique fait naturellement
appara\^itre des $\cal D$-modules sur certains champs de fibr\'es
vectoriels (~\cite{gin}~). L'\'etude de ces objets (~par exemple leurs
indices~) n\'ecessite alors une g\'en\'eralisation des r\'esultats de
\cite{jos} au cadre plus g\'en\'eral, et bien plus maniable, des champs
alg\'ebriques. \\

Il va sans dire que l'utilisation de m\'ethodes analytiques en
g\'eom\'etrie alg\'ebrique complexe est un outil extr\^emement puissant.
Leur pertinence est d'autant plus grande qu'il existe des th\'eor\`emes
de comparaison entre la g\'eom\'etrie alg\'ebrique et analytique (~par
exemple les "th\'eor\`emes GAGA"~). Il nous semble alors int\'eressant de
poss\'eder de tels r\'esultats pour
les champs alg\'ebriques complexes.
\vspace{10mm}

Il se trouve qu'une partie non n\'egligeable des travaux consacr\'es aux
champs alg\'ebriques traite de la th\'eorie des intersections (~
\cite{mu,g2,vi2,eg,kr}~), et plus particuli\`erement de la d\'efinition de
groupes (~ou anneaux~) de Chow pour des champs alg\'ebriques. Le point
commun \`a toutes ces d\'efinitions est la remarquable propri\'et\'e que
pour la projection d'un champ de Deligne-Mumford sur son espace de
modules $p : F \longrightarrow M$, le morphisme d'images directes
$p_{*} : CH(F)\otimes \bf Q\mit \longrightarrow CH(M)\otimes \bf Q$
est un isomorphisme. Ainsi, ces groupes ne contiennent pas (~ou pas
assez~) l'information \'equivariante suppl\'ementaire que poss\`ede $F$ par rapport \`a
$M$. Ceci implique en particulier, que m\^eme pour les champs les plus
simples (~les champs classifiant d'un groupe fini par exemple~), il
ne peut exister de formule du type Hirzebruch-Riemann-Roch \`a valeurs
dans ces groupes de Chow.

C'est cette derni\`ere remarque qui nous a naturellement amen\'es a nous
int\'eresser au probl\`eme majeur trait\'e dans cette th\`ese, qui est
d'\'elargir les pr\'ec\'edentes d\'efinitions des groupes de Chow, de fa\c{c}on \`a
ce que des formules d'Hirzebruch-Riemann-Roch puissent exister. Nous
r\'esolvons ce probl\`eme en introduisant la notion de
"cohomologie \`a coefficients dans les repr\'esentations", ce qui nous
permet de d\'emontrer des formules g\'en\'erales de
Grothendieck-Riemann-Roch dans le cadre des champs alg\'ebriques. \\

Comme il a \'et\'e expliqu\'e ci-dessus, le point de d\'epart
est la constatation que les groupes de Chow d\'ej\`a existant (~
\cite{mu,g2,vi2,eg,kr}~) ne sont pas adapt\'es aux formules de
Riemann-Roch. Pour r\'esoudre ce probl\`eme nous avons d\'ecid\'e de
remonter "\`a la source", et d'\'etudier directement les propri\'et\'es de la
$K$-th\'eorie des champs alg\'ebriques. Comme on sait, la $K$-th\'eorie
alg\'ebrique joue un r\^ole de "cohomologie universelle", et donc, \`a
travers son \'etude, ce sont les propri\'et\'es cohomologiques (~voire m\^eme
"motiviques"~) que nous \'etudions. Les r\'esultats cl\'es que nous d\'emontrons
sont les th\'eor\`emes de d\'evissage. Pour simplifier, supposons que $F$
soit un champ de Deligne-Mumford lisse sur le corps des nombres complexes.
Dans ce cas, le th\'eor\`eme de d\'evissage assure l'existence d'un isomorphisme
d'anneaux, fonctoriel pour les images r\'eciproques
$$\phi_{F} : \bf G\mit_{*}(F)\otimes \bf C\mit \longrightarrow
H^{-*}(I_{F},\underline{G})\otimes \bf C$$
o\`u $I_{F}$ est le champ des ramifications de $F$, et le membre de
droite est la cohomologie de $I_{F}$ \`a valeurs dans le pr\'efaisceau en
spectres de \\
$G$-th\'eorie (~remarquons d\`es maintenant que nous sommes oblig\'es
ici d'utiliser les techniques d'alg\`ebre homotopique pour lui donner
un sens~).

Il est \`a noter, que lorsque $F=[X/H]$ est le champ quotient d'une vari\'et\'e
par un groupe fini, la formule pr\'ec\'edente est \'equivalente \`a une
formule connue depuis longtemps, d\'ecrivant la $K$-th\'eorie
\'equivariante de $X$ en fonction de la $K$-th\'eorie des points fixes (~
\cite{as,vi}~)
$$\bf K\mit_{*}(X,H)\otimes \bf C\mit \simeq \bigoplus_{h \in c(H)}
\bf K\mit_{*}(X^{h})^{Z(h)}$$
o\`u $c(H)$ est l'ensemble des classes de conjugaisons de $H$, $X^{h}$ le
sous-sch\'ema des points fixes de $X$ par $h$, et $Z(h)$
le centralisateur de $h$ dans $H$.
Comme un champ de Deligne-Mumford est localement un quotient par un
groupe fini, on peut l\'egitimement penser que notre th\'eor\`eme de
d\'evissage est obtenu en "recollant" les isomorphismes ci-dessus. C'est
pr\'ecis\'ement cette op\'eration de recollement qui fait appara\^itre la
cohomologie g\'en\'eralis\'ee dans la formule de d\'evissage.

Notons enfin, que sans l'utilisation du formalisme des champs, la
construction du morphisme $\phi_{F}$ pourrait para\^itre un peu technique,
alors qu'il s'agit d'une construction tout \`a fait naturelle, consistant \`a
diagonaliser chaque fibr\'e vectoriel suivant les actions des
automorphismes du champ $F$ (~construction qui existe au niveau des
cat\'egories~).  \\

En gardant \`a l'esprit l'isomorphisme $\phi_{F}$, il est alors tr\`es
naturel de penser que la cohomologie du champ $I_{F}$ est le bon objet \`a
consid\'erer. Nous avons choisi d'appeler cette cohomologie la
"cohomologie \`a coefficients dans les repr\'esentations". Le choix de
cette terminologie r\'eside dans le fait qu'elle s'interpr\`ete
effectivement comme la cohomologie de l'espace de modules de $F$, mais
o\`u les coefficients des cycles sont pris dans les fonctions centrales
des groupes d'automorphismes de ses points.

Un fois cette id\'ee d\'egag\'ee, on peut dire, en exag\'erant tout de m\^eme un
peu, que la d\'emonstration du th\'eor\`eme de Grothendieck-Riemann-Roch
n'est plus qu'un gros exercice technique. La m\'ethode que nous avons
choisie pour r\'esoudre cet "exercice" est la suivante. Dans un premier
temps nous nous inspirons des d\'emonstrations des formules de
Lefschetz-Riemann-Roch (~\cite{bfm,th3}~) pour d\'emontrer le th\'eor\`eme
dans le cas d'un morphisme fortement projectif. Ensuite, \`a l'aide de
techniques de descente par quasi-enveloppes de Chow (~notion qui
remplace celle des enveloppes de Chow de \cite{g3}~), et du cas
pr\'ec\'edent, nous r\'eduisons le probl\`eme au cas des champs classifiants
de groupes finis. Le th\'eor\`eme devient alors \'equivalent \`a des formules
classiques de la th\'eorie des repr\'esentations des groupes finis. \\

Avant de d\'ecrire plus en d\'etail les chapitres \`a venir,
notons que la "cohomologie \`a coefficients dans les
repr\'esentations" permet aussi d'interpr\'eter de nombreux travaux ant\'erieurs.

\begin{itemize}
\item La premi\`ere r\'ef\'erence que nous connaissons dans laquelle on peut
reconna\^itre cette notion, est l'article de P. Baum, W. Fulton, et R.
MacPherson sur la formule de Lefschetz-Riemann-Roch \cite{bfm}. En
effet, si l'on consid\`ere le champ quotient d'une vari\'et\'e par un
automorphisme d'ordre fini, sa "cohomologie \`a coefficients dans les
repr\'esentations" contient naturellement le r\'eceptacle du caract\`ere de
Chern \'equivariant construit dans \cite{bfm}.
La formule de Grothendieck-Riemann-Roch que l'on d\'emontre, appliqu\'ee
au morphisme naturel d'un tel champ vers le classifiant du groupe
engendr\'e par l'automorphisme, redonne alors la formule de
Lefschetz-Riemann-Roch d\'emontr\'ee dans \cite{bfm}. Plus g\'en\'eralement, les
formules de Riemann-Roch \'equivariantes se
retrouvent en appliquant notre th\'eor\`eme de Grothendieck-Riemann-Roch \`a la
projection d'un champ quotient par un sch\'ema en groupe vers le
classifiant de ce groupe.

\item Dans ses deux articles \cite{ka,ka2} T. Kawasaki d\'emontre des
formules d'indices et de Riemann-Roch pour des $V$-vari\'et\'es, qui, dans le
cadre des champs alg\'ebriques, s'interpr\`etent comme
des champs de Deligne-Mumford lisses sur le corps des nombres
complexes,
g\'en\'eriquement non-ramifi\'es (~i.e. des orbifolds complexes~).
Pour cela, il associe \`a toute $V$-vari\'et\'e $X$, une nouvelle $V$-vari\'et\'e
$\Sigma X$. Cette derni\`ere correspond exactement au champ des
ramifications (~\ref{d1.8}~) du champ $X$. Par ces identifications, la
formule d'Hirzebruch-Riemann-Roch qu'il d\'emontre est \'equivalente \`a
celle que nous d\'eduisons du th\'eor\`eme de Grothendieck-Riemann-Roch,
lorsqu'on l'applique au morphisme structural.

\item En g\'eom\'etrie non commutative, on trouve une description de
l'homologie p\'eriodique d'un groupoide \'etale diff\'erentiel
$X_{\bullet}$, \`a l'aide de la cohomologie du groupoide des lacets
$\Omega X_{\bullet}$ (~\cite[$6.12$]{cm}~). Or, si $F$ est le champ
diff\'erentiel associ\'e au groupoide $X_{\bullet}$, alors le champ associ\'e \`a
$\Omega X_{\bullet}$ s'identifie naturellement au champ des
ramifications de $F$. Ainsi, par ces identifications, le caract\`ere de
Chern non-commutatif pour les groupoides diff\'erentiels \'etales est une
version $\cal
C\mit^{\infty}$ du caract\`ere de Chern \`a "coefficients dans les
repr\'esentations".
\end{itemize}

\vspace{10mm}

Le corps de la th\`ese est constitu\'e de cinq chapitres et un appendice.
Le premier est un r\'esum\'e de notations et de d\'efinitions. Les chapitres
deux et trois vont ensembles, et traitent de la $K$-th\'eorie des champs
alg\'ebriques et des formules de Riemann-Roch. Dans le chapitre quatre,
on applique ces formules aux $\cal D$-modules. Le cinqui\`eme chapitre
est ind\'ependant de tous les autres, et rassemble des r\'esultats de
types "GAGA". Enfin, l'appendice comporte quatre parties qui
rassemblent des faits (~quelques fois avec d\'emonstrations~) qui nous
sont utiles tout au cours de ce texte. \\

Dans le premier chapitre, le lecteur trouvera deux paragraphes que
nous invitons \`a consid\'erer comme des annexes de notations et de
d\'efinitions (~bien que l'on y trouve des d\'emonstrations~). On s'y
r\'ef\'erera tout au long du texte, surtout au cours des chapitres $2$, $3$
et $4$.

Dans le premier paragraphe nous donnons essentiellement deux
d\'efinitions : les foncteurs de $K$-th\'eorie et de $K$-cohomologie. Nous
avons choisi de les d\'efinir dans le cadre le plus g\'en\'eral possible, \`a
savoir celui d'un site muni d'une cat\'egorie cofibr\'ee en cat\'egories
exactes (~dont l'exemple standard est celui du grand site lisse des
sch\'emas, muni de son champ des fibr\'es vectoriels~). Ce degr\'e de
g\'en\'eralit\'e nous permettra de traiter, de fa\c{c}on uniforme et sans r\'ep\'etitions,
aussi bien le cas des fibr\'es vectoriels, que celui
des faisceaux coh\'erents, des $\cal D$-modules, des fibr\'es avec
connexion ... Le seul r\'esultat d\'emontr\'e dans ce paragraphe est le
fait qu'il existe toujours une transformation naturelle de la
$K$-th\'eorie vers la $K$-cohomologie.

Le second paragraphe est consacr\'e aux champs alg\'ebriques. On y rappelle
les principales d\'efinitions, ainsi que quelques r\'esultats concernant
les espace de modules. Enfin, nous introduisons les
quasi-enveloppes de Chow (~qui remplaceront pour les champs alg\'ebriques
la notion d'enveloppe de Chow de \cite{g3}~), et on d\'emontre des
r\'esultats d'existence. \\

Dans le second chapitre on \'etudie la $G$-th\'eorie
et la $G$-cohomologie des champs alg\'ebriques. Nous n'insistons pas
beaucoup sur les propri\'et\'es g\'en\'erales (~fonctorialit\'es, localisation,
homotopie\dots~), mais nous concentrons plut\^ot nos efforts sur les
r\'esultats de descente et de d\'evissage.

Les th\'eor\`emes de descente sont de deux types : contravariants et
covariants. Dans le premier cas on montre que le foncteur de
$G$-th\'eorie v\'erifie la propri\'et\'e de descente \'etale au-dessus d'un espace
alg\'ebrique. C'est un r\'esultat tr\`es commode qui permet par exemple
d'utiliser la descente galoisienne, ainsi que de r\'eduire certains
\'enonc\'es au cas de champs quotients. Notons aussi que l'on montre que
le morphisme de la $G$-th\'eorie vers la $G$-cohomologie n'est pas un
isomorphisme (~m\^eme rationnellement~), ce qui est un ph\'enom\`ene
nouveau quand on compare avec le cas des sch\'emas. En r\'ealit\'e, ce
ph\'enom\`ene explique d\'ej\`a pourquoi il ne peut y avoir de formule de
Riemann-Roch \`a valeurs dans les groupes de Chow usuels (~ou m\^eme dans
la cohomologie usuelle~).

Les th\'eor\`emes de descente covariante
permettent \`a l'aide des quasi-enveloppe de Chow de ramener certains
calculs au cas des gerbes (~voir m\^eme des gerbes triviales~). La
moralit\'e de ces th\'eor\`emes est que l'on conna\^it la $G$-th\'eorie (~resp.
la $G$-cohomologie~) des champs alg\'ebriques, si l'on conna\^it celle des gerbes
(~resp. celle des sch\'emas~).

Comme nous l'avons expliqu\'e plus haut dans cette introduction, les
r\'esultats cl\'es sont les th\'eor\`emes de d\'evissage. Bien qu'inspir\'es de
constructions existantes en g\'eom\'etrie \'equivariante (~\cite{bfm,th3,vi}~),
ces r\'esultats sont nouveaux. De fa\c{c}on intuitive,
il s'agit de diagonaliser les fibr\'es vectoriels sur un champ
alg\'ebrique. Dans le cas des champs de Deligne-Mumford, tous les
automorphismes sont d'ordre fini, et donc diagonalisables s'ils sont
d'ordre premier aux caract\'eristiques. Ainsi, le th\'eor\`eme de
d\'evissage dans ce cas, compare la $K$-th\'eorie d'un champ $F$, \`a la
$K$-cohomologie de son champ des ramifications mod\'er\'e $I_{F}^{t}$
(~classifiant les automorphismes d'ordre premier aux caract\'eristiques~).
En contre-partie, dans le cas des champs alg\'ebriques g\'en\'eraux, il
nous faut utiliser une construction un peu plus technique, s'appuyant
sur un th\'eor\`eme de repr\'esentabilit\'e d'Alexandre Grothendieck. On
construit le "champ des sous-groupes de type multiplicatif" $D_{F}$,
qui remplacera le champ $I_{F}^{t}$. Le th\'eor\`eme de d\'evissage compare
alors la $K$-th\'eorie d'un champ $F$, avec la $K$-cohomologie du champ
$D_{F}$, tordue par le "faisceau des caract\`eres". Le formalisme
n\'ecessaire pour donner un sens \`a cette expression est bri\`evement
rappel\'e en appendice.

Notons aussi, qu'il existe une autre fa\c{c}on de traiter le cas g\'en\'eral.
Cela consiste \`a ne consid\'erer que les automorphismes d'ordre fini et
premier aux caract\'eristiques. A la fin du troisi\`eme chapitre, on d\'emontrera
donc aussi un th\'eor\`eme de
d\'evissage comparant la $K$-th\'eorie de $F$, avec la $K$-cohomologie du
champ des automorphisme d'ordre fini et mod\'er\'es $I_{F}^{t,f}$. La
justification de la pertinence d'une telle d\'efinition, est \`a chercher
dans le fait que pour un sch\'ema en groupes affine, les \'el\'ements d'ordre fini
forment un ensemble sch\'ematiquement dense dans l'ensemble des
\'el\'ements semi-simples. L'information contenue dans la $K$-cohomologie
de $I_{F}^{t,f}$ doit donc s'av\'erer suffisante pour d\'emontrer des
formules de Riemann-Roch, ce que nous montrerons \^etre le cas sous
certaines hypoth\`eses.

A la fin de ce chapitre nous nous int\'eressons \`a la description de la
$G$-th\'eorie des gerbes born\'ees par des groupes r\'eductifs. Ce r\'esultat
ne sera pas utilis\'e par la suite. \\

Le troisi\`eme chapitre est le c\oe ur de la th\`ese. On commence par y
\'etudier les propri\'et\'es g\'en\'erales de la cohomologie \`a coefficients dans
les repr\'esentations (~ou encore dans les caract\`eres~). On montre en
particulier l'existence du caract\`ere de Chern, et on d\'efinit la
transformation de Riemann-Roch pour des champs "bien ramifi\'es" (~qui
comprennent les champs lisses qui sont localement des quotients par des
groupes affines et lisses sur un corps~).

Ensuite, on d\'emontre les
diff\'erentes formules de Riemann-Roch (~Lefschetz-Riemann-Roch et
Grothendieck-Riemann-Roch~). Tous ces th\'eor\`emes sont nouveaux, et
g\'en\'eralisent les diff\'erents th\'eor\`emes de Riemann-Roch \'equivariants
d\'ej\`a existants (~\cite{bfm,th3}~), ainsi que la formule de
Riemann-Roch pour les $V$-vari\'et\'es de T. Kawasaki (~\cite{ka}~).

Les d\'emonstrations ont toutes deux \'etapes. La premi\`ere traite le cas
de morphismes fortement projectifs. Les d\'emonstrations propos\'ees ici
sont alors tr\`es proches de celles utilis\'ees en g\'eom\'etrie \'equivariante
(~d\'eformation vers le c\^one normal, et utilisation du calcul de la
$K$-th\'eorie des fibr\'es projectifs~). La seconde \'etape traite du cas
des morphismes propres g\'en\'eraux (~\'eventuellement non
repr\'esentables~). La m\'ethode consiste \`a utiliser la
premi\`ere \'etape et les quasi-enveloppes de Chow, afin de se ramener au
cas des champs classifiants de sch\'emas en groupes. Les th\'eor\`emes se
traduisent alors par des formules connues en th\'eorie des
repr\'esentations.

La fin de ce chapitre est consacr\'e \`a quelques exemples
d'applications de telles formules dans le cas des champs de
Deligne-Mumford (~formule d'Hirzebruch-Riemann-Roch, cas particulier
des courbes, formules de Gauss-Bonnet, de signature\dots~), ainsi qu'\`a
la comparaison de la cohomologie \`a coefficients dans les
repr\'esentations avec deux autres th\'eories cohomologiques : les groupes
de Chow du complexe de Gersten, et la cohomologie p\'eriodique. La
seconde comparaison nous semble int\'eressante, car elle confirme les
relations entre la th\'eorie des champs alg\'ebriques et la g\'eom\'etrie
non-commutative. \\

Depuis l'article de G. Laumon \cite{l}, on sait que le th\'eor\`eme de
Riemann-Roch permet de calculer des indices de $\cal D$-modules
alg\'ebriques. Dans ce quatri\`eme chapitre nous suivons cette id\'ee afin
de d\'emontrer des formules de Riemann-Roch pour les $\cal D$-modules
sur des champs de Deligne-Mumford. Nous n'avons pas trait\'e le cas des
champs alg\'ebriques g\'en\'eraux que par souci de simplicit\'e. Il semble
clair que les r\'esultats s'\'etendent \`a ce cas, et donnent
ainsi des th\'eor\`emes g\'en\'eralisant les th\'eor\`emes de R. Joshua
\cite{jos}, bien que les techniques utilis\'ees pour prendre en compte
la $K$-th\'eorie sup\'erieure, soient diff\'erentes.

Un autre r\'esultat qui nous semble nouveau (~m\^eme dans le cas des
sch\'emas~) est la description de la
$K$-th\'eorie des $\cal D$-modules holon\^omes, \`a l'aide de "fonctions
constructibles". \\

Dans le cinqui\`eme chapitre nous \'etudions les relations entre champs
alg\'ebriques complexes (~de Deligne-Mumford~), et champs analytiques.
Nous commen\c{c}ons par y d\'emontrer que les th\'eor\`emes "GAGA" restent valables.
Nous nous int\'eressons ensuite aux probl\`emes d'alg\'ebrisation des
champs analytiques. A ce sujet nous posons la question de savoir si
un champ analytique propre est alg\'ebrique si et seulement si son
espace de modules l'est. On d\'emontrera que ceci est vrai apr\`es
\'eclatement, et on proposera une m\'ethode, bas\'ee sur les id\'ees
de M. Artin (~\cite{a}~), pour d\'emontrer que ceci est
suffisant pour r\'epondre par l'affirmative \`a la question pr\'ec\'edente. \\

Enfin, nous avons rassembl\'e en appendice des r\'esultats concernant les
spectres, la
descente cohomologique, l'extension des coefficients, la
strictification des pseudo-foncteurs, et la th\'eorie de Hodge pour les
champs complexes. Le premier de ces appendices \`a pour but d'expliquer
tr\`es rapidement le fonctionnement des spectres aux lecteurs peu
habitu\'es \`a ce language. Les deux appendices suivants rassemblent des
r\'esultats qui seront utilis\'es tout au long du chapitre deux, dans l'\'etude
des foncteurs de $K$-th\'eorie. En ce qui concerne la strictification,
nous rappelons juste que la th\'eorie des pseudo-foncteurs est,
d'un point de vue de la th\'eorie de l'homotopie,
\'equivalente \`a celle des foncteurs stricts. Ceci nous permet de
contourner des difficult\'es (~dues au fait que la cat\'egorie des champs
est une $2$-cat\'egorie~) pour d\'efinir certains
objets (~spectres de $G$-th\'eorie des champs simpliciaux augment\'es \ref{chsimp},
construction du morphisme $\psi_{F}$ dans la preuve de \ref{th3.3},
\dots~). Dans la derni\`ere partie de cet appendice on montre tr\`es
bri\`evement comment la th\'eorie de Hodge reste valable pour des champs
de Deligne-Mumford complexes. Ces r\'esultats sont, semble-t-il, des
faits connus, bien que n'apparaissant pas, ou peu (~\cite{tel}~), dans la
litt\'erature.

\end{section}

\newpage

\underline{\bf Notations et Conventions:} \rm \\

Pour une cat\'egorie $C$, on notera $Ob(C)$ son ensemble d'objets,
$Fl(C)$ son ensemble de morphismes, et $\pi_{0}(C)$ l'ensemble des
classes d'isomorphie d'objets de $C$. \\

La cat\'egorie simpliciale standard sera not\'ee $\Delta$. Pour
toute cat\'egorie $C$, la cat\'egorie des objets simpliciaux de $C$ est
$$SC:=Hom_{Cat}(\Delta^{op},C)$$

La cat\'egorie des ensembles sera not\'ee $Ens$, celle des groupes
$Gp$, et celle des groupes ab\'eliens $Ab$.

Pour tous groupes ab\'eliens $M$ et $A$, on notera
$M_{A}:=M\otimes_{\bf Z}A$. \\

Un groupoide est une cat\'egorie pour laquelle tout morphisme admet un
inverse. La $2$-cat\'egorie des groupoides est celle dont les objets sont
les groupoides, les $1$-morphismes sont les foncteurs, et les $2$-morphismes
les transformations naturelles entre foncteurs. \\

Si $C$ est une cat\'egorie de mod\`eles ferm\'ee au sens de \cite{q2}, nous
noterons $HoC$ la cat\'egorie homotopique associ\'ee.

Une r\'esolution injective d'un objet $X \in Ob(C)$, est une
cofibration triviale $X \hookrightarrow X'$, avec $X'$ fibrant. \\

La cat\'egorie des ensembles simpliciaux est not\'ee $SEns$. Par
convention, on appliquera syst\'ematiquement la construction de
Kan d\'ecrite dans \cite{j2}, et on supposera donc qu'ils sont tous
fibrants.

Si $I$ est une cat\'egorie, et $H : I \longrightarrow SEns$ un
pr\'efaisceau, on notera sa limite homotopique et sa colimite
homotopique (~\cite[$XI$, $XII$]{bk}~) par
$$holim_{I}H \qquad hocolim_{I}H$$

La cat\'egorie des spectres est not\'ee $Sp$. Pour tout objet $E$ de
$Sp$, nous noterons $E_{[n]}\in Ob(SEns)$ son "$n$-\`eme \'etage".

Si $I$ est une cat\'egorie,
et $H : I \longrightarrow Sp$ un pr\'efaisceau, sa limite homotopique et
sa colimite homotopique (~\cite[$5$]{th4}~) seront not\'ees par
$$holim_{I}H \qquad hocolim_{I}H$$

Si $X$ est un espace alg\'ebrique, $X_{et}$ (~resp. $X_{li}$~) d\'esignera
le site des espaces alg\'ebriques \'etales et de type fini sur $X$ (~resp.
lisses et de
type fini sur $X$~) muni de la topologie
\'etale (~resp. lisse~).

\newpage
\begin{part}*{\Huge Premi\`ere Partie :\\
\vspace{25mm}Th\'eor\`emes de Grothendieck-Riemann-Roch \large}
\newpage
\begin{section}{Chapitre $1$ : G\'en\'eralit\'es sur les champs}
\hspace{5mm}
Dans ce chapitre nous fixerons les notations et les d\'efinitions dont
nous aurons besoin par la suite. Comme il ne contient aucun r\'esultat
vraiment nouveau, nous invitons le lecteur \`a le consid\'erer comme une
annexe de notations.

Dans un premier temps nous rappelons les relations entre champs et
pr\'efaisceaux simpliciaux. Cela nous permettra par la suite de donner
un cadre naturel pour la cohomologie d'un champ \`a valeurs dans un
pr\'efaisceau en spectres, formalisme qu'il est indispensable de
poss\'eder pour \'etudier la $K$-th\'eorie de tels objets. Comme
nous aurons \`a faire de nombreuses constructions directement au niveau
des spectres, en particulier pour appliquer des techniques de
descente, nous avons
choisi de travailler dans le cadre des cat\'egories homotopiques.

Dans la seconde partie du chapitre nous fixerons quelques d\'efinitions
concernant les champs alg\'ebriques. Certaines se trouvent dans la
litt\'erature classique (~\cite{dm,lm}~), d'autres pas (~\ref{d1.12}~). Le
lecteur y trouvera aussi des \'enonc\'es sur les (~quasi~) enveloppes de Chow dans
le cadre des champs alg\'ebriques. Cette notion nous sera tr\`es utile
pour "approximer" certains champs alg\'ebriques par des gerbes.

\begin{subsection}{Champs, pr\'efaisceaux simpliciaux et pr\'efaisceaux
en spectres}
\hspace{5mm}
Tout au long de ce paragraphe, $C$ d\'esignera un site de Grothendieck,
dans lequel les produits fibr\'es existent.

\begin{subsubsection}{Champs et pr\'efaisceaux simpliciaux}

\hspace{5mm}
Soit $SPr(C)$ la cat\'egorie des pr\'efaisceaux
simpliciaux sur $C$. D'apr\`es \cite{j2}, c'est une cat\'egorie de mod\`eles
ferm\'ee simpliciale. Rappelons que lorsque $C$ poss\`ede suffisamment de
points (~par exemple lorsque $C$ est le grand site \'etale des
sch\'emas~), un morphisme $f : F \longrightarrow F'$ entre deux
pr\'efaisceaux simpliciaux est une \'equivalence faible, si pour tout
point $x$, le morphisme induit sur les fibres $f_{x} : F_{x}
\longrightarrow F_{x}'$ est une \'equivalence faible.

Si $F$ et $G$ sont deux objets de $SPr(C)$, nous noterons
$Hom_{s}(F,G)$ l'ensemble simplicial des morphismes de $F$ dans $G$,
et
$$\bf R\mit Hom_{s}(F,G):=Hom_{s}(F,HG)$$
o\`u l'on a choisi une r\'esolution injective $G \hookrightarrow HG$.
Comme les r\'esolutions injectives sont essentiellement uniques,
$\bf R\mit Hom_{s}(F,G)$ est un objet d\'etermin\'e \`a isomorphisme unique
pr\`es dans la cat\'egorie homotopique $HoSEns$. \\

Soit $U \longrightarrow X$ un morphisme de $C$, et
$\cal N\mit(U/X)$ son nerf. C'est l'objet simplicial de $C$ d\'efini par
$$\begin{array}{cccc}
\cal N\mit(U/X) : & \Delta^{op} & \longrightarrow & C \\
                  & [p]         & \mapsto         &
                  \underbrace{U\times_{X} U
                  \times_{X} \dots \times_{X} U}_{p \; fois}
\end{array}$$
Ainsi, pour tout pr\'efaisceau simplicial $F$ sur $C$, on dispose du
foncteur compos\'e
$$\begin{array}{cccc}
F\circ \cal N\mit(U/X) : & \Delta & \longrightarrow & SEns \\
                         & [p]         & \mapsto    & F(\cal
                         N\mit(U/X)([p]))
\end{array}$$
Nous noterons alors
$$\bf H\mit(U/X,F):=holim_{\Delta}(F\circ \cal N\mit(U/X))$$
Dans le cas o\`u $U \longrightarrow X$ est un morphisme couvrant,
$\bf H\mit(U/X,F)$ est l'espace de cohomologie de $\Check{C}$ech du
recouvrement
$U \longrightarrow X$ \`a coefficients dans $F$. Par la propri\'et\'e
universelle des limites homotopiques, il existe un morphisme naturel
d'ensembles simpliciaux
$$F(X) \longrightarrow \bf H\mit(U/X,F)$$

\begin{df}\label{d1.1}
Un objet $F$ de $SPr(C)$ est appel\'e flasque, si pour tout morphisme
couvrant $U \longrightarrow X$ de $C$, le morphisme naturel
$$F(X) \longrightarrow \bf H\mit(U/X,F)$$
est une \'equivalence faible.
\end{df}

La relation entre pr\'efaisceaux flasques et fibrants est donn\'ee par le
th\'eor\`eme de descente, dont une partie de la d\'emonstration figure en appendice
(~\ref{desc}~).

\begin{thm}\label{th1.1}
Un pr\'efaisceau simplicial $F$ sur $C$ est flasque, si et seulement si
pour toute r\'esolution injective
$$F \hookrightarrow HF$$
et tout objet $X \in Ob(C)$, le morphisme
$$F(X) \longrightarrow HF(X)$$
est une \'equivalence faible.
\end{thm}

\underline{Remarque:} Dans la terminologie de \cite{j2}, notre notion
de flasque se traduit par "flasque par rapport \`a tout objet $X$ de
$C$".

Il existe aussi une notion analogue pour les pr\'efaisceaux en
spectres, et le th\'eor\`eme pr\'ec\'edent reste encore valable.\\

Rappelons (~\cite{lm}~) qu'une cat\'egorie fibr\'ee en groupoides sur $C$
est la donn\'ee d'un cat\'egorie $\cal C$ et d'un foncteur
$$\pi : \cal C\mit \longrightarrow C$$
v\'erifiant les deux conditions suivantes

\begin{enumerate}
\item Pour tout morphisme $f :Y \longrightarrow X$ dans $C$, et tout
objet $x \in Ob(\cal C\mit)$ tel que $\pi(x)=X$, il existe un morphisme
$u : y \longrightarrow x$ dans $\cal C$ tel que $\pi(u)=f$.
\item Pour toute paire de morphismes $u : y \longrightarrow x$ et
$v : z \longrightarrow x$ dans $\cal C$, et tout morphisme
$f : \pi(y) \longrightarrow \pi(z)$ dans $C$ tel que $\pi(v)\circ
\pi(f)=\pi(u)$, il existe un unique morphisme
$w : y \longrightarrow z$ dans $\cal C$ tel que $\pi(w)=f$.
\end{enumerate}

Pour tout pr\'efaisceau d'ensembles $E$ sur $C$, on d\'efinit
la cat\'egorie fibr\'ee en groupoides $\pi : \widetilde{E} \longrightarrow
C$ de la fa\c{c}on suivante

\begin{itemize}
\item les objets de $\widetilde{E}$ sont les couples $(X,s)$ o\`u $X\in
Ob(C)$, et $s \in E(X)$
\item un morphisme $u : (Y,t) \longrightarrow (X,s)$ est la donn\'ee d'un
morphisme
$f : Y \longrightarrow X$ dans $C$, tel que $u^{*}(s)=t$
\item le foncteur $\pi$ est d\'efini par $\pi(X,s)=X$, et $\pi(u)=f$.
\end{itemize}

Si $\pi : \cal C\mit \longrightarrow C$ et $\pi' : \cal C\mit'
\longrightarrow C$
sont deux cat\'egories fibr\'ees en groupoides sur $C$, nous noterons
$Hom_{C}(\cal C\mit,\cal C\mit')$
la cat\'egorie des foncteurs de $\cal C$ vers $\cal C\mit'$ qui
commutent avec $\pi$ et $\pi'$. Remarquons que cette cat\'egorie est en
r\'ealit\'e un groupoide (~\cite{lm}~).

\begin{df}\label{d1.2}
Soit $\pi : \cal C\mit \longrightarrow C$ une cat\'egorie fibr\'ee en
groupoides sur $C$. Le pr\'efaisceau en groupoides associ\'e est d\'efini par
$$\begin{array}{cccc}
F_{\cal C\mit} : & C & \longrightarrow & Gpd \\
                 & X & \mapsto         &
                 Hom_{C}(\widetilde{X},\cal C\mit)
\end{array}$$
o\`u $\widetilde{X}$ est le pr\'efaisceau repr\'esent\'e par $X$.

Le pr\'efaisceau simplicial d\'eduit de $F_{\cal C}$ par le foncteur qui \`a
un groupoide associe son ensemble simplicial classifiant, sera not\'e
$BF_{\cal C}$.
\end{df}

Rappelons qu'un cat\'egorie fibr\'ee en groupoides $\cal C$ est un champ
en groupoides,
si toutes les donn\'ees de descente sont effectives (~\cite{lm}~).
Par abus de langage le mot "champ" signifiera toujours "champ en
groupoides", sauf mention explicite du contraire.
Il est alors facile de voir que $\cal C$ est un champ si et seulement
si $BF_{\cal C}$ est flasque. La cat\'egorie des champs en groupoides sur $C$
sera not\'ee $Ch(C)$.
Nous noterons aussi $HoCh(C)$ sa cat\'egorie homotopique. C'est
la cat\'egorie qui poss\`ede les m\^emes objets que $Ch(C)$, et dont
l'ensemble des morphismes entre deux champs $\cal C$ et $\cal C\mit'$
est
$$Hom_{HoCh(C)}(\cal C\mit,\cal C\mit'):=
\pi_{0}Hom_{C}(\cal C\mit,\cal C\mit')$$

En r\'ealit\'e le foncteur $\cal C\mit \mapsto BF_{\cal C}$ induit une
\'equivalence de la cat\'egorie homotopique des champs sur $C$ avec celle
des pr\'efaisceaux simpliciaux flasques $1$-tronqu\'es et morphismes
flexibles (~\cite{s}~). Remarquons que cette derni\`ere est
elle-m\^eme \'equivalente \`a la sous-cat\'egorie pleine de $HoSPr(C)$,
form\'ee des objets $1$-tronqu\'es.

Ainsi, il nous est permis de voir un champ (~ou une cat\'egorie fibr\'ee
en groupoides~) comme un pr\'efaisceau simplicial. De cette fa\c{c}on, pour
tout champ $\cal C$ sur $C$, et tout $F \in SPr(C)$, nous pouvons
d\'efinir
$$\bf R\mit Hom_{s}(\cal C\mit,F):=\bf R\mit Hom_{s}(BF_{\cal C},F) \in HoSEns$$

\end{subsubsection}

\begin{subsubsection}{Cohomologie g\'en\'eralis\'ee des pr\'efaisceaux
Simpliciaux}

\hspace{5mm}
Soit $Sp(C)$ la cat\'egorie des pr\'efaisceaux en spectres sur $C$. Nous
savons d'apr\`es \cite[$2.53$]{j}, que c'est une cat\'egorie de mod\`eles ferm\'ee
munie de "Hom" internes que l'on notera $\underline{Hom}_{sp}(.,.)$.
Le spectre des morphismes entre deux objets $F,G \in Sp(C)$, est
d\'efini par
$$Hom_{sp}(F,G):=lim_{C}(\underline{Hom}_{sp}(F,G))$$
On d\'efinit alors
$$\bf R\mit Hom_{sp}(F,G):=Hom_{sp}(F,HG)$$
o\`u l'on a choisi une r\'esolution injective $G \hookrightarrow HG$.
C'est un objet d\'etermin\'e \`a isomorphisme unique pr\`es dans $HoSp$.\\

Si $F$ est un pr\'efaisceau simplicial, et $K$ un pr\'efaisceau en
spectres, on peut d\'efinir exactement comme dans \ref{spectre} le
pr\'efaisceau en spectres des morphismes
$$\underline{Hom}_{sp}(F,K).$$
On dipose alors du spectre des morphismes de $F$ vers $K$, d\'efini par
$$Hom_{sp}(F,K):=lim_{C}(\underline{Hom}_{sp}(F,G)).$$

\begin{df}\label{d1.3}
Soit $F$ un pr\'efaisceau simplicial, et $K$ un pr\'efaisceau en spectres
sur $C$. Alors le spectre de cohomologie de $F$ \`a coefficients dans $K$ est
d\'efini par
$$\bf H\mit (F,K):=\bf R\mit Hom_{sp}(F,K)$$
Si $\cal C \mit \longrightarrow C$ est une cat\'egorie fibr\'ee en
groupoides sur $C$, son spectre de cohomologie \`a coefficients dans
$K$ est d\'efini par
$$\bf H\mit(\cal C\mit,K):=\bf H\mit(F_{\cal C\mit},K)$$
\end{df}

Remarquons que $\bf H$ d\'efinit des bifoncteurs
$$\bf H\mit : HoSpr(C)\times HoSp(C) \longrightarrow HoSp$$
$$\bf H\mit : HoCh(C)\times HoSp(C) \longrightarrow HoSp$$

\end{subsubsection}

\end{subsection}

\begin{subsection}{Champs et spectres de $K$-th\'eorie}\label{s1}

\hspace{5mm}
Pour ce paragraphe, on supposera de
plus qu'il existe une cat\'egorie cofibr\'ee en cat\'egories exactes
$$p : \cal E\mit \longrightarrow C$$
C'est \`a dire que $\cal E$ est une cat\'egorie exacte,
et $p$ est un foncteur tel que

\begin{itemize}
\item Pour tout objet $x \in Ob\cal E$, et tout morphisme de $C$ $f :
Y \longrightarrow X$, il existe une "image r\'eciproque de $x$ par
$f$". C'est \`a dire, il existe un morphisme $u : y \longrightarrow x$
dans $\cal E$, tel que $p(u)=f$, et tel que
pour tout morphisme de $\cal E$ $v : z \longrightarrow x$ avec
$p(v)=f$, il existe un unique morphisme $w : y \longrightarrow z$ tel
que $p(w)=Id$.
\item Soit $\xymatrix{Z \ar[r]^{g} & Y \ar[r]^{f} & X}$ sont deux
morphismes de $C$, $x \in Ob \cal E$,
et $y~\in~Ob\cal E$ une image r\'eciproque de $x$ par $f$. Alors toute
image r\'eciproque de $y$ par $g$ est une image r\'eciproque de $x$ par
$f\circ g$.
\item Si $\xymatrix{x \ar[r]^{u} & y \ar[r]^{v} & z}$ est une suite
exacte dans $\cal E$, alors $p(u)~=~p(v)~=~Id$.
\item Si $E : \xymatrix{x \ar[r]^{u} & y \ar[r]^{v} & z}$ est une suite
exacte avec $p(x)=X$, alors pour tout morphisme de $C$ $f : Y
\longrightarrow X$, toute image r\'eciproque de $E$ par $f$ est
encore une suite exacte.
\end{itemize}

Nous lui associons le pr\'efaisceau en cat\'egories exactes suivant
$$\begin{array}{cccc}
F_{\cal E} : & C & \longrightarrow & CatEx \\
             & X & \mapsto         & Hom_{Cart}(X,\cal E\mit)
\end{array}$$
o\`u $CatEx$ est la cat\'egorie des cat\'egories exactes et foncteurs
exacts, et $Hom_{Cart}(X,\cal E\mit)$ la
sous-cat\'egorie de $Hom_{C}(\widetilde{X},\cal E\mit)$ des sections
cart\'esiennes (~\cite[$1.1.1$]{gi}~).

L'exemple standard que l'on utilisera est celui o\`u $C=(Sch/S)_{li}$ est le
site des
sch\'emas muni de la topologie lisse, et $\cal E$ la cat\'egorie cofibr\'ee
des fibr\'es vectoriels sur $C$. Ses objets sont les couples $(X,V)$,
avec $X$ un sch\'ema sur $S$, et $V$ un fibr\'e vectoriel sur $X$, et un
morphisme entre $(Y,W)$ et $(X,V)$ est la donn\'ee d'un morphisme de
sch\'emas $f : Y \longrightarrow X$ et d'un morphisme de fibr\'es
vectoriels sur $Y$, $W \longrightarrow f^{*}(V)$.\\

Si $\cal C$ est un champ, on posera
$$\int_{\cal C}\cal E\mit:=Hom_{Cart}(\cal C\mit,\cal E\mit)$$
la cat\'egorie des morphismes cart\'esiens de champs sur $C$
(~\cite[$1.1.1$]{gi}~). \\

Si on note $\cal C\mit(X)$ la cat\'egorie des
fl\`eches de $\cal C$ au-dessus de l'identit\'e de $X$,
alors la cat\'egorie $\int_{\cal C}\cal E$ est
\'equivalente \`a la cat\'egorie suivante \\

\underline{un objet :} est d\'efini par la donn\'ee suivante~:
\begin{itemize}
\item
Pour tout objet $X \in Ob(C)$, et tout objet $s \in Ob(\cal C\mit(X))$
la donn\'ee d'un objet $V_{(s)} \in Ob\cal E\mit(X)$
\item
Pour tout morphisme de $C$, $f : Y \longrightarrow X$, et
toute paire d'objets
$s\in Ob(\cal C\mit(X))$, $s' \in Ob(\cal C\mit(Y))$, et
tout isomorphisme dans $\cal C\mit(Y)$, $h~:~f^{*}(s)~\simeq~t$,
un isomorphisme dans $\cal E\mit(Y)$
$$\phi_{s,f,h} : f^{*}V_{(s)} \simeq V_{(t)}$$
\item
Pour toute paire de morphismes de $C$
$$\xymatrix{Z \ar[r]^{g} & Y \ar[r]^{f} & X}$$
tout triplets d'objets $s \in Ob(\cal C\mit(X))$, $t \in Ob(\cal
C\mit(Y))$,
$u \in Ob(\cal C\mit(Z))$, et toute paire d'isomorphismes dans
$\cal E\mit(Y)$ et $\cal E\mit(Z)$
$$h : f^{*}(s) \simeq  t$$
$$j : g^{*}(t) \simeq u$$
une \'egalit\'e dans $\cal E\mit(Z)$
$$\phi_{t,g,j}\circ g^{*}\phi_{s,f,h} = \phi_{u,f\circ g,h\circ j}$$
\end{itemize}

\underline{un morphisme :} entre $V$ et $W$ est d\'efini par la donn\'ee
suivante~:
\begin{itemize}
\item
Pour tout objet
$s \in Ob(\cal C\mit(X))$,
un morphisme dans $\cal E\mit(E)$
$$a_{s} : V_{(s)} \longrightarrow W_{(s)}$$
\item
Pour tout morphisme de $C$, $f : Y \longrightarrow X$,
toute paire d'objets $s~\in~Ob\cal C\mit(X)$ et
$t~\in~\cal C\mit(Y)$, et tout isomorphisme dans $\cal E\mit(Y)$,
$h~:~f^{*}(s)~\simeq~t$, une \'egalit\'e dans $\cal E\mit(Y)$
$$\phi_{s,f,h}^{W} \circ f^{*}(a_{s})=a_{t} \circ \phi_{s,f,h}^{V}$$
\end{itemize}

Ainsi, $\int_{\cal C}\cal E$ est \'equivalent \`a la cat\'egorie des
pseudo-transformations naturelles entre les pseudo-foncteurs
(~\ref{strict}~)
$$\begin{array}{ccc}
C & \longrightarrow & Cat \\
X & \mapsto & \cal C\mit(X)
\end{array}$$
et
$$\begin{array}{ccc}
C & \longrightarrow & Cat \\
X & \mapsto & \cal E\mit(X)
\end{array}$$

Nous noterons aussi $\int_{C}\cal E$ la cat\'egorie des sections
cart\'esiennes globales de $\cal E$ sur $C$. En clair
$$\int_{C}\cal E\mit := Hom_{Cart}(\widetilde{*},\cal E\mit)$$
o\`u $*$ est le pr\'efaisceau d'ensembles constant associ\'e \`a un ensemble \`a
un \'el\'ement.

\begin{df}\label{d1.4}
Le pr\'efaisceau en spectres de $K$-th\'eorie associ\'e au couple $(C,\cal
E\mit)$ est d\'efini par
$$\begin{array}{cccc}
\underline{K} : & C & \longrightarrow & Sp \\
                & X & \mapsto & K(F_{\cal E}(X))
\end{array}$$
o\`u
$$K : CatEx \longrightarrow Sp$$
est le foncteur de $K$-th\'eorie d\'efini dans \cite[$1.3$]{wal} par exemple.
On d\'efinit alors la $K$-cohomologie d'un champ $\cal C$
\`a coefficients dans $\cal E$ par
$$\underline{\bf K}(\cal C\mit):=\bf H\mit(\cal C\mit,\underline{K})$$
Le spectre de $K$-th\'eorie d'un champ $\cal C$ \`a coefficients
dans $\cal E$ est d\'efini
$$\bf K\mit(\cal C\mit):=K(\int_{\cal C}\cal E\mit)$$
Les groupes de $K$-cohomologie et de $K$-th\'eorie  de $\cal C$
\`a coefficients dans $\cal E$ sont d\'efinis respectivement par
$$\underline{\bf K}_{m}(\cal C\mit):=\pi_{m}\underline{\bf K}(\cal
C\mit)$$
$$\bf K\mit_{m}(\cal C\mit):=\pi_{m}\bf K\mit(\cal C\mit)$$
\end{df}

Les correspondances $\cal C \mit \mapsto \underline{\bf K}(\cal C\mit)$ et
$\cal C\mit \mapsto \bf K\mit(\cal C\mit)$, d\'efinissent des foncteurs
$$\underline{\bf K} : HoSPr(C) \longrightarrow HoSp$$
$$\bf K\mit : HoCh(C) \longrightarrow HoSp$$
En effet, pour $\underline{\bf K}$ cela provient directement de sa
d\'efinition. Pour le second, il suffit de garder \`a l'esprit que
$$\cal C\mit \mapsto \int_{\cal C}\cal E$$
transforme \'equivalences de champs en \'equivalences de cat\'egories, et
d\'etermine donc un foncteur
$$HoCh(C) \longrightarrow HoCatEx$$

\begin{prop}\label{p1.1}
Il existe une transformation naturelle de foncteur
$$can : \bf K\mit \longrightarrow \underline{\bf K}$$
\end{prop}

\underline{\bf Preuve:} \rm Soit $\cal C$ un champ en groupoides sur
$C$. Alors, pour chaque objet $X$ de $C$, et chaque morphisme de champs
$s : \widetilde{X} \longrightarrow \cal C$, on dispose du foncteur
image r\'eciproque
$$s^{*} : \int_{\cal C}\cal E\mit:=Hom_{Cart}(\cal C\mit,\cal E\mit)
\longrightarrow Hom_{Cart}(\widetilde{X},\cal E\mit)=:F_{\cal E}(X)$$
Ces images r\'eciproques v\'erifient de plus $(s\circ s')^{*}=s^{*}\circ
(s')^{*}$. Ainsi, pour $X$ variable dans $Ob(C)$, les foncteurs
$$Hom_{Cart}(\widetilde{X},\cal C\mit) \longrightarrow F_{\cal E}(X)$$
d\'efinissent un foncteur exact de cat\'egories exactes
$$\int_{\cal C}\cal E\mit \longrightarrow
Hom_{C}(F_{\cal C},F_{\cal E})$$
o\`u $Hom_{C}(F_{\cal C},F_{\cal E})$ est la cat\'egorie exacte des
morphismes de pr\'efaisceaux en cat\'egories sur $C$. Or, il existe un morphisme
naturel de cat\'egories
$$WHom_{C}(F_{\cal C},F_{\cal E}) \longrightarrow
Hom_{C}(F_{\cal C},WF_{\cal E})$$
o\`u $W$ d\'esigne la construction de Waldhausen (~\cite[$1.3$]{wal}~) qui
a une cat\'egorie exacte associe son spectre de $K$-th\'eorie. Ainsi, on
trouve un foncteur naturel en $\cal C$
$$W\int_{\cal C}\cal E\mit \longrightarrow
Hom_{C}(F_{\cal C},WF_{\cal E})$$
que l'on compose avec le foncteur classifiant
$$BW\int_{\cal C}\cal E\mit \longrightarrow
BHom_{C}(F_{\cal C},WF_{\cal E})$$
Comme il existe un morphisme canonique
$$BHom_{Cat}(A,B) \longrightarrow Hom_{SEns}(BA,BB)$$
on en d\'eduit un morphisme d'ensembles simpliciaux
$$\bf K\mit(\cal C\mit)_{[0]}:=BW\int_{\cal C}\cal E\mit
\longrightarrow Hom_{s}(BF_{\cal C},\underline{K}_{[0]})$$
Par la naturalit\'e de cette construction, ce morphisme s'\'etend en un
morphisme de spectres
$$\bf K\mit(\cal C\mit) \longrightarrow Hom_{sp}(S(BF_{\cal
C}),\underline{K})$$
On peut alors composer ce morphisme avec une r\'esolution injective
$$\underline{K} \hookrightarrow H\underline{K}$$
pour obtenir le morphisme cherch\'e
$$\bf K\mit(\cal C\mit) \longrightarrow \bf H\mit(\cal
C\mit,\underline{K})$$
Une fois que cette r\'esolution injective a \'et\'e choisie, ce morphisme
est fonctoriel en $\cal C$. $\Box$ \\

Voici deux exemples qui montrent que le morphisme de $HoSp$
$$\bf K\mit(\cal C\mit) \longrightarrow \underline{\bf K}(\cal C\mit)$$
peut \^etre, ou ne pas \^etre, un isomorphisme.

Consid\'erons $C=(QProj/k)_{Zar}$, le gros site des sch\'emas quasi-projectifs
sur un corps $k$, muni de la topologie de Zariski. Prenons
$\cal E$ la cat\'egorie cofibr\'ee des faisceaux coh\'erents localement libres
et de rang fini sur $C$, et $\bf K$ et $\underline{\bf K}$ les foncteurs
de $K$-th\'eorie et de $K$-cohomologie associ\'es.

Alors, pour $\cal C\mit=X$ un sch\'ema de $C$, $\bf K\mit(F)$ est le
spectre de $K$-th\'eorie de la cat\'egorie des faisceaux localement libres
et de rang fini sur $X$, et $\underline{\bf K}(\cal C\mit)$ est le spectre
de cohomologie de $X_{Zar}$ \`a coefficients dans $\bf K$. On sait
alors que le morphisme canonique
$$\bf K\mit(X) \longrightarrow \bf H\mit(X_{Zar},\bf K\mit)$$
est une \'equivalence faible (~\cite[$10.5$]{th}~).

Prenons maintenant $C=(QProj/k)_{et}$, le gros site des sch\'emas quasi-projectifs
sur un corps $k$, muni de la topologie de \'etale. Dans ce cas $\underline{\bf K}(X)$ 
est le spectre de cohomologie de $X_{et}$ \`a coefficients dans $\bf K$. On sait alors
que le morphisme naturel 
$$\bf K\mit(X) \longrightarrow \bf H\mit(X_{Zar},\bf K\mit)$$
n'est une \'equivalence que rationnellement. \\

Supposons pour terminer que l'on dispose d'une autre cat\'egorie cofibr\'ee
en cat\'egories exactes $\cal E\mit'$ sur $C$,  qu'il existe un
produit tensoriel exact dans les deux variables
$$\otimes : \cal E\mit\times_{C} \cal E\mit \longrightarrow \cal E\mit$$
ainsi qu'une structure de $\cal E$-module sur $\cal E\mit'$
$$\otimes : \cal E\mit\times_{C} \cal E\mit' \longrightarrow \cal E\mit'$$
qui est exacte en la premi\`ere variable. Notons $\bf K\mit'$ et
$\underline{\bf K}'$ les foncteurs de $K$-th\'eorie et de
$K$-cohomologie \`a coefficients dans $\cal E\mit'$. On sait alors que l'on peut
construire des produits (~\cite[$5.3$]{j}~)
$$\bf K\mit \wedge \bf K\mit \longrightarrow \bf K$$
$$\bf K\mit \wedge \bf K\mit' \longrightarrow \bf K\mit'$$
$$\underline{\bf K}\wedge \underline{\bf K} \longrightarrow
\underline{\bf K}$$
$$\underline{\bf K}\wedge \underline{\bf K}' \longrightarrow
\underline{\bf K}'$$
qui font de $\bf K$ (~resp. $\underline{\bf K}$~) un foncteur en
spectres en anneaux, et de $\bf K\mit'$ (~resp. $\underline{\bf K}'$~)
un foncteur en spectres en $\bf K$-modules (~resp. $\underline{\bf
K}$-modules~).

\end{subsection}

\begin{subsection}{Champs alg\'ebriques}

\hspace{5mm}
Nous noterons $S$ un sch\'ema noeth\'erien int\`egre de dimension finie et
universellement japonais (~i.e. toutes les normalisations de sch\'emas
de type fini sur $S$ sont des morphismes finis~), et $(Esp/S)_{li}$ le gros
site des $S$-espaces alg\'ebriques essentiellement de type fini et s\'epar\'es
sur $S$,
o\`u un morphisme est couvrant s'il est lisse et surjectif.

A l'aide du lemme de Yoneda, nous identifierons la cat\'egorie
$(Esp/S)$ \`a une sous-cat\'egorie pleine des pr\'efaisceaux simpliciaux sur
$(Esp/S)_{li}$. Nous dirons alors qu'un objet $F \in
SPr((Esp/S)_{li})$ est repr\'esentable s'il est isomorphe dans
$HoSPr((Esp/S)_{li})$
\`a un objet provenant
de $(Esp/S)$. Par extension, nous appellerons tout objet
repr\'esentable un "$S$-espace alg\'ebrique".

Comme nous l'avons fait remarquer au $2.1.1$, nous pouvons d\'efinir les
champs en groupoides comme des pr\'efaisceaux simpliciaux. Il nous arrivera
cependant de les d\'efinir comme des cat\'egorie fibr\'ees en groupoides.
Comme nous les consid\'ererons dans la cat\'egorie homotopique, le point
de vu adopt\'e importe peu.

\begin{subsubsection}{Quelques d\'efinitions et propri\'et\'es}\label{s2}

\begin{df}\label{d1.5}
\begin{enumerate}
\item
Un morphisme de pr\'efaisceaux simpliciaux $f~:~F~\longrightarrow~F'$
sur $(Esp/S)_{li}$ est repr\'esentable, si pour tout $S$-espace
alg\'ebrique $X$, et tout morphisme $s : X \longrightarrow F'$ de
pr\'efaisceaux simpliciaux, le
pr\'efaisceau simplicial $f^{-1}(X):=F\times_{F'}X$ est repr\'esentable.
\item
Soit $\bf P$ un type de morphismes de $S$-espaces alg\'ebriques,
stable par changement de base (~e.g. lisse, \'etale, plat, immersion
ferm\'ee, immersion ouverte, surjectif, de type fini\dots~). Un
morphisme repr\'esentable de pr\'efaisceaux simpliciaux sur $(Esp/S)_{li}$
$f : F \longrightarrow F'$ est de type $\bf P$, si
pour chaque $S$-espace alg\'ebrique $X$, et chaque morphisme
$s : X \longrightarrow F'$ de pr\'efaisceaux simpliciaux,
le morphisme d'espaces alg\'ebriques
$$f_{X} : f^{-1}(X) \longrightarrow X$$
est de type $\bf P$.
\item
Soit $\bf P$ un type de morphismes de $S$-espaces alg\'ebriques,
local pour la topologie lisse (~e.g. localement d'intersection
compl\`ete, immersion r\'eguli\`ere\dots~). Un
morphisme repr\'esentable de pr\'efaisceaux simpliciaux sur $(Esp/S)_{li}$
$f : F \longrightarrow F'$ poss\`ede la propri\'et\'e $\bf P$, si
pour chaque $S$-espace alg\'ebrique $X$, et chaque morphisme
$s : X \longrightarrow F'$ repr\'esentable et lisse, le morphisme
d'espaces alg\'ebriques
$$f_{X} : f^{-1}(X) \longrightarrow X$$
est de type $\bf P$.
\item
Un pr\'efaisceau simplicial $1$-tronqu\'e sur $(Esp/S)_{li}$ $F$, est
alg\'ebrique (~quasi-s\'epar\'e et localement de type fini~), si
\begin{itemize}
\item le morphisme diagonal
$$\Delta : F \longrightarrow F\times_{S}F$$
est repr\'esentable, quasi-compact et s\'epar\'e
\item il existe un $S$-espace alg\'ebrique $X$, et un morphisme
(~automatiquement repr\'esentable par le premier point~), lisse et surjectif
$$f : X \longrightarrow F$$
\end{itemize}
On dira alors que $F$ est lisse (~resp. r\'egulier, normal, de type fini~)
sur $S$,
si on peut prendre $X$ lisse (~resp. r\'egulier, normal, de type fini~) sur $S$.
\item
une cat\'egorie fibr\'ee en groupoides $\cal C \mit$ sur $(Esp/S)_{li}$
est un champ alg\'ebrique, si
le pr\'efaisceau simplicial $F_{\cal C}$ est flasque et
alg\'ebrique.
\end{enumerate}
\end{df}

Remarquons qu'une cat\'egorie fibr\'ee en groupoides $\cal C \mit$ sur
$(Esp/S)_{li}$
est un champ alg\'ebrique dans la terminologie de \cite{lm}, si et
seulement si c'est un champ alg\'ebrique pour la d\'efinition pr\'ec\'edente.

\begin{df}\label{d1.6}
La $2$-cat\'egorie des champs alg\'ebriques sur $S$ sera not\'ee
$ChAlg(S)$. Le groupoide des $1$-morphismes de $\cal C$ vers
$\cal C\mit'$ sera not\'e $Hom_{Ch}(\cal C\mit,\cal C\mit')$. \\

La cat\'egorie homotopique des champs alg\'ebriques sur $S$, not\'ee
$HoChAlg(S)$, est l'image essentielle dans $HoSPr((Esp/S)_{li})$ du
foncteur
$$\begin{array}{ccc}
ChAlg(S) & \longrightarrow & HoSPr((Esp/S)_{li}) \\
\cal C\mit & \mapsto & F_{\cal C}
\end{array}$$
\end{df}

\underline{\bf Notations et Terminologie:} \rm
Le mot "morphisme" fera toujours r\'ef\'erence \`a un morphisme dans
$HoChAlg(S)$, alors que nous pr\'eciserons \\
"$1$-morphisme" pour ceux
dans $ChAlg(S)$. De m\^eme nous parlerons de "diagrammes commutatifs"
pour les diagrammes commutatifs de $HoChAlg(S)$, et de "diagrammes
$1$-commutatifs" pour ceux de $ChAlg(S)$.
Nous dirons que deux champs sont \'equivalents s'ils
sont isomorphes dans $HoChAlg(S)$.

Un $1$-morphisme $s : X \longrightarrow F$, avec $X$ un espace
alg\'ebrique sera appel\'e une section de $F$ au-dessus de $X$. Le lemme de
Yoneda permettant d'identifier canoniquement le groupoide des sections
au-dessus de $X$ avec $F(X)$, nous parlerons alors de l'objet $s \in
ObF(X)$ associ\'e \`a $s$. Ce n'est que l'image par $s$ de l'identit\'e.

Par la suite un "champ alg\'ebrique" sera toujours un champ alg\'ebrique
de type fini sur $S$. Dans les quelques cas o\`u les champs ne seront
pas de type fini, nous pr\'eciserons "champs alg\'ebriques localement de
type fini". En particulier, comme $S$ est noeth\'erien, tout champ
alg\'ebrique est aussi noeth\'erien.

Bien que nous travaillerons essentiellement dans $HoChAlg(S)$, nous
aurons besoin quelque-fois de revenir \`a $ChAlg(S)$ pour d\'efinir
certains objets.\\

\begin{df}\label{d1.7}
Un $1$-morphisme de champs alg\'ebriques $f : F \longrightarrow
F'$ est propre, si pour tout $S$-espace alg\'ebrique $X$, et
tout $1$-morphisme $s : X \longrightarrow F'$, il existe un
$S$-espace alg\'ebrique $Y$, et un diagramme commutatif
$$\xymatrix{Y \ar[d]_{p} \ar[dr]^{q} & \\
f^{-1}(X) \ar[r]_{f_{X}} & X }$$
avec $q$ propre, et $p$ surjective.
\end{df}

\underline{Remarque:} Comme la propri\'et\'e d'\^etre repr\'esentable pour
un $1$-morphisme est invariante par \'equivalence,
la notion de morphismes repr\'esentables dans $HoChAlg(S)$ est bien
d\'efinie (~comme morphismes isomorphes \`a des images de
$1$-morphismes repr\'esentables~). De m\^eme, pour un morphisme
repr\'esentable, la propri\'et\'e d'\^etre \'etale, surjectif, lisse, plat,
localement d'intersection compl\`ete, une immersion ferm\'ee, une
immersion ouverte ... poss\`ede un sens.

Il en est de m\^eme pour la notion de morphisme propre.

\begin{df}\label{d1.8}
Le champ des ramifications $I_{F}$ d'un champ alg\'ebrique $F$ est
d\'efini par
$$I_{F}:=F\times_{F\times_{S}F}F$$
On notera
$$\pi_{F} : I_{F} \longrightarrow F$$
le morphisme naturel.\\

Soit $F$ un champ alg\'ebrique.
\begin{enumerate}
\item
On dira que $F$ est s\'epar\'e, si le morphisme
$$\Delta : F \longrightarrow F\times_{S}F$$
est propre.
\item
Nous dirons que $F$ est de Deligne-Mumford, si
le morphisme
$$\Delta : F \longrightarrow F\times_{S}F$$
est non-ramifi\'e.
\item
Nous dirons que $F$ est $\Delta$-affine, si le morphisme
$$\Delta : F \longrightarrow F\times_{S}F$$
est affine.
\end{enumerate}
\end{df}

\textit{Par la suite, tous les champs de Deligne-Mumford que l'on rencontrera
seront suppos\'es s\'epar\'es.} \\

\underline{Remarques:} \begin{itemize}
\item
Soit $F$ un champ alg\'ebrique. Pour chaque objet
$X \in (Esp/S)$, et chaque objet $s \in ObF(X)$,
le faisceau des automorphismes de $s$
$$\begin{array}{cccc}
\underline{Aut}_{X}(s) : & (Esp/X)_{li} & \longrightarrow & Gp \\
                         & (~u : Y \rightarrow X~) & \mapsto &
                         Aut_{F(Y)}(u^{*}(s))
\end{array}$$
est repr\'esentable par un $X$-espace alg\'ebrique en groupes (~\cite{lm}~), not\'e
$Aut_{X}(s)$. Dire alors que $F$ est $\Delta$-affine, est
\'equivalent \`a dire que pour tout $X$ et $s$ comme ci-dessus,
$Aut_{X}(s)$ est affine sur $X$.
\item
D'apr\`es \cite{lm}, la d\'efinition $(2)$ pr\'ec\'edente est \'equivalente \`a la
d\'efinition donn\'ee dans \cite[$4.6$]{dm}, ce qui explique le choix de la
terminologie.
\item
Un champ alg\'ebrique s\'epar\'e de Deligne-Mumford est $\Delta$-affine. En
effet, $\Delta$ \'etant quasi-compact et non-ramifi\'e, il est
quasi-fini. Ainsi, $\Delta$ est propre et quasi-fini, donc fini, et en
particulier affine.
\item
Si $S$ est de caract\'eristique nulle, tout champ alg\'ebrique
$\Delta$-affine et s\'epar\'e est de Deligne-Mumford. En effet,
$\Delta$ \'etant affine et propre, il est fini, et donc quasi-fini.
C'est donc un champ de Deligne-Mumford d'apr\`es \cite[$7.17$]{vi2}.
Ceci n'est plus vrai pour $S$ g\'en\'eral.
\item
Remarquons aussi, qu'un $S$-espace alg\'ebrique est un champ
alg\'ebrique $F$ tel que
$$\Delta : F \longrightarrow F\times_{S}F$$
soit un monomorphisme. (~\cite{lm}~).
\end{itemize}

\underline{Exemple:} Soit $X$ un $S$-espace alg\'ebrique, et
$H \longrightarrow S$ un $X$-espace alg\'ebrique en groupes, lisse
sur $S$. On suppose que $H$ op\`ere sur $X$ au-dessus de $S$
$$a : X\times_{S} H \longrightarrow X$$
On d\'efinit
le champ classifiant $[X/H]$, comme la cat\'egorie fibr\'ee sur
$(Esp/S)_{li}$, dont la cat\'egorie fibre au-dessus d'un objet
$Y \in (Esp/S)$ est le groupoide des diagrammes
$$\xymatrix{ Y & \ar[l]_{p}  P \ar[r]^{f} & X }$$
o\`u $p$ est un $H$-fibr\'e principal, et $f$ un morphisme $H$-\'equivariant.
On sait que $[X/H]$ est un champ alg\'ebrique (~\cite{lm}~).

De plus, dans le cas o\`u $H$ est affine sur $S$, $[X/H]$ est
$\Delta$-affine. En effet, le champ $I_{F}$ est alors \'equivalent au
champ quotient $[\overline{X}/H]$, o\`u $\overline{X}$ est le
sous-espace alg\'ebrique en groupes de $X\times_{S}H$ des couples
$(x,h)$ tels que $h.x=x$. \\

\begin{df}\label{d1.9}
Soit $F$ un champ alg\'ebrique.
\begin{itemize}
\item Un espace de modules pour $F$ est un $S$-espace alg\'ebrique
$M$, muni d'un morphisme
$$p : F \longrightarrow M$$
tel que
\begin{enumerate}
\item pour tout $S$-corps s\'eparablement clos $K$, le morphisme induit
$$p_{*} : \pi_{0}F(Spec K) \longrightarrow \pi_{0}M(Spec K)$$
est une bijection
\item le morphisme $p$ est universel (~dans $HoChAlg(S)$~) pour les
morphismes vers les espaces alg\'ebriques.
\end{enumerate}
\item Un quotient g\'eom\'etrique uniforme pour $F$ est un espace de
modules $M$, tel que
\begin{enumerate}
\item la projection $p : F \longrightarrow M$ est un morphisme
submersif (~i.e. un sous-champ $F' \hookrightarrow F$
est ouvert si et seulement si $p(F')$ est ouvert dans $M$~)
\item le morphisme naturel $F \longrightarrow F\times_{M}F$ est
surjectif
\item pour tout morphisme plat d'espaces alg\'ebriques $f : M'
\longrightarrow M$, la projection
$p' : F':=F\times_{M}M' \longrightarrow M'$ est un morphisme
submersif, qui fait de $M'$ un espace de modules pour $F'$, et
v\'erifiant la propri\'et\'e $(2)$ ci-dessus
\end{enumerate}
\end{itemize}
\end{df}

\underline{Remarque:} Si $F=[X/H]$, est un champ quotient d'une
action d'un $S$-sch\'ema en groupes sur un $S$-sch\'ema $X$, alors $M$ est
un quotient g\'eom\'etrique uniforme pour $F$ si et seulement s'il est un
quotient g\'eom\'etrique uniforme de $X$ par $H$ au sens de
\cite[$0.6$]{mu2}. \\

Rappelons les deux principaux r\'esultats d'existence.

\begin{thm}\label{th1.2}{\cite{km}}
Si $F$ est un champ alg\'ebrique tel que
$$\Delta : F \longrightarrow F\times_{S}F$$
soit fini, alors $F$ poss\`ede un quotient g\'eom\'etrique uniforme.

En particulier, tout champ de Deligne-Mumford poss\`ede un quotient
g\'eom\'etrique uniforme.
\end{thm}

\underline{Question:} Le th\'eor\`eme \ref{th1.2} reste-t-il vrai si on
remplace "fini" par "\'equidimensionnel" ?

Dans le cas o\`u $S$ est de caract\'eristique nulle, on poss\`ede le
lemme suivant, r\'epondant partiellement \`a la question pr\'ec\'edente.

\begin{lem}\label{l1.2}
Soit $F$ un champ alg\'ebrique sur $S$, tel que le morphisme
diagonal $\Delta : F \longrightarrow F\times_{S} F$ est
\'equidimensionnel. Alors, si $S$ est de caract\'eristique nulle,
il existe une factorisation unique \`a homotopie pr\`es
$$\xymatrix{F \ar[r]^{f} \ar[d]_{p} & F_{0} \ar[ld]_{q} \\
M &}$$
o\`u $F_{0}$ est un champ de Deligne-Mumford ( \'eventuellement non quasi-s\'epar\'e, i.e. avec 
un morphisme diagonal \'eventuellement non s\'epar\'e ), et $f$ fait de
$F$ une gerbe born\'ee par des espaces alg\'ebriques en groupes lisses sur $F_{0}$.
\end{lem}

\underline{\bf Preuve:} \rm
Pour chaque paire de sections $s,t : X \longrightarrow F$, on dispose
de l'espace alg\'ebrique
des isomorphismes $Isom_{X}(s,t) \longrightarrow X$. C'est un torseur
sous le $X$-espace alg\'ebrique en groupes $Aut_{X}(s)$. Par hypoth\`ese,
le morphisme de projection
$Aut_{X}(s) \longrightarrow X$ est \'equidimensionnelle.
Comme la caract\'eristique de $S$ est nulle, on conclut donc par
\cite[Exp. $IV_{B}$ Cor. $4.4$]{sga3I} que le sch\'ema des composantes
connexes $K=\pi_{0}Aut_{X}(s)
\longrightarrow X$ existe, mais peut \^etre non s\'epar\'e sur $X$.
Le torseur $Isom_{X}(s,t)$ induit donc un
$K$-torseur sur $X$, $\pi_{0}Isom_{X}(s,t) \longrightarrow X$.

Ceci nous permet alors de d\'efinir $F_{0}(X)$ comme \'etant le groupoide
poss\'edant les m\^emes objets que $F(X)$, et avec l'ensemble des sections
de $\pi_{0}Isom_{X}(s,t) \longrightarrow X$ comme morphisme de
$s$ vers $t$. \\

De cette fa\c{c}on, $F_{0}$ est clairement un champ tel que le morphisme
diagonal $F_{0} \longrightarrow F_{0}\times F_{0}$ est quasi-fini.
Comme nous sommes en caract\'eristique nulle, c'est un champ de
Deligne-Mumford. Enfin, le morphisme naturel $F \longrightarrow F_{0}$
est localement sur $s : X \longrightarrow F_{0}$ de la forme
$$BAut_{X}(s) \longrightarrow B\pi_{0}Aut_{X}(s).$$
C'est donc une gerbe born\'e par l'espace alg\'ebrique en groupes lisse
sur $X$ repr\'esentant la composante neutre de $Aut_{X}(s)$. $\Box$\\

\begin{df}\label{ps}
Supposons que $S$ est de caract\'eristique nulle.
Un champ alg\'ebrique $F$, tel que le morphisme diagonal est
\'equidimensionnel
est appel\'e $\Delta$-\'equidimensionnel.

Dans ce cas, si le champ $F^{0}$ du lemme \ref{l1.2} est s\'epar\'e, on
dira que $F$ est pseudo-s\'epar\'e.
\end{df}

Ainsi, on voit que si $S$ est de caract\'erisitique nulle, tout champ
$\Delta$-\'equidimensionnel pseudo-s\'epar\'e poss\`ede un quotient
g\'eom\'etrique uniforme.

Inversemment, si $F$ est normal et poss\`ede un quotient g\'eom\'etrique
uniforme, alors
le morphisme diagonal est forcemment \'equidimensionnel.

\begin{thm}\label{th1.3}{\cite{lm}}
Si $F$ est un champ alg\'ebrique, tel que le morphisme
$$\pi_{F} : I_{F} \longrightarrow F$$
soit plat, alors $F$ poss\`ede un quotient g\'eom\'etrique uniforme
$M$. De plus le morphisme naturel $F \longrightarrow M$, fait de
$F$ une gerbe sur $M$, born\'ee par des $M$-espaces alg\'ebriques en groupes
plats sur $M$.
\end{thm}

Par la suite, nous nous int\'eresserons particuli\`erement aux champs
qui sont localement des quotients par des groupes affines.

\begin{df}\label{d1.10}
Un champ alg\'ebrique $F$ est localement un quotient (~resp.
localement un quotient affine~), s'il existe
un morphisme
$$p : F \longrightarrow X$$
o\`u $X$ est un $S$-espace alg\'ebrique, tel qu'il existe un recouvrement
\'etale $\{U_{i}\}_{i \in I}$ de $X$, des $S$-espaces alg\'ebriques en groupes
$H_{i}$ lisses (~resp. lisses et affines~) sur $S$, op\'erant sur des $S$-espaces
alg\'ebriques $X_{i}$, et des \'equivalences
$$F_{U_{i}}:=F\times_{X}U_{i} \simeq [X_{i}/H_{i}]$$
Si de plus, $F$ poss\`ede un quotient g\'eom\'etrique uniforme, on dira
que $F$ est localement un quotient g\'eom\'etrique uniforme (~resp.
quotient g\'eom\'etrique uniforme affine~).
\end{df}

Pour terminer nous rappelons un fait bien connu,
mais pour lequel nous n'avons pas trouv\'e de r\'ef\'erence sous cette
forme.

\begin{prop}\label{p1.2}
Soit $F$ un champ de Deligne-Mumford, et $p~:~F~\longrightarrow~M$ son
espace de modules.
Alors, il existe un recouvrement \'etale $\{U_{i}\}_{i \in I}$ de $M$,
des groupes finis $H_{i}$, des espaces alg\'ebriques $X_{i}$ et une
action de $H_{i}$ sur $X_{i}$, tel que pour tout $i \in I$,
le champ $F_{U_{i}}:=p^{-1}(U_{i})$ soit \'equivalent au champ
$[X_{i}/H_{i}]$.
\end{prop}

\underline{\bf Preuve:} \rm Voir la premi\`ere partie de la preuve de
\cite[$2.8$]{vi2}. $\Box$\\

\begin{cor}\label{c1.1}
Tout champ $F$ de Deligne-Mumford est localement un quotient
g\'eom\'etrique uniforme affine.

Si $S$ est de caract\'eristique nulle, tout champ alg\'ebrique
$\Delta$-\'equidimensionnel et pseudo-s\'epar\'e est localement un quotient
g\'eom\'etrique uniforme.
\end{cor}

\underline{Question:} Si $F$ est un champ alg\'ebrique $\Delta$-affine
poss\'edant un quotient g\'eom\'etrique uniforme, $F$ est-il localement un quotient
g\'eom\'etrique uniforme ?

\end{subsubsection}

\begin{subsubsection}{Quasi-enveloppes de Chow}

\hspace{5mm}
Pour un champ alg\'ebrique $F$, nous noterons $|F|$ l'ensemble de ses
sous-champs ferm\'es int\`egres. Les \'el\'ements de $|F|$ seront appel\'es les
points de $F$.

Pour chaque point $x \in |F|$, le sous-champ correspondant sera not\'e
$\overline{\{x\}}$. Comme ce champ est int\`egre, il existe un sous-champ
ouvert $U$ de $\overline{\{x\}}$ qui est une gerbe sur un espace alg\'ebrique
int\`egre $M$. Notons \\
$i : Spec K(M) \longrightarrow M$, le point
g\'en\'erique de $M$. Alors \mbox{$i^{*}F:=U\times_{M}Spec K(M)$} est une gerbe
sur $Spec K(M)$.

\begin{df}\label{d1.11}
La gerbe $i^{*}F$ d\'efinie ci-dessus est appel\'ee la gerbe r\'esiduelle de
$F$ au point $x$. Elle sera not\'ee $\widetilde{x}$.

La classe du groupe d'isotropie d'un point $x \in |F|$, est la classe
de conjugaison du groupe alg\'ebrique sur $Spec k(x)^{sp}$ qui borne la gerbe
$\widetilde{x}$. On la notera $|H_{x}|$.

L'ordre de ramification de $F$ en un point $x$ est par d\'efinition
l'ordre de $|H_{x}|$, s'il existe.
\end{df}

Remarquons que pour chaque point $x$ de $F$, on dispose d'un morphisme
repr\'esentable canonique
$$i_{x} : \widetilde{x} \longrightarrow F$$

\begin{df}\label{d1.12}
Un morphisme propre de champs alg\'ebriques
$$f : F \longrightarrow F'$$
est une quasi-enveloppe de Chow, si
pour tout point $x \in |F'|$, le morphisme induit
$$f_{\widetilde{x}} : f^{-1}(\widetilde{x}):=F\times_{F'}\widetilde{x}
\longrightarrow \widetilde{x}$$
admet une section apr\`es un changement de base fini de $k(x)$.
\end{df}

Notons que les quasi-enveloppes de Chow sont stables par changements de
base quelconques, ainsi que par composition.\\

\underline{Exemple:} Soit $F$ un champ alg\'ebrique s\'epar\'e, et
$X \longrightarrow F$ une quasi-enveloppe de Chow, avec $X$ un
espace alg\'ebrique. Alors $F$ est automatiquement un espace
alg\'ebrique. En effet, lorsqu'un morphisme repr\'esentable
$f : F \longrightarrow F'$
est une quasi-enveloppe de Chow, pour tout point $x \in |F'|$, il
existe un point $y \in |F|$ avec $f(y)=x$, tel que le morphisme induit
$$|H_{y}|\otimes Spec \overline{k(y)^{sp}} \longrightarrow |H_{x}|\otimes
Spec \overline{k(y)^{sp}}$$
soit un isomorphisme. Comme $X$ est tel que $|H_{y}|$ est trivial pour
chaque $y \in |X|$, on en d\'eduit que $F$ est un champ alg\'ebrique tel
que
$$\Delta : F \longrightarrow F\times_{S} F$$
est une immersion ferm\'ee. C'est donc un espace alg\'ebrique.\\

\begin{thm}\label{th1.4}
\begin{enumerate}
\item
Soit $F$ un champ alg\'ebrique de Deligne-Mumford r\'eduit. Alors, il
existe un nombre fini de $S$-espaces alg\'ebriques $X_{i}$, des groupes
finis $H_{i}$ (~op\'erant trivialement sur $X_{i}$~), et un morphisme
repr\'esentable fini
$$f : \coprod_{i}[X_{i}/H_{i}] \longrightarrow F$$
qui est une quasi-enveloppe de Chow.
\item
Supposons que $S$ est de caract\'eristique nulle.
Soit $F$ un champ alg\'ebrique normal,  $\Delta$-\'equidimensionnel
et pseudo-s\'epar\'e (\ref{ps}). Alors, il existe un
nombre finis d'espaces alg\'ebriques $X_{i}$, et un morphisme repr\'esentable fini
$$f : \coprod_{i}F_{i}  \longrightarrow F$$
qui est une quasi-enveloppe de Chow, avec $F_{i}$ une gerbe sur
$X_{i}$ born\'ee par des espaces alg\'ebriques en groupes lisses sur $X_{i}$.
\end{enumerate}
\end{thm}

\underline{\bf Preuve:} \rm
$(1)$ Si on construit une quasi-enveloppe de Chow v\'erifiant les
conclusions du th\'eor\`eme pour chacune
des composantes irr\'eductibles de $F$, leur union disjointe satisfera
aux conditions demand\'ees. On peut donc se restreindre au cas o\`u
$F$ est int\`egre.
Soit $M$ l'espace de modules de $F$. D'apr\`es \cite[$2.6$]{vi2}, il existe un
$S$-espace alg\'ebrique normal $X$ et un morphisme repr\'esentable fini
$$f : X \longrightarrow F$$
Notons $F_{0}$ la normalisation de
$F\times_{M}X$, et $X_{0}$ l'espace de modules de $F_{0}$. Alors le
morphisme naturel $X_{0} \longrightarrow X$ est un morphisme fini
et birationnel entre deux espaces alg\'ebriques normaux, c'est donc un
isomorphisme. D\'emontrons alors que $F_{0}$ est une gerbe triviale sur son
espace de modules $X$.

\begin{lem}\label{l1.1}
Soit $F$ un champ de Deligne-Mumford normal tel que la projection
naturelle $p : F \longrightarrow M$ sur son espace de modules admette
une section. Alors $F$ est \'equivalent \`a une gerbe triviale sur $M$.
\end{lem}

\underline{\bf Preuve:} \rm Il suffit de montrer que la projection
$F \longrightarrow M$ fait de $F$ une gerbe sur $M$. Comme ceci
est local sur $M_{et}$, on peut supposer par \ref{p1.2}, que $F=[X/H]$
est un champ quotient d'un groupe fini $H$ op\'erant sur un $S$-sch\'ema
normal et irr\'eductible $X$.

Dans ce cas $M=X/H$, et la section $M \longrightarrow F$ est d\'efinie
par un diagramme commutatif
$$\xymatrix{
X/H \ar[r]^{Id} & X/H \\
Y \ar[u]^{q} \ar[r]_{f} & X \ar[u]_{p} }$$
o\`u $p$ est la projection canonique, $q$ un $H$-fibr\'e principal, et
$f$ un morphisme $H$-\'equivariant. Comme $q$ est \'etale, $f$ est
non-ramifi\'e. Or $Y$ et $X$ sont normaux et de m\^eme dimension, donc
$f$ est \'etale. Ce qui implique que $p$ est aussi \'etale. Ainsi, si
$H_{0}=Ker(H \longrightarrow Aut(X))$, l'action de $H/H_{0}$ est libre ,
et $F$ est \'equivalent \`a $[(X/H)/H_{0}]$. $\Box$\\

Ainsi, $F_{0}\simeq [X_{0}/H_{0}]$ pour un groupe fini op\'erant
trivialement sur $X_{0}$. De plus, le morphisme
$$f_{0} : F_{0} \longrightarrow F$$
est repr\'esentable fini, et est g\'en\'eriquement une quasi-enveloppe de
Chow. Il existe donc un sous-champ ouvert dense $U \hookrightarrow F$
tel que $F_{0}\times_{F}U \longrightarrow U$ soit une quasi-enveloppe
de Chow. Notons $F'$ le ferm\'e compl\'ementaire r\'eduit de $U$ dans $F$.
Par r\'ecurrence noeth\'erienne, la proposition est vraie pour $F'$. Soit
$X_{i}'$, $H'_{i}$ et
$$f' : \coprod_{i} [X'_{i}/H'_{i}] \longrightarrow F'$$
une quasi-enveloppe de Chow pour $F'$. Alors
$$f=f_{0}\coprod f' : [X_{0}/H_{0}] \coprod_{i}[X'_{i}/H'_{i}]
\longrightarrow F$$
est une quasi-enveloppe de Chow pour $F$ qui v\'erifie les conditions
demand\'ees.\\

$(2)$ Le r\'esultat se d\'eduit imm\'ediatemment du cas $(1)$ et du lemme
\ref{l1.2}. $\Box$\\

\end{subsubsection}

\end{subsection}

\end{section}

\newpage

\begin{section}{Chapitre $2$ : $K$-th\'eorie des champs alg\'ebriques}
\hspace{5mm}
Dans ce chapitre nous allons \'etudier les spectres de $K$-th\'eorie des
champs alg\'ebriques. Il s'agit d'essayer de d\'ecrire ces spectres
en fonctions de "choses connues", \`a savoir les spectres de
$K$-th\'eorie des sch\'emas ou des espaces alg\'ebriques, ou encore la
$K$-cohomologie des champs alg\'ebriques.

Les premiers r\'esultats dans cette direction sont les th\'eor\`emes de
descente (~\ref{th2.1}, \ref{th2.2}~). Cependant, ils ne sont pas r\'eellement
utilisables pour d\'ecrire la $K$-th\'eorie, mais sont plut\^ot
des outils techniques permettant de ramener les calculs \`a des cas
connus (~quotients par des groupes finis par exemple~). Nous en feront
un usage intensif dans le chapitre suivant, lors de la preuve des
formules de Riemann-Roch.

Bien que poss\'edant de nombreuses propri\'et\'es analogues \`a celles
de la $K$-th\'eorie des sch\'emas (~localisation, homotopie, axiome du
fibr\'e projectif ...~), la $K$-th\'eorie des champs alg\'ebriques diff\`ere
de celle-ci par le fait qu'elle ne poss\`ede plus la propri\'et\'e
de descente \'etale (~\cite[$11.10$]{th}~). Ceci provient de la nature mixte
des faisceaux coh\'erents sur les champs alg\'ebriques, dans le sens o\`u ils
font intervenir d'une part des faisceaux coh\'erents sur des sch\'emas, et
d'autre part des repr\'esentations de groupes alg\'ebriques. Dans le cas
des champs de Deligne-Mumford par exemple, la $K$-cohomologie rationnelle
ne peut
pas retenir l'information sur ces repr\'esentations, car la cohomologie
d'un groupe fini est de torsion. Bien que je n'aie pas v\'erifi\'e tous
les d\'etails, il est m\^eme probable que la condition de descente \'etale
pour la $K$-th\'eorie caract\'erise les espaces alg\'ebriques parmi les
champs de Deligne-Mumford. C'est
alors le but des th\'eor\`emes de d\'evissage (~\ref{th2.3}, \ref{th2.4},
\ref{th2.5}~) de d\'ecrire la partie de la
$K$-th\'eorie qui dispara\^it dans la $K$-cohomologie. Ces th\'eor\`emes sont
d'un certain point de vue orthogonaux aux th\'eor\`emes de descente. En
effet, il est difficile de les utiliser dans les calculs, mais en
contre-partie, leur caract\`ere descriptif permet de d\'efinir,
de mani\`ere assez \'evidente, le caract\`ere de Chern qui sera utilis\'e dans
les formules de Riemann-Roch. \\

Le site $(Esp/S)_{li}$ est muni d'un faisceau d'anneaux coh\'erent
$$\cal O\mit : X \mapsto \cal O\mit_{X}(X)$$
On lui associe le champ en cat\'egories $\bf Vect\mit
\longrightarrow (Esp/S)_{li}$
(~resp. \\
$\bf Coh\mit~\longrightarrow~(Esp/S)_{li}$~), dont la cat\'egorie
des sections au-dessus d'un espace alg\'ebrique $X$ est la cat\'egorie
des faisceaux de $\cal O\mit$-modules
localement libres et de rang fini (~resp. localement de pr\'esentation
finie~) sur le site restreint $(Esp/X)_{li}$.
Remarquons que $\bf Vect$ est un champ en cat\'egories exactes. Il n'en
n'est plus de m\^eme de $\bf Coh$.

Si $F$ est un champ alg\'ebrique, on peut d\'efinir son petit site lisse
$F_{li}$. Ses objets sont les $1$-morphismes lisses $s : X \longrightarrow
F$, avec $X$ un espace alg\'ebrique. Un morphisme entre $s :
X \longrightarrow F$ et $t : Y \longrightarrow F$, est la donn\'ee d'un morphisme
$f : X \longrightarrow Y$, et d'un $2$-morphisme $h$ entre
$s$ et $t\circ f$.

Ce site est muni du faisceau d'anneaux
$$\cal O\mit_{F} : (~X \longrightarrow F~) \mapsto \cal O\mit_{X}(X)$$
On lui associe le champ en cat\'egories exactes $\bf Coh\mit_{F}
\longrightarrow F_{li}$, dont la cat\'egorie des sections au-dessus de
l'objet $X \longrightarrow F$ est la cat\'egorie des faisceaux
de $\cal O\mit_{F}$-modules coh\'erents sur $X$. C'est la cat\'egorie cofibr\'ee
$\bf Coh$ restreinte \`a $F_{li}$.\\

Pour un champ alg\'ebrique $F$, on pose alors
$$\bf Vect\mit(F):=\int_{F}\bf Vect$$
$$\bf Coh\mit(F):=\int_{F}\bf Coh\mit$$

D'apr\`es la d\'efinition des sections cart\'esiennes globales d'un
champ (~\cite[$1.1.1$]{gi}~), la cat\'egorie $\bf Vect\mit (F)$ (~resp. $\bf Coh\mit(F)$~) 
peut-\^etre d\'efinie de la fa\c{c}on suivante \\

\underline{un objet :} est d\'efini par la donn\'ee suivante~:
\begin{itemize}
\item
Pour toute section
$s : X \longrightarrow F$, avec $X$ un espace alg\'ebrique,
la donn\'ee d'un fibr\'e vectoriel (~resp. faisceau coh\'erent~)
$V_{(s)}$ sur $X$
\item
pour tout morphisme d'espaces alg\'ebriques $f : Y \longrightarrow X$, et
toute paire de sections
$$s : X \longrightarrow F$$
$$t : Y \longrightarrow F$$
et tout $2$-morphisme $h : s\circ f \Rightarrow t$, un isomorphisme
de fibr\'es vectoriels (~resp. faisceaux coh\'erents~)
$$\phi_{s,f,h} : f^{*}V_{(s)} \simeq V_{(t)}$$
\item
Pour toute paire de morphismes d'espaces alg\'ebriques
$$\xymatrix{Z \ar[r]^{g} & Y \ar[r]^{f} & X}$$
tout triplet de $1$-morphismes $s : X \longrightarrow F$, $t : Y
\longrightarrow F$, $u~:~Z~\longrightarrow~F$, et toute paire de
$2$-morphismes
$$h : s \circ f \Rightarrow t$$
$$j : t\circ g \Rightarrow u$$
une \'egalit\'e
$$\phi_{t,g,j}\circ g^{*}\phi_{s,f,h} = \phi_{u,f\circ g,h\circ j}$$
\end{itemize}

\underline{un morphisme :} entre $V$ et $W$ est d\'efini par la donn\'ee
suivante~:
\begin{itemize}
\item
Pour toute section
$s : X \longrightarrow F$, avec $X$ un espace alg\'ebrique,
un morphisme de fibr\'es vectoriels (~resp. faisceaux coh\'erents~) sur $X$
$$a_{s} : V_{(s)} \longrightarrow W_{(s)}$$
\item
Pour tout morphisme d'espaces alg\'ebrique $f : Y \longrightarrow X$, et
toute paire de sections
$$s : X \longrightarrow F$$
$$t : Y \longrightarrow F$$
et tout $2$-morphisme $h : s\circ f \Rightarrow t$, une \'egalit\'e
$$\phi_{s,f,h}^{W} \circ f^{*}(a_{s})=a_{t} \circ \phi_{s,f,h}^{V}$$
\end{itemize}

Plus g\'en\'eralement, pour $\cal E$ une cat\'egorie cofibr\'ee sur
$(Esp/S)_{li}$, $F$ un champ, et $V$ une section globale cart\'esienne de
$\cal E$ sur $F$,  nous noterons $V_{(s)}$ la section de $\cal E$
sur l'espace alg\'ebrique $X$ d\'efinie par un $1$-morphisme
$s : X \longrightarrow F$. \\

Supposons que $f : F \longrightarrow F'$ soit un morphisme propre de
champs alg\'ebriques. Alors, on peut d\'efinir une pseudo-transformation
naturelle entre les
pseudo-foncteurs sur $F'_{li}$
$$U \mapsto \bf Coh\mit(f^{-1}(U))$$
$$U \mapsto \bf Coh\mit(U)$$
qui \`a un faisceau coh\'erent $\cal F$ sur $f^{-1}(U)$, associe le
faisceau $f_{*}(\cal F\mit)$ sur $U$. Le fait que ceci d\'efinit bien
une pseudo-transformation naturelle est une cons\'equence de la formule
de transfert pour les morphismes lisses. Ainsi, $f_{*}$ d\'efinit un
morphisme de champs sur $F'_{li}$
$$f_{*} : f_{*}\bf Coh\mit \longrightarrow \bf Coh\mit$$
Par la m\^eme construction, si $C\bf Mod\mit_{qcoh}$ d\'esigne la
cat\'egorie cofibr\'ee
des complexes de $\cal O$-modules, \`a cohomologie quasi-coh\'erente et
born\'ee, on d\'efinit une image directe
$$f_{*} : f_{*}C\bf Mod\mit_{qcoh} \longrightarrow C\bf Mod\mit_{qcoh}$$

\begin{subsection}{Premi\`eres propri\'et\'es}

\begin{df}\label{d2.1}
Le pr\'efaisceau en spectres de $K$-th\'eorie \`a coefficients dans
$\bf Vect$ est not\'e $\underline{K}$. La $K$-cohomologie d'un champ
alg\'ebrique $F$ est
$$\underline{\bf K}(F):=\bf H\mit(F,\underline{K}_{\mathbf{Q}})$$

Le foncteur de $K$-th\'eorie \`a coefficients dans $\bf Vect$ est not\'e
$$\begin{array}{cccc}
\bf K\mit : & Ch(S) & \longrightarrow & Sp \\
            &        F          & \mapsto & \bf K\mit(F):=K(\bf
            Vect\mit(F))
\end{array}$$

Si $F$ est un champ alg\'ebrique, son pr\'efaisceau en spectres de $K$-th\'eorie
\`a coefficients dans $\bf Coh\mit_{F}$ sera not\'e $\underline{G}$.
On d\'efinit la $G$-cohomologie de $F$ par
$$\underline{\bf G}(F):=\bf H\mit(F_{li},\underline{G}_{\mathbf{Q}})$$

Le spectre de $G$-th\'eorie d'un champ alg\'ebrique $F$ est
d\'efini par
$$\bf G\mit(F):=K(\bf Coh\mit(F))$$
Nous noterons
$$can_{F} : \bf K\mit(F) \longrightarrow \underline{\bf K}(F)$$
$$can_{F} : \bf G\mit(F) \longrightarrow \underline{\bf G}(F)$$
les morphismes canoniques (~\ref{p1.1}~).
\end{df}

\underline{Remarque:} Pour le site lisse $F_{li}$ les morphismes de transitions ne sont pas forcement
plats. Ainsi, $\underline{G}$ n'est pas r\'eelement d\'efini comme pr\'efaisceau en spectres sur $F_{li}$. 
Pour r\'esoudre cette difficult\'e il suffit de travailler avec les d\'efinitions de \cite{th}. \\

Notons que si $f : F \longrightarrow F'$ est un
$1$-morphisme plat entre deux champs alg\'ebriques, le foncteur
$$f^{*} : f^{*}\bf Coh\mit_{F'} \longrightarrow (\bf Coh\mit_{F})$$
est exact, et induit donc des morphismes
$$f^{*} : \bf G\mit(F) \longrightarrow \bf G\mit(F')$$
$$f^{*} : \underline{\bf G}(F) \longrightarrow
\underline{\bf G}(F')$$
De cette fa\c{c}on, $F \mapsto \bf G\mit(F)$ et
$F \mapsto \underline{\bf G}(F)$
sont des foncteurs stricts de la $2$-cat\'egorie $(ChAlg(S),fl)$
des champs alg\'ebriques et $1$-morphismes plats, vers celle des
spectres, morphismes de spectres et classe d'homotopie d'homotopie
entre morphismes.

Si $f : F \longrightarrow F'$ est un morphisme
propre et de dimension cohomologique finie. Alors, on peut d\'efinir un
morphisme dans $HoSp$
$$f_{*} : \bf G\mit(F) \longrightarrow \bf G\mit(F')$$
De cette fa\c{c}on, $F \mapsto \bf G\mit(F)$ est un foncteur covariant
de la cat\'egorie $(HoChAlg(S),pr<\infty)$ des champs alg\'ebriques et morphismes
propres de dimension cohomologique finie, vers $HoSp$.

Nous aurons besoin \`a certains moments d'avoir une version
fonctorielle de ces images directes. Pour cela, nous utiliserons les
constructions de Thomason (~\cite{th}~).

Par exemple supposons que l'on dispose d'un $I$-pseudo-diagramme de champs
alg\'ebriques
$$F : I \longrightarrow ChAlg(S)$$
tel que pour chaque $u : i \longrightarrow j$, le morphisme $F(u) :
F(i) \longrightarrow F(j)$ soit propre et de dimension cohomologique
finie.
Notons, pour chaque $i \in I$, $A(i)$ la cat\'egorie bi-compliciale de
Waldhausen des complexes de $\cal O\mit_{SF(i)}$-modules acycliques, et
\`a cohomologie coh\'erente et born\'ee. Alors, la correspondance
$$i \mapsto A(i)$$
$$(u : i \rightarrow j) \mapsto F(u)_{*} : A(i) \longrightarrow
A(j)$$
d\'efinit un $I$-pseudo-diagramme de cat\'egories bi-compliciales de
Waldhausen. Et par le proc\'ed\'e de strictification \ref{strict}, un
$I$-diagramme de
cat\'egories bi-compliciales de Waldhausen
$$SA : i \mapsto SA(i)$$
En prenant l'image par le foncteur covariant $\bf G$, on obtient donc
un \\
$I$-diagramme dans $Sp$
$$\bf G\mit : i \mapsto K(SA(i))$$
Ainsi, pour chaque diagramme de champs alg\'ebriques, avec des morphismes
de transition propres et de dimension cohomologique finie, on disposera
d'un diagramme de spectres de $G$-th\'eorie
associ\'e. Le cas que nous utiliserons le plus est celui o\`u
$I=\Delta^{op}$. \\

Supposons maintenant que $f : F \longrightarrow F'$ est propre et
repr\'esentable. On dispose alors d'un morphisme de champs
sur $F'_{li}$
$$f_{*} : f_{*}\bf Coh\mit \longrightarrow \bf Coh\mit$$
Si on note $C\bf Mod\mit_{coh}$ (~resp. $C\bf Mod\mit^{ac}_{coh}$~)
la cat\'egorie cofibr\'ee en cat\'egories
bi-compliciales de Waldhausen des complexes de $\cal O$-modules
(~resp. des $\cal O$-modules acycliques et \`a cohomologie coh\'erente et
born\'ee~),
\`a cohomologie born\'ee et coh\'erente, ce morphisme induit un morphisme
exact de cat\'egories cofibr\'ees en cat\'egories bi-compliciales de Waldhausen
$$f_{*} : f_{*}C\bf Mod\mit^{ac}_{coh} \longrightarrow C\bf
Mod\mit_{coh}$$
Ainsi, en passant aux spectres de $K$-th\'eorie,
on a construit un morphisme de pr\'efaisceaux en spectres sur
$F'_{li}$
$$f_{*} : f_{*}\underline{G} \longrightarrow \underline{G}$$

Les produits tensoriels
$$\otimes : \bf Vect\mit\times_{(Esp/S)_{li}}\bf Vect\mit
\longrightarrow \bf Vect$$
$$\otimes : \bf Vect\mit_{F}\times_{F_{li}}\bf Vect\mit_{F} \longrightarrow
\bf Vect\mit_{F}$$
$$\otimes : \bf Vect\mit_{F}\times_{F_{li}}\bf Coh\mit_{F} \longrightarrow
\bf Coh\mit_{F}$$
induisent des produits dans $HoSp$ (~\ref{s1}~)
$$\otimes : \bf K\mit \wedge \bf K\mit \longrightarrow \bf K$$
$$\otimes : \bf K\mit(F) \wedge \bf G\mit(F)\longrightarrow \bf
G\mit(F)$$
$$\otimes : \underline{\bf K} \wedge \underline{\bf K} \longrightarrow
\underline{\bf K}$$
$$\otimes : \underline{\bf K}(F) \wedge \underline{\bf G}(F)\longrightarrow
\underline{\bf G}(F)$$

Rappelons les principales propri\'et\'es de ce foncteurs.

\begin{prop}\label{p2.1}
\begin{enumerate}
\item (~axiome du fibr\'e projectif~)
Soit $\pi : \bf P\mit(V) \longrightarrow F$ un fibr\'e projectif associ\'e \`a un
fibr\'e vectoriel $V$ de rang $r+1$ sur un champ alg\'ebrique $F$, et
$x=\cal O\mit_{P}(1)$ le fibr\'e inversible canonique. Alors il existe
des isomorphismes dans $HoSp$
$$\begin{array}{ccc}
\bigvee_{i=0}^{i=r}\bf K\mit(F) & \longrightarrow & \bf K\mit(\bf
P\mit(V)) \\
\vee_{i}a_{i} & \mapsto & \sum_{i}x^{i}\otimes\pi^{*}(a_{i})
\end{array}$$
et de m\^eme avec $\underline{\bf K}$, $\bf G$ et $\underline{\bf G}$.
\item (~localisation~)
Si $j : F' \hookrightarrow F$ est une immersion ferm\'ee de champs
alg\'ebriques, et $i : U \hookrightarrow F$ l'immersion ouverte
compl\'ementaire. Alors il existe des triangles fonctoriels pour les
images r\'eciproques de morphismes plats
$$\xymatrix{
    & \bf G\mit(F') \ar[rd]^{j_{*}} &  &   & & \underline{\bf G}(F')
  \ar[rd]^{j_{*}} & \\
\bf G\mit(U) \ar[ru]^{-1} & & \ar[ll]_{i^{*}} \bf G\mit(F) & &
\underline{\bf G}(U) \ar[ru]^{-1} & & \ar[ll]_{i^{*}} \underline{\bf
G}(F) }$$
\item (~homotopie~)
Si $p : V \longrightarrow F$ est un morphisme de champs alg\'ebriques qui
est un torseur affine sur $F_{li}$, alors
les morphismes naturels dans $HoSp$
$$p^{*} : \bf G\mit(F) \longrightarrow \bf G\mit(V)$$
$$p^{*} : \underline{\bf G}(F) \longrightarrow \underline{\bf G}(V)$$
sont des isomorphismes.
\item (~dualit\'e de Poincar\'e~)
Si $F$ est un champ alg\'ebrique r\'egulier, le morphisme naturel
dans $HoSp(F_{li})$
$$\underline{\bf K} \longrightarrow \underline{\bf G}$$
est un isomorphisme.
\item (~descente~)
Pour tout champ alg\'ebrique $F$, le morphisme canonique
$$\underline{\bf K}(F) \longrightarrow \bf
H\mit(F_{li},\underline{K}_{\mathbf{Q}})$$
est un isomorphisme d'anneaux dans $HoSp$, compatible avec les images
r\'eciproques.

Si $X$ est un espace alg\'ebrique, alors le morphisme canonique
$$\bf G\mit_{*}(X)\otimes \bf Q\mit \longrightarrow
\underline{\bf G}_{*}(X)$$
est un isomorphisme.
\item (~invariance topologique~)
Si $j : F_{red} \hookrightarrow F$ est l'immersion canonique du
sous-champ alg\'ebrique r\'eduit d'un champ alg\'ebrique, les morphismes
$$j_{*} : \bf G\mit(F_{red}) \longrightarrow \bf G\mit(F)$$
$$j_{*} : \underline{\bf G}(F_{red}) \longrightarrow \underline{\bf
G}(F)$$
sont des isomorphismes de $HoSp$.
\item (~continuit\'e~)
Soit $\{F_{i}\}_{i \in I}$ un syst\`eme inductif filtrant de champs
alg\'ebriques avec $1$-morphismes de transition plats. Alors, si \\
$F=Colim_{I}F_{i}$ est une $1$-limite inductive dans $ChAlg(S)$,
le morphisme naturel
$$\bf G\mit(F) \longrightarrow lim_{I}\bf G\mit(F_{i})$$
est un isomorphisme dans $HoSp$.
\item (~transfert~)
Si
$$\xymatrix{F'_{0} \ar[d]_{f'} \ar[r]^{u_{0}} & F_{0} \ar[d]^{f}\\
F' \ar[r]_{u} & F }$$
est un carr\'e $2$-cart\'esien dans $ChAlg(S)$, avec $f$ propre et
repr\'esentable, et $u$ plat et repr\'esentable, alors les diagrammes
suivants commutent \`a homotopie naturelle pr\`es dans $Sp$
$$\xymatrix{
\bf G\mit(F'_{0}) \ar[d]_{f_{*}'} & \ar[l]_{u_{0}^{*}}  \bf
G\mit(F_{0}) \ar[d]^{f_{*}} & &
\underline{\bf G}(F'_{0}) \ar[d]_{f_{*}'} & \ar[l]_{u_{0}^{*}}
\underline{\bf G}(F_{0}) \ar[d]^{f_{*}}\\
\bf G\mit(F') &  \ar[l]^{u^{*}} \bf G\mit(F) & &
\underline{\bf G}(F') &  \ar[l]^{u^{*}}
\underline{\bf G}(F) }$$
\item (~formule de projection~)
Soit $f : F' \longrightarrow F$ un morphisme propre
et repr\'esentable de champs alg\'ebriques. Alors les diagrammes suivants
commutent \`a homotopie naturelle pr\`es dans $Sp$
$$\xymatrix{
\bf K\mit(F)\wedge \bf G\mit(F') \ar[r]^-{Id\otimes
f_{*}} \ar[d]_-{f^{*}\otimes Id} & \bf G\mit(F) & &
\underline{\bf K}(F)\wedge \underline{\bf G}(F') \ar[r]^-{Id\otimes
f_{*}} \ar[d]_-{f^{*}\otimes Id} & \underline{\bf G}(F) \\
\bf K\mit(F')\wedge \bf G\mit(F')
\ar[r]^-{\otimes} & \bf G\mit(F') \ar[u]_-{f_{*}} & &
\underline{\bf K}(F')\wedge \underline{\bf G}(F')
\ar[r]^-{\otimes} & \underline{\bf G}(F') \ar[u]_-{f_{*}} }$$
\end{enumerate}
\end{prop}

\underline{\bf Preuve:} \rm Les points $(1)$, $(2)$, $(3)$, $(6)$,
$(7)$, $(8)$ et $(9)$ se d\'emontrent exactement comme dans le cas d'un sch\'ema
(~\cite{q}~). Le point $(4)$ se d\'emontre comme dans \cite{th2}.

Le point $(5)$ provient directement du th\'eor\`eme de
descente de la \\
$G$-th\'eorie \'etale (~\cite[$11.10$]{th}~), et du fait
qu'un pr\'efaisceau en spectres est flasque pour la topologie \'etale si
et seulement il l'est pour la topologie lisse (~car tout morphisme
lisse poss\`ede une section apr\`es un recouvrement \'etale~). $\Box$\\

\underline{Remarques:}
\begin{itemize}
\item
Le probl\`eme de savoir si pour un champ r\'egulier $F$, le morphisme
naturel $\bf K\mit(F) \longrightarrow \bf G\mit(F)$ est
un isomorphisme semble difficile,
m\^eme \`a coefficients rationnels. Les seuls cas o\`u nous connaissons une
r\'eponse partielle est celui des champs quotients par des groupes
affines (~\cite{th2}~).
\item
La propri\'et\'e de descente implique que l'on a un isomorphisme naturel
dans $HoSp(F_{li})$
$$\underline{G}_{\bf Q} \longrightarrow \underline{\bf G}$$
Ainsi, si $f : F \longrightarrow F'$ est un morphisme propre et
repr\'esentable, on a $\bf R\mit f_{*}\underline{G}_{\bf Q}
\simeq f_{*}\underline{G}_{\bf Q}$ dans $HoSp(F'_{li})$, et donc
on peut construire une image directe
$$f_{*} : \underline{\bf G}(F)\simeq \bf
H\mit(F'_{li},\bf R\mit f_{*}\underline{G})\simeq
\bf H\mit(F'_{li},\bf R\mit f_{*}\underline{G}_{\bf Q})
\stackrel{\bf H\mit(f_{*})}{\longrightarrow}
\bf H\mit(F'_{li},\underline{G}_{\bf Q})\simeq \underline{\bf G}(F')$$
dans $HoSp$. En utilisant les images directes fonctorielles pour
$\bf G$, et cette identification, on peut montrer que
$$F \mapsto \underline{\bf G}$$
est un foncteur covariant de la $2$-cat\'egorie des champs alg\'ebriques
repr\'esentables sur un champ de base fixe $F'$, et
$1$-morphismes propres et repr\'esentables, vers celle des spectres,
morphismes de spectres et classes d'homotopie d'homotopie entre morphismes.
\end{itemize}

\begin{cor}\label{c2.1}
Soit $f : F \longrightarrow F'$ un $1$-morphisme propre et
repr\'esentable de champs alg\'ebriques. Alors le diagramme suivant commute
dans $HoSp$
$$\xymatrix{
\bf G\mit(F) \ar[r]^{f_{*}} \ar[d]_-{can}& \bf G\mit(F') \ar[d]^-{can} \\
\underline{\bf G}(F) \ar[r]_{f_{*}} & \underline{\bf
G}(F')}$$
\end{cor}

\underline{\bf Preuve:} \rm C'est imm\'ediat d'apr\`es la fonctorialit\'e
du morphisme $can$. $\Box$\\

\end{subsection}

\begin{subsection}{Descente de la $G$-th\'eorie rationnelle}
\hspace{5mm}

\begin{subsubsection}{Descente au-dessus d'un espace alg\'ebrique}

\begin{thm}\label{th2.1}
Soit $F$ un champ alg\'ebrique, et $p : F  \longrightarrow X$
un \\
$1$-morphisme, avec $X$ un espace alg\'ebrique. Alors
$$\begin{array}{cccc}
p_{*}\bf G\mit_{\bf Q} : & X_{li} & \longrightarrow & Sp \\
                                  & U      & \mapsto & \bf
                                  G\mit(p^{-1}U)_{\bf Q}
\end{array}$$
est un pr\'efaisceau flasque sur $X_{li}$.
\end{thm}

\underline{\bf Preuve:} Comme un morphisme lisse poss\`ede des sections
apr\`es un changement de base \'etale et surjectif, un pr\'efaisceau est flasque
sur $X_{li}$ si et seulement si sa restriction l'est sur $X_{et}$. On
travaillera donc avec la topologie \'etale sur $X$.\\

Remarquons aussi, que s'il existe un triangle
$$\xymatrix{ & E \ar[rd] &  \\
F \ar[ru]^{-1} \ar[rr] & & G}$$
dans $Sp(X_{et})$, $G$ est flasque si et seulement si $E$ et $F$ le
sont. Ainsi, en utilisant la localisation (~\ref{p2.1}~), la continuit\'e
(~\ref{p2.1}~) et une r\'ecurrence noeth\'erienne on peut supposer que
$X=Spec K$, est le spectre d'un
corps.

\begin{lem}\label{l2.1}
Soit $K \hookrightarrow L$ une extension galoisienne de corps, et $H$
son groupe de galois. Soit $F$ un champ alg\'ebrique sur $K$, et
\mbox{$F_{L}:=F\times_{Spec K}Spec L$}, muni de l'action de $H$ induite.
Alors le morphisme naturel
$$\bf G\mit(F)_{\bf Q} \longrightarrow
holim_{H}\bf G\mit(F_{L})_{\bf Q}$$
est un isomorphisme dans $HoSp$.
\end{lem}

\underline{\bf Preuve:} \rm Comme $H$ est un groupe fini, et que les
spectres sont \`a coefficients rationnels, le lemme est \'equivalent au
fait que le morphisme naturel
$$q^{*} : \bf G\mit_{*}(F)_{\bf Q} \longrightarrow \bf G\mit_{*}(F_{L})_{\bf
Q}^{H}$$
est un isomorphisme. Mais d'apr\`es la formule de projection
(~\ref{p2.1} $(9)$~) pour la 
projection $q : F_{L} \longrightarrow F$, un inverse de $q^{*}$ est
$\frac{1}{m}q_{*}$, o\`u $m$ est l'ordre de $H$. $\Box$\\

D'apr\`es le lemme, $p_{*}\bf G\mit_{\bf Q}$ poss\`ede la propri\'et\'e de
descente pour tout recouvrement dans $(Spec K)_{et}$ qui est de la
forme $Spec L \longrightarrow Spec K$, pour une extension galoisienne
$L/K$. Mais comme tout morphisme couvrant de $(Spec K)_{et}$ poss\`ede
une section apr\`es un changement de base de cette forme,
le pr\'efaisceau en spectres $p_{*}\bf G\mit_{\bf Q}$ est flasque sur
$(Spec K)_{et}$. $\Box$\\

\begin{cor}\label{c2.2}
\begin{enumerate}
\item
Soit $F$ un champ alg\'ebrique, $F \longrightarrow X$ un $1$-morphisme
vers un espace alg\'ebrique, et $Y \longrightarrow X$ un fibr\'e principal
homog\`ene sous un groupe fini $H$. Notons $p : F_{Y}:=F\times_{X}Y
\longrightarrow F$ son image r\'eciproque sur $F$. Alors, il existe un
diagramme commutatif d'isomorphismes dans $HoSp$
$$\xymatrix{
hocolim_{H}\bf G\mit(F_{Y})_{\bf Q} \ar[r]^{can} \ar[d]_{p_{*}} &
holim_{H}\bf G\mit(F_{Y})_{\bf Q} \ar[ld]^{p_{*}}\\
\bf G\mit(F)_{\bf Q} & }$$
\item
Si $F$ est un champ de Deligne-Mumford r\'egulier, et
$p : F \longrightarrow M$ la projection sur son espace de modules,
alors il existe un isomorphisme naturel dans $HoSp$
$$\bf H\mit(M_{et},p_{*}\bf K\mit\otimes \bf Q\mit)
\simeq \bf G\mit(F)_{\bf Q}$$
\end{enumerate}
\end{cor}

\underline{\bf Preuve:} \rm $(1)$
Comme $H$ est fini, on sait que pour toute
action de $H$ sur un $\bf Q$-espace vectoriel $V$, le morphisme naturel
$V_{H} \longrightarrow V^{H}$, des co-invariants vers les invariants,
est un isomorphisme. Ceci implique que le morphisme $can$ du
corollaire est un isomorphisme.

Il suffit ensuite de remarquer que $p^{*}\circ p_{*}=\times m$, avec
$m$ l'ordre de $H$. Or une application du th\'eor\`eme \ref{th2.1} au morphisme
$F \longrightarrow X$, et au recouvrement \'etale $Y \longrightarrow
X$, montre que
$$p^{*} : \bf G\mit(F)_{\bf Q} \longrightarrow \bf H\mit(Y/X,f_{*}\bf
G\mit_{\bf Q})$$
est un isomorphisme.
Mais le membre de droite est canoniquement isomorphe dans $HoSp$ \`a
$holim_{H}\bf G\mit(F_{Y})_{\bf Q}$. \\

$(2)$ Consid\'erons les morphismes naturels de $HoSp$
$$\xymatrix{
\bf G\mit(F)_{\bf Q} \ar[r]^-{a} & \bf H\mit(M_{et},p_{*}\bf
G\mit_{\bf Q}) & \ar[l]_-{b} \bf H\mit(M_{et},p_{*}\bf
K\mit_{\bf Q})}$$
Le th\'eor\`eme \ref{th2.1} implique que $a$ est un isomorphisme.
Localement sur $M_{et}$, on a $F\simeq [X/H]$, avec $H$ un groupe fini
op\'erant sur un sch\'ema r\'egulier $X$ (~\ref{p1.2}~). Mais dans ce cas, on
sait que le morphisme naturel
$$\bf K\mit([X/H]) \longrightarrow \bf G\mit([X/H])$$
est un isomorphisme dans $HoSp$ (~\cite[$5.3$]{th2}~). Ce qui montre que le
morphisme
$$p_{*}\bf K\mit_{\bf Q} \longrightarrow p_{*}\bf G\mit_{\bf Q}$$
est une \'equivalence faible sur $M_{et}$, et donc que $b$ est un
isomorphisme dans $HoSp$. $\Box$\\

\underline{Remarques:}
\begin{itemize}
\item
Il est en g\'en\'eral faux que $\bf G\mit_{\bf Q}$ soit flasque sur $F_{li}$. Un
exemple simple est celui de $F=[Spec K/H]$, avec $H$ un groupe fini
(~\ref{s1}~).
\item
Un cas d'application important du pr\'ec\'edent th\'eor\`eme est celui o\`u $F$
est localement un quotient affine g\'eom\'etrique uniforme
(~par exemple un champ de Deligne-Mumford~), et $p : F
\longrightarrow M$ la projection naturelle sur l'espace de modules. En
effet dans ce cas le th\'eor\`eme nous permet de localiser sur $M_{et}$,
et donc de ramener certaines situations au cas o\`u $F=[X/H]$
est un quotient par un groupe affine.
\end{itemize}

\end{subsubsection}

\begin{subsubsection}{Descente covariante}
\hspace{5mm}
Dans ce paragraphe on consid\'erera des "champs simpliciaux augment\'es".
Comme ces objets ne sont pas r\'eellement des objets simpliciaux
strictement augment\'es, mais seulement "augment\'es \`a isomorphismes pr\`es",
nous allons commencer par fixer le vocabulaire pour \'eviter les
confusions.

\begin{df}\label{chsimp}
\begin{enumerate}
\item
Soit $F$ un champ alg\'ebrique. On d\'efinit la $1$-cat\'egorie quotient
des champs sur $F$, $Ch/F$, par
\begin{itemize}
\item les objets de $Ch/F$ sont les $1$-morphismes
$$u : F' \longrightarrow F$$
\item un morphisme entre $u : F' \longrightarrow F$ et
$v : F'' \longrightarrow F$, est un couple $(f,h)$, o\`u $f$ est un
$1$-morphisme entre $F''$ et $F'$, et
$h$ un $2$-morphisme entre $u\circ f$ et $v$.
\end{itemize}
\item
Un champ simplicial augment\'e vers un champ alg\'ebrique $F$, est un
objet simplicial de $Ch/F$.
\end{enumerate}
\end{df}

\underline{Remarque:} Si on consid\`ere les champs comme des cat\'egories
fibr\'ees (~donc, en particulier, comme des cat\'egories~), un champ
simplicial augment\'e vers $F$ donne lieu \`a un pseudo-foncteur covariant
$$F_{\bullet} : \Delta^{op} \longrightarrow Cat$$
En particulier, si chaque $1$-morphisme de transition est propre et de
dimension
cohomologique finie sur $F$, on peut en prendre l'image par
le pseudo-foncteur covariant $\bf Coh$, et obtenir un
$\Delta^{op}$-pseudo-diagramme de cat\'egories
$$\begin{array}{cccc}
\bf Coh\mit(F_{\bullet}) : & \Delta^{op} & \longrightarrow & Cat \\
& [n] & \mapsto & \bf Coh\mit(F_{n})
\end{array}$$
En consid\'erant les cat\'egories de complexes de $\cal O$-modules
acycliques et \`a cohomologie coh\'erente et born\'ee, on en d\'eduit un
$\Delta^{op}$-pseudo-diagramme (~covariant~) de
cat\'egories bi-compliciales de Waldhausen (~\ref{strict}~). Par les
proc\'ed\'es de strictification \ref{strict}, on en d\'eduit un
$\Delta^{op}$-diagramme
dans $Sp$
$$\begin{array}{cccc}
\bf G\mit : & \Delta^{op} & \longrightarrow & Sp \\
  & [m] & \mapsto & \bf G\mit(F_{m})
\end{array}$$
On posera alors
$$\bf G\mit(F_{\bullet}):=hocolim_{[n] \in \Delta^{op}}\bf
G\mit(F_{n})$$
De plus, l'augmentation de $F_{\bullet}$ vers $F$ induit un
pseudo-morphisme de $F_{\bullet}$ vers le
champ simplicial augment\'e constant $F$. Et donc, un morphisme dans
$HoSp$ (~\ref{strict}~)
$$\bf G\mit(F_{\bullet})\longrightarrow \bf G\mit(F)$$
De m\^eme, si $f : F \longrightarrow F'$
est un $1$-morphisme de champs, $q : F_{\bullet} \longrightarrow F$ et
$q' : F'_{\bullet} \longrightarrow F'$ deux champs simpliciaux
augment\'es, et $f_{\bullet} : F_{\bullet} \longrightarrow F'_{\bullet}$
un morphisme de champs simpliciaux augment\'es sur $F'$,
propre sur chaque $F_{m}$, alors on trouve un diagramme commutatif dans $HoSp$
$$\xymatrix{
\bf G\mit(F_{\bullet}) \ar[r]^{(f_{\bullet})_{*}} \ar[d]_{q_{*}} &
\bf G\mit(F'_{\bullet}) \ar[d]^{q'_{*}} \\
\bf G\mit(F) \ar[r]_{f_{*}} & \bf G\mit(F')}$$

Soit $f :  F_{0} \longrightarrow F$ un $1$-morphisme
propre de champs
alg\'ebriques. Son nerf est le champ alg\'ebrique simplicial augment\'e
sur $F$,
$\cal N\mit(F_{0}/F)$, d\'efini par
$$\begin{array}{cccc}
\cal N\mit(F_{0}/F) : & \Delta &
\longrightarrow & ChAlg(S) \\
 & [m] & \mapsto & \underbrace{F_{0}\times_{F}F_{0} \dots
\times_{F}F_{0}}_{m \;
 fois }
\end{array}$$

\begin{df}\label{d2.2}
Le spectre $\bf G\mit(\cal N\mit(F_{0}/F))$ sera not\'e
$\bf G\mit(F_{0}/F)$.

De la m\^eme fa\c{c}on, si $f$ est repr\'esentable,
on d\'efinit les spectres de $G$-cohomologie relatifs
$$\underline{\bf G}(F_{0}/F):=
hocolim_{[n]\in \Delta^{op}}\underline{\bf G}(\cal
N\mit(F_{0}/F)_{n})$$
\end{df}

Comme la construction pr\'ec\'edente est fonctorielle, elle
garde un sens au niveau des cat\'egories homotopiques. Ainsi, si
$f : F_{0} \longrightarrow F$ est un morphisme propre et de dimension
cohomologique finie dans $HoChAlg(S)$, on dispose des spectres
$\bf G\mit(F_{0}/F)$, et $\underline{\bf G}(F_{0}/F)$ dans le
cas o\`u $f$ est repr\'esentable, munis de
morphismes dans $HoSp$
$$f_{*} : \bf G\mit(F_{0}/F) \longrightarrow \bf G\mit(F)$$
$$f_{*} : \underline{\bf G}(F_{0}/F) \longrightarrow
\underline{\bf G}(F)$$

\begin{thm}\label{th2.2}
\begin{enumerate}
\item
Soit $f : X \longrightarrow F$ un morphisme propre et surjectif
de champs alg\'ebriques, avec $F$ de Deligne-Mumford et $X$ un espace
alg\'ebrique. Alors
le morphisme naturel
$$f_{*} : \underline{\bf G}(X/F) \longrightarrow
\underline{\bf G}(F)$$
est un isomorphisme.
\item
Soit $f : F_{0} \longrightarrow F$ une quasi-enveloppe de Chow
de dimension cohomologique finie. Alors
le morphisme naturel
$$f_{*} : \bf G\mit(F_{0}/F)_{\bf Q} \longrightarrow \bf
G\mit(F)_{\bf Q}$$
est un isomorphisme.
\end{enumerate}
\end{thm}

\underline{\bf Preuve:} \rm $(1)$ En utilisant la localisation
(~\ref{p2.1}~), la continuit\'e (~\ref{p2.1}~), et une r\'ecurrence
noeth\'erienne, on
peut se ramener au cas o\`u $F$ est une gerbe sur un corps $K$. De plus,
comme les colimites homotopiques commutent entre elles,  une utilisation du
corollaire \ref{c2.1} nous ram\`ene au cas o\`u $K$ est s\'eparablement clos,
et donc au cas o\`u $F$ est une gerbe triviale, born\'ee par un groupe fini $H$. \\

Formons le diagramme cart\'esien
$$\xymatrix{
Spec K \ar[r] & F \\
Y \ar[u]^{g} \ar[r] & X \ar[u]_{f}}$$
Il nous permet de construire un diagramme commutatif
$$\xymatrix{
hocolim_{H}\underline{\bf G}(Spec K) \ar[r] & \underline{\bf
G}(F) \\
hocolim_{H}\underline{\bf G}(Y/Spec K) \ar[r]
\ar[u]^{g_{*}} & \underline{\bf G}(X/F) \ar[u]_{f_{*}}}$$
La formule de projection (~\ref{p2.1}~) implique que les fl\`eches
horizontales sont des isomorphismes. Il reste donc \`a montrer que
$g_{*}$ est un isomorphisme. Comme les colimites homotopiques
pr\'eservent les \'equivalences faibles, il suffit de montrer que
$$g_{*} : \underline{\bf G}(Y/Spec K) \longrightarrow \underline{\bf
G}(Spec K)$$
est un isomorphisme.

Soit $L/K$ une extension purement ins\'eparable, telle que $Y$ poss\`ede
un point rationnel sur $L$. Posons $Y_{L}=Y\times_{Spec K}Spec L$.
Alors, on sait que les morphismes
$$\underline{\bf G}(Y_{L}) \longrightarrow \underline{\bf
G}(Y)$$
$$\underline{\bf G}(Spec L) \longrightarrow \underline{\bf
G}(Spec K)$$
sont des isomorphismes (~\cite[$4.7$]{q}~). On peut donc supposer que $Y$
poss\`ede un point rationnel sur $K$. Mais dans ce cas, la section
$s : Spec K \longrightarrow Y$, induit un inverse homotopique
de $g_{*}$ (~\cite[$4.1$ $(I)$]{g3}~). \\

$(2)$ Comme pour le point $(1)$, on peut se ramener au cas o\`u $F$ est
une gerbe sur un corps s\'eparablement clos $K$. Comme, par d\'efinition
des quasi-enveloppes de Chow, le morphisme $F_{0} \longrightarrow F$
dispose alors d'une section apr\`es un changement de base par une
extension purement ins\'eparable de $K$, on peut supposer que
$$F_{0}=F\times_{Spec K}Spec L \longrightarrow F$$
o\`u $L/K$ est purement ins\'eparable.

Dans ce cas on a un isomorphisme
$$\cal N\mit(F_{0}/F)\simeq F\times_{Spec K}\cal N\mit(Spec L/Spec K)$$
De plus, $\cal N\mit(Spec L/Spec K)_{red}$ est isomorphe au champ
simplicial constant $Spec L$. Par invariance topologique, il faut donc
montrer que
$$f_{*} : \bf G\mit(F_{0})_{\bf Q} \longrightarrow \bf G\mit(F)_{\bf Q}$$
est un isomorphisme.

Or la formule de projection (~\ref{p2.1}~) implique que $f_{*}\circ
f^{*}=\times m$, o\`u $m$ est le degr\'e de $L$ sur $K$.

Consid\'erons le carr\'e cart\'esien
$$\xymatrix{
F_{0} \ar[r]^-{f} & F \\
F_{0}\times_{F}F_{0} \ar[r]_-{p} \ar[u]^{q} & F_{0} \ar[u]_{f}}$$
Comme $L/K$ est purement ins\'eparable, la section diagonale
\mbox{$a~:~F_{0}~\longrightarrow~F_{0}\times_{F}F_{0}$}, est une immersion
nilpotente. En particulier
\mbox{$a_{*}~:~\bf G\mit(F_{0})_{\bf Q}~\longrightarrow~\bf G\mit
(F_{0}\times_{F}F_{0})_{\bf Q}$} est un isomorphisme
(\ref{p2.1}~). Ainsi, $p_{*}=q_{*}=(a_{*})^{-1}$.
Un autre application de la formule de projection implique
alors que $p^{*}=q^{*}=m\times a_{*}$.

Le transfert implique alors que
$$\begin{array}{cl}
f^{*} \circ f_{*} & =q_{*} \circ p^{*} \\
& =m\times (a_{*})^{-1}\circ a_{*} \\
& =\times m
\end{array}$$
Ainsi, $f_{*}\circ f^{*}=\times m$, et $f^{*}\circ f_{*}=\times m$.
Donc $f_{*} :  \bf G\mit(F_{0})_{\bf Q} \longrightarrow \bf G\mit(F)$
est un isomorphisme. $\Box$\\

\begin{cor}\label{c2.3}
\begin{enumerate}
\item
Pour toute hyper-quasi-enveloppe de Chow (~\ref{d1.11}~)
$$p : F_{\bullet} \longrightarrow F$$
le morphisme naturel
$$\bf G\mit(F_{\bullet})_{\bf Q}\longrightarrow \bf G\mit(F)_{\bf Q}$$
est un isomorphisme dans $HoSp$.
\item
Le foncteur covariant
$$\begin{array}{cccc}
\underline{\bf G} : & (HoChAlg(S),pr.rep.) & \longrightarrow & HoSp \\
                            & F   & \mapsto & \underline{\bf G}(F)
\end{array}$$
s'\'etend en un foncteur covariant
$$\begin{array}{cccc}
\underline{\bf G} : & (HoChAlgDM(S),pr.) & \longrightarrow & HoSp \\
                            & F   & \mapsto & \underline{\bf G}(F)
\end{array}$$
o\`u $(HoChAlgDM(S),pr.)$ est la sous-cat\'egorie des champs alg\'ebriques
de Deligne-Mumford et morphismes propres. Ce foncteur v\'erifie encore
les formule de transfert et de projection pour des morphismes
non-n\'ecessairement repr\'esentables.

De plus, si $F$ est de Deligne-Mumford, et $p : F
\longrightarrow M$ la projection sur son espace de modules, alors
$$p_{*} : \underline{\bf G}(F) \longrightarrow
\underline{\bf G}(M) \simeq \bf G\mit(M)_{\bf Q}$$
est un isomorphisme.
\item
Si $S=Spec K$, avec $K$ un corps de caract\'eristique nulle, alors le
foncteur pr\'ec\'edent s'\'etend en un foncteur covariant
$$\begin{array}{cccc}
\underline{\bf G} : & (HoChAlg^{aff}(S),pr.) & \longrightarrow & HoSp \\
                            & F   & \mapsto & \underline{\bf G}(F)
\end{array}$$
o\`u $(HoChAlg^{aff}(S),pr.)$ est la sous-cat\'egorie des champs
alg\'ebriques $\Delta$-affines, et morphismes propres.
\end{enumerate}
\end{cor}

\underline{\bf Preuve:} \rm $(1)$ Il suffit d'appliquer un
raisonnement analogue \`a celui fait dans \cite[$4.1$]{g3}.\\

$(2)$ Soit $f : F \longrightarrow F'$ un morphisme propre de champs
de Deligne-Mumford. D'apr\`es \cite[$2.6$]{vi2}, il existe un espace
alg\'ebrique $X$, et un
morphisme propre et surjectif $p : X \longrightarrow F$. On d\'efinit
alors $f_{*}$ par le diagramme commutatif dans $HoSp$ suivant
$$\xymatrix{
\underline{\bf G}(X/F) \ar[dr]^-{(f\circ p)_{*}}
 & \\
\underline{\bf G}(F) \ar[r]_{f_{*}} \ar[u]^-{(p_{*})^{-1}} &
 \underline{\bf G}(F')}$$
Comme le produit fibr\'e de deux morphismes propres et surjectifs est
encore un morphisme propre et surjectif, on v\'erifie ais\'ement que
$f_{*}$ ne d\'epend pas du choix du morphisme $p : X \longrightarrow F$,
v\'erifie la formule de transfert et de projection, et d\'efini bien
un foncteur covariant.
Pour un choix fix\'e de $X \longrightarrow F$, on dispose en r\'ealit\'e
d'un morphisme dans $HoSp(F'_{li})$
$$f_{*} : \bf R\mit f_{*}\underline{G}_{\bf Q} \longrightarrow
\underline{G}_{\bf Q}$$

Soit $p : F \longrightarrow M$ la projection d'un champ de
Deligne-Mumford sur
son espace de modules, et montrons que $p_{*}$ est un isomorphisme.

Par invariance topologique, on peut supposer que $F$ et $M$ sont
r\'eduits. Consid\'erons le morphisme
$$f_{*} : \bf R\mit f_{*}\underline{G}_{\bf Q} \longrightarrow
\underline{G}_{\bf Q}$$
sur $M_{li}$, et montrons que c'est une \'equivalence faible.
Comme localement sur $M_{et}$, $F$ est un quotient par un groupe fini,
on peut se restreindre au cas o\`u $F=[X/H]$, avec $H$ un groupe fini
op\'erant sur un sch\'ema r\'eduit $X$.

On consid\`ere alors le morphisme propre et surjectif $q : X \longrightarrow
F$. On a alors un isomorphisme
$$q_{*} : hocolim_{H}\underline{\bf G}(X) \simeq \underline{\bf
G}(X/F)$$
Il reste \`a montrer que le morphisme induit par la projection
$r~:~X~\longrightarrow~X/H$
$$r_{*} : hocolim_{H}\underline{\bf G}(X) \longrightarrow
\underline{\bf G}(X/H)$$
est un isomorphisme. Ce qui est \'equivalent au lemme suivant

\begin{lem}\label{l2.2}
Soit $X$ un sch\'ema, et $H$ un groupe fini op\'erant sur $X$. Alors le
morphisme naturel $p : X \longrightarrow X/H$ induit un isomorphisme
$$p_{*} : (\underline{\bf G}_{*}(X))_{H} \longrightarrow
\underline{\bf G}_{*}(X/H)$$
\end{lem}

\underline{\bf Preuve:} \rm On peut clairement supposer que $H$ op\`ere
fid\`element, et que $X$ est r\'eduit. En utilisant alors la
localisation (~\ref{p2.1}~), on se ram\`ene au cas o\`u l'action de $H$ est
libre sur $X$. Le lemme provient alors de \ref{c2.1}. $\Box$\\

On vient de voir que $f_{*} : \bf R\mit f_{*}\underline{G}_{\bf Q}
\longrightarrow \underline{G}_{\bf Q}$ est un isomorphisme dans
$HoSp(M_{li})$. Il induit donc un isomorphisme dans $HoSp$
$$f_{*} : \bf H\mit(M_{li},\bf R\mit f_{*}\underline{G}_{\bf Q}) \simeq
\underline{\bf G}(F)
\longrightarrow
\bf H\mit(M_{li},\underline{G}_{\bf Q}) \simeq \underline{\bf G}(M)$$
$\Box$\\

$(3)$ Soit $f : F \longrightarrow F'$ un morphisme propre de champs
alg\'ebriques $\Delta$-affines. Pour chaque $U \longrightarrow F'$
lisse, avec $U$ un espace alg\'ebrique, $f^{-1}(U)$ est un champ alg\'ebrique
$\Delta$-affine, et propre sur $U$. En particulier, il est s\'epar\'e,
donc de Deligne-Mumford (~\ref{s2}~). Ainsi, $f_{*}F$ d\'etermine un
champ en champs de Deligne-Mumford sur $F_{li}$, avec morphismes
de restriction lisses et repr\'esentables. \\

Rappelons que $(Ch/F')$ d\'esigne la cat\'egorie des $1$-morphismes vers
$F'$, dans laquelle un morphisme de $f : F_{1} \longrightarrow F'$
vers $g : F_{2} \longrightarrow F'$ est une paire compos\'ee d'un
$1$-morphisme $u : F_{1}
\longrightarrow F_{2}$ et d'un $2$-morphisme de $g\circ u$ vers $f$.
De m\^eme, nous noterons $(Esp/F')$ la sous-cat\'egorie de $(Ch/F')$ des espaces
alg\'ebriques sur $F'$.

Soit $C\bf Mod\mit^{ac}_{coh}$ la cat\'egorie cofibr\'ee des complexes de $\cal
O$-modules
acycliques, \`a cohomologie born\'ee et coh\'erente,
sur le site $(Esp/F',li)_{li}$ des espaces alg\'ebriques sur $F'$ et
morphismes lisses.
Par strictification,
on peut associer \`a cette cat\'egorie cofibr\'ee un pr\'efaisceau en
cat\'egories bi-compliciales de Waldhausen sur $(Esp/F',li)_{li}$
(~qui sera encore
not\'e $C\bf Mod\mit^{ac}_{coh}$~), et v\'erifiant les hypoth\`eses
suivantes~:
\begin{itemize}
\item Pour tout morphisme propre de $Esp/F'$, $f : X \longrightarrow
Y$, il existe un foncteur "exact"
$$f_{*} : C\bf Mod\mit^{ac}_{coh}(X) \longrightarrow C\bf
Mod\mit^{ac}_{coh}(Y)$$
\item Pour toute paire de morphismes propres $f : X \longrightarrow Y$
et $g~:~Y~\longrightarrow~Z$ d'espaces
alg\'ebriques sur $F'$, on a une \'egalit\'e de foncteurs
$$(g\circ f)_{*}=g_{*}\circ f_{*} : C\bf
Mod\mit^{ac}_{coh}(X) \longrightarrow C\bf Mod\mit^{ac}_{coh}(Z)$$
\item Pour tout diagramme cart\'esien d'espaces alg\'ebriques sur $F'$
$$\xymatrix{Y' \ar[r]^{f'} \ar[d]_{u'} & X' \ar[d]^{u} \\
Y \ar[r]_{f} & X}$$
avec $f$ propre et $u$ lisse, une \'egalit\'e de foncteurs
$$u^{*}\circ f_{*}=f'_{*}\circ (u')^{*} : C\bf
Mod\mit^{ac}_{coh}(Y) \longrightarrow C\bf
Mod\mit^{ac}_{coh}(X')$$
\item Pour tout objet $X$ de $Esp/F'$, il existe une \'equivalence
faible dans $Sp$, compatible avec les images r\'eciproques, et les images
directes
$$K(C\bf Mod\mit^{ac}_{coh}(X))\simeq \bf G\mit(X)$$
\end{itemize}

On pourra donc supposer que $X \mapsto \bf G\mit(X)_{\bf Q}$, o\`u
$X \in Ob(Esp/F')$, est un foncteur strict pour les images r\'eciproques de
morphismes lisses et pour les images directes de morphismes propres, et
qui v\'erifie de plus la formule de transfert strictement.

En clair, si on note $(Esp/F',pr.)$ la cat\'egorie des espaces
alg\'ebriques sur $F'$ et morphismes propres, on a construit un
foncteur covariant
$$\begin{array}{ccc}
(Esp/F',pr.) & \mapsto & Sp(F'_{li}) \\
(p : X \rightarrow F') & \mapsto & p_{*}\bf G\mit
\end{array}$$

Pour tout champ de Deligne-Mumford $F_{0}$ sur $F'$, nous noterons
$Pr(F_{0})$ la cat\'egorie des $1$-morphismes propres et surjectifs
$p : X \longrightarrow F_{0}$ avec $X$ un sch\'ema (~un morphisme entre
$p : X \longrightarrow F_{0}$ et $q : Y \longrightarrow F_{0}$ \'etant
la donn\'ee d'un morphisme $f : X \longrightarrow Y$ et d'un
$2$-morphisme entre $q\circ f$ et $p$~). Le th\'eor\`eme \ref{th2.2}
$(1)$, et la descente \ref{p2.1} $(5)$ impliquent qu'il existe un
isomorphisme naturel dans $HoSp$
$$hocolim_{X \in Pr(F_{0})}\bf G\mit(X/F_{0})_{\bf Q} \longrightarrow
\underline{\bf
G}(F_{0})_{\bf Q}$$
De plus, la propri\'et\'e universelle des colimites homotopique, et
la formule de transfert (~stricte~) pour $\bf G$ montrent que
si $u : F_{0}' \longrightarrow F_{0}$ est un \\
$1$-morphisme lisse et
repr\'esentable de champs sur $F'$, il existe un morphisme naturel
$$u^{*} : hocolim_{X \in Pr(F_{0})}\bf G\mit(X/F_{0})_{\bf Q}
\longrightarrow hocolim_{X' \in Pr(F'_{0})}\bf G\mit(X'/F'_{0})_{\bf Q}$$
qui fait de $U \mapsto hocolim_{X \in Pr(f^{-1}(U)}\bf
G\mit(X/f^{-1}(U))_{\bf Q}$
un foncteur contravariant de $F'_{li}$ vers $Sp$. L'isomorphisme pr\'ec\'edent
montre alors
qu'il existe un isomorphisme dans $HoSp(F'_{li})$ entre
$f_{*}\underline{\bf G}$ et
$U~\mapsto~hocolim_{X \in Pr(f^{-1}(U))}\bf G\mit(X/f^{-1}(U))_{\bf Q}$. De
plus, ce foncteur est naturellement
muni d'une augmentation
$$hocolim_{X \in Pr(f^{-1}(U))}\bf G\mit(X/f^{-1}(U))_{\bf Q}
\longrightarrow \bf G\mit(U)_{\bf Q}$$
donn\'e par les images directes. Ainsi, cette augmentation d\'efinit
un morphisme dans $HoSp(F'_{li})$
$$f_{*} : \bf R\mit f_{*}\underline{G}_{\bf Q} \longrightarrow
\underline{G}_{\bf Q}$$
En passant \`a la cohomologie, on obtient le morphisme cherch\'e dans
$HoSp$
$$f_{*} : \underline{\bf G}(F) \longrightarrow \underline{\bf
G}(F')$$
On v\'erifie ais\'ement que ces images directes sont fonctorielles,
v\'erifient les propri\'et\'es de transfert et de projection, et
\'etendent les images directes pour les champs de Deligne-Mumford.
$\Box$\\

\underline{Remarque:} Au cours de la preuve du point $(3)$, nous
avons d\'efini un foncteur covariant
$$\begin{array}{ccc}
(Esp/F',pr.) & \mapsto & Sp(F'_{li}) \\
(p : X \rightarrow F') & \mapsto & p_{*}\bf G\mit
\end{array}$$
Il est \`a noter, que c'est cette forme de covariance renforc\'ee que nous
imposerons par la suite aux th\'eories cohomologiques utilis\'ees
(~\ref{d3.1}, $3$~).

\end{subsubsection}

\end{subsection}

\begin{subsection}{Th\'eor\`emes de d\'evissage de la $G$-th\'eorie}

\begin{subsubsection}{Cas des champs de Deligne-Mumford}\label{dev}
\hspace{5mm}
Rappelons que pour un champ $F$, nous avons d\'efini son
champ des ramifications $I_{F}$, muni de sa projection canonique
(~\ref{d1.8}~)
$$\pi_{F} : I_{F} \longrightarrow F$$
Remarquons que si $F$ est repr\'esent\'e par un pr\'efaisceau en groupoides,
alors $I_{F}$ est canoniquement \'equivalent au pr\'efaisceau en
groupoide $\Omega F$, dont le groupoide des sections au-dessus d'un
$S$-espace alg\'ebrique $X$ est d\'efini par~:
\begin{itemize}
\item Les objets de $\Omega F$ sont les couples $(s,h)$, avec $s~\in~ObF(X)$,
et $h \in Hom_{F(X)}(s,s)$.
\item Les fl\`eches entre $(s,h)$ et $(s',h')$ sont les fl\`eches $u \in
Hom_{F(X)}(s,s')$, telle que $u^{-1}.h'.u=h$.
\end{itemize}

\begin{df}\label{d2.3}
Soit $F$ un $S$-champ de Deligne-Mumford. On d\'efinit son champ des
ramifications mod\'er\'ees $I_{F}^{t}$, comme le sous-champ de $I_{F}$
form\'e des objets $(s,h) \in ObF(X)\times Aut(s)$, avec $h$ d'ordre premier aux
caract\'eristiques de $S$.

Nous dirons alors que $F$ et mod\'er\'e sur $S$ si $I_{F}^{t}=I_{F}$.
\end{df}

Remarquons que si $F$ est un champ de Deligne-Mumford, alors tous les
automorphismes des objets de $F$ sont d'ordre fini. Ainsi, la
d\'efinition pr\'ec\'edente a un sens. De plus, $I_{F}^{t}$ est un
sous-champ ouvert et ferm\'e de $I_{F}$.

Un champ $F$ de Deligne-Mumford est mod\'er\'e sur son espace de modules, si et seulement si pour
tout point $x \in |F|$, l'ordre de ramification de $F$ en $x$
(~\ref{d1.11}~) est
premier avec la caract\'eristique du corps r\'esiduel $k(x)$.

Par d\'efinition m\^eme, le faisceau sur $(I_{F}^{t})_{et}$ repr\'esent\'e par
$I_{I_{F}^{t}}
\longrightarrow I_{F}^{t}$, poss\`ede une section canonique
$$\cal S\mit : I_{F}^{t} \longrightarrow I_{I_{F}^{t}}$$
qui \`a un objet $(s,h) \in I_{F}^{t}(X)$ au-dessus d'un espace
alg\'ebrique $X$, associe l'automorphisme $h : s \simeq s$, qui v\'erifie
bien $h^{-1}.h.h=h$. \\

Nous noterons $\mu_{\infty}^{t}$ le faisceau sur $(Esp/S)_{et}$, qui \`a
un $S$-espace alg\'ebrique $X$ associe le sous-groupe de
$\cal O\mit_{X}^{*}(X)$ des \'el\'ements de torsion, dont l'ordre est
premier aux caract\'eristiques de $S$. On dispose alors d'un faisceau en
$\bf Q$-alg\`ebres
$$\begin{array}{cccc}
\bf Q\mit[\mu_{\infty}^{t}] : & (Esp/S)_{et} & \longrightarrow &
\bf Q\mit-alg \\
& X & \mapsto & \bf Q\mit[\mu_{\infty}^{t}(X)]
\end{array}$$
Si $H$ est un groupe cyclique d'ordre fini, on peut choisir un
isomorphisme
$$\bf Q\mit[H] \simeq \frac{\bf Q\mit[X]}{X^{m}-1}$$
o\`u $m$ est l'ordre de $H$. On \'ecrit alors
$$X^{m}-1=\prod_{d | m}\phi_{d}(X)$$
o\`u $\phi_{d}(X)$ est le $d$-\`eme polyn\^ome cyclotomique. Le noyau du quotient 
de $\bf Q\mit[H]$ correspondant \`a la projection naturelle
$$\frac{\bf Q\mit[X]}{X^{n}-1}
\longrightarrow \frac{\bf Q\mit[X]}{\phi_{m}(X)},$$
sera not\'ee $I_{H}$. 
Il est ind\'ependant de l'isomorphisme choisi entre
$\bf Q\mit[H]$ et $\frac{\bf Q\mit[X]}{X^{m}-1}$. On notera alors
$$\bf Q\mit(H):=\frac{\bf Q\mit[H]}{I_{H}}.$$
De plus, si $H
\hookrightarrow H'$ est un homomorphisme injectif de groupes cycliques, alors
le morphisme naturel
$$\bf Q\mit[H] \longrightarrow \bf Q\mit[H']$$
induit un morphisme de $\bf Q\mit$-alg\`ebres 
entre $\bf Q\mit(H)$ et $\bf Q\mit(H')$.

En appliquant cette construction aux faisceaux $\mu_{m}$ des racines\\
$m$-\`eme de l'unit\'e sur $(Esp/S)_{et}$, on peut d\'efinir un faisceau en
\mbox{$\bf Q$-alg\`ebres} $\bf Q\mit(\mu_{m})$. On passe alors \`a la limite
inductive sur les entiers $m$ premiers avec les caract\'eristiques de
$S$, et on obtient un faisceau de \mbox{$\bf Q$-alg\`ebres} sur $(Esp/S)_{et}$
$$\begin{array}{cccc}
\Lambda : & (Esp/S)_{et} & \longrightarrow & \bf Q\mit-alg \\
& X & \mapsto & colim_{m}\bf Q\mit(\mu_{m}(X))
\end{array}$$
C'est une alg\`ebre quotient  de l'alg\`ebre de groupes $\bf Q\mit[\mu_{\infty}^{t}]$. 
On dispose donc d'une projection naturelle
$$\bf Q\mit[\mu_{\infty}^{t}] \longrightarrow
\Lambda$$

Dans le cas o\`u $S=Spec k$, avec $k$ un corps contenant les racines de
l'unit\'e, on peut choisir un plongement (~non canonique~)
$$\mu_{\infty}^{t}(Spec k) \hookrightarrow \mu_{\infty}(\bf C\mit)$$
Ceci nous permet alors de d\'efinir un plongement (~non canonique~)
$$\Lambda=\bf Q\mit(\mu_{\infty}^{t}(k)) \hookrightarrow \bf Q\mit^{ab}
\hookrightarrow \bf C$$
o\`u $\bf Q\mit^{ab}$ est la cl\^oture ab\'elienne de $\bf Q$, vue comme un
faisceau constant sur $(Esp/S)_{et}$. \\

Comme il est expliqu\'e dans l'appendice, il existe des objets
$\underline{K}\otimes \Lambda$ et
$\underline{G}\otimes \Lambda$ dans $HoSp((I_{F}^{t})_{et})$. Le premier
est un objet en $\bf Q$-alg\`ebres dans $HoSp((I_{F}^{t})_{et})$, et le
second est un objet en $\underline{K}\otimes \Lambda$-modules. \\

A l'aide de ces notations on va d\'efinir la $K$-cohomologie "\`a
coefficients dans les repr\'esentations". Le choix de la terminologie
sera expliqu\'e au cours du chapitre suivant.

\begin{df}\label{katrep}
Soit $F$ un champ de Deligne-Mumford.
La $K$-cohomologie et la $G$-cohomologie \`a coefficients dans les
repr\'esentations sont d\'efinies par
$$\underline{\bf K}^{rep}(F):=\bf
H\mit((I_{F}^{t})_{et},\underline{K}\otimes \Lambda)$$
$$\underline{\bf G}^{rep}(F):=\bf
H\mit((I_{F}^{t})_{et},\underline{G}\otimes \Lambda)$$
\end{df}

Des propri\'et\'es des foncteurs $\underline{K}$ et $\underline{G}$
(~\ref{p2.1}~) on
tire imm\'ediatement la proposition suivante.

\begin{prop}
La correspondance $F \mapsto \underline{\bf K}^{rep}(F)$ est un foncteur
contravariant de la sous-cat\'egorie de $HoChAlg(S)$ des champs de
Deligne-Mumford vers les $\bf Q$-alg\`ebres de $HoSp$.

La correspondance $F \mapsto \underline{\bf G}^{rep}(F)$ est un
foncteur covariant de la sous-cat\'egorie de $HoChAlg(S)$ des
champs de Deligne-Mumford et morphismes propres, vers $HoSp$.
C'est aussi un foncteur contravariant pour les
morphismes \'etales repr\'esentables.
De plus, on a les propri\'et\'es suivantes

\begin{enumerate}
\item
pour tout champ $F$ de Deligne-Mumford,
$\underline{\bf G}^{rep}(F)$ est un \\
$\underline{\bf K}^{rep}(F)$-module.
Si $f : F' \longrightarrow F$ est un morphisme propre,
$x~\in~\underline{\bf K}_{*}^{rep}(F)$
et $y~\in~\underline{\bf G}_{*}^{rep}(F')$, alors
$$f_{*}(f^{*}(x).y)=x.f_{*}(y)$$
\item
si le carr\'e suivant est cart\'esien
$$\xymatrix{
G' \ar[r]^{q} \ar[d]_{v} & F' \ar[d]^{u} \\
G \ar[r]_{p} & F }$$
avec $p$ propre, et $u$ \'etale, alors
$$q_{*} \circ v^{*} = u^{*} \circ p_{*}$$
\item
 si $j : F' \hookrightarrow F$ est une immersion ferm\'ee, et
$i : U \hookrightarrow F$ l'immersion compl\'ementaire, alors il
existe un triangle fonctoriel dans $HoSp$
$$\xymatrix{
& \underline{\bf G}^{rep}(F') \ar[rd]^{j_{*}} & \\
\underline{\bf G}^{rep}(U) \ar[ru]^{-1} & & \ar[ll]_{i^{*}}
\underline{\bf G}^{rep}(F) & &  }$$
\item
si $p : T \rightarrow F$ est un torseur sous un fibr\'e vectoriel $V$
sur $F$, alors le
morphisme naturel
$$p^{*} : \underline{\bf G}^{\chi}(F) \longrightarrow
\underline{\bf G}^{\chi}(V)$$
est un isomorphisme.
\end{enumerate}
\end{prop}

\underline{\bf Preuve} \rm Il suffit d'appliquer les propri\'et\'es des
foncteurs $\underline{K}$ et $\underline{G}$ au cas des champs
$I_{F}^{t}$. $\Box$\\

\underline{Remarque:} D'apr\`es les propri\'et\'es rappel\'ees en appendice
(~\ref{extcoeff}~),
si $\mu_{\infty}$ est un faisceau constant sur $I_{F}^{t}$
(~on dira dans ce cas que "$S$ contient les racines de l'unit\'e"~),
on a un isomorphisme canonique
$$\underline{\bf K}^{rep}(F)\simeq
\underline{\bf K}(I_{F}^{t})\otimes \Lambda(I_{F}^{t})$$
Comme, $\Lambda$ est le sous-corps de $\bf Q\mit^{ab}$
engendr\'e par les racines de l'unit\'e d'ordre premier aux
caract\'eristiques de $S$, on a
$$\underline{\bf K}_{*}^{rep}(F) \hookrightarrow \underline{\bf
K}_{*}(I_{F}^{t})\otimes \bf Q\mit^{ab}$$

\begin{thm}\label{th2.3}
Soit $F$ un champ de Deligne-Mumford.
\begin{enumerate}
\item
Il existe un morphisme de spectres en anneaux
$$\phi_{F} : \bf K\mit(F) \longrightarrow
\underline{\bf K}^{rep}(F)$$
qui est fonctoriel pour les images r\'eciproques.
\item
Si de plus, $S$ est d'\'egales caract\'eristiques et contient les racines
de l'unit\'e, et $F$ est r\'egulier
alors le morphisme induit
$$\phi_{F} : \bf G\mit(F)\otimes \Lambda(S)
\otimes \bf Q\mit \longrightarrow
\underline{\bf G}^{rep}(F)$$
est un isomorphisme.
\end{enumerate}
\end{thm}

\underline{\bf Preuve:} \rm $(1)$ Soit $\zeta \in
\mu_{\infty}^{t}(I_{F}^{t})$ une section globale.
Commen\c{c}ons par d\'efinir un foncteur exact
$$F_{\zeta} : Vect(I_{F}^{t}) \longrightarrow Vect(I_{F}^{t})$$
Soit $V \in Vect(I_{F}^{t})$, un fibr\'e vectoriel sur $I_{F}^{t}$. Il
est d\'efini par les donn\'ees suivantes~:
\begin{itemize}
\item Pour tout $S$-espace alg\'ebrique $X$, et tout objet $(s,h)$ de
$I_{F}^{t}(X)$, un fibr\'e vectoriel $V_{(s,h)}$ sur $X$.
\item Pour tout couple $(s,h) \in Ob(I_{F}^{t}(X))$, $(s',h') \in
Ob(I_{F}^{t}(Y))$, tout morphisme $f : Y \longrightarrow X$ de
$S$-espaces alg\'ebriques, et tout isomorphisme $H : f^{*}(s,h)
\simeq (s',h')$, un isomorphisme de fibr\'es vectoriels sur $Y$
$$\phi_{f,H} : f^{*}V_{(s,h)} \simeq V_{(s',h')}$$
\item Pour tout couple de morphismes de $S$-espaces alg\'ebriques
$$\xymatrix{ Z \ar[r]^{g} & Y \ar[r]^{f} & X}$$
tout objet $(s,h)~\in~Ob(I_{F}^{t}(X))$, $(s',h')~\in~Ob(I_{F}^{t}(Y))$, et
$(s'',h'')~\in~Ob(I_{F}^{t}(Z))$, et tout isomorphisme
$H_{1} : f^{*}(s,h) \simeq (s',h')$, et
$H_{2}~:~g^{*}(s',h')~\simeq~(s'',h'')$, une \'egalit\'e
$$g^{*}\phi_{f,H_{1}}\circ \phi_{g,H_{2}}=\phi_{f\circ
g,g^{*}H_{1}\circ H_{2}}$$
\end{itemize}
En particulier, comme pour tout objet $(s,h) \in I_{F}^{t}(X)$,
$h \in Hom_{F(X)}(s,s)$, d\'efinit un isomorphisme $h : (s,h)
\longrightarrow (s,h)$ dans $I_{F}^{t}(X)$, le fibr\'e vectoriel
$V_{(s,h)}$ sur $X$ est muni d'une action du groupe cyclique
$<h>$, engendr\'e par $h$ dans $Hom_{F(X)}(s,s)$. Par hypoth\`ese sur
l'ordre de $h$, cette action se diagonalise canoniquement
$$V_{(s,h)} \simeq  V_{(s,h)}^{(\zeta)}\bigoplus W_{(s,h)}$$
o\`u $V_{(s,h)}^{(\zeta)}$, est le sous-fibr\'e o\`u $h$
op\`ere par multiplication par $\zeta$.

De plus, pour tout couple $(s,h) \in Ob(I_{F}^{t}(X))$, $(s',h') \in
Ob(I_{F}^{t}(Y))$, tout morphisme $f : Y \longrightarrow X$ de
$S$-espaces alg\'ebriques, et tout isomorphisme $H : f^{*}(s,h)
\longrightarrow (s',h')$, $\phi_{f,H}$ commute avec les actions de $<h>$ et
$<h'>$.
Il induit donc des isomorphismes
$$\phi_{f,H}^{\zeta} : f^{*}V_{(s,h)}^{(\zeta)} \longrightarrow
V_{(s',h')}^{(\zeta)}$$
La condition de cocycle pour $\phi$ induit alors une condition de cocycle
pour chaque $\phi^{\zeta}$. Ainsi, la donn\'ee de $V_{(s,h)}^{(\zeta)}$, pour
chaque $(s,h) \in Ob(I_{F}^{t})$, et des cocycles $\phi^{\zeta}$,
d\'efinissent un fibr\'e vectoriel $V^{(\zeta)}$ sur $I_{F}^{t}$.

Ce fibr\'e vectoriel $V^{(\zeta)}$, est aussi le sous-fibr\'e de $V$
sur lequel la section canonique $\cal S\mit : I_{F}^{t} \longrightarrow
I_{I_{F}^{t}}$ op\`ere par multiplication par $\zeta$.

Posons alors
$$\begin{array}{cccc}
F_{\zeta} : & Vect(I_{F}^{t}) & \longrightarrow & Vect(I_{F}^{t}) \\
            & V               & \mapsto & V^{(\zeta)}
\end{array}$$
Pour chaque $\zeta \in \mu_{\infty}^{t}$, $F_{\zeta}$ est
un foncteur exact. On peut consid\'erer la somme directe de ces foncteurs
$$\begin{array}{cccc}
F : & \bf Vect\mit(I_{F}^{t}) & \longrightarrow &
\bigoplus_{\mu_{\infty}^{t}(I_{F}^{t})}\bf Vect\mit(I_{F}^{t}) \\
            & V               & \mapsto &
            \oplus_{\zeta \in \mu_{\infty}^{t}(I_{F}^{t})}V^{(\zeta)}
\end{array}$$
Nous noterons $\bigoplus_{\mu_{\infty}}\bf Vect$ le champ associ\'e \`a la
cat\'egorie fibr\'ee
$$\begin{array}{ccc}
(I_{F}^{t})_{et} & \longrightarrow & Cat \\
U & \mapsto & \bigoplus_{\mu_{\infty}(U)}\bf Vect\mit(U)
\end{array}$$
On peut alors v\'erifier que le foncteur $F$ se localise sur
$(I_{F}^{t})_{et}$, et induit un morphisme sur les champs en cat\'egories
exactes sur $(I_{F}^{t})_{et}$
$$F : \bf Vect\mit \longrightarrow \bigoplus_{\mu_{\infty}^{t}}\bf
Vect\mit$$
Nous faisons ici remarquer au lecteur que la cat\'egorie cofibr\'ee
$\bigoplus_{\mu_{\infty}^{t}}\bf
Vect\mit$ n'est g\'en\'eralement pas un champ.
Ce morphisme induit donc un morphisme dans $HoSp((I_{F}^{t})_{et})$
$$F : \underline{K} \longrightarrow K(\bigoplus_{\mu_{\infty}^{t}}\bf
Vect\mit)$$
Or, d'apr\`es \ref{a13} $(4)$, on sait que
$$K(\bigoplus_{\mu_{\infty}^{t}}\bf Vect\mit)_{\bf Q}\simeq \underline{K}\otimes
\bf Q\mit[\mu_{\infty}^{t}]$$
dans $HoSp((I_{F}^{t})_{et})$.

On compose alors avec les deux morphismes canoniques
$$\bf Q\mit[\mu_{\infty}^{t}] \longrightarrow \Lambda$$
$$can : \bf K\mit(I_{F}^{t}) \longrightarrow
\underline{\bf K}(I_{F}^{t})$$
pour obtenir le morphisme cherch\'e
$$F : \bf K\mit(I_{F}^{t}) \longrightarrow \underline{\bf
K}^{rep}(F)$$
Posons alors
$$\phi_{F}=F \circ \pi_{F}^{*} : \bf K\mit(F)
\longrightarrow \underline{\bf K}^{rep}(F)$$
Ce morphisme est clairement fonctoriel en $F$ pour les images
r\'eciproques. La compatibilit\'e avec les produits, provient du fait que
$F$ est un morphisme de champs en cat\'egories tensorielles. Ce qui
peut se v\'erifier localement sur $(I_{F})_{et}$ \`a l'aide de
l'isomorphisme naturel de foncteurs
$$F_{\zeta}(-\otimes -) \simeq
\bigoplus_{\eta.\eta'=\zeta}F_{\eta}(-)\otimes F_{\eta'}(-)$$
Remarquons que le morphisme au niveaux des cat\'egories cofibr\'ees
$$\bf Vect\mit~\longrightarrow~\bigoplus_{\mu_{\infty}^{t}}~\bf
Vect\mit$$
ne pr\'eserve pas le produit tensoriel en g\'en\'eral (~prendre par
exemple $S=Spec \bf R$, et $F=B\bf Z\mit/3$~). \\

$(2)$ Comme $S$ "contient les racines de l'unit\'e", on a (~\ref{a13}
$(5)$~)
$$\underline{\bf G}^{rep}(F)\simeq \underline{\bf
G}(I_{F}^{t})\otimes\Lambda(S)$$
Soit $p : F \longrightarrow M$ l'espace de modules de $F$.
D'apr\`es \ref{th2.1}, pour montrer que le morphisme $\phi_{F}$ induit une
\'equivalence faible en $G$th\'eorie, on peut remplacer $M$ par un recouvrement
\'etale.
Ainsi, \ref{p1.2} nous permet de supposer que $F$ est un champ quotient
$[X/H]$, avec $X$ un $S$-sch\'ema r\'egulier d'\'egales caract\'eristiques,
et $H$ un groupe fini op\'erant sur $X$. On sait alors que $I_{F}^{t}$
est \'equivalent \`a $[\overline{X}/H]$, o\`u
$$\overline{X}:=\{ (x,h) \in X\times H /h.x=x \; et \; h \; d'ordre \;
premier \; avec \; carS \}$$
Ainsi
$$\underline{\bf G}_{*}(I_{F}^{t})\simeq \bigoplus_{h \in H \; o(h)
\; premier \; avec \; carS}\bf G\mit_{*}(X^{h})^{Z(h)}$$
avec $X^{h}$ le sous-sch\'ema des points fixes de $h$, et $Z(h)$ le
centralisateur de $h$ dans $H$.

Mais alors le morphisme
$$\phi_{F} : \bf G\mit_{*}(F)\otimes \Lambda(S) \longrightarrow
\bigoplus_{h \in c(H) \; o(h)
\; premier \; avec \; carS}\bf G\mit_{*}(X^{h})^{Z(h)}\otimes
\Lambda(S)$$
est celui d\'efini dans \cite[Th. $1$]{vi}. On sait alors que c'est un
isomorphisme. $\Box$\\

Notons que  le foncteur
$F$ d\'efini au cours de la preuve du point $(1)$,
existe plus g\'en\'eralement sur la cat\'egorie des
$\cal O\mit_{I_{F}^{t}}$-modules. De plus, il est exact et pr\'eserve la
coh\'erence, la quasi-coh\'erence, et la propri\'et\'e d'\^etre acyclique. La
construction est rigoureusement la m\^eme. Il suffit de se rappeler que
l'action d'un groupe diagonalisable sur un $\cal O$-module se
diagonalise canoniquement.

\begin{cor}\label{c2.4}
Soit $F$ un champ alg\'ebrique de Deligne-Mumford r\'egulier, avec $S$ d'\'egales
caract\'eristiques et contenant les racines de l'unit\'e.
Soit $IM$ l'espace de modules de $I_{F}^{t}$. Alors, il existe un isomorphisme
de $\Lambda(I_{F}^{t})$-modules
$$\bf G\mit(F)\otimes \Lambda(S) \longrightarrow
\bf G\mit(IM)\otimes \Lambda(S)$$
\end{cor}

\underline{\bf Preuve:} \rm Il suffit de composer le \ref{th2.3} $(2)$ avec
\ref{c2.3}. $\Box$\\

\underline{Remarque:}
Notons que cet isomorphisme est compos\'e de deux isomorphismes, dont
l'un pr\'eserve les images r\'eciproques, et l'autre les images directes.
Ainsi, l'isomorphisme pr\'ec\'edent n'est pas fonctoriel en g\'en\'eral. Une
version covariante, ainsi qu'une extension aux cas des champs avec
singularit\'es, seront donn\'ees comme cons\'equences des th\'eor\`emes de
Grothendieck-Riemann-Roch.

Dans le cas o\`u le faisceau $\mu_{\infty}^{t}$ est constant sur $S$
(~i.e. "$S$ contient les racines de l'unit\'e"~), on peut
choisir un plongement
$$\mu_{\infty}^{t} \hookrightarrow \mu_{\infty}(\bf C\mit)$$
Dans ce cas la $\bf Q$-alg\`ebre $\Lambda(S)$
s'identifie \`a un sous-corps de $\bf Q\mit^{ab}$, la cl\^oture ab\'elienne de
$\bf Q$. Ainsi, le corollaire pr\'ec\'edent donne un isomorphisme
$$\bf G\mit_{*}(F)\otimes \bf Q\mit^{ab} \longrightarrow
\bf G\mit_{*}(IM)\otimes \bf Q\mit^{ab}$$

\end{subsubsection}

\begin{subsubsection}{$K$-th\'eorie \`a coefficients dans les caract\`eres}
\hspace{5mm}
Nous venons de voir que pour un champ de
Deligne-Mumford $F$, les spectres de $G$-th\'eorie sont approximables
(~et m\^eme calculables, apr\`es extension des coefficients et dans le cas
r\'egulier~) par la $G$-cohomologie du champ $I_{F}^{t}$. Une
g\'en\'eralisation directe de ses r\'esultats pose quelques difficult\'es.
Nous proposons dans ce paragraphe une approche dont les id\'ees sont
proches de celles du paragraphe pr\'ec\'edent. Cependant, m\^eme dans le cas des
champs qui poss\`edent des quotients g\'eom\'etriques uniformes affines sur un
corps de
caract\'eristique nulle, cette derni\`ere m\'ethode
ne permet pas d'obtenir des r\'esultats analogues \`a \ref{c2.4}. \\

Soit $X$ un $S$-espace alg\'ebrique, et $H \longrightarrow X$ un $X$-espace
alg\'ebrique en groupes. Par abus de notations, nous dirons que $H$ est de
type multiplicatif, s'il est de type multiplicatif et de type fini sur $X$
(~donc quasi-isotrivial d'apr\`es \cite[$X$ $4.5$]{sga3II}~). Avant tout
rappelons que si
$H \longrightarrow X$ est un $X$-espace alg\'ebrique en groupes affine
et lisse sur $X$, alors le foncteur
$$\begin{array}{cccc}
\cal T\mit_{H} : & Esp/X & \longrightarrow & Ens \\
                 &  Y    & \mapsto & \cal T\mit_{H}(Y)
\end{array}$$
o\`u $\cal T\mit_{H}(Y)$ est l'ensemble des sous-groupes de type
multiplicatif
de \\
$H_{Y}:=H\times_{X}Y$, est repr\'esentable (~localement de type fini~),
et lisse sur
$X$ (~\cite[$XI$ $4.1$]{sga3II}~). De plus, si $H' \longrightarrow H$ est
un sous-groupe ferm\'e, alors
$\cal T\mit_{H} \longrightarrow \cal T\mit_{H'}$ repr\'esentable par
une immersion ferm\'ee (~\cite[$XI$ $4.3$]{sga3II}~).

\begin{df}\label{d2.4}
\begin{enumerate}
\item
Soit $F$ un champ alg\'ebrique. Nous d\'efinissons le champ des sous-groupes
de type multiplicatif de $F$
$$\begin{array}{cccc}
D_{F} : & Esp/S & \longrightarrow & Gpd \\
              &  X    & \mapsto         & D_{F}(X)
\end{array}$$
o\`u $D_{F}(X)$ est le groupoide dont les objets sont les couples
$(s,D)$, avec $s \in ObF(X)$, et $D$ un sous-groupe de type
multiplicatif
ferm\'e de $Aut_{X}(s)$, et un isomorphisme entre $(s,D)$ et $(s',D')$
est donn\'e par un isomorphisme $u : s \simeq s'$ dans $F(X)$, tel que
$u^{-1}.D'.u=D$.
\item
Nous d\'efinissons le champ des
caract\`eres de $F$
$$\begin{array}{cccc}
X^{*}_{F} : & (Esp/S) & \longrightarrow & Gpd \\
           &  X    & \mapsto         & X^{*}_{F}(X)
\end{array}$$
o\`u $X^{*}_{F}(X)$ est le groupoide des triplets $(s,D,\chi)$, avec
$(s,D)~\in~ObD_{F}(X)$, et $\chi \in Hom_{Gp/X}(D,\bf G\mit_{m})$ un
caract\`ere de $D$. Un morphisme de $(s,D,\chi)$ vers $(s',D',\chi')$
est un morphisme $u~:~(s,D)~\longrightarrow~(s',D')$ dans
$D_{F}(X)$, tel que $u^{-1}.\chi'.u=\chi$.
\end{enumerate}
\end{df}

D'apr\`es \cite[$IX$ $6.8$]{sga3II}, si $f : H \longrightarrow H'$ est un
morphisme de
$X$-espaces alg\'ebriques en groupes affines sur $X$, et $D
\hookrightarrow H$ un sous-groupe de type multiplicatif, le morphisme induit
$D \longrightarrow H'$ se factorise canoniquement en
$$\xymatrix{ D \ar[r] & D' \ar[r]^{j} & H}$$
o\`u $D'$ est de type multiplicatif, et $j$ une immersion ferm\'ee. Ainsi, on
dispose d'une application $f_{*}$ qui envoie sous-groupes de type
multiplicatif
de $H$ vers sous-groupes de type multiplicatif de $H'$. Ceci permet de voir
que $F \mapsto D_{F}$ est fonctoriel.

Il est clair que ces d\'efinitions passent \`a la cat\'egorie homotopique.
On peut donc d\'efinir le foncteur suivant
$$\begin{array}{cccc}
D_{F} : & HoChAlg(S) & \longrightarrow & HoCh(S) \\
        &   F        & \mapsto         & D_{F}
\end{array}$$
De plus, si $f : F \longrightarrow F'$ est un morphisme repr\'esentable de
champs,
il existe un morphisme naturel
$$X^{*}(f) : X^{*}_{F} \longrightarrow X^{*}_{F'}$$
En effet, si $f$ est repr\'esentable, pour chaque section
$s \in ObF(X)$, et $s'=f(s)\in ObF'(X)$, le morphisme induit
$$Aut_{X}(s) \longrightarrow Aut_{X}(s')$$
est une immersion. Ainsi, si $D$ est un sous-groupe de type
multiplicatif de
$Aut_{X}(s)$, et $\chi \in X^{*}(D)$, l'image de $(s,D,\chi)$ par $f$
est par d\'efinition $(f(s),D,\chi)$, o\`u $D$ est vu comme sous-groupe de
$Aut_{X}(f(s))$.

\begin{prop}\label{p2.2}
Soit $F$ un champ alg\'ebrique tel que, pour chaque objet
$s \in ObF(X)$ au-dessus d'un $S$-espace alg\'ebrique $X$, le
$X$-espace alg\'ebrique en groupes $Aut_{X}(s)$ soit localement sur
$X_{et}$, un sous-groupe
ferm\'e d'un $X$-espace alg\'ebrique en groupes affine et lisse sur $X$
(~par exemple $F$ est un localement un quotient affine~).
Alors le champ
$D_{F}$ est un champ alg\'ebrique, localement de type fini sur $F$.
\end{prop}

\underline{\bf Preuve:} \rm Il nous suffit de d\'emontrer que
$$D_{F} \longrightarrow F$$
est repr\'esentable par une r\'eunion disjointe de champs alg\`ebriques de type fini sur $F$. Comme $D_{F}$
est r\'eunion disjointe des sous-champs $D_{F}^{M}$, param\'etrisants
les sous-groupes de types constants \'egaux \`a $M$, o\`u
$M$ parcourt l'ensemble d'isomorphie des groupes ab\'eliens de type fini,
il suffit de montrer que chaque $D_{F}^{M}$ est alg\'ebrique et de type fini sur $F$.\\

Il faut donc montrer que pour chaque $1$-morphisme $s : X \longrightarrow
F$, avec $X$ repr\'esentable, le champ $D_{F}^{M}\times_{F}X$ est
un espace alg\'ebrique de type fini sur $X$. Consid\'erons ce champ comme un champ sur $X$.
Son groupoide des sections au-dessus de $f : Y \longrightarrow X$, est
form\'e des triplets $(t,D,h)$, o\`u $(t,D)~\in~ObD_{F}^{M}(Y)$, et
$h$ est un $2$-morphisme entre $s \circ f$ et $t$. Par le morphisme de
groupoides $(t,D,h) \mapsto h^{-1}.D.h$, il est \'equivalent au
groupoide discret des sous-groupes de type multiplicatif de type $M$ de
$f^{*}Aut_{X}(s) \longrightarrow Y$. Ainsi, $D_{F}^{M}\times_{F}X$ est
\'equivalent au foncteur $\cal T\mit^{M}(Aut_{X}(s))$
des sous-groupes de type multiplicatif de $Aut_{X}(s)$, de types
constants \'egaux \`a $M$ (~\cite[$IX$ $1.4$]{sga3II}~).
Mais d'apr\`es l'hypoth\`ese, on sait que ce foncteur est repr\'esentable
par un $X$-espace alg\'ebrique de type fini sur $X$. $\Box$\\

\begin{cor}\label{c2.5}
Si $F$ est un champ alg\'ebrique v\'erifiant l'hypoth\`ese de la
proposition pr\'ec\'edente. Alors
le morphisme naturel
$$X^{*}_{F} \longrightarrow F$$
est repr\'esentable.
\end{cor}

\underline{\bf Preuve:} \rm D'apr\`es la proposition, il suffit de
montrer que
$$X^{*}_{F} \longrightarrow D_{F}$$
est repr\'esentable. Montrons que $X^{*}_{F} \longrightarrow D_{F}$ est un
champ constant tordu quasi-isotrivial (~au sens de \cite[$X$
$7$]{sga3II}~).\\

On peut se restreindre
au-dessus d'une composante $D_{F}^{M}$, et supposer que les types
des sous-groupes de type multiplicatif sont constants \'egaux \`a $M$. \\

Soit $X$ un espace alg\'ebrique, et $X \longrightarrow D_{F}^{M}$
une section, correspondant \`a l'objet $(s,D) \in ObD_{F}^{M}$.
Ainsi, $D$ est un sous-groupe de type multiplicatif de $Aut_{X}(s)$. Son
faisceau des caract\`ere $Hom_{Gp/X}(D,\bf G\mit_{m})$, est un faisceau
constant tordu quasi-isotrivial sur $X$, dont les fibres sont
isomorphes au sch\'ema discret $M$. Il est donc repr\'esentable par un
espace alg\'ebrique constant tordu quasi-isotrivial $X^{*}(D) \longrightarrow X$.
Comme cet espace alg\'ebrique repr\'esente le champ $X^{*}_{F}\times_{D_{F}}X$,
il vient que $X^{*}_{F} \longrightarrow D_{F}$ est un champ
constant tordu quasi-isotrivial. $\Box$\\

La proposition et le corollaire pr\'ec\'edents permettent de donner
un sens aux foncteurs suivants
$$\begin{array}{cccc}
D : & HoChAlg'(S) & \longrightarrow & \overline{HoChAlg}(S) \\
          &   F        & \mapsto         & D_{F}
\end{array}$$
$$\begin{array}{cccc}
X^{*} : & (HoChAlg'(S),rep.) & \longrightarrow & \overline{HoChAlg}(S) \\
          &   F        & \mapsto         & X^{*}_{F}
\end{array}$$
o\`u $HoChAlg'(S)$ (~resp. $(HoChAlg'(S),rep.)$~) est la sous-cat\'egorie pleine de
$HoAlgCh(S)$ des champs v\'erifiant l'hypoth\`ese de la proposition (~resp.
et morphismes repr\'esentables.~), et
$\overline{HoChAlg}(S)$ celle des champs alg\'ebriques locamenent de type fini sur $S$.\\

Remarquons aussi, que si $F' \hookrightarrow F$ est une immersion
localement ferm\'ee, alors les diagrammes suivants sont cart\'esiens
$$\xymatrix{
D_{F'} \ar[r] \ar[d] & D_{F} \ar[d]  & &
X^{*}_{F'} \ar[r] \ar[d] & X^{*}_{F} \ar[d] \\
F' \ar[r] & F & &
F' \ar[r] & F }$$

Il est assez difficile de d\'ecrire $D_{F}$ en fonction de $F$ en toute
g\'en\'eralit\'e.
Par contre, la proposition suivante donne une description locale du
morphisme induit $Df : D_{F} \longrightarrow D_{F'}$, lorsque
$f$ est un morphisme repr\'esentable.

\begin{prop}\label{p2.3}
Soit $f : F \longrightarrow F'$ un morphisme repr\'esentable de champs
alg\'ebriques, et $Df : D_{F} \longrightarrow D_{F'}$ le morphisme induit.

Soit $(s,D) : X \longrightarrow D_{F'}$ un $1$-morphisme, avec $X$ un
espace alg\'ebrique, correspondant \`a un sous-groupe de type
multiplicatif
$$D \hookrightarrow Aut_{X}(s)$$
et \`a un $1$-morphisme $s : X \longrightarrow F$.
Notons $Z=X\times_{D_{F'}}D_{F}$, et
\mbox{$X'=X\times_{F'}F$}. Alors il existe un isomorphisme
canonique entre $Z$ et $(X')^{D}$, le sous-espace des points fixes de
$X'$ par l'action de $D$, induite par celle de $Aut_{X}(s)$ sur $X'$.
\end{prop}

\underline{\bf Preuve:} \rm C'est une v\'erification imm\'ediate \`a l'aide
des d\'efinitions. $\Box$\\

Par construction, on dispose du groupe de type multiplicatif universel,
$\cal D\mit_{F} \longrightarrow D_{F}$. C'est un champ en groupes
de type multiplicatif sur $D_{F}$, tel que pour toute section
au-dessus d'un espace alg\'ebrique
$(s,D)~:~X \longrightarrow~D_{F}$, $X\times_{D_{F}}\cal D\mit_{F}\simeq D$
comme $X$-espace alg\'ebrique en groupes. On peut alors d\'efinir le faisceau des
caract\`eres universels
$$\begin{array}{cccc}
\cal X\mit^{*}_{F} : & (Esp/D_{F})_{li} & \longrightarrow & Ab \\
                    & X \rightarrow D_{F} & \mapsto & Hom_{Gp/X}(\cal
                    D\mit_{F}\times_{D_{F}}X,\bf G\mit_{m})
\end{array}$$
De fa\c{c}on plus explicite, si $(s,D) : X \longrightarrow D_{F}$ est une
section au-dessus de $X$, correspondant au sous-groupe $D \hookrightarrow
Aut_{X}(s)$, alors
$$\cal X\mit^{*}_{F}(X)=Hom_{Gp/X}(D,\bf G\mit_{m})$$
C'est aussi le faisceau sur $(Esp/D_{F})_{li}$ repr\'esent\'e par le
morphisme $X^{*}_{F}~\longrightarrow~D_{F}$.

Pour tout groupe ab\'elien de type fini $M$, nous noterons $K[[M]]$
l'anneau $K[M]$ localis\'e le long du syst\`eme multiplicatif engendr\'e par les
\'el\'ements de la forme $1-m$, o\`u $m \in M$.

\begin{df}\label{d2.5}
\begin{itemize}
\item
La $\bf Q$-alg\`ebre $\bf Q\mit[[\cal X\mit^{*}(F)]]$ est not\'ee
$\cal A\mit_{F}$.
\item
La $K$-cohomologie d'un champ $F$, \`a
coefficients dans les caract\`eres, est d\'efinie par
$$\underline{\bf K}^{\chi}(F):=\bf H\mit (D_{F},\underline{K}\otimes
\cal A\mit_{F})$$
o\`u l'on prend la cohomologie dans la cat\'egorie $Sp((Esp/D_{F})_{li})$.

La $G$-cohomologie d'un champ $F$, \`a
coefficients dans les caract\`eres, est d\'efinie par
$$\underline{\bf G}^{\chi}(F):=\bf H\mit ((D_{F})_{li},\underline{G}\otimes
\cal A\mit_{F})$$
o\`u l'on prend la cohomologie dans la cat\'egorie $Sp((D_{F})_{li})$.
\end{itemize}
\end{df}

Notons que $\underline{\bf K}^{\chi}(F)$ est muni naturellement d'une structure
d'alg\`ebre, qui est le produit tensoriel des produits sur
$\underline{K}$ et sur $\bf Q\mit[[\cal X\mit^{*}_{F}]]$
(~\ref{extcoeff}~). De m\^eme,
$\underline{\bf G}^{\chi}(F)$ est un $\underline{\bf
K}^{\chi}(F)$-module.

\begin{prop}\label{p2.4}
La correspondance $F \mapsto \underline{\bf K}^{\chi}(F)$ est un foncteur
contravariant de $HoChAlg'(S)$ vers les anneaux de $HoSp$.

La correspondance $F \mapsto \underline{\bf G}^{\chi}(F)$ est un
foncteur covariant de $(HoChAlg'(S),pr.rep.)$, la sous-cat\'egorie de
champs alg\'ebriques et morphismes propres repr\'esentables, vers celle des
groupes ab\'eliens. C'est aussi un foncteur contravariant pour les
morphismes \'etales repr\'esentables.
De plus, on a les propri\'et\'es suivantes

\begin{enumerate}
\item
pour tout champ $F$, $\underline{\bf G}^{\chi}(F)$ est un
$\underline{\bf K}^{\chi}(F)$-module. Si $f : F' \longrightarrow
F$ est un morphisme propre,
$x \in \underline{\bf K}_{*}^{\chi}(F)$
et $y \in \underline{\bf G}_{*}^{\chi}(F')$, alors
$$f_{*}(f^{*}(x).y)=x.f_{*}(y)$$
\item
si le carr\'e suivant est cart\'esien
$$\xymatrix{
G' \ar[r]^{q} \ar[d]_{v} & F' \ar[d]^{u} \\
G \ar[r]_{p} & F }$$
avec $p$ propre, et $u$ \'etale, alors
$$q_{*} \circ v^{*} = u^{*} \circ p_{*}$$
\item
 si $j : F' \hookrightarrow F$ est une immersion ferm\'ee, et
$i : U \hookrightarrow F$ l'immersion compl\'ementaire, alors il
existe un triangle fonctoriel dans $HoSp$
$$\xymatrix{
& \underline{\bf G}^{\chi}(F') \ar[rd]^{j_{*}} & \\
\underline{\bf G}^{\chi}(U) \ar[ru]^{-1} & & \ar[ll]_{i^{*}}
\underline{\bf G}^{\chi}(F) & &  }$$
\item
si $p : T \rightarrow F$ est un torseur sous un fibr\'e vectoriel $V$
sur $F$, alors le
morphisme naturel
$$p^{*} : \underline{\bf G}^{\chi}(F) \longrightarrow
\underline{\bf G}^{\chi}(V)$$
est un isomorphisme.
\end{enumerate}
\end{prop}

\underline{\bf Preuve:} \rm Commen\c{c}ons par les fonctorialit\'es.

Soit $f : F \longrightarrow F'$ un morphisme de champs de
$HoChAlg'(S)$. Il existe un morphisme de restriction des caract\`eres
$$Res_{f} : Df^{*}\cal X\mit^{*}_{F'} \longrightarrow \cal X\mit^{*}_{F}$$
o\`u $Df : D_{F} \longrightarrow D_{F'}$ est le morphisme induit par
$f$. Il est d\'efini de la fa\c{c}on suivante. Une section de
$Df^{*}\cal X\mit^{*}(F')$ au-dessus d'un espace alg\'ebrique $X$,
est la donn\'ee d'une section $(s,D) \in ObD_{F}(X)$,
d'une section $(s',D',\chi')~\in~Ob\cal X\mit_{F'}(X)$,
et d'un isomorphisme
$u~\in~Isom_{X}(f(s),s')$, tel que
$u^{-1}.Df(D).u=D'$. On pose alors
$$Res_{f}(s,D,\chi'):=u.\chi'.u^{-1}\circ f : D \hookrightarrow Df(D)
\longrightarrow \bf G\mit_{m}$$
Ce morphisme induit donc un morphisme d'alg\`ebres
$$Res_{f} : \bf Q\mit[[Df^{*}\cal X\mit^{*}_{F'}]]=
Df^{*}\cal A\mit_{F} \longrightarrow \cal A\mit_{F}$$
On dispose ainsi d'un morphisme dans $HoSp(Esp/D_{F})$
$$Id\otimes Res_{f} : \underline{K}\otimes Df^{*}\cal A\mit_{F'}
\longrightarrow \underline{K}\otimes \cal A\mit_{F}$$
En prenant l'image par $\bf R\mit Df_{*} : HoSp(Esp/D_{F})
\longrightarrow HoSp(Esp/D_{F'})$, on obtient un morphisme dans
$HoSp(Esp/D_{F'})$
$$\bf R\mit Df_{*}(\underline{K})\otimes \cal A\mit_{F'} \longrightarrow
\bf R\mit Df_{*}(\underline{K}\otimes \cal A\mit_{F})$$
que l'on compose avec les images r\'eciproques pour $\underline{K}$
$$\xymatrix{
\underline{K}\otimes \cal A\mit_{F'} \ar[r]^-{Df^{*}\otimes Id} &
\bf R\mit Df_{*}(\underline{K})\otimes \cal A\mit_{F'} \ar[r] &
\bf R\mit Df_{*}(\underline{K}\otimes \cal A\mit_{F})}$$
En prenant les sections globales sur $D_{F'}$, on trouve
$$f^{*} : \underline{\bf K}^{\chi}(F') \longrightarrow
\underline{\bf K}^{\chi}(F)$$
Il est clair que cette construction est naturelle pour la composition,
et pr\'eserve la structure d'anneau.\\

Soit $f : F \longrightarrow F'$ un morphisme propre
et repr\'esentable. Il induit alors un morphisme naturel
$Ind_{f}~:~\cal X\mit_{F}~\longrightarrow~\cal X\mit_{F}$, et
donc un morphisme d'alg\`ebres
$$Ind_{f} : \cal A\mit_{F} \longrightarrow Df^{*}\cal A\mit_{F'}$$
et donc
$$Id \otimes Ind_{f} : \underline{G} \otimes \cal A\mit_{F} \longrightarrow
\underline{G} \otimes Df^{*}\cal A\mit_{F'}$$
En en prenant l'image directe par $Df$
$$\bf R\mit Df_{*}(\underline{G} \otimes \cal A\mit_{F})
\longrightarrow \bf R\mit Df_{*}\underline{G} \otimes \cal A\mit_{F'}$$
et en composant avec les images directes pour $\underline{G}_{\bf
Q}$
$$\bf R\mit Df_{*} \underline{G}_{\bf Q} \longrightarrow
\underline{G}_{\bf Q}$$
on obtient un morphisme dans $HoSp((D_{F'})_{li})$
$$Df_{*} : \bf R\mit Df_{*}(\underline{G} \otimes \cal A\mit_{F})
\longrightarrow \underline{G} \otimes \cal A\mit_{F'}$$
et donc le morphisme cherch\'e sur les sections globales
$$f_{*} : \underline{\bf G}^{\chi}(F) \longrightarrow
\underline{\bf G}^{\chi}(F')$$
Il est aussi clair que cette construction est naturelle pour la
composition. \\

Passons \`a la d\'emonstration des points $(1)$ \`a $(4)$.\\

Une remarque g\'en\'erale qui nous sera utile, est que lorsque
$f~:~F~\longrightarrow~F'$ est repr\'esentable, le morphisme de
restriction
$$Res_{f} : Df^{*}\cal A\mit_{F'} \longrightarrow \cal A\mit_{F}$$
est un isomorphisme, et $Ind_{f}=(Res_{f})^{-1}$. \\

Les points $(1)$ et $(2)$ se d\'eduisent directement de la formule de
projection \ref{p2.1}, et du fait que $Res_{f}=Ind_{f}^{-1}$ pour un
morphisme repr\'esentable.\\

$(3)$ Comme le diagramme suivant est cart\'esien
$$\xymatrix{
D_{F'} \ar[r]^-{Dj} \ar[d] & D_{F} \ar[d] & D_{U} \ar[l]_{Di} \ar[d] \\
F' \ar[r]_{j} & F & U \ar[l]_{i} }$$
on a un triangle dans $HoSp((D_{F})_{li})$
$$\xymatrix{
& j_{*}\underline{G} \ar[rd]^{j_{*}} & \\
i^{*}\underline{G} \ar[ru]^{-1} & & \ar[ll]_{i^{*}}
\underline{G} & &  }$$
qui induit un triangle apr\`es tensorisation par $\cal A\mit_{F}$ (~car
$\cal A\mit_{F}$ est plat sur $\bf Z\mit$~)
$$\xymatrix{
& j_{*}\underline{G}\otimes \cal A\mit_{F} \ar[rd]^{j_{*}} & \\
i^{*}\underline{G}\otimes \cal A\mit_{F} \ar[ru]^{-1} & & \ar[ll]_{i^{*}}
\underline{G}\otimes \cal A\mit_{F} & &  }$$
Or
$$Ind_{j} : \bf R\mit j_{*}(\underline{G}\otimes \cal A\mit_{F'})
\simeq j_{*}\underline{G}\otimes \cal A\mit_{F}$$
$$Res_{i} : i^{*}\underline{G}\otimes \cal A\mit_{F} \simeq
i^{*}(\underline{\bf G}\otimes \cal A\mit_{U})$$
Ainsi, le triangle pr\'ec\'edent induit un triangle sur la cohomologie
$$\xymatrix{
& \underline{\bf G}^{\chi}(F') \ar[rd]^{j_{*}} & \\
\underline{\bf G}^{\chi}(U) \ar[ru]^{-1} & & \ar[ll]_{i^{*}}
\underline{\bf G}^{\chi}(F) & &  }$$

$(4)$ Commen\c{c}ons par montrer que $Dp : D_{T} \longrightarrow D_{F}$ est
encore un torseur sous un fibr\'e vectoriel.
Soit $(s,D) : X \longrightarrow D_{F}$ une
section, avec $X$ un espace alg\'ebrique, correspondant au
sous-groupe $D \hookrightarrow Aut_{X}(s)$. Alors
$D_{V}\times_{D_{F}}X$ est \'equivalent \`a l'espace $V_{X}^{D}$ des points
fixes de l'action de $D$ sur la restriction $V_{X}$, de $V$ sur $X$
(~\ref{p2.3}~).
Or $V_{X}^{D}$ est un sous-fibr\'e vectoriel de $V_{X}$. Ceci montre que
$D_{V}$ est un sous-fibr\'e vectoriel de $V\times_{F}D_{F}$. De plus,
$D_{T}$ \'etant un torseur sous $D_{V}$, c'est encore un torseur affine
sur $D_{F}$.

On sait alors que (~\ref{p2.1}~)
$$Dp^{*} : \underline{G} \longrightarrow \bf R\mit Dp_{*}\underline{G}$$
est une \'equivalence faible, et donc
$$Dp^{*}\otimes Id : \underline{G}\otimes \cal A\mit_{F'} \longrightarrow
\bf R\mit Dp_{*} \underline{G}\otimes \cal A\mit_{F'}$$
aussi. Mais
$$\bf R\mit Dp_{*}\underline{G}\otimes \cal A\mit_{F'}\simeq \bf R\mit
Dp_{*}(\underline{G}\otimes Dp^{*}\cal A\mit_{F'})$$
et comme la restriction $Res_{f} : Dp^{*}\cal A\mit_{F'} \longrightarrow
\cal A\mit_{F}$ est un isomorphisme, on obtient bien que
$$f^{*}=Dp^{*}\otimes Res_{f} :
\underline{G}\otimes \cal A\mit_{F'} \longrightarrow \bf R\mit
Dp_{*}(\underline{G}\otimes Dp^{*}\cal A\mit_{F'})\simeq
\bf R\mit Dp_{*}(\underline{G}\otimes \cal A\mit_{F})$$
est un isomorphisme. $\Box$\\

\begin{thm}\label{th2.4}
Soit $F$ un champ alg\'ebrique de $HoChAlg'(S)$.
Alors, il existe un morphisme d'alg\`ebres dans $HoSp$
$$\chi_{F} : \bf K\mit(F) \longrightarrow  \underline{\bf K}^{\chi}(F) $$
qui est compatible avec les images r\'eciproques.
\end{thm}

\underline{\bf Preuve:} \rm Comme il existe un morphisme naturel
d'alg\`ebres $\bf Q\mit[\cal X\mit^{*}_{F}]~\longrightarrow~\cal A\mit_{F}$,
il suffit de construire un morphisme
$$\bf K\mit(F) \longrightarrow \bf H\mit(D_{F},\underline{K}_{\bf Q}\otimes
\bf Q\mit[\cal X\mit^{*}_{F}])$$
Construisons un foncteur exact
$$\rho : \bf Vect\mit(D_{F}) \longrightarrow \bigoplus_{\cal
X\mit^{*}_{F}(D_{F})}
\bf Vect\mit(D_{F})$$
Si $V$ est un fibr\'e vectoriel sur $D_{F}$, il est canoniquement muni
d'une action du groupe de type multiplicatif universel
$\cal D\mit_{F} \hookrightarrow I_{D_{F}}$. Pour cette action, $V$ se
d\'ecompose en somme de sous-fibr\'es propres
$$V\simeq \bigoplus_{\chi \in \cal
X\mit^{*}_{F}(D_{F})}V^{(\chi)}\bigoplus W$$
o\`u $V^{(\chi)}$ est le sous-fibr\'e de $V$ sur lequel $\cal D\mit_{F}$
op\`ere par multiplication par $\chi \in Hom_{Gp/D_{F}}(\cal
D\mit_{F},\bf G\mit_{m})$. Le foncteur exact $V \mapsto
\bigoplus_{\chi \in \cal X\mit^{*}_{F}(D_{F})}V^{(\chi)}$ se localise sur
$(D_{F})_{li}$, et induit un morphisme de champs en cat\'egories exactes sur
$(D_{F})_{li}$
$$\rho : \bf Vect\mit \longrightarrow \bigoplus_{\cal
X\mit^{*}_{F}}\bf Vect\mit$$
Remarquons qu'il se passe ici le m\^eme ph\'enom\`ene que dans la preuve de
\ref{th2.4}, \`a savoir que la cat\'egorie cofibr\'ee
$$U \mapsto \bigoplus_{\cal
X\mit^{*}_{F}(U)}\bf Vect\mit(U)$$
n'est pas un champ sur $(D_{F})_{li}$. Nous mettons en garde le
lecteur que la notation
$\bigoplus_{\cal X\mit^{*}_{F}}\bf Vect\mit$ fait r\'ef\'erence au champ associ\'e \`a
cette cat\'egorie cofibr\'ee.

En termes de cocycles, les fibr\'es $V^{(\chi)}$ sont construits de la
fa\c{c}on suivante.

Le fibr\'e $V$ sur $D_{F}$ est d\'efini par les
donn\'ee suivantes~:
\begin{itemize}
\item Pour chaque espace alg\'ebrique $X$, et chaque section
$(s,D) \in ObD_{F}(X)$, un fibr\'e vectoriel $V_{(s,D)}$
sur $X$.
\item Pour chaque morphisme d'espaces alg\'ebriques $f : Y
\longrightarrow X$, chaque couple de sections
$(s,D)\in ObD_{F}(X)$, $(s',D') \in ObD_{F}(Y)$, et
chaque isomorphisme $u : f^{*}(s,D) \simeq (s',D')$ dans
$D_{F}(Y)$, un isomorphisme de fibr\'es vectoriels
$$\phi_{f,u} : f^{*}V_{(s,D)} \simeq V_{(s',D')}$$
\item Pour tout couple de morphismes d'espaces alg\'ebriques
$$\xymatrix{ Z \ar[r]^{g} & Y \ar[r]^{f} & X}$$
tout objet $(s,D)~\in~ObD_{F}(X)$, $(s',D')~\in~ObD_{F}(Y)$, et
$(s'',D'')~\in~ObD_{F}(Z)$, et tous isomorphismes
$u_{1} : f^{*}(s,D) \simeq (s',D')$, et
$u_{2}~:~g^{*}(s',D')~\simeq~(s'',D'')$,
une \'egalit\'e
$$g^{*}\phi_{f,u_{1}}\circ \phi_{g,u_{2}}=\phi_{f\circ
g,g^{*}u_{1}\circ u_{2}}$$
\end{itemize}
Ainsi, chaque $V_{(s,D)}$ sur $X$, est muni d'une action
lin\'eaire du sous-groupe $D$. Comme $D$ est de type multiplicatif,
cette action se diagonalise
$$V_{(s,D)}\simeq \bigoplus_{\chi \in X^{*}(D)}V_{(s,D)}^{(\chi)}
\bigoplus W_{(s,D)}$$
o\`u $X^{*}(D)=Hom_{Gp/X}(D,\bf G\mit_{m})$ est le groupe des caract\`eres de $D$,
et $V^{(\chi)}$ le sous-fibr\'e de $V$ sur lequel $D$ op\`ere par multiplication
par $\chi$. Alors, pour
chaque morphisme $f : Y \longrightarrow X$, chaque couple de sections
$(s,D)\in~ObD_{F}(X)$, $(s',D')~\in ObD_{F}(Y)$, et
chaque isomorphisme $u~:~f^{*}(s,D)~\simeq~(s',D')$ dans
$D_{F}(Y)$, l'isomorphisme $\phi_{f,u}$ pr\'eserve cette
d\'ecomposition, et induit donc un isomorphisme
$$\phi_{f,u} : f^{*}V_{(s,D)}^{(\chi)} \simeq V_{(s',D')}^{(\chi)}$$
Ces isomorphismes v\'erifient \'evidemment la condition de cocycles, et
permettent donc de recoller les $V_{(s,D)}^{(\chi)}$ en un fibr\'e
vectoriel sur $D_{F}$, not\'e $V^{(\chi)}$. \\

Par le morphisme canonique \ref{p1.1} appliqu\'e \`a
la $K$-th\'eorie \`a coefficients dans $\bigoplus_{\cal X\mit^{*}_{F}}\bf
Vect\mit$, le morphisme $\rho$ induit un morphisme de spectres
$$\rho : \bf K\mit(D_{F}) \longrightarrow \bf
H\mit(D_{F},\underline{K}) \longrightarrow
\bf H\mit(D_{F},K(\bigoplus_{\cal X\mit^{*}_{F}}\b Vect\mit) \simeq
\bf H\mit(D_{F},\bigoplus_{\cal X\mit^{*}_{F}}\underline{K}) \longrightarrow
\bf H\mit(D_{F},\underline{K}\otimes \bf Q\mit[\cal X\mit^{*}_{F}])$$
et en composant avec le morphisme naturel $\bf Q\mit[\cal X\mit^{*}_{F}] \hookrightarrow
\cal A\mit_{F}$
$$\rho : \bf K\mit(D_{F}) \longrightarrow \underline{\bf K}^{\chi}(F)$$
On d\'efinit alors
$$\chi_{F} : \xymatrix{
\bf K\mit(F) \ar[r]^-{d_{F}^{*}} & \bf K\mit(D_{F}) \ar[r]^{\rho} &
\underline{\bf K}^{\chi}(F)}$$
o\`u $d_{F} : D_{F} \longrightarrow F$ est la projection.\\

Pour finir, la compatibilit\'e avec le produit provient
du fait que $\rho$ est un foncteur de champs en cat\'egories
tensorielles. Ce qui se v\'erifie localement sur $(D_{F})_{li}$,
\`a l'aide de la formule
$$(V\otimes W)^{(\chi)}\simeq \bigoplus_{\chi_{1}.\chi_{2}=\chi}
V^{(\chi_{1})}\otimes W^{(\chi_{2})}$$
Remarquons qu'ici aussi, le morphisme $\rho$ preserve la
structure tensorielle car nous avons consid\'er\'e le champ
$\bigoplus_{\cal X\mit^{*}_{F}}\bf Vect\mit$, et non la cat\'egorie cofibr\'ee
$U \mapsto \bigoplus_{\cal X\mit^{*}_{F}(U)}\bf Vect\mit(U)$. $\Box$\\

\underline{Remarque:} Lorsque le morphisme $p_{F} : D_{F}
\longrightarrow F$ est de $Tor$-dimension finie, on peut d\'efinir une
image r\'eciproque
$$d_{F}^{*} : \bf G\mit(F) \longrightarrow \bf G\mit(D_{F})$$
La m\^eme construction que pr\'ec\'edemment, donne \'egalement un morphisme
$$\rho : \bf G\mit(D_{F}) \longrightarrow \underline{\bf G}^{\chi}(F)$$
Ce morphisme sera utilis\'e dans le cas des "champs bien ramifi\'es"
(~\ref{d3.4}~), pour d\'emontrer une formule de Riemann-Roch.
Notons alors que les deux constructions sont compatibles pour le
morphisme naturel $\bf Vect\mit~\longrightarrow~\bf Coh$.\\

Une fois de plus, la construction du foncteur $\rho$ garde un sens
dans le cadre plus g\'en\'eral des $\cal O\mit_{D_{F}}$-modules.\\

\underline{Exemple:} Soit $k$ un corps alg\'ebriquement clos, $S=Spec k$,
$H$ un groupe
alg\'ebrique affine et lisse op\'erant sur un sch\'ema $X$, et
\mbox{$F=[X/H]$}. Notons $A$
l'ensemble des classes de conjugaison de sous-groupes de type
multiplicatif
de $H$, et fixons pour chaque $a \in A$, un repr\'esentant
$D_{a}~\hookrightarrow~H$. Notons $M_{a}$ le groupe des caract\`eres de $H$, et
$\cal N\mit_{a}$ le normalisateur de $D_{a}$ dans $H$. Alors $\cal
N\mit_{a}$ op\`ere par conjugaison sur $M_{a}$, et par l'action induite
sur $X^{a}$, le sous-sch\'ema des points fixes de $D_{a}$. Alors on a
des \'equivalences
$$X^{*}(F)\simeq \coprod_{a \in A} [X^{a}\times M_{a}/\cal N\mit_{a}]$$
$$D_{F} \simeq \coprod_{a \in A} [X^{a}/ \cal N\mit_{a}]$$
Ceci provient du fait que le sch\'ema des sous-groupes de type
multiplicatif
de $H$, est la r\'eunion disjointes de ses orbites sous l'action de
$H$ par conjugaison (~voir la preuve de \cite[$XII$ $5.5$]{sga3II}~). Ainsi
$$\underline{\bf K}^{\chi}_{0}(F)\simeq \bigoplus_{a \in A}
\bf H\mit^{0}([X^{a}/\cal N\mit_{a}],
\underline{K}\otimes \bf Q\mit [[M_{a}]])$$
et le morphisme
$$\chi_{F} : \bf K\mit_{0}(F) \longrightarrow \bigoplus_{a \in A}
\bf H\mit^{0}([X^{a}/\cal N\mit_{a}],\underline{K}\otimes \bf Q\mit [[M_{a}]])$$
induit un morphisme
$$\bf K\mit_{0}(F) \longrightarrow (\bigoplus_{a \in A}
\bf K\mit_{0}(X^{a})\otimes \bf Q\mit[[M_{a}]])^{\cal
N\mit_{a}}$$
qui est la somme des morphismes $res_{a}$, d\'ecrits dans \cite[$5.11$]{cg}.
Ainsi, le th\'eor\`eme pr\'ec\'edent est une globalisation du morphisme de localisation
en $K$-th\'eorie \'equivariante. \\

Pour finir, nous allons montrer que l'on peut aussi g\'en\'eraliser le
th\'eor\`eme de d\'evissage des champs de Deligne-Mumford au cas des champs
alg\'ebriques v\'erifiant
l'hypoth\`ese de la proposition \ref{p2.2}. Comme les id\'ees et les
preuves sont tout \`a fait analogues, nous nous contenterons d'une
br\`eve description. \\

\begin{df}
Soit $F$ un champ alg\'ebrique. On d\'efinit le champ $I_{F}^{t,f}$, des
automorphismes d'ordre fini et non-ramifi\'es, comme le sous-champ de
$I_{F}$ form\'e des couples $(s,h)$, o\`u $h$ est d'ordre fini et
premier aux caract\'eristiques de $S$.
\end{df}

En clair, le champ $I_{F}^{t,f}$ est d\'efini par le pr\'efaisceau en
groupoides, dont la valeur sur un espace alg\'ebrique $X$ est le
groupoide des couples $(s,h)$, avec $s\in ObF(X)$, et $h \in
Aut_{F(X)}(s)$ est d'ordre fini, premier aux caract\'eristiques de $S$.

\begin{lem}
Si $F$ est dans $HoChAlg(S)'$ (~i.e. satisfait \`a l'hypoth\`ese de la
proposition \ref{p2.2}~), alors le champ
$I_{F}^{t,f}$ est alg\'ebrique, et la projection naturelle
$$\pi_{F} : I_{F}^{t,f} \longrightarrow F$$
est repr\'esentable.
\end{lem}

\underline{\bf Preuve:} \rm On consid\`ere le morphisme naturel
$$p : I_{F}^{t,f} \longrightarrow D_{F}$$
qui \`a une section $(s,h)$ au-dessus de $X$, associe la section
$(s,<h>)$ de $D_{F}$, o\`u $<h>$ est le sous-sch\'ema en groupe de
$Aut_{X}(s)$ engendr\'e par $h$. Il est clair que ce morphisme est
repr\'esentable, fini et \'etale. En effet, si $(s,D)$ est une section de
$D_{F}$ au-dessus de $X$, $I_{F}^{t,f}\times_{D_{F}}X$ est vide si $D$
n'est pas fini sur $X$, et un $D$-torseur sinon.

Comme d'apr\`es \ref{p2.2}, $D_{F}$ est repr\'esentable, $I_{F}^{f,t}$
aussi. $\Box$\\

\underline{Remarque:} On peut aussi d\'emontrer ce lemme en remarquant que
$I_{F}^{t,f}$ est
un sous-champ ouvert et ferm\'e de $\cal D\mit_{F}$. \\

\begin{df}
Soit $F$ un champ de $HoChAlg(S)'$. On d\'efinit sa \\
$K$-cohomologie et
sa $G$-cohomologie \`a
coefficients dans les repr\'esentations par
$$\underline{\bf K}^{rep}(F):=\bf
H\mit(I_{F}^{t,f},\underline{K}\otimes \Lambda)$$
$$\underline{\bf G}^{rep}(F):=\bf
H\mit((I_{F}^{t,f})_{li},\underline{G}\otimes \Lambda)$$
\end{df}

La m\^eme construction que celle d\'ecrite dans \ref{th2.3} donne alors le
th\'eor\`eme suivant.

\begin{thm}\label{th2.3'}
Soit $F$ un champ de $HoChAlg(S)'$. Alors il existe un morphisme
d'anneaux dans $HoSp$, et fonctoriels pour les images r\'eciproques
$$\phi_{F} : \bf K\mit(F) \longrightarrow
\underline{\bf K}^{rep}(F)$$
\end{thm}

\end{subsubsection}

\begin{subsubsection}{D\'evissage des gerbes r\'eductives}
\hspace{5mm}
Bien que l'on ne sache pas d\'emontrer de r\'esultats analogues \`a
\ref{th2.3}
pour le cas des champs d'Artin, il existe un r\'esultat pour les gerbes
born\'ee par des groupes r\'eductifs. \\

Rappelons que tout sch\'ema constant tordu
quasi-isotrivial sur un sch\'ema normal, est isotrivial
(~c'est un corollaire de \cite[$X$ $5.13$]{sga3II}~) (~i.e. trivial
apr\`es un changement de base \'etale et fini~). \\

Soit $X$ un espace alg\'ebrique normal, et $H \longrightarrow X$ un espace
alg\'ebrique en groupes r\'eductifs (~\cite[$XIX$ $2.7$]{sga3III}~).
On se propose de
"calculer" $\bf G\mit(F)_{\bf Q}$, o\`u $F$ est une gerbe sur
$X$ born\'ee par $H$.

Pour cela nous construisons l'espace des caract\`eres de $F$. Sans perte
de g\'en\'eralit\'e on pourra supposer $X$  connexe. Soit $r$ le rang
r\'eductif de $H$. Notons $\cal T\mit^{max}(F)$ le sous-champ ouvert et ferm\'e
de $D_{F}$ des sous-groupes de type multiplicatif de type $\bf
Z\mit^{r}$. Ainsi, $\cal T\mit^{max}(F)$ est le "champ des tores maximaux
de $F$". Nous noterons $X^{*,max}_{F}=X^{*}_{F}\times_{D_{F}}\cal
T\mit^{max}(F)$.

\begin{lem}
Le champ $X^{*,max}_{F}$ est une gerbe born\'ee par $H$, sur un espace
alg\'ebrique $MX^{*}_{F}$, tel que la projection naturelle
$$p : MX^{*}_{F} \longrightarrow X$$
fasse de $MX^{*}_{F}$ un espace constant tordu et quasi-isotrivial
sur $X$ (~\cite[$X$ $7$]{sga3II}~).

De plus, pour tout point g\'eom\'etrique $\overline{x}$ de $X$, la fibre
de $p$ au-dessus de $\overline{x}$ s'identifie \`a l'ensemble des
caract\`eres invariants par conjugaison, d'un tore maximal de
$H_{x}\otimes k(x)^{sp}$.
\end{lem}

\underline{\bf Preuve:} \rm Comme ceci est local sur $X_{et}$, il
suffit de traiter le cas o\`u $F$ est une gerbe triviale born\'ee par
$H$. On peut aussi supposer que $H$ admet un tore maximal $\bf T\mit
\hookrightarrow H$, diagonalisable sur $X$. Alors $X^{*,max}_{F}$ est \'equivalent
au quotient $[(H/\cal N\mit(\bf T\mit))\times \bf Z\mit^{r}/H]$
(~o\`u $\cal N\mit(\bf T\mit)$ est le normalisateur de $\bf T$ dans $H$~)
ou encore au quotient $[X\times \bf Z\mit^{r}/H]$, o\`u $H$
op\`ere sur $\bf Z\mit^{r}$ par conjugaison \`a travers l'identification
$$Hom_{Gp/X}(\bf T\mit,\bf G\mit_{m})\simeq \bf Z\mit^{r}$$
$\Box$\\

\begin{thm}\label{th2.5}
Soit $Z(H)$ le centre de $H$.
Si la classe d\'efinie par $F$ dans
$H^{2}_{et}(X,Z(H))$ est de torsion, alors
il existe un isomorphisme dans $HoSp$
$$\psi_{F} : \bf G\mit(F)_{\bf Q} \simeq \bf G\mit(MX^{*}_{F})_{\bf Q}$$
\end{thm}

\underline{\bf Preuve:} \rm Le principe de la d\'emonstration
consiste a construire $f : Y
\longrightarrow X$, une normalisation de $X$ dans une extension
galoisienne de $K(X)$. Alors d'apr\`es le
th\'eor\`eme de descente on a
$$f^{*} : \bf G\mit(F)_{\bf Q} \simeq \bf G\mit(F_{Y})_{\bf Q}^{Gal}$$
o\`u $Gal$ est le groupe de galois de $Y$ sur $X$, et
$F_{Y}=F\times_{X}Y$. De plus, comme le carr\'e suivant est cart\'esien
$$\xymatrix{
Y \ar[r]^{f} & X \\
T^{max}(F_{Y}) \ar[u] \ar[r] & T^{max}(F) \ar[u] }$$
on a aussi
$$f^{*} : \bf G\mit(MX^{*}_{F})_{\bf Q} \simeq
\bf G\mit(MX^{*}_{(F_{Y})})_{\bf Q}^{Gal}$$
On construira alors $\psi_{F}$ \`a l'aide de $\psi_{F_{Y}}$ et du carr\'e
commutatif
$$\xymatrix{
\bf G\mit(F)_{\bf Q} \ar[d]_{f^{*}} \ar[r]^-{\psi_{F}} & \bf
G\mit(MX^{*}_{F})_{\bf Q} \ar[d]^{f^{*}} \\
\bf G\mit(F_{Y})_{\bf Q} \ar[r]_-{\psi_{F_{Y}}} & \bf
G\mit(MX^{*}(F_{Y}))_{\bf Q} }$$
Il suffira pour cela de v\'erifier que la construction est compatible
avec les changements de bases \'etales.\\

Comme $X$ est normal, il existe un rev\^etement
galoisien $Y \longrightarrow X$, tel que le sch\'ema en groupes $Out_{X}(H)$,
soit constant sur $Y$ (~\cite[$XXIV$ $4.16$]{sga3III} et
\cite[$XXIV$ $1.3$ $(ii)$]{sga3III}~). Par l'argument ci-dessus, on
peut donc supposer que le sch\'ema des automorphismes ext\'erieurs,
$Out_{X}(H)$ est constant sur $X$.

Soit $Z(H)$ le centre de $H$. Comme $H$ est r\'eductif, alors $Z(H)$ est
un groupe de type multiplicatif de type fini. Or comme $X$ est normal,
il existe un rev\^etement galoisien $Y \longrightarrow X$, tel que la
restriction de $Z(H)$ sur $Y$ soit le produit
direct d'un tore d\'eploy\'e et d'un groupe constant fini. On peut donc
aussi supposer que $Z(H)\simeq (\bf G\mit_{m}/X)^{r}\times Z$, avec
$Z$ un groupe fini.

Soit $a \in H^{1}(X_{et},Out_{X}(H))$ le
lien de $F$ (~\cite{s2}~). Comme $Out_{X}(H)$ est constant et $X$
normal, toute classe dans
$H^{1}(X_{et},Out_{X}(H))$ devient triviale apr\`es un rev\^etement
galoisien $Y \longrightarrow X$.
Ainsi, on peut supposer que $F$ est d\'etermin\'ee par sa
classe $b \in H^{2}(X_{et},Z(H))$ (~\cite{s2}~), qui est de torsion
par hypoth\`ese. Elle provient donc d'un \'el\'ement dans
$H^{2}_{et}(X_{et},\mu_{n}^{r}\times Z)$, correspondant \`a une gerbe
de groupe $\mu_{n}^{n}\times Z$ sur $X_{et}$. Etant un champ de
Deligne Mumford, on sait qu'il existe un rev\^etement ramifi\'e galoisien
$Y \longrightarrow X$, qui la trivialise. Ainsi, on peut supposer que
$F\simeq X\times BH$.\\

Enfin, comme $H$ est r\'eductif sur $X$, qui est normal, le th\'eor\`eme
d'isotrivialit\'e (~\cite[$XXIV$ $4.16$]{sga3III}~), permet de se ramener au cas
o\`u $H$ admet un tore maximal d\'eploy\'e sur $H$,  que nous noterons
alors $\bf T\mit \hookrightarrow H$, et $M$ son groupe des caract\`eres.\\

Construisons un morphisme
$$\chi : \bf G\mit(F) \longrightarrow \bf G\mit(X\times M)^{W}$$
o\`u $W$ est le groupe de Weil de $\bf T$. On proc\`ede de la fa\c{c}on
suivante. En choisissant une section de $F \longrightarrow M$, on peut
identifier $F$ \`a un quotient $[M/H]$, et donc $\bf Coh\mit(F)$ \`a la
cat\'egorie des faisceaux coh\'erents sur $M$ munis d'une action de $H$.
Nous noterons $\bf Coh\mit(M,H)$ cette cat\'egorie. Par restriction, on
dispose d'un foncteur exact
$$Res : \bf Coh\mit(M,H) \longrightarrow \bf Coh\mit(M,\bf T\mit)$$
De plus le groupe $W$ op\`ere par automorphismes int\'erieurs sur $\bf T$,
et donc sur la cat\'egorie $\bf Coh\mit(M,\bf T\mit)$. Le foncteur de
restriction $Res$ devient alors une pseudo-transformation naturelle de
$W$-pseudo-foncteurs, $\bf Coh\mit(M,H)$ \'etant vu comme un diagramme
trivial. Par les proc\'ed\'es de strictification \ref{strict}, on en
d\'eduit un morphisme de spectres de $K$-th\'eorie
$$\bf G\mit(F) \longrightarrow holim_{W}\bf G\mit(B\bf T\mit\times X)=:
\bf G\mit(B\bf T\mit\times X)^{W}$$
Or, on sait que la cat\'egorie $\bf Coh\mit(B\bf T\mit\times X)$ est
canoniquement \'equivalent \`a la cat\'egorie $\bigoplus_{M}\bf Coh\mit(X)$.
Ainsi, on a des isomorphismes canoniques dans $HoSp$
$$\bf G\mit(B\bf T\mit\times X)^{W}\simeq (\bigvee_{M}\bf
G\mit(X))^{W}\simeq \bf Coh\mit(X\times M)^{W}$$
Mais, on sait que
le morphisme $r : X\times M \longrightarrow (X\times M)/W=MX^{*}_{F}$
induit un isomorphisme
$$r^{*} : \bf G\mit(MX^{*}_{F})_{\bf Q} \simeq \bf G\mit(X\times
M)^{W}_{\bf Q}$$
Ainsi, on a construit le morphisme cherch\'e
$$\psi_{F} : \bf G\mit(F)_{\bf Q} \longrightarrow \bf
G\mit(MX^{*}_{F})_{\bf Q}$$
Par construction, ce morphisme est clairement compatible aux
changements de base \'etales sur $X$. Ainsi, l'existence de $\psi_{F}$ est
d\'emontr\'ee. Pour voir que c'est une \'equivalence faible, une
localisation sur $X_{et}$ (~\ref{th2.1}~), permet de se ramener au cas o\`u
$F$ est triviale, et $H$ poss\`ede un tore maximal d\'eploy\'e. On utilise
alors une localisation sur $X$ (~\ref{p2.1}~), et un raisonnement par
r\'ecurrence
noeth\'erienne, pour ce ramener au cas o\`u $X=Spec k$, le spectre d'un
corps. Par une autre application de la descente (~\ref{th2.1}~), on peut
aussi supposer que $k$ est s\'eparablement clos. Et comme une extension
purement ins\'eparable induit un isomorphisme en $G$-th\'eorie
rationnelle, on peut m\^eme supposer que $k$ est alg\'ebriquement clos.

Alors, la cat\'egorie $\bf Coh\mit(F)$ est semi-simple, car $H$ est
r\'eductif. De plus, comme $k$ est alg\'ebriquement clos, l'anneau des
endomorphismes d'un objet simple est isomorphe \`a $k$.
On sait alors que
$$\bf G\mit(F)_{\bf Q} \simeq  \bigvee_{S(\bf Coh\mit(F))}\bf
G\mit(Spec k)_{\bf Q}$$
o\`u $S(\bf Coh\mit(F))$ est l'ensemble des classes d'isomorphie
d'objets simples de $\bf Coh\mit(F)$. Ceci implique, que le
morphisme de Kunneth
$$\bf G\mit_{0}(F) \otimes \bf G\mit_{q}(Spec k) \longrightarrow
\bf G\mit_{q}(F)$$
est un isomorphisme. On se ram\`ene donc \`a d\'emontrer que
$$\begin{array}{cccc}
ch : & \bf G\mit_{0}(F)_{\bf Q} & \longrightarrow & \bf Q\mit[M]^{W}\\
     &     x                    & \mapsto         & ch(x)
\end{array}$$
qui \`a une repr\'esentation associe son caract\`ere, est un isomorphisme. Ce
qui est bien connu. $\Box$\\

Remarquons que l'hypoth\`ese concernant la torsion de l'\'el\'ement dans
$H^{2}_{et}(X,Z(H))$, est
de torsion est v\'erifi\'ee lorsque $X$ est lisse, on encore lorsque $H$
est semi-simple. \\

\end{subsubsection}

\end{subsection}

\end{section}

\newpage
\begin{section}{Chapitre $3$ : Cohomologie \`a coefficients dans les
caract\`eres et th\'eor\`emes de Grothendieck-Riemann-Roch}
\hspace{5mm}
Ce chapitre est consacr\'e aux th\'eor\`emes de Riemann-Roch. On y d\'emontre
essentiellement deux types de formules, les formules de type
Lefschetz-Riemann-Roch, et celles de type Grothendieck-Riemann-Roch.

Dans le premier cas (~\ref{th3.1}, \ref{th3.3}~), ces formules explicitent
le comportement des
morphismes de d\'evissage (~\ref{th2.3}, \ref{th2.4}, \ref{th2.3'}~), par
rapport aux images directes.
Lorsque le champ est un quotient d'un sch\'ema par un automorphisme
d'ordre fini, nous retrouvons la formule de Lefschetz-Riemann-Roch
de \cite{bfm}. Plus g\'en\'eralement, si le champ est un quotient par un
sch\'ema en groupes affine, nous retrouvons la formule
des traces de Lefschetz pour les faisceaux coh\'erents d\'emontr\'ee dans
\cite[$6.4$]{th3}. Remarquons que lorsque l'on applique ces formules \`a des
sch\'emas, on trouve que l'identit\'e commute avec les images directes, ce
qui n'est pas une grande d\'ecouverte.

Les formules de Grothendieck-Riemann-Roch (~\ref{p3.6}, \ref{th3.4}~),
quant \`a elles, explicitent
le comportement du caract\`ere de Chern en $K$-cohomologie par rapport
aux images directes. Dans le cas o\`u on les applique \`a des sch\'emas,
on retrouve les formules de Grothendieck-Riemann-Roch d\'emontr\'ees dans
\cite[$4.1$]{g}. Cependant, appliqu\'ees au morphisme structural d'un champ
alg\'ebrique sur un corps, elles ne donnent pas de formule de
Hirzebruch-Riemann-Roch.

On peut alors dire que le cas le plus int\'eressant est lorsque l'on
compose les formules de Grothendieck-Riemann-Roch, avec celles de
Lefschetz-Riemann-Roch. On obtient, dans ce cas, le th\'eor\`eme de
Grothendieck-Riemann-Roch sous sa forme finale (~\ref{c3.3},
\ref{th3.5}, \ref{th3.5'}~), qui permet de
calculer des caract\'eristiques d'Euler de faisceaux coh\'erents, au moins
dans le cas des champs de Deligne-Mumford.

Notons que nous n'avons r\'eussi \`a d\'emontrer les th\'eor\`emes de
Riemann-Roch pour des morphismes non-repr\'esentables de champs
alg\'ebriques d'Artin, que dans le cas des champs qui admettent des quotients
g\'eom\'etriques uniformes. Cette restriction nous est impos\'ee par le
manque de r\'esultats concernant les quasi-enveloppes de Chow de
morphismes propres, analogues \`a \ref{th1.4} (~mais dans un cadre
relatif~).

Notons aussi que le th\'eor\`eme de Riemann-Roch est d\'emontr\'e en toute
g\'en\'eralit\'e en utilisant la cohomologie \`a coefficients dans les
repr\'esentations, et non celle \`a coefficients dans les caract\`eres. Il
se trouve que cette derni\`ere n'est pas bien adapt\'ee au cas des
morphismes non repr\'esentables, et que le th\'eor\`eme serait faux, m\^eme
pour les cas les plus simples. \\

Bien que cela multiplie les notations et les \'enonc\'es nous avons tenu
\`a donner les formules pr\'ec\'edentes s\'epar\'ement avant de les composer,
car les deux formules nous semblent s\'epar\'ement int\'eressantes. Le
lecteur pourra se convaincre par exemple de l'utilit\'e de la formule
de Grothendieck-Riemann-Roch en $K$-cohomologie dans l'\'etude de la
topologie orbifold des champs de Deligne-Mumford (~\ref{gb1}~). \\

\begin{subsection}{Cohomologie des champs alg\'ebriques}
\hspace{5mm} Nous commencerons ce chapitre par des pr\'eliminaires sur la
cohomologie
des champs alg\'ebriques. Cela nous semble n\'ecessaire, car nous ne
savons d\'efinir des images directes que sous certaines hypoth\`eses
concernant la th\'eorie cohomologique utilis\'ee. Nous expliciterons donc
de quelles propri\'et\'es nous avons besoin, et donnerons deux exemples
de th\'eorie cohomologique v\'erifiant ces hypoth\`eses. \\

On supposera, sauf exception au paragraphe \ref{s3}, que
$S=Spec k$, avec $k$ un corps.

\begin{subsubsection}{Th\'eorie cohomologique avec images
directes}\label{s4}

\begin{df}\label{d3.1}
Une th\'eorie cohomologique avec images directes est la donn\'ee de~:
\begin{enumerate}
\item
Pour chaque entier $i$, un pr\'efaisceau en spectres sur $(Esp/S)_{li}$
$\cal H\mit^{i}$, muni d'une structure de groupe ab\'elien dans $HoSp(S)$.
\item
Des morphismes dans $HoSp(S)$
$$\cal H\mit^{i}\wedge \cal H\mit^{j} \longrightarrow \cal H\mit^{i+j}$$
qui font de $\cal H\mit:=\bigvee_{i}\cal H\mit^{i}$ un anneau gradu\'e
commutatif dans $HoSp(S)$.
\item
Pour tout espace alg\'ebrique lisse $X$, une structure de foncteur covariant
$$\begin{array}{ccc}
(Li/X,pr.) & \longrightarrow & (\cal H\mit_{X_{li}})-mod \\
f : Y\rightarrow X & \mapsto & \bf R\mit f_{*}\cal H\mit \\
\end{array}$$
o\`u $(Li/X,pr.)$ est la cat\'egorie des $X$-espaces alg\'ebriques lisses,
et fortement quasi-projectifs sur $X$ avec morphismes propres,
et \mbox{$(\cal H\mit_{X_{li}})$-mod} celle des objets en
$\cal H\mit$-modules dans $Sp(X_{li})$.
L'image d'un \\
$X$-morphisme
$$\xymatrix{
Z \ar[r]^{u} \ar[rd]_{g} & Y \ar[d]^{f} \\
 & X}$$
sera not\'ee
$$u_{*} : \bf R\mit g_{*}(\cal H\mit) \longrightarrow \bf R\mit
f_{*}(\cal H\mit)$$
Si de plus, $u : Y \longrightarrow X$ est un morphisme propre
d'espaces alg\'ebriques lisses et irr\'eductibles, alors
$u_{*}$ est gradu\'e de degr\'e $d.p$, o\`u $d$ est un entier \'egal \`a $1$
o\`u $2$, et $p=Dim X-Dim Y$.
\item
Pour toute immersion ferm\'ee $j : Y \hookrightarrow X$ d'espaces
alg\'ebriques lisses, un triangle fonctoriel dans $HoSp(X_{li})$
$$\xymatrix{
& \bf R\mit j_{*} (\cal H\mit) \ar[rd]^{j_{*}} & \\
i^{*}\cal H\mit \ar[ur]^{-1} & & \ar[ll]^{i^{*}} \cal H\mit }$$
o\`u $i : X-Y \hookrightarrow X$ est l'immersion ouverte compl\'ementaire.
\item
Un morphisme dans $HoSp(S)$
$$C_{1} : BPic \longrightarrow \cal H\mit^{d}_{[0]}$$
tel que si $p : \bf P\mit(V) \longrightarrow X$ est un fibr\'e projectif
associ\'e \`a un fibr\'e vectoriel $V$ de rang $r+1$ sur $X$, et
$x=C_{1}(\cal O\mit_{1}(\bf P\mit(V))$, alors le morphisme
$$\begin{array}{ccc}
\bigvee_{i=0}^{i=r} \cal H\mit_{X_{li}} & \longrightarrow & \bf
R\mit p_{*}\cal H\mit \\
\vee a_{i} & \mapsto & \sum_{i} p^{*}(a_{i}).x^{i}
\end{array}$$
est un isomorphisme dans $HoSp(X_{li})$.
\end{enumerate}
\end{df}

Remarquons, qu'en dehors du point $(3)$, ces axiomes sont ceux donn\'es
dans \cite[$2.1$]{g}. Le point $(3)$ quant \`a lui, est une covariance
renforc\'ee, qui n'est pas v\'erifi\'ee, \`a priori, pour toutes les
th\'eories cohomologiques habituelles cit\'ees dans \cite[$1.4$]{g}. Nous
connaissons essentiellement deux exemples v\'erifiant ces axiomes. \\

\underline{Exemples:}
\begin{itemize}
\item \underline{La th\'eorie de Gersten (~\cite[$1.4$]{g}~):}\\

Soit $\cal K\mit_{i} : X \mapsto \bf K\mit_{i}(X)$ le pr\'efaisceau du
$i$-\`eme
groupe de $K$-th\'eorie sur $(Esp/S)_{li}$. On lui associe par la
construction de Dold-Puppe un pr\'efaisceau en spectres (~\cite{q2}~)
$$\cal H\mit^{i} := K(\cal K\mit_{i}\otimes_{\bf Z}\bf Q\mit,i)$$
Pour montrer que cette th\'eorie v\'erifie $(3)$, on utilise les
r\'esolutions de Gersten (~\cite[$7.17$]{g}~).

Notons $\cal K\mit_{i} \longrightarrow \cal R\mit^{i}$ la r\'esolution de
Gersten de $\cal K\mit_{i}$ sur $(Esp/S,fl)_{li}$, le site des $S$-espace
alg\'ebriques et morphismes plats. Soit
$$\cal H\mit'_{i}:= K(\cal R\mit^{i}\otimes_{\bf Z}\bf Q\mit,i)$$
Alors, pour tout espace alg\'ebrique lisse $X$, le morphisme induit
$$\cal H\mit^{i} \longrightarrow (\cal H\mit'_{i})$$
est une \'equivalence faible sur $X_{li}$. De plus, comme on a pris le
complexe $\cal R\mit^{i}$ \`a coefficients rationnels, le pr\'efaisceau
$\cal H\mit'_{i}$ est flasque sur $X_{li}$. Ainsi, pour tout
morphisme $f : Y \longrightarrow X$, avec $Y$ lisse, $\bf R\mit f_{*}\cal
H\mit$
est canoniquement isomorphe dans $HoSp(X_{li})$ \`a $f_{*}\cal H\mit'$.
Mais on sait alors, que les images directes de morphismes propres sur
$\cal H\mit'$ v\'erifient les formules de transfert et de projection,
elles permettent donc de d\'efinir le foncteur cherch\'e dans $(3)$.\\

\item \underline{La cohomologie de De Rham (~\cite{h}~):}\\

Supposons que $k$ soit de caract\'eristique nulle. Posons alors
$$\cal H\mit^{i}:=K(\Omega,2i)$$
o\`u $\Omega$ est le complexe de De Rham sur $(Esp/S)_{li}$ d\'efini
dans \cite{h}. Pour montrer que cette th\'eorie v\'erifie $(3)$, on
utilise les r\'esolutions canoniques de \cite[$2$]{h}. On conclut alors par
la m\^eme m\'ethode que pr\'ec\'edemment.\\

\end{itemize}

Pour la suite, on se fixe une th\'eorie cohomologique avec images
directes. Si $X$ est un espace alg\'ebrique, la restriction de $\cal H$
sur le petit site lisse $X_{li}$, sera not\'ee $\cal H\mit_{X}$.

Si $X$ est un espace alg\'ebrique irr\'eductible, on peut trouver d'apr\`es
\cite{jo}, une
hyper-quasi-enveloppe de Chow
$$p : Z_{\bullet} \longrightarrow X$$
telle que chaque $Z_{m}$ soit lisse irr\'eductible et quasi-projectif sur
$Spec k$.
En appliquant $(3)$ au-dessus de $Z_{0}$, on peut d\'efinir un spectre
simplicial
$$\bf H\mit(Z_{\bullet},\cal H\mit^{i}) : [m] \mapsto \bf
H\mit((Z_{m})_{li},\cal H\mit^{i-d_{m}})$$
o\`u $d_{m}=Dim Z_{0}-Dim Z_{m}$.
Notons alors
$$\bf H\mit(Z/X,\cal H\mit^{i}):=hocolim_{\Delta^{op}}\bf H\mit(Z_{\bullet},
\cal H\mit^{i})$$
et
$$\cal H\mit_{i}'(X):=hocolim_{Z_{\bullet}\in HE(X)}\bf H\mit(Z/X,\cal
H\mit^{i})$$
o\`u la limite est prise sur la cat\'egorie $HE(X)$ des hyper-enveloppes de
Chow lisses et quasi-projectives sur $Spec k$. Notons que les images
directes sur $\cal H$ (~\ref{d3.1} $3$~) induisent des morphismes pour
tout entier $i$
$$\cal H\mit'_{i}(X) \longrightarrow \cal H\mit^{i}(X)$$

\begin{lem}\label{l3.1}
Si $X$ est lisse, alors le morphisme naturel
$$\cal H\mit_{i}'(X)_{\bf Q} \longrightarrow
\bf H\mit^{i}(X_{li},\cal H\mit)_{\bf Q}$$
est un isomorphisme dans $HoSp$.

De plus, la correspondance $X \mapsto \cal H\mit_{i}'(X)$ est un
pr\'efaisceau en spectres sur $(Esp/S,li)$, la cat\'egorie des espaces
alg\'ebriques et morphismes lisses.

Si $f : X \longrightarrow Y$ est un morphisme propre d'espaces
alg\'ebriques, alors il existe un morphisme dans $HoSp(Y_{li})$
$$f_{*} : f_{*}\cal H\mit'_{X} \longrightarrow \cal
H\mit'_{Y}$$
compatible avec l'\'equivalence pr\'ec\'edente lorsque $X$ et $Y$ sont
lisses, et $f$ fortement quasi-projectif.
\end{lem}

\underline{\bf Preuve:} \rm La premi\`ere partie se d\'emontre comme les
th\'eor\`emes de descente \ref{th2.1}. La seconde assertion est une
cons\'equence directe de $(3)$, et de la propri\'et\'e universelle des
colimites homotopiques.

Soit $f : X \longrightarrow Y$ un morphisme propre. On consid\`ere
la cat\'egorie $f^{*}HE(Y)$ form\'ee des triplets
$(Z_{\bullet},Z'_{\bullet},u)$, o\`u $Z_{\bullet}$ est un
objet de $HE(X)$, $Z'_{\bullet}$ un objet de $HE(Y)$, et $u$ un morphisme de
$Z_{\bullet}$ vers $f^{-1}Z'_{\bullet}:=Z'_{\bullet}\times_{Y}X$ dans $HE(X)$.
Alors, par $(3)$ on a un morphisme naturel
$$f_{*} : hocolim_{(Z_{\bullet},Z'_{\bullet},u)\in f^{*}HE(Y)}
\bf H\mit(Z/X,\cal H\mit)
\longrightarrow hocolim_{Z'_{\bullet}\in HE(Y)}
\bf H\mit(Z'/Y,\cal H\mit)=\cal H\mit'(Y)$$
Or, comme toute objet $f^{-1}Z'_{\bullet}$, o\`u $Z'_{\bullet} \in HE(Y)$,
est domin\'e par un objet de $HE(X)$, le morphisme naturel
$$a : hocolim_{(Z_{\bullet},Z'_{\bullet},u)\in f^{*}HE(Y)}
\bf H\mit(Z/X,\cal H\mit)\longrightarrow
 hocolim_{Z_{\bullet}\in HE(X)}
\bf H\mit(Z/X,\cal H\mit)=\cal H\mit'(X)$$
est une \'equivalence faible. Comme cette construction est compatible
avec le changement de base par des morphismes lisses, on a donc construit
un diagramme dans $HoSp(Y_{li})$
$$\xymatrix{
hocolim_{(Z_{\bullet},Z'_{\bullet},u)\in f^{*}HE(Y)}
\bf H\mit(Z/X,\cal H\mit)_{X}
\ar[d]_{a} \ar[r]^-{f_{*}} &  \cal H\mit'_{Y} \\
f_{*}\cal H\mit'_{X} & }$$
avec $u$ un isomorphisme. Ceci d\'efinit donc de fa\c{c}on unique un
morphisme dans $HoSp(Y_{li})$
$$f_{*} : f_{*}\cal H\mit'_{X} \longrightarrow \cal
H\mit'_{Y}$$
$\Box$\\

\begin{df}\label{d3.2}
Soit $F$ un champ alg\'ebrique. On d\'efinit sa cohomologie et son
homologie par
$$H^{p}(F,q):=\pi_{dq-p}\bf H\mit(F_{li},\cal H\mit^{q}\otimes \bf
Q\mit)$$
$$H_{p}(F,q):=\pi_{dq-p}\bf H\mit(F_{li},\cal H\mit'_{q}\otimes \bf
Q\mit)$$
\end{df}

Les principales propri\'et\'es sont r\'epertori\'ees dans la proposition suivante.

\begin{prop}\label{p3.1}
La correspondance $F \mapsto H^{\bullet}(F,*)$ est un foncteur
contravariant de $HoChAlg(S)$ vers les
$\bf Q$-alg\`ebres commutatives bi-gradu\'ees.

La correspondance $F \mapsto H_{\bullet}(F,*)$ est un
foncteur covariant de $(HoChAlg(S),pr.rep.)$, la sous-cat\'egorie de
champs alg\'ebriques et morphismes propres repr\'esentables, vers celle des
groupes ab\'eliens. C'est aussi un foncteur contravariant pour les
morphismes lisses et repr\'esentables.
De plus, on a les propri\'et\'es suivantes~:

\begin{enumerate}
\item
Pour tout champ $F$, $H_{\bullet}(F,*)$ est un
$H^{\bullet}(F,*)$-module bi-gradu\'e. Si $f : F' \longrightarrow
F$ est un morphisme propre repr\'esentable, $x \in H^{\bullet}(F,*)$
et $y \in H_{\bullet}(F',*)$, alors
$$f_{*}(f^{*}(x).y)=x.f_{*}(y)$$
\item
Si $F$ est lisse, il existe un isomorphisme
de $H^{\bullet}(F,*)$-modules
$$p_{F} : H^{\bullet}(F,*) \simeq
H_{\bullet}(F,*)$$
compatible avec les images r\'eciproques et les produits.
\item
Si le carr\'e suivant est cart\'esien
$$\xymatrix{
G' \ar[r]^{q} \ar[d]_{v} & F' \ar[d]^{u} \\
G \ar[r]_{p} & F }$$
avec $p$ propre repr\'esentable, et $u$ lisse et repr\'esentable, alors
$$q_{*} \circ v^{*} = u^{*} \circ p_{*}$$
\item
Si $j : F' \hookrightarrow F$ est une immersion ferm\'ee, et
$i : U \hookrightarrow F$ l'immersion compl\'ementaire, alors il
existe une suite exacte fonctorielle
$$\xymatrix{
\dots \ar[r] & H_{\bullet}(F',*) \ar[r]^{j_{*}} & H_{\bullet}(F,*)
\ar[r]^{i^{*}} & H_{\bullet}(U,*) \ar[r] & \dots }$$
\item
Si $p : V \rightarrow F$ est un fibr\'e vectoriel, alors le
morphisme naturel
$$p^{*} : H_{\bullet}(F,*) \longrightarrow
H_{\bullet}(V,*)$$
est un isomorphisme.
\item
Si $p : \bf P\mit(V) \longrightarrow F$ est la projection
d'un fibr\'e projectif associ\'e \`a un fibr\'e vectoriel $V$ de rang $r+1$,
et si $x=C_{1}(\cal O\mit_{P}(1)) \in H^{d}(F,1)$, alors le
morphisme naturel
$$\begin{array}{ccc}
\bigoplus_{i=0}^{i=r} H_{\bullet}(F,*) & \longrightarrow &
H_{\bullet}(P,*) \\
\sum_{i}x_{i} & \mapsto & \sum x^{i}.p^{*}(x_{i})
\end{array}$$
est un isomorphisme.
\end{enumerate}
\end{prop}

\underline{\bf Preuve:} \rm L'existence des images directes pour les morphismes
propres repr\'esentables est une cons\'equence directe du lemme \ref{l3.1}.
Toutes les
autres propri\'et\'es se d\'eduisent alors ais\'ement des axiomes \ref{d3.1}.
$\Box$\\

Il nous reste \`a traiter le cas des images directes pour des morphismes
propres non-n\'ecessairement repr\'esentables. Pour cela, on utilise les
m\^emes arguments que \ref{c2.2}, et on d\'emontre la proposition suivante.

\begin{prop}\label{p3.2}
\begin{enumerate}
\item
Le foncteur covariant
$$\begin{array}{cccc}
H_{\bullet}(-,*) : & (HoChAlg(S),pr.rep.) & \longrightarrow & Ab \\
                            & F   & \mapsto & H_{\bullet}(F,*)
\end{array}$$
s'\'etend en un foncteur covariant
$$\begin{array}{cccc}
H_{\bullet}(-,*) : & (HoChAlgDM(S),pr.) & \longrightarrow & Ab \\
                            & F   & \mapsto & H_{\bullet}(F,*)
\end{array}$$
o\`u $(HoChAlgDM(S),pr.)$ est la sous-cat\'egorie des champs alg\'ebriques
de Deligne-Mumford et morphismes propres. Ce foncteur v\'erifie encore
les formules de transfert et de projection pour des morphismes
non-n\'ecessairement repr\'esentables.

De plus, si $F$ est de Deligne-Mumford, et $p : F
\longrightarrow M$ la projection sur son espace de modules, alors
$$p_{*} : H_{\bullet}(F,*) \longrightarrow
H_{\bullet}(M,*)$$
est un isomorphisme.
\item
Si $S=Spec K$, avec $K$ un corps de caract\'eristique nulle, alors le
foncteur pr\'ec\'edent s'\'etend en un foncteur covariant
$$\begin{array}{cccc}
H_{\bullet}(-,*) : & (HoChAlg^{aff}(S),pr.) & \longrightarrow & Ab \\
                            & F   & \mapsto & H_{\bullet}(F,*)
\end{array}$$
o\`u $(HoChAlg^{aff}(S),pr.)$ est la sous-cat\'egorie des champs
alg\'ebriques $\Delta$-affines, et morphismes propres.
\end{enumerate}
\end{prop}

Terminons par la d\'efinition des classes caract\'eristiques.\\

Dans \cite[$2.2$]{g} sont construits des morphismes dans
$HoSpr((Esp/S)_{li})$
$$C_{i} : \underline{K}_{[0]} \longrightarrow \cal H\mit^{i}_{[0]}$$
Ces morphismes induisent donc, pour tout champ alg\'ebrique $F$, la $i$-\`eme
classe de Chern
$$C_{i} : \underline{\bf K}_{p}(F)\simeq\pi_{p}\bf
H\mit(F_{li},\underline{K}_{[0]}) \longrightarrow H^{di-p}(F,i)$$
Ces classes caract\'eristiques v\'erifient \'evidemment tous les axiomes de \cite{g}.
Elles permettent aussi de d\'efinir le caract\`ere de Chern
$$Ch : \underline{\bf K}_{*}(F) \longrightarrow H^{\bullet}(F,*)$$
qui est un morphisme d'alg\`ebres, fonctoriel pour les images
r\'eciproques.

Enfin, on dispose aussi d'une classe de Todd
$$Td : \underline{\bf K}_{0}(F) \longrightarrow H^{\bullet}(F,*)$$
qui est multiplicatif, et fonctoriel pour les images r\'eciproques.\\

A l'aide de ces d\'efinitions et du "splitting principle", on d\'emontre
l'\'equation suivante (~\cite[$5.3$]{fl}~), qui nous sera utile par la suite.

\begin{lem}\label{l3.2}
Soit $x \in \bf K\mit_{0}(F)$, $r$ son rang, que l'on suppose positif,
et $x^{\vee}$ son dual. Soit
$$\lambda_{i} : \bf K\mit_{0}(F) \longrightarrow \bf K\mit_{0}(F)$$
la $i$-\`eme $\lambda$-op\'eration (~\cite[$V$, $1$]{fl}~), et
$\lambda_{-1}(x):=\sum_{i}(-1)^{i}.\lambda_{i}(x)$. Alors
$$Ch(can(\lambda_{-1}(x))).Td(can(x^{\vee}))=C_{r}(can(x^{\vee}))$$
\end{lem}

\underline{Remarque :} Il nous arrivera d'\'ecrire $C_{i}(x)$ (~de m\^eme
pour $Ch$ et $Td$~), pour signifier $C_{i}(can(x))$,
lorsque $x \in \bf K\mit_{*}(F)$. \\

Pour terminer ce paragraphe, nous allons montrer que ces th\'eories
cohomologiques ne permettent pas de d\'emontrer un th\'eor\`eme de
Riemann-Roch. Le contre-exemple est le suivant.\\

Soit $F=[Spec \bf C\mit/H]$ avec $H$ un groupe fini ab\'elien. Dans ce
cas le fibr\'e tangent est trivial, et la transformation de Riemann-Roch
associ\'ee aux d\'efinitions de $Ch$ et $Td$ pr\'ec\'edentes
$$\tau_{F} : \bf K\mit_{0}(F) \rightarrow H^{*}(F,0)$$
est donc le caract\`ere de Chern. Ainsi, c'est un morphisme d'anneaux. La
propri\'et\'e \ref{p3.2} implique que
$$p_{*} : H^{*}(F,0) \longrightarrow H^{*}(Spec \bf C\mit,0)$$
est un isomorphisme.
Supposons que la th\'eorie cohomologique soit telle que
$H^{*}(Spec \bf C\mit,0)$ soit un corps $K$ de caract\'eristique nulle.
C'est le cas pour les deux exemples que nous avons cit\'e.
Si la formule d'Hirzebruch-Riemann-Roch \'etait
v\'erifi\'ee, on aurait un diagramme commutatif
$$\xymatrix{
\bf K\mit_{0}(F) \ar[r]^-{Ch} \ar[d]_{p_{*}} &
H^{*}(F,*) \ar[d]^{\wr p_{*}}\\
K \ar[r]_{Id} & K }$$
En identifiant $\bf K\mit_{0}(F)$ avec le groupe de Grothendieck des
repr\'esentations lin\'eaires de $H$ dans $\bf C$, on aurait
pour tout $\bf C\mit[H]$-module
$V$ de dimension finie
$$p_{*}Ch(V)=Dim(V^{H})$$
Prenons $V$ de dimension $1$ non triviale. Alors $V^{H}=(0)$. Or si
$m$ est l'ordre de $H$, alors $V^{\otimes m}=1$. On aurait donc
$$Ch(V^{\otimes m})=Ch(V)^{m}=Ch(1)=1 \Rightarrow Ch(V)\not = 0$$
Ce qui est absurde. \\

\end{subsubsection}

\begin{subsubsection}{Cohomologie \`a coefficients dans les caract\`eres}
\hspace{5mm} Fixons-nous une th\'eorie cohomologique avec images
directes $\cal H$. Les groupes de cohomologie et d'homologie
associ\'es \`a cette th\'eorie, seront not\'es comme dans le paragraphe
pr\'ec\'edent, $H^{\bullet}(-,*)$, et $H_{\bullet}(-,*)$. \\

\begin{df}\label{d3.3}
Soit $F$ un champ alg\'ebrique de $HoChAlg'(S)$.
\begin{enumerate}
\item
Sa cohomologie et son homologie, \`a
coefficients dans les caract\`eres, sont d\'efinis par
$$H^{p}_{\chi}(F,q):=\pi_{dq-p}\bf H\mit((D_{F})_{li},\cal
H\mit^{q}\otimes \cal A\mit_{F}\otimes \bf Q\mit)$$
$$H_{p}^{\chi}(F,q):=\pi_{dq-p}\bf H\mit((D_{F})_{li},\cal
H\mit'_{q}\otimes \cal A\mit_{F}\otimes \bf Q\mit)$$
On notera aussi
$$H^{\bullet}_{\chi}(F,*):=\bigoplus_{p,q}H^{p}_{\chi}(F,q)$$
$$H_{\bullet}^{\chi}(F,*):=\bigoplus_{p,q}H^{\chi}_{p}(F,q)$$
\item
Les classes caract\'eristiques \`a coefficients dans les caract\`eres, sont
d\'efinies par \\

\hspace{-25mm}
$\displaystyle{C_{i}^{\chi} : \xymatrix{
\bf K\mit_{p}(F) \ar[r]^-{\chi_{F}} & \bf
H\mit^{-p}((D_{F})_{li},\underline{K}\otimes \cal A\mit_{F}\otimes \bf Q\mit)
\ar[r]^-{C_{i}\otimes Id} & H^{di-p}((D_{F})_{li},\cal H\mit^{i}\otimes
\cal A\mit_{F}\otimes \bf Q\mit)=H^{di-p}_{\chi}(F,i)}}$
\item
Le caract\`ere de Chern et la classe de Todd sont d\'efinis de mani\`ere
analogue
$$Ch^{\chi} :  \xymatrix{
\bf K\mit_{*}(F) \ar[r]^{\chi_{F}} & \underline{\bf K}_{*}^{\chi}(F)
\ar[r]^{Ch} & H^{\bullet}_{\chi}(F,*)}$$
$$Td^{\chi} :  \xymatrix{
\bf K\mit_{0}(F) \ar[r]^{\chi_{F}} & \underline{\bf K}^{\chi}_{0}(F)
\ar[r]^{Td} & H^{\bullet}_{\chi}(F,*)}$$
\end{enumerate}
\end{df}

\underline{Remarques:}
\begin{itemize}
\item
Les classes de Chern d\'efinies ci-dessus ne
v\'erifient par l'axiome de normalisation (~\cite[$2.1$]{g}~). En effet, prenons
le cas o\`u $F=BH$ est le champ classifiant d'un groupe alg\'ebrique sur
un corps alg\'ebriquement clos $k$, $\cal H$ la th\'eorie de Gersten,
et $V$ une repr\'esentation lin\'eaire de $H$. Alors
$$C_{0}(V) \in H^{0}_{\chi}(F,0)\simeq \bigoplus_{D \in
\cal T\mit_{H}(k)}\bf Q\mit[[M_{D}]]$$
est la somme des caract\`eres des repr\'esentations obtenues par
restriction de $D$ sur $V$. Ainsi, $C_{0}(V)\neq 1$ en g\'en\'eral.

\item
Pour tout champ alg\'ebrique $F$, il existe une section canonique
$F~\hookrightarrow~D_{F}$, correspondant au sous-groupe trivial.
Ainsi, on dispose d'isomorphismes fonctoriels
$$H^{\bullet}_{\chi}(F,*)\simeq H^{\bullet}(F,*)\oplus
H^{\bullet}_{\chi\neq 1}(F,*)$$
$$H_{\bullet}^{\chi}(F,*)\simeq H_{\bullet}(F,*)\oplus
H_{\bullet}^{\chi\neq 1}(F,*)$$
Par ces isomorphismes, il nous arrivera d'identifier
$H^{\bullet}(F,*)$ (~resp. $H_{\bullet}(F,*)$~)
\`a son image dans $H^{\bullet}_{\chi}(F,*)$ (~resp.
$H_{\bullet}^{\chi}(F,*)$~).
De cette fa\c{c}on, si $x \in \bf K\mit_{p}(F)$, les classes de Chern
$C_{i}(x)\in H^{d.i-p}(F,i)$ seront vues comme des \'el\'ements de
$H^{\bullet}_{\chi}(F,*)$. De plus, il est clair que
la projection de $C^{\chi}_{i}(x)$ dans $H^{\bullet}(F,*)$ est
$C_{i}(x)$.
\end{itemize}

Les produits sur $\cal H$ et $\bf Z\mit[[\cal X\mit^{*}_{F}]]$ induisent
des produits
sur $\cal H\mit\otimes \cal A\mit \otimes \bf Q$.
De cette fa\c{c}on, les groupes
$H^{\bullet}_{\chi}(F,*)$ sont des anneaux bi-gradu\'es. De m\^eme,
$H_{\bullet}^{\chi}(F,*)$ est un $H^{\bullet}_{\chi}(F,*)$-modules
bi-gradu\'es. Ces structures sont compatibles avec la d\'ecomposition
pr\'ec\'edente. \\

Rappelons que $HoChAlg'(S)$ est la sous-cat\'egorie de $HoChAlg(S)$ des
champs qui v\'erifient l'hypoth\`ese de \ref{p2.2}.

\begin{prop}\label{p3.3}
La correspondance $F \mapsto H^{\bullet}_{\chi}(F,*)$ est un foncteur
contravariant de $HoChAlg'(S)$ vers les
alg\`ebres commutatives bi-gradu\'ees.

La correspondance $F \mapsto H^{\chi}_{\bullet}(F,*)$ est un
foncteur covariant de $(HoChAlg'(S),pr.)$, la sous-cat\'egorie des
champs alg\'ebriques et morphismes propres et repr\'esentables, vers celle des
 groupes ab\'eliens. C'est aussi un foncteur contravariant pour les
morphismes \'etales repr\'esentables.
De plus, on a les propri\'et\'es suivantes

\begin{enumerate}
\item
pour tout champ $F$, $H_{\bullet}^{\chi}(F,*)$ est un
$H^{\bullet}_{\chi}(F,*)$-module bi-gradu\'e. Si $f : F' \longrightarrow
F$ est un morphisme propre repr\'esentable, $x \in H^{\bullet}_{\chi}(F,*)$
et $y \in H^{\chi}_{\bullet}(F',*)$, alors
$$f_{*}(f^{*}(x).y)=x.f_{*}(y)$$
\item
si le carr\'e suivant est cart\'esien
$$\xymatrix{
G' \ar[r]^{q} \ar[d]_{v} & F' \ar[d]^{u} \\
G \ar[r]_{p} & F }$$
avec $p$ propre repr\'esentable, et $u$ \'etale repr\'esentable, alors
$$q_{*} \circ v^{*} = u^{*} \circ p_{*}$$
\item
 si $j : F' \hookrightarrow F$ est une immersion ferm\'ee, et
$i : U \hookrightarrow F$ l'immersion compl\'ementaire, alors il
existe une suite exacte naturelle
$$\xymatrix{
\dots H_{\bullet}^{\chi}(F',*) \ar[r]^{j_{*}} &
H_{\bullet}^{\chi}(F,*) \ar[r]^{i^{*}} & H_{\bullet}^{\chi}(U,*)
\ar[r] & H_{\bullet}^{\chi}(F',*) \dots }$$
\item
si $p : V \rightarrow F$ est un fibr\'e vectoriel, alors le
morphisme naturel
$$p^{*} : H_{\bullet}^{\chi}(F,*) \longrightarrow
H_{\bullet}^{\chi}(V,*)$$
est un isomorphisme.
\end{enumerate}
\end{prop}

\underline{\bf Preuve:} \rm C'est la m\^eme que pour \ref{p2.4}. $\Box$\\

Revenons un moment au cas d'une base g\'en\'erale $S$.\\

Lorsqu'un morphisme $f : F \longrightarrow F'$,
est localement une intersection compl\`ete, on peut d\'efinir le fibr\'e
conormal $\cal N\mit_{f}^{\vee}:=\cal N\mit(F/F')^{\vee}$ de $f$. Nous
noterons alors
$$\lambda_{f}=\lambda_{-1}(\cal N\mit(F/F')^{\vee}) \in \bf
K\mit_{0}(F)$$

Rappelons que nous avons d\'efini un morphisme
$$\rho : \bf K\mit(D_{F}) \longrightarrow \underline{\bf
K}^{\chi}(F)$$
pour tout champ alg\'ebrique $F$ (~voir la preuve de \ref{th2.4}~).

\begin{df}\label{d3.4}
\begin{itemize}
\item
Soit $F$ un champ alg\'ebrique tel que la projection
$$d_{F} : D_{F} \longrightarrow F$$
soit localement d'intersection compl\`ete (~repr\'esentable~).
La classe de ramification de $F$ est d\'efinie par
$$\alpha_{F}:=\rho (\lambda_{d_{F}}) \in \underline{\bf K}^{\chi}_{0}(F)$$
o\`u $\rho$ est d\'efini dans \ref{th2.4}.
\item
Un champ alg\'ebrique de $HoChAlg'(S)$ est dit bien ramifi\'e, si le
morphisme
$$d_{F} : D_{F} \longrightarrow F$$
est localement d'intersection compl\`ete,
et $\alpha_{F}$ est inversible dans $\underline{\bf
K}^{\chi}_{0}(F)$.
\end{itemize}
\end{df}

\begin{prop}\label{p3.4}
Dans tous les cas suivants, le champ $F$ de $HoChAlg'(S)$ est bien ramifi\'e.
\begin{enumerate}
\item Le champ $F$ est une gerbe sur $X$, born\'ee par un $X$-espace
alg\'ebrique en groupes affine et lisse sur $X$.
\item Le champ $F$ est un champ de Deligne-Mumford r\'egulier.
\item Le champ $F$ est localement un quotient affine, lisse sur $S=Spec k$,
avec $k$ un corps.
\item Le champ $F$ est localement un quotient affine, r\'egulier, et pour toute
section $s \in ObF(X)$ au-dessus d'un $S$-espace alg\'ebrique,
$Aut_{X}(s)$ est ab\'elien.
\end{enumerate}
\end{prop}

\underline{\bf Preuve:} Pour montrer qu'un \'el\'ement de $\underline{\bf
K}^{\chi}_{0}(F)$ est inversible, on utilisera le lemme suivant,
appliqu\'e \`a $C=(D_{F})_{li}$, et $K=\underline{K}\otimes\cal
A\mit_{F}$.

\begin{lem}\label{l3.3}
Soit $C$ un site, et $K \in HoSp(C)$, un objet en anneaux. Un \'el\'ement
$x \in \bf H\mit^{0}(C,K)$ est inversible, si et seulement pour tout
objet $X \in C$, il existe un
recouvrement $i : U \longrightarrow X$, tel que
$i^{*}(x) \in \bf H\mit^{0}(U,K)$ est inversible.
\end{lem}

\underline{\bf Preuve:} \rm L'\'el\'ement $x$ est
repr\'esent\'e par un morphisme de pr\'efaisceaux en spectres
$$x : * \longrightarrow HK$$
o\`u $K \hookrightarrow HK$ est une r\'esolution injective. Alors, dire
que $x$ est inversible est \'equivalent \`a dire que le morphisme de
$HoSp(C)$
$$\begin{array}{cccc}
-\wedge x : & HK & \longrightarrow & HK \\
            &  y & \mapsto         & y\wedge x
\end{array}$$
est un isomorphisme. Mais ceci est local sur $C$. $\Box$\\

$(1)$ Comme ceci est local sur $X_{et}$,
on peut supposer que $F=[X/H]$ est une gerbe triviale de groupe $H$. Alors
$$D_{F} \simeq [\cal T\mit_{H}/H]$$
o\`u $\cal T\mit_{H}$ est le sch\'ema des sous-groupes de type
multiplicatif de $H$.
Or, comme $\cal T\mit_{H} \longrightarrow X$ est lisse d'apr\`es
\cite[$XI$ $4.1$]{sga3II}, le morphisme $D_{F} \longrightarrow F$ est lisse.
Ainsi, le fibr\'e conormal de $D_{F} \longrightarrow F$ est trivial, et
donc $\alpha_{F}=1$. \\

$(2)$ En localisant sur l'espace de modules de $F$, on peut supposer
que $F=[X/H]$, o\`u $X$ est un sch\'ema r\'egulier, et $H$ un groupe fini
op\'erant sur $X$. Notons $T(H)$ l'ensemble des sous-groupes
de type multiplicatif de $H\times X \longrightarrow X$.

Soit $\overline{X}=\coprod_{D \in T(H)}X^{D}$, o\`u $X^{D}$ est le sous-sch\'ema
de $X$ des points fixes de $D$. On sait alors que $\overline{X}$ est
r\'egulier (~\cite[$6.2$]{th3}~). Or, $D_{F}\simeq [\overline{X}/H]$. Donc
$D_{F}$
est r\'egulier. Ainsi, le morphisme $D_{F} \longrightarrow F$ est
localement d'intersection compl\`ete.

Supposons toujours que $F=[X/H]$, et notons $M(D)$ le groupe des caract\`eres
de $D$.
Comme $\coprod_{D \in T(H)}X^{D} \longrightarrow D_{F}$
est lisse et surjectif, le lemme \ref{l3.3}, nous permet de nous ramener \`a
d\'emontrer que la restriction de $\alpha_{F}$ \`a
chaque $X^{D}$ est inversible.

Pour chaque $D \in T(H)$, soit $\cal N\mit_{D}^{\vee}$ le fibr\'e conormal de
$X_{D}$ dans $X$. Alors, la
composante de $\alpha_{F}$ support\'ee par $X^{D}$, est de
la forme
$$\alpha_{F}(D)=\prod_{t \in M(D)}
\lambda_{-t}((\cal N\mit_{D}^{\vee})^{(t)}) \in
\bf K\mit_{0}(X^{D})[[M(D)]]_{\bf Q}$$
o\`u $(\cal N\mit_{D}^{\vee})^{(t)}$ est le sous-fibr\'e de $\cal
N\mit_{D}^{\vee}$ o\`u $D$ op\`ere par
multiplication par $t$. Ainsi, le rang de $\alpha_{F}(D)$ est
$$\prod_{t \in M(D)}(1-t)^{rg (\cal N\mit^{\vee}_{D})^{(t)}} \in
\bf Q\mit[[M(D)]]$$
Or, comme $D$ fixe exactement $X^{D}$,
$(\cal N\mit^{\vee}_{D})^{(1)}=0$, et ainsi le rang de $\alpha_{F}(D)$ est
inversible.
Donc $\alpha_{F}(D)$ aussi. \\

$(3)$ On peut clairement supposer que $k$ est alg\'ebriquement clos.\\

Soit $f : F \longrightarrow M$ un morphisme, avec $M$ un espace
alg\'ebrique, tel que localement sur $M_{et}$, $F$ soit un quotient par
un sch\'ema en groupes lisse et affine sur $k$. \\

En localisant sur $M_{et}$, on peut supposer
que $F=[X/H]$, o\`u $X$ est un sch\'ema lisse, et $H$ un sch\'ema en groupes
lisse et affine op\'erant sur $X$.

Alors $D_{F}$ est \'equivalent au champ
quotient
$$D_{F} \simeq \coprod_{a \in A}[X^{D_{a}}/\cal N\mit(D_{a})]$$
o\`u $A$ est l'ensemble des classes de conjugaisons de sous-groupes
de type multiplicatif de $H$, $D_{a}$ un repr\'esentant de $a \in A$,
$\cal N\mit(D_{a})$ le normalisateur de $D_{a}$ dans $H$, et
$X^{D_{a}}$ le sous-sch\'ema des points fixes de $D_{a}$.
Alors, comme $X^{D_{a}} \hookrightarrow X$ est une immersion
r\'eguli\`ere (~\cite[$6.2$]{th3}~), ceci montre que $D_{F} \longrightarrow F$ est
localement d'intersection compl\`ete.

En localisant sur $\coprod_{a \in A} X^{D_{a}} \longrightarrow D_{F}$,
qui est lisse et surjectif, on d\'emontre que $\alpha_{F}$ est
inversible comme dans le $(2)$. \\

$(4)$ Tout comme dans le point $(3)$ on peut supposer que
$F=[X/H]$, avec
$H$ un $S$-groupe lisse et affine sur $S$, op\'erant sur un $S$-espace
alg\'ebrique r\'egulier $X$. En localisant sur $X \longrightarrow F$,
on se ram\`ene \`a d\'emontrer le lemme suivant.

\begin{lem}\label{l3.4}
Soit $H$ un $S$-sch\'ema en groupes ab\'eliens lisse et affine sur $S$, op\'erant
sur un $S$-sch\'ema r\'egulier $X$. Soit $H'$ le $X$-sch\'ema en groupes
$$H'=\{ (x,h) \in X\times_{S}H / h.x=x\}$$
Alors $\cal T\mit_{H'}$ est r\'egulier.
\end{lem}

\underline{\bf Preuve:} \rm Comme ceci est local sur $X_{et}$, on
peut supposer que $H$ est d\'eploy\'e sur $X$ (~i.e. que son sous-groupe de
type multiplicatif maximal est diagonalisable sur $X$~). \\

Soit $T(H)$ l'ensemble des sous-groupes
de type multiplicatif de $H_{X}$. Alors, comme $H$ est
ab\'elien, $\cal T\mit_{H_{X}}=\coprod_{D \in T(H)}X$. Ainsi
$$\cal T\mit_{H'}=\coprod_{D \in T(H)}X^{D}$$
o\`u $X^{D}$ est le sous-sch\'ema de $X$ des points fix\'es par $D$. Or on
sait que $X^{D}$ est r\'egulier (~\cite[$6.2$]{th3}~). $\Box$\\

Le lemme ci-dessus montre que $D_{F}
\longrightarrow F$ est localement d'intersection compl\`ete.

La d\'emonstration du fait que $\alpha_{F}$ est inversible se fait de la
m\^eme fa\c{c}on que dans $(2)$.

$\Box$\\

\underline{Remarque:} Tous les exemples de la proposition sont des
champs qui sont localement des quotients. Il serait tr\`es int\'eressant
de savoir si de fa\c{c}on plus g\'en\'erale un champ alg\'ebrique lisse sur un corps
est bien ramifi\'e. Comme la plupart des champs alg\'ebriques connus sont
des quotients par des groupes affines, la recherche d'un contre
exemple \'eventuel est assez d\'elicate. \\

Pour la fin de ce paragraphe, on revient au cas o\`u $S=Spec k$.

\begin{df}\label{d3.5}
Soit $F$ un champ  alg\'ebrique lisse et bien ramifi\'e. La classe de
Todd de $F$ est d\'efinie par
$$Td^{\chi}(F):=Ch(\alpha_{F}^{-1}).Td(T_{D_{F}}) \in
H^{\bullet}_{\chi}(F,*)$$
o\`u $T_{D_{F}} \in \underline{\bf K}_{0}(F)$ est le fibr\'e virtuel tangent de $F$.

La transformation de Riemann-Roch \`a coefficients dans les caract\`eres est alors
d\'efinie par
$$\begin{array}{cccc}
\tau^{\chi}_{F} : & \bf G\mit_{*}(F) & \longrightarrow &
H^{\bullet}_{\chi}(F,*) \\
           & x                & \mapsto & Ch^{\chi}(x).Td^{\chi}(F)
\end{array}$$
\end{df}

\underline{Remarque:} Dans \ref{d3.3}, nous n'avons d\'efini $Ch^{\chi}$ que
pour les \'el\'ements de $\bf K\mit_{*}(F)$. Cependant, lorsque $F$ est
lisse et bien ramifi\'e, on peut refaire le th\'eor\`eme \ref{th2.4} en
$G$-th\'eorie. On d\'efinit $Ch^{\chi}$ de fa\c{c}on analogue \`a
\ref{d3.3}, en gardant \`a l'esprit que $\underline{\bf K}^{\chi}(F)\simeq
\underline{\bf G}^{\chi}(F)$, car $D_{F}$ est lisse (~\ref{p2.1}~).
Ainsi, la d\'efinition pr\'ec\'edente poss\`ede un sens.\\

\end{subsubsection}

\end{subsection}

\begin{subsection}{Formules de Riemann-Roch}
\hspace{5mm}
Nous poss\'edons maintenant tous les outils pour d\'emontrer une formule
de Grothendieck-Riemann-Roch \`a valeur dans la cohomologie \`a
coefficients dans les caract\`eres. Cependant, nous n'avons pas r\'eussi \`a
d\'emontrer le th\'eor\`eme pour le cas g\'en\'eral des morphismes propres
et repr\'esentables de
champs localement quotients affines sur un corps, alors que ce degr\'e
de g\'en\'eralit\'e semble pourtant accessible. En majeure partie, c'est
le manque de r\'esultats concernant les (~quasi~) enveloppes de Chow des
morphismes propres de champs alg\'ebriques, ainsi que la r\'esolution des
faisceaux coh\'erents par des fibr\'es vectoriels, qui limite le cadre
d'application du th\'eor\`eme. D'autre part le cas des morphismes propres
non-repr\'esentables ne peut pas \^etre trait\'e avec les d\'efinitions que
nous avons pour le moment. Nous montrerons cependant en fin de
chapitre, qu'il est possible de g\'en\'eraliser les r\'esultats obtenus pour
les champs de Deligne-Mumford aux champs qui poss\`edent des
quotients g\'eom\'etriques uniformes (~tout au moins en caract\'eristique
nulle~). Je n'ai pas r\'eussi \`a d\'emontrer le th\'eor\`eme dans le cas
plus g\'en\'eral des champs $\Delta$-affines bien ramifi\'es,
alors qu'il me semble que la formule de Riemann-Roch reste vraie.  \\

Un autre fa\c{c}on d'\'eliminer les hypoth\`eses superflues serait de trouver
un moyen de ramener le th\'eor\`eme de Riemann-Roch \`a un probl\`eme "local
en bas". C'est \`a dire que pour d\'emonter la formule pour un morphisme
$f : F\longrightarrow F'$, il suffirait de la d\'emontrer apr\`es un
changement de base par un morphisme lisse est surjectif $X
\longrightarrow F'$. Mais pour l'instant, les preuves existantes de la formule
de Riemann-Roch ne permettent pas ce genre de r\'eduction.
Peut-\^etre qu'une utilisation syst\'ematique des m\'ethodes homotopiques
de \cite{s3}, permettrait d'avancer dans cette direction.

\begin{subsubsection}{Formule de Lefschetz-Riemann-Roch}
\hspace{5mm}

Le point crucial pour la formule de Lefschetz-Riemann-Roch est le cas
particulier de la formule d'intersection suivante.

\begin{prop}\label{p3.5}
Soit $j : F \longrightarrow F'$ une immersion r\'eguli\`ere de champs
alg\'ebriques bien ramifi\'es. Notons
$$\alpha_{j}:=\lambda_{-1}(Dj^{*}\cal N\mit_{D_{F'}/F'}^{\vee} - \cal
N\mit_{D_{F}/F}^{\vee}) \in \bf K\mit_{0}(D_{F})$$
Alors le diagramme suivant commute dans $HoSp$
$$\xymatrix{
\bf G\mit(F) \ar[d]_{\alpha_{j}.d_{F}^{*}}\ar[r]^{j_{*}} & \bf G\mit(F')
\ar[d]^{d_{F'}^{*}}\\
\bf G\mit(D_{F}) \ar[r]_{Dj_{*}} & \bf G\mit(D_{F'}) }$$
\end{prop}

\underline{\bf Preuve:} \rm Une application de la d\'eformation vers le
c\^one normal (~\cite[$4.1$]{g}~), nous permet de ne consid\'erer que le cas o\`u
$F'=\bf
P\mit(V\oplus 1)$ est le compl\'et\'e projectif d'un fibr\'e vectoriel $V$
sur $F$, et $j$ est la section nulle $e : F \hookrightarrow V
\hookrightarrow \bf P\mit(V\oplus 1)$. Notons $P=\bf P\mit(V\oplus
1)$, et
$$p : P \longrightarrow F$$
$$Dp : D_{P} \longrightarrow D_{F}$$
$$De : D_{F} \longrightarrow D_{P}$$
les morphismes induits.
Remarquons qu'il suffit de d\'emontrer que
$$De_{*}(\alpha_{e})=d_{P}^{*}(e_{*}(1)) \in \bf K\mit_{0}(D_{P})$$
En effet, si la formule ci-dessus est vraie, on a les \'egalit\'es
suivantes dans $HoSp$
$$\begin{array}{cl}
De_{*}(\alpha_{e}.d_{F}^{*}) & = De_{*}(\alpha_{e}.De^{*}Dp^{*}d_{F}^{*}) \\
& = De_{*}(\alpha_{e}).Dp^{*}d_{F}^{*} \\
& = De_{*}(\alpha_{e}).d_{P}^{*}p^{*} \\
& = d_{P}^{*}(e_{*}(1)).d_{P}^{*}p^{*}\\
& = d_{P}^{*}(e_{*}(1).p^{*}) \\
& = d_{P}^{*}(e_{*}(1.e^{*}p^{*})) \\
& = d_{P}^{*}(e_{*})
\end{array}$$

Soit $i : D_{P} \longrightarrow D_{F}\times_{F}P$ le morphisme
canonique. C'est l'immersion canonique des points fixes du fibr\'e
projectif $D_{F}\times_{F}P$, pour l'action lin\'eaire du groupe de
type multiplicatif universel $\cal D\mit_{F}$ (~\ref{d2.5}~). Ainsi,
localement sur
$(D_{F})_{li}$, si $W=V\oplus 1\simeq \bigoplus_{\chi \in
\cal X\mit^{*}_{F}(F)}W^{(\chi)}$, on a
$$D_{P}\simeq \coprod_{\chi \in \cal X\mit^{*}_{F}(F)}P^{\chi}$$
o\`u $P^{\chi}$ est le fibr\'e projectif associ\'e au fibr\'e vectoriel
$W^{(\chi)}$.

\begin{lem}
On a $\alpha_{e}=\lambda_{-1}(De^{*}\cal N\mit^{\vee}_{i})$, o\`u
$\cal N\mit_{i}^{\vee}$ est le fibr\'e conormal de l'immersion r\'eguli\`ere
$i : D_{P} \hookrightarrow D_{F}\times_{F}P$.
\end{lem}

\underline{\bf Preuve:} \rm Consid\'erons le diagramme commutatif
suivant
$$\xymatrix{
D_{P} \ar[r]^-{i} \ar[dr]_{Dp} & D_{F}\times_{F}P \ar[r]^-{a}
\ar[d]^{q} & P \ar[d]^{p} \\
& D_{F} \ar[r]_{d_{F}} & F}$$
Notons $f : D_{F} \hookrightarrow D_{F}\times_{F}P$ la section
induite par $e$. Alors
$$\cal N\mit_{d_{P}}^{\vee}=\cal N\mit_{i}^{\vee}+Di^{*}\cal
N\mit_{a}^{\vee}$$
donc
$$De^{*}(\cal N\mit_{d_{P}}^{\vee})=De^{*}Di^{*}\cal N\mit_{a}^{\vee}+
De^{*}\cal N\mit_{i}^{\vee}
=Df^{*}\cal N\mit_{a}^{\vee}+De^{*}\cal N\mit_{i}^{\vee}$$
et
$$q^{*}\cal N\mit_{d_{F}}^{\vee}=\cal N\mit_{a}^{\vee}$$
donc
$$Df^{*}\cal N\mit_{a}^{\vee}=Df^{*}q^{*}\cal N\mit_{d_{F}}^{\vee}=\cal
N\mit_{d_{F}}^{\vee}$$
Ainsi
$$De^{*}(\cal N\mit_{d_{P}}^{\vee}-Dp^{*}\cal N\mit_{d_{F}}^{\vee})=
Df^{*}\cal N\mit_{a}^{\vee}+De^{*}\cal N\mit_{i}^{\vee}-
\cal N\mit^{\vee}_{d_{F}}=De^{*}\cal N\mit_{i}^{\vee}$$
$\Box$\\

On a donc
$$De_{*}(\alpha_{e})=De_{*}(De^{*}\lambda_{-1}(\cal N\mit_{i}^{\vee}))
=\lambda_{-1}(\cal N\mit_{i}^{\vee}).De_{*}(1)$$
Or comme la section $De_{*}$ envoi $D_{F}$ dans $P^{1}=\bf
P\mit(W^{(1)})$, $De_{*}(\alpha_{e})$ est support\'e par
$P^{1}$.

Notons $\cal N\mit^{\vee}_{1}$ le fibr\'e conormal de $P^{1}$ dans $P$.
Il faut donc montrer que
$$De_{*}(1).\lambda_{-1}(\cal N\mit_{1}^{\vee})=d_{P}^{*}(e_{*}(1))$$
dans $\bf K\mit_{0}(P^{1})$.

Soit
$$0 \longrightarrow E_{P} \longrightarrow p^{*}V\oplus 1
\longrightarrow \cal O\mit_{P}(1)\longrightarrow 0$$
$$0 \longrightarrow E_{P^{1}} \longrightarrow Dp_{1}^{*}W^{(1)}
\longrightarrow \cal O\mit_{P^{1}}(1)\longrightarrow 0$$
les suites exactes canoniques sur $P$ et $P^{1}$, o\`u $Dp_{1} : P^{1}
\longrightarrow D_{F}$ est la projection. Alors, on sait que
(~\cite[$4.3$]{fl}~)
$$e_{*}(1)=\lambda_{-1}(E_{P})$$
$$De_{*}(1)=\lambda_{-1}(E_{P^{1}})$$
Notons $W\simeq W^{(1)}\oplus W^{(\neq 1)}$, o\`u $W^{(1)}$ est le
sous-fibr\'e o\`u $\cal D\mit_{F}$ op\`ere trivialement. Remarquons que par
d\'efinition, $W^{(1)}\simeq V^{(1)}\oplus 1$.
Ainsi, $\cal N\mit_{1}^{\vee}=Dp_{1}^{*}W^{\neq 1}$,
et donc,
$$\begin{array}{cl}
De_{*}(1).\lambda_{-1}(\cal N\mit_{1}^{\vee}) &
= \lambda_{-1}(E_{P^{1}}).\lambda_{-1}(Dp_{1}^{*}W^{(\neq 1)}) \\
& = \lambda_{-1}(Dp_{1}^{*}(V^{(1)})+1-\cal O\mit_{P^{1}}(1)+
Dp_{1}^{*}W^{(\neq 1)}) \\
& = \lambda_{-1}(Dp_{1}^{*}(V+1)-\cal O\mit_{P^{1}}(1))
\end{array}$$
Par ailleurs
$$\begin{array}{cl}
d_{P}^{*}(e_{*}(1)) & = i^{*}(f_{*}(1)) \\
& = i_{1}^{*}(f_{*}(1))
\end{array}$$
o\`u $i_{1} : P^{1} \hookrightarrow P$ est induit par $i$. Donc
$$\begin{array}{cl}
d_{P}^{*}(e_{*}(1)) & = i_{1}^{*}\lambda_{-1}(a^{*}(E_{P})) \\
& = \lambda_{-1}(Dp_{1}^{*}V + 1 -\cal O\mit_{P^{1}}(1))
\end{array}$$
Ce qui prouve que
$$De_{*}(1).\lambda_{-1}(\cal N\mit_{1}^{\vee})=d_{P}^{*}(e_{*}(1))$$
et donc la proposition. $\Box$\\

\begin{thm}\label{th3.1}{(~Lefschetz-Riemann-Roch~)}
Soit $f : F \longrightarrow F'$ un morphisme localement d'intersection
compl\`ete et fortement projectif de
champs alg\'ebriques bien ramifi\'es (~\ref{d3.4}~). Alors le diagramme
suivant commute
$$\xymatrix{
\bf G\mit_{*}(F) \ar[d]_{f_{*}}\ar[r]^-{\chi_{F}.\alpha_{F}^{-1}} &
\underline{\bf
G}^{\chi}_{*}(F) \ar[d]_{f_{*}}\\
\bf G\mit_{*}(F') \ar[r]_-{\chi_{F'}.\alpha_{F'}^{-1}} & \underline{\bf
G}^{\chi}_{*}(F')}$$
\end{thm}

\underline{\bf Preuve:} \rm En factorisant $f$ en une immersion
ferm\'ee r\'eguli\`ere, suivie de la projection d'un fibr\'e projectif associ\'e \`a un
fibr\'e vectoriel, et \`a l'aide du lemme suivant, on peut ne traiter que ces
deux cas.

\begin{lem}\label{l3.5}
Soit $F$ un champ bien ramifi\'e, et $P \longrightarrow F$ un
fibr\'e projectif associ\'e \`a un fibr\'e vectoriel $V$ sur $F$. Alors
$P$ est encore bien ramifi\'e.
\end{lem}

\underline{\bf Preuve:} \rm On a d\'ej\`a vu que le morphisme naturel
$$D_{V} \longrightarrow D_{F}\times_{F}V$$
faisait de $D_{V}$ un sous-fibr\'e vectoriel de $D_{F}\times_{F}V$ sur
$D_{F}$. De plus, \\
$D_{P} \hookrightarrow D_{F}\times_{F}P$
est une immersion r\'eguli\`ere. Ainsi, la projection
$D_{P}~\longrightarrow~P$ se
factorise par $D_{P}~\hookrightarrow~D_{F}\times_{F}P~\longrightarrow~P$,
qui comme $P~\longrightarrow~F$ est lisse,
est encore un morphisme localement d'intersection compl\`ete.

Pour montrer que $\alpha_{P}$ est inversible, on va montrer que
$$Dp^{*}(\alpha_{F}^{-1}).\alpha_{P} \in \underline{\bf K}^{\chi}_{0}(P)$$
est inversible. Par le lemme \ref{l3.3}, ceci est local sur
$(D_{F})_{li}$, ce qui nous permet de faire les calculs au-dessus d'un
morphisme lisse et surjectif $(s,D) : X \longrightarrow D_{F}$,
qui correspond \`a un sous-groupe diagonalisable
$$D \hookrightarrow Aut_{X}(s)$$
Le fibr\'e vectoriel $V_{X}$ \'etant muni d'une action de $D$, le fibr\'e
projectif $P_{X}$ aussi. Alors $D_{P}\times_{D_{F}}X$ est repr\'esent\'e
par $P_{X}^{D}$, le sous-espace de $P_{X}$ des points fixes de $D$.
Or la restriction de $Dp^{*}(\alpha_{F}^{-1}).\alpha_{P}$ dans
$\bf K\mit_{0}(P^{D}_{X})[[M(D)]]_{\bf Q}$, est \'egale \`a
$$\prod_{t \in M(D)}\lambda_{-t}(\cal N\mit_{D}^{\vee})$$
o\`u $M(D)$ est le groupe des caract\`eres de $D$, et
$\cal N\mit_{D}^{\vee}$ le fibr\'e conormal de $P_{X}^{D}$ dans $P_{X}$.
C'est donc un \'el\'ement inversible. $\Box$\\

\underline{Cas d'une immersion ferm\'ee:} \\

Soit $j : F \hookrightarrow F'$ une immersion ferm\'ee, r\'eguli\`ere. Alors le
diagramme suivant est cart\'esien
$$\xymatrix{
F \ar[r]^{j} & F' \\
D_{F} \ar[u] \ar[r]^{Dj} & D_{F'} \ar[u]}$$
Notons $\lambda_{F}=\lambda_{-1}(\cal N\mit_{d_{F}}^{\vee})$
(~resp. $\lambda_{F'}=\lambda_{-1}(\cal N\mit_{d_{F}}^{\vee})$~), o\`u $\cal
N\mit_{d_{F}}^{\vee}$
est le fibr\'e conormal de $d_{F} : D_{F} \longrightarrow F$
(~resp. $d_{F'} : D_{F'} \longrightarrow D_{F}$~). Alors, la formule
d'auto-intersection (~\ref{p3.5}~) implique, que pour tout $x \in \bf
G\mit_{*}(F)$
$$Dj_{*}(\lambda_{-1}(Dj^{*}\cal N\mit_{d_{F'}}^{\vee}-\cal
N\mit_{d_{F}}^{\vee}).d_{F}^{*}(x))
=d_{F'}^{*}(j_{*}(x))$$
En composant avec le morphisme $\rho$ (~\ref{th2.4}~), qui commute avec
les images directes et r\'eciproques pour $j$, et les produits sur $D_{F}$,
on obtient
$$j_{*}(\rho(\lambda_{-1}(Dj^{*}\cal N\mit_{d_{F'}}^{\vee}-\cal
N\mit_{d_{F}}^{\vee})).\rho(d_{F}^{*}(x)))=\rho(d_{F'}^{*}(j_{*}(x))$$
Or, $\rho(d_{F}^{*}(x))=\chi_{F}(x)$, et
$\rho(d_{F'}^{*}(j_{*}(x))=\chi_{F'}(x)$, par d\'efinition. D'autre part
$$\begin{array}{cl}
\rho(\lambda_{-1}(Dj^{*}\cal N\mit_{d_{F'}}^{\vee}-\cal
N\mit_{d_{F}}^{\vee})).\alpha_{F} &
 = \rho(\lambda_{-1}(Dj^{*}\cal N\mit_{d_{F'}}^{\vee}-\cal
N\mit_{d_{F}}^{\vee})).\rho(\lambda_{F}) \\
& = \rho((\lambda_{-1}(Dj^{*}\cal N\mit_{d_{F'}}^{\vee}-\cal
N\mit_{d_{F}}^{\vee})).\lambda_{-1}(\cal N\mit_{d_{F}}^{\vee})) \\
& = \rho(\lambda_{-1}(Dj^{*}\cal N\mit_{d_{F'}}^{\vee})) \\
& = j^{*}\rho(\lambda_{F'}) \\
& = j^{*}\alpha_{F'}
\end{array}$$
Ce qui montre que
$$\rho(\lambda_{-1}(Dj^{*}\cal N\mit_{d_{F'}}^{\vee}-\cal
N\mit_{d_{F}}^{\vee})) = j^{*}\alpha_{F'}.\alpha_{F}^{-1}$$
et donc
$$j_{*}(j^{*}\alpha_{F'}.\alpha{F}^{-1}.\chi_{F}(x))=
\chi_{F'}(j_{*}(x))$$
ce qui, par la formule de projection est \'equivalent \`a
$$j_{*}(\alpha_{F}^{-1}.\chi_{F}(x))=\alpha_{F'}^{-1}.j_{*}(x)$$

Remarquons que la d\'emonstration pr\'ec\'edente montre que l'on a une
\'egalit\'e dans $HoSp$
$$j_{*}(\alpha_{F}^{-1}.\chi_{F})=\alpha_{F'}^{-1}.j_{*}$$

\underline{Cas d'un fibr\'e projectif:} \\

Soit $V \longrightarrow F$ un fibr\'e vectoriel de rang $r+1$, et
$p : P=\bf P\mit(V) \longrightarrow F$ le fibr\'e projectif associ\'e.
Remarquons tout d'abord qu'il suffit de montrer la  formule pour les \'el\'ements
$x^{i}=\cal O\mit_{P}(i) \in \bf G\mit_{0}(P)$.
En effet, d'apr\`es \ref{p2.1}, tout \'el\'ement $y \in \bf
G\mit_{*}(P)$ s'\'ecrit de fa\c{c}on unique
$$y = \sum_{i=0}^{i=r}p^{*}(a_{i}).x^{i}$$
$a_{i} \in \bf G\mit_{*}(F)$.
Supposons que la formule soit d\'emontr\'ee pour les \'el\'ements $x^{i}$, alors
$$\begin{array}{cl}
p_{*}(\alpha_{P}^{-1}.\chi_{P}(y)) &
= \sum_{i}p_{*}(\alpha_{P}^{-1}.p^{*}(\chi_{F}(a_{i})).x^{i}) \\
& = \sum_{i}p_{*}(\alpha_{P}^{-1}.x^{i}).\chi_{F}(a_{i}) \\
& = \sum_{i}\alpha_{F}^{-1}.\chi_{F}(p_{*}(x^{i})).\chi_{F}(a_{i}) \\
& = \alpha_{F}^{-1}.\sum_{i}\chi_{F}(p_{*}(x^{i}).a_{i}) \\
& = \alpha_{F}^{-1}.\sum_{i}\chi_{F}(p_{*}(y))
\end{array}$$

Commen\c{c}ons par montrer que l'on peut supposer que $V$ est une somme
directe de fibr\'es inversibles.

\begin{lem}\label{l3.6}
Soit $p : P \longrightarrow F$ un fibr\'e projectif, avec $F$ un champ
lisse et bien ramifi\'e. Alors le morphisme
$$p^{*} : \underline{\bf G}^{\chi}_{*}(F) \longrightarrow
\underline{\bf G}^{\chi}_{*}(P)$$
est injectif.
\end{lem}

\underline{\bf Preuve:} \rm Comme $D_{P} \longrightarrow D_{F}$ est
un morphisme propre et lisse, il existe
des images r\'eciproques
$$p^{*} : \underline{\bf G}^{\chi}(F) \longrightarrow
\underline{\bf G}^{\chi}(P)$$
et des images directes
$$p^{*} : \underline{\bf K}^{\chi}(P) \longrightarrow
\underline{\bf K}^{\chi}(F)$$
On peut alors appliquer la formule de
projection \ref{p2.1} pour $x \in \underline{\bf G}^{\chi}_{*}(F)$
$$p_{*}p^{*}(x)=x.p_{*}(1)$$
Or, $p_{*}(1)$ est inversible dans $\underline{\bf K}^{\chi}(F)$. En
effet, ceci est local sur $(D_{F})_{li}$, et localement sur
$(D_{F})_{li}$, $p : P \longrightarrow F$ est la projection d'une
r\'eunion disjointe de fibr\'es projectifs (non vides). $\Box$\\

Soit alors $f : G \longrightarrow
F$ un morphisme projectif et lisse, tel que $g^{*}(G)$ admette une
filtration par de sous-fibr\'es de quotients successifs de rang $1$. On
sait qu'un tel morphisme existe, et de plus, d'apr\`es le
lemme \ref{l3.6}, le morphisme
$$f^{*} : \underline{\bf G}^{\chi}_{0}(F) \longrightarrow
\underline{\bf G}^{\chi}_{0}(G)$$
est injectif.

Notons  $g : G_{P}:=G\times_{F}P \longrightarrow P$ et
$q : G_{P} \longrightarrow P$ les morphismes induits, et $x_{G}=q^{*}(x)$.
Comme $G \longrightarrow F$ est lisse, on a
$$q^{*}(\alpha_{P}^{-1}.p^{*}\alpha_{F})=
\alpha_{G_{P}}^{-1}.q^{*}\alpha_{G}$$
Avec ces notations, si la formule est d\'emontr\'ee pour $q$ et $x_{G}^{i}$, on a
$$\begin{array}{cl}
f^{*}(p_{*}(\lambda_{P}^{-1}.p^{*}\lambda_{F}.\chi_{P}(x^{i}))) &
 = q_{*}(g^{*}(\lambda_{P}.p^{*}\lambda_{F}.\chi_{P}(x^{i})) \\
& = q_{*}(\lambda_{G_{P}}.q^{*}\lambda_{G}.g^{*}(\chi_{P}(x^{i}))) \\
& = q_{*}(\lambda_{G_{P}}.q^{*}\lambda_{G}.\chi_{G}(g^{*}(x)^{i})) \\
& = q_{*}(\lambda_{G_{P}}.q^{*}\lambda_{G}.\chi_{G}(x_{G}^{i})) \\
& = \chi_{G}(q_{*}(x_{G}^{i})) \\
& = \chi_{G}(q_{*}(g^{*}(x^{i}))) \\
& = \chi_{G}(f^{*}(p_{*}(x^{i}))) \\
& = f^{*}(\chi_{F,f}(p_{*}(x^{i})))
\end{array}$$
Or comme $f^{*}$ est injectif, on en d\'eduit la formule pour $p$ et
$x^{i}$. \\

On est donc ramen\'e au cas o\`u $V$ est sous forme triangulaire. Il
existe alors un morphisme lisse
$$r : W \longrightarrow F$$
qui est une composition de torseurs sous des fibr\'es vectoriel, tel
que $r^{*}(V)$ soit une somme directe de fibr\'es inversibles. En
utilisant \ref{p2.1}, et un calcul analogue \`a celui fait ci-dessus,
on peut supposer que $V$ est diagonalisable.\\

Soit $V\simeq V_{1}\oplus \cal L$ une d\'ecomposition de $V$, avec
$\cal L$ de rang $1$. On note $j : P_{1}:=\bf P\mit(V_{1})
\hookrightarrow P$ l'immersion canonique. Comme tout \'el\'ement
$x \in \bf G\mit_{0}(P)$ s'\'ecrit comme
$$x=j_{*}(y)+p^{*}(z)$$
le cas d'une immersion ferm\'ee permet de nous ramener au cas o\`u
\mbox{$x=p^{*}(z)$}, et par la formule de projection au cas o\`u $x=1$. \\

On va utiliser l'argument de d\'eformation de \cite[App. $3$]{bfm}, pour
montrer que
$$p_{*}(\alpha_{P}^{-1}.p^{*}(\alpha_{F}))=1$$
Comme $V$ est diagonalisable, la construction de l'espace de
d\'eformation de \cite[App. $3$]{bfm}, se globalise sur $F_{li}$, et donne un
sous-champ
alg\'ebrique ferm\'e
$$J : \cal X\mit \hookrightarrow P\times_{F}P\times_{S}\bf A\mit^{1}$$
tel que la projection induite
$$\pi :\cal X\mit \longrightarrow F\times_{S}\bf A\mit^{1}$$
soit un morphisme plat, et que
$$J_{1} : \cal X\mit_{1}:=\pi^{-1}(F\times\{1\}) \hookrightarrow P\times_{F}P$$
soit l'immersion diagonale, et
$$J_{0} : \cal X\mit_{0}:=\pi^{-1}(F\times\{0\})=
\cup_{a+b=r}P^{a}\times_{F}P^{b}$$
o\`u $P^{n}$ d\'esigne un sous-fibr\'e projectif de $P$ de rang $n$.

Pour un champ bien ramifi\'e, notons
$$\lambda_{F}=\alpha_{F}^{-1}.\chi_{F} : \bf G\mit_{0}(F)
\longrightarrow \underline{\bf G}^{\chi}_{0}(F)$$
Pour tout sous-champ ferm\'e $\cal Y\mit$ de $P\times_{F}P\times_{S} \bf
A\mit^{1}$ de
compl\'ementaire $U$, on notera $\lambda_{\cal Y\mit}$
le morphisme induit sur la fibre des suites exactes suivantes
$$\xymatrix{
\bf G\mit(\cal Y\mit) \ar[r] \ar[d]_{\lambda_{\cal Y\mit}} &
\bf K\mit(P\times_{F}P\times \bf A\mit^{1}) \ar[r]
\ar[d]_{\lambda_{P\times_{F}P\times \bf A\mit^{1}}} & \bf K\mit(U)
\ar[d]^{\lambda_{U}} \\
\underline{\bf K}^{\chi}(\cal Y\mit) \ar[r] &
\underline{\bf K}^{\chi}(P\times_{F}P\times \bf
A\mit^{1}) \ar[r] & \underline{\bf K}^{\chi}(U) }$$
Remarquons que le th\'eor\`eme dans le cas d'une immersion ferm\'ee r\'eguli\`ere
de champs bien ramifi\'es permet de conclure que les deux d\'efinitions
ci-dessus de $\lambda_{\cal Y\mit}$ coincident si $\cal Y\mit
\hookrightarrow P\times_{F}P\times_{S}\bf A\mit^{1}$ est une
immersion r\'eguli\`ere, et $\cal Y$ est bien ramifi\'e
(~exactement comme le corollaire \cite[$3.7$, $(ii)$]{g} se d\'eduit
du th\'eor\`eme \cite[$3.1$]{g}~). En particulier pour
$$\cal X\mit_{1}\simeq P \hookrightarrow P\times_{F}P\times_{S}\bf
A\mit^{1}$$

D\'efinissons alors, pour $\cal Y$ un sous-champ ferm\'e de
$P\times_{F}\times_{S}\bf A\mit^{1}$
$$\chi(\cal Y\mit):=q_{*}(\lambda_{\cal Y\mit}(1))$$
o\`u $q : \cal Y\mit \longrightarrow F$ est la projection induite.

Si on montre que
$\chi(\cal X\mit_{t})$ est ind\'ependant de $t \in \bf A\mit^{1}$, alors
$$p_{*}(\alpha_{P}^{-1}.p^{*}\alpha_{F})=\chi(\cal X\mit_{1})
=\chi(\cal X\mit_{0})$$
De la m\^eme fa\c{c}on que dans
\cite[App. $3$]{bfm}, on conclut alors par r\'ecurrence sur le rang
de $V$ que
$$p_{*}(\alpha_{P}^{-1}.p^{*}\alpha_{F})=1$$

Il nous reste donc \`a montrer que $\chi(\cal X\mit_{t})$ est
ind\'ependant de $t \in \bf A\mit^{1}$.

Soit $L:=\lambda_{\cal X}(1) \in \underline{\bf G}^{\chi}_{0}(\cal
X\mit)$. Comme dans le diagramme cart\'esien suivant
$$\xymatrix{
\cal X\mit_{t} \ar[r]^{J_{t}} \ar[d]_{q_{t}} &
\cal X\mit \ar[d]^{\pi} \\
F \ar[r]_-{i_{t}} & F\times_{S}\bf A\mit^{1}}$$
le morphisme $\pi$ est plat, les fibr\'e conormaux de $i_{t}$ et de
$J_{t}$ sont triviaux. Ainsi, on a
$$(q_{t})_{*}\circ J_{t}^{*}(L)=i_{t}^{*}\circ \pi_{*}(L)$$
Or $i_{t}^{*}=i_{t'}^{*}$ pout tout $t$ et $t' \in \bf A\mit^{1}$,
ceci montre que $(q_{t})_{*}\circ J_{t}^{*}(L)$ est ind\'ependant de $t$.
Or $J_{t}^{*}(L)=\lambda_{\cal X\mit_{t}}(1)$, et donc
$\chi(\cal X\mit_{t})=(q_{t})_{*}\circ J_{t}^{*}(L)$ est ind\'ependant
de $t$. $\Box$\\

\underline{Remarque:} La formule pr\'ec\'edente reste de toute \'evidence
valable pour un morphisme propre repr\'esentable
de champs bien ramifi\'es. Cependant, faute de
r\'esultats concernant les enveloppes (~ou quasi-enveloppes~) de Chow
d'un morphisme propre, nous ne savons pas la d\'emontrer dans cette
g\'en\'eralit\'e.

Cela pose aussi la question de savoir quand un morphisme propre
repr\'esentable est fortement projectif. On sait que cela est vrai
lorsque les champs satisfont deux hypoth\`eses suppl\'ementaires. A savoir
l'existence d'un fibr\'e ample, et le fait que tout faisceau coh\'erent
est quotient d'un fibr\'e vectoriel. Ainsi, on a le corollaire suivant.

\begin{cor}\label{c3.1}
Soit $HoChAlg''(S)$ la sous-cat\'egorie de $HoChAlg'(S)$ des champs alg\'ebriques
bien ramifi\'es et v\'erifiant les hypoth\`eses suivantes~:
\begin{itemize}
\item
Il existe un fibr\'e inversible ample sur $F$. C'est \`a dire $\cal L$,
tel que pour tout faisceau coh\'erent $\cal F$ sur $F$, il existe un
entier $n$, tel que
$H^{i}(F,\cal F\mit\otimes \cal L\mit^{\otimes n})=0$ pour $i>n$.
\item
Tout faisceau coh\'erent est quotient d'un fibr\'e vectoriel.
\end{itemize}
Alors, pour tout morphisme propre, repr\'esentable et
localement d'intersection compl\`ete de $HoChAlg''(S)$
$f : F \longrightarrow F'$, le carr\'e suivant commute
$$\xymatrix{
\bf G\mit_{*}(F) \ar[d]_{f_{*}}\ar[r]^-{\chi_{F}.\alpha_{F}^{-1}} &
\underline{\bf
G}^{\chi}_{*}(F) \ar[d]^{f_{*}}\\
\bf G\mit_{*}(F') \ar[r]_-{\chi_{F'}.\alpha_{F'}^{-1}} & \underline{\bf
G}^{\chi}_{*}(F')}$$
\end{cor}

\begin{cor}\label{c3.2}
Soit $S=Spec k$, et $F$ et $F'$ deux champs alg\'ebriques, qui sont
\'equivalents \`a des champs quotients de sch\'emas quasi-projectifs par des
sch\'emas en groupes r\'eductifs. Alors, pour tout morphisme
propre, repr\'esentable et localement d'intersection compl\`ete
$f : F \longrightarrow F'$, le carr\'e suivant
commute
$$\xymatrix{
\bf G\mit_{*}(F) \ar[d]_{f_{*}}\ar[r]^-{\chi_{F}.\alpha_{F}^{-1}} &
\underline{\bf
G}^{\chi}_{*}(F) \ar[d]^{f_{*}}\\
\bf G\mit_{*}(F') \ar[r]_-{\chi_{F'}.\alpha_{F'}^{-1}} & \underline{\bf
G}^{\chi}_{*}(F')}$$
\end{cor}

\underline{\bf Preuve:} \rm On sait que ces deux champs v\'erifient les
hypoth\`eses du corollaire pr\'ec\'edent (~\cite{th2}~). $\Box$\\

\end{subsubsection}

\begin{subsubsection}{Formule de Grothendieck-Riemann-Roch}
\hspace{5mm}Pour tout ce paragraphe, $S=Spec k$ est le spectre d'un
corps, et $\cal H$ d\'esigne une th\'eorie cohomologique avec images
directes.

\begin{prop}\label{p3.6}
Soit $f : F \longrightarrow F'$ un morphisme fortement projectif de
champs alg\'ebriques lisses. Alors, le diagramme suivant commute
$$\xymatrix{
\underline{\bf G}_{*}(F) \ar[r]^{f_{*}} \ar[d]_-{Td(T_{F}).Ch} &
\underline{\bf G}_{*}(F') \ar[d]^-{Td(T_{F'}).Ch} \\
H^{\bullet}(F,*) \ar[r]_{f_{*}} & H^{\bullet}(F',*)}$$
o\`u $T_{F}\in \underline{\bf K}_{0}(F)$ et
$T_{F'}\in \underline{\bf K}_{0}(F')$ sont les fibr\'es tangents virtuels de
$F$ et $F'$.
\end{prop}

\underline{\bf Preuve:} \rm C'est mot pour mot la m\^eme que celle pour
des sch\'emas (~\cite[$4.1$]{g}~). Le cas d'une immersion ferm\'ee provient
de la d\'eformation vers le c\^one normal et des propri\'et\'es des
classes caract\'eristiques, et celui d'un fibr\'e projectif,
de la formule du projectif \ref{p2.1} en $G$-cohomologie. $\Box$\\

\begin{prop}\label{p3.7}
Soit $f : F \longrightarrow F'$ un morphisme fortement projectif de
champ alg\'ebriques lisses et bien ramifi\'es. Alors le diagramme suivant
commute
$$\xymatrix{
\bf H\mit^{*}(D_{F},\underline{K}\otimes\cal A\mit_{F})
\ar[rr]^-{Df_{*}\otimes Ind_{f}} \ar[d]_-{Td(T_{D_{F}}).Ch\otimes
Id} & &
\bf H\mit^{*}(D_{F'},\underline{K}\otimes \cal A\mit_{F'})
\ar[d]^-{Td(T_{D_{F'}}).Ch\otimes Id} \\
\bf H\mit^{*}(D_{F},\cal H\mit\otimes \cal A\mit_{F})
\ar[rr]_{Df_{*}\otimes Ind_{f}} & &
\bf H\mit^{*}(D_{F'},\cal H\mit\otimes \cal A\mit_{F'})}$$
\end{prop}

\underline{\bf Preuve:} \rm En effet, si $F$ et $F'$ sont lisses et
bien ramifi\'es, $D_{F}$ et $D_{F'}$ sont lisses. De plus, le morphisme
induit $Df : D_{F} \longrightarrow D_{F'}$ est encore fortement
projectif. En effet, soit
$$\xymatrix{F \ar[r]^{j} &  P \ar[r]^{p} &  F'}$$
une factorisation de $f$, avec $j$ une immersion ferm\'ee et $p$ la
projection d'un fibr\'e projectif. Alors, $Df$ se factorise par
$$Df :  \xymatrix{D_{F} \ar[r]^{Dj} &  D_{P} \ar[r]^{Dp} &  D_{F}}$$
Or on sait que $Dj$ est une immersion ferm\'ee, et de plus $Dp$ se
factorise par
$$Dp : \xymatrix{D_{P} \ar[r]^-{i} &  D_{F'}\times_{F'}P \ar[r]^-{q} &
D_{F'}}$$
o\`u $i$ est une immersion ferm\'ee.

On peut alors refaire la preuve de la proposition pr\'ec\'edente avec
$\underline{K}\otimes \cal A\mit$ \`a la place de
$\underline{K}$, et
$\cal H\mit \otimes \cal A$ \`a la place de $\cal H$. $\Box$\\

Pour le corollaire suivant, on a not\'e $H^{\chi,f}_{\bullet}(F',*)$ la composante de $H^{\chi}_{\bullet}(F',*)$
support\'ee par $D_{F',f}$. 

\begin{cor}\label{c3.3}{(~Grothendieck-Riemann-Roch~)}
Soit $f : F \longrightarrow F'$ un morphisme fortement projectif de
champs alg\'ebriques lisses et bien ramifi\'es. Alors le diagramme
suivant commute
$$\xymatrix{
\bf G\mit_{*}(F) \ar[r]^-{\tau_{F}^{\chi}} \ar[d]_-{f_{*}} &
H^{\chi}_{\bullet}(F,*) \ar[d]^-{f_{*}} \\
\bf G\mit_{*}(F') \ar[r]_-{\tau_{F'}^{\chi}} &
H^{\chi,f}_{\bullet}(F',*)}$$
\end{cor}

\underline{\bf Preuve:} \rm Il suffit de composer la formule \ref{th3.1} et
le corollaire pr\'ec\'edent. $\Box$\\

\end{subsubsection}

\begin{subsubsection}{Cas des champs de Deligne-Mumford : Cohomologie \`a
coefficients dans les repr\'esentations}\label{s3}
\hspace{5mm}
Dans le cas des champs de Deligne-Mumford, le th\'eor\`eme de d\'evissage
\ref{th2.3}
est plus pr\'ecis, et permet d'\'etendre les formules de Riemann-Roch au
cas des morphismes non repr\'esentables.

Nous revenons un moment au cas o\`u $S$ est un sch\'ema r\'egulier
universellement japonais. \\

Pour simplifier les \'enonc\'es et les notations, on supposera que $S$
contient les racines de l'unit\'e. On fixera alors un plongement
$$\mu_{\infty}^{t}(S) \hookrightarrow \mu_{\infty}(\bf C\mit)$$
qui nous permettra d'identifier par la suite le faisceau
$\Lambda$ avec un sous-faisceau constant $K
\hookrightarrow \bf Q\mit^{ab}$. Comme $K$ est un corps, on a
$$\underline{\bf K}^{rep}_{*}(F)\simeq \underline{\bf
K}_{*}(I_{F}^{t})\otimes K \hookrightarrow
\underline{\bf K}_{*}(I_{F}^{t})\otimes \bf Q\mit^{ab}$$

Nous faisons remarquer au lecteur que cette hypoth\`ese suppl\'ementaire
n'est pas tr\`es restrictive. En effet, si $F$ est un champ de
Deligne-Mumford, il existe un changement de base galoisien de $S$, tel
que le faisceau $\mu_{m}$ soit constant sur $S_{et}$, avec $m$ un
multiple de tous les ordres de ramifications de $F$. On peut alors
redescendre sur $S$ par le proc\'ed\'e de descente galoisienne
\ref{c2.2} $(1)$. \\

Rappelons que l'on a construit un morphisme
$$F : \bf K\mit(I_{F}^{t}) \longrightarrow \underline{\bf
K}^{rep}(F)$$
Nous noterons \'egalement, pour $\zeta \in \mu_{\infty}(S)$
$$F_{\zeta} : \bf Vect\mit(I_{F}^{t}) \longrightarrow
\bf Vect\mit(I_{F}^{t})$$
le foncteur qui a un fibr\'e vectoriel $V$ sur $I_{F}^{t}$, associe le sous-fibr\'e
vectoriel $V^{(\zeta)}$ sur lequel la section canonique $\cal S$
(~\ref{dev}~) op\`ere par multiplication par $\zeta$.

\begin{df}\label{d3.6}
Soit $F$ un champ de Deligne-Mumford.
\begin{enumerate}
\item
Si $F$ est r\'egulier, sa classe de ramification est d\'efinie par
$$\alpha^{rep}_{F}:=F(\lambda_{\pi_{F}}) \in \underline{\bf
K}^{rep}(F)$$
o\`u $\pi_{F} : I_{F}^{t} \longrightarrow F$ est la projection, et
$F : \bf K\mit_{0}(I_{F}^{t}) \longrightarrow
\underline{\bf K}^{rep}_{0}(F)$
est d\'efini dans \ref{th2.3}.
\item
Soit $S=Spec k$, et $\cal H$ est une des deux th\'eories cohomologiques
avec images directes cit\'ee dans \ref{s4}.
On d\'efinit la cohomologie et
l'homologie \`a coefficients dans les repr\'esentations par
\begin{itemize}
\item
$$H^{p}_{rep}(F,q):=\pi_{dq-p}\bf H\mit((I_{F}^{t})_{li},\cal
H\mit^{q}\otimes K)$$
$$H_{p}^{rep}(F,q):=\pi_{dq-p}\bf H\mit((I_{F}^{t})_{li},\cal
H\mit'_{q}\otimes K)$$
si $\cal H$ est la th\'eorie de Gersten.
\item
$$H^{p}_{rep}(F,q):=\pi_{dq-p}\bf H\mit((I_{F}^{t})_{li},\cal
H\mit^{q})$$
$$H_{p}^{rep}(F,q):=\pi_{dq-p}\bf H\mit((I_{F}^{t})_{li},\cal
H\mit'_{q})$$
si $k$ est de caract\'eristique nulle, et $\cal H$ est la cohomologie
de De Rham
\end{itemize}
On notera aussi
$$H^{\bullet}_{rep}(F,*):=\bigoplus_{p,q}H^{p}_{rep}(F,q)$$
$$H_{\bullet}^{rep}(F,*):=\bigoplus_{p,q}H^{rep}_{p}(F,q)$$
\item
Si $F$ est lisse sur $S=Spec k$, sa classe de Todd est d\'efinie par
$$Td^{rep}(F):=Ch(\alpha_{F}^{rep})^{-1}.Td(T_{I_{F}^{t}}) \in
H_{rep}{\bullet}(F,*)$$
\item
Le caract\`ere de Chern est d\'efini par
$$Ch^{rep} :  \xymatrix{
\bf {K}_{*}(F) \ar[r]^{\phi_{F}} & \underline{\bf K}_{*}^{rep}(F)
\ar[r]^{Ch} & H^{\bullet}_{rep}(F,*)}$$
\item
Si $F$ est lisse sur $S=Spec k$ et que son espace de modules
est quasi-projectif, alors
la transformation de Riemann-Roch est d\'efinie par
$$\begin{array}{cccc}
\tau^{rep}_{F} : & \bf G\mit_{*}(F) & \longrightarrow &
H^{\bullet}_{rep}(F,*) \\
           & x               & \mapsto & Ch^{rep}(x).Td^{rep}(F)
\end{array}$$
\end{enumerate}
\end{df}

Remarquons que le m\^eme argument que \ref{p3.4} montre que, si $F$ est un
champ de
Deligne-Mumford r\'egulier, alors $I_{F}^{t}$ est aussi r\'egulier,
et de plus $\alpha_{F}^{rep}$ est inversible. Ainsi, la d\'efinition $(2)$
poss\`ede un sens. \\

Pour tout champ alg\'ebrique de Deligne-Mumford, il existe une section
canonique $F \hookrightarrow I_{F}^{t}$ correspondant \`a l'identit\'e.
Cette section induit donc des isomorphismes naturels
$$H^{\bullet}_{rep}(F,*)\simeq H^{\bullet}(F,*)\oplus
H^{\bullet}_{rep\neq 1}(F,*)$$
$$H_{\bullet}^{rep}(F,*)\simeq H_{\bullet}(F,*)\oplus
H_{\bullet}^{rep\neq 1}(F,*)$$
compatibles avec les produits et les classes caract\'eristiques. \\

Nous ne red\'emontrerons pas les propri\'et\'es des th\'eories
$H_{rep}^{\bullet}(F,*)$, et $H^{rep}_{\bullet}(F,*)$. Elles sont en
tout point analogues \`a celles de la cohomologie \`a coefficients dans
les caract\`eres. \\

\underline{\bf Remarque:} \rm Supposons que $k=\bf C$, et que l'on
prenne la cohomologie de De Rham.
Si $F$ est un champ de Deligne-Mumford lisse, on peut lui associer un champ
analytique $F^{an}$ (~\ref{d5.5}~).
Le th\'eor\`eme de comparaison de Grothendieck (~\cite[Thm. $1'$]{gr}~)
donne alors un isomorphisme de \\
$\bf C$-alg\`ebres
$$H_{rep}^{\bullet}(F,0)\simeq H^{\bullet}((I_{F})^{top},\bf C\mit)$$
o\`u $(I_{F})^{top}$ est le site topologique sur $(I_{F})^{an}$, et $\bf
C$ le faisceau constant de fibre $\bf C$. Soit $\pi : I_{F}
\rightarrow F$ la projection, et  $p : F \rightarrow M$ la projection sur
l'espace de
modules. On pose $R=p_{*}\circ \pi_{*}(\bf C\mit)$; c'est un
faisceau en \\
$\bf C$-alg\`ebres sur $M^{top}$, dont la fibre au point
$x \in M$ est isomorphe \`a $\bf C\mit(H_{x})$, la
$\bf C$-alg\`ebre des fonctions centrales sur le groupe d'isotropie $H_{x}$
de $x$. Ainsi
$$H_{rep}^{\bullet}(F,0)\simeq H^{\bullet}(M^{top},R)$$
Cet isomorphisme explique le nom de "cohomologie \`a coefficients dans
les repr\'esentations", en gardant \`a l'esprit que les \'el\'ements de $\bf
C\mit(H_{x})$ s'identifient \`a des repr\'esentations virtuelles de
$H_{x}$ \`a coefficients complexes. \\

Le th\'eor\`eme suivant est la version \'etendue de \ref{th3.1}

\begin{thm}\label{th3.3}{(~Lefschetz-Riemann-Roch~)}
Supposons que $S$ est r\'egulier.
Soit $\cal QDM$ (~resp. $\cal DM$~) la sous-cat\'egorie de $HoChAlg(S)$ des
champs de Deligne-Mumford sur $S$, et
dont l'espace de modules est quasi-projectif sur $S$ (~resp. des champs de
Deligne-Mumford sur $S$~).
Alors, pour chaque $F$ objet de $\cal DM$, il existe un unique morphisme
$$\psi_{F} : \bf G\mit_{*}(F) \longrightarrow \underline{\bf
G}^{rep}_{*}(F)$$
tel que~:
\begin{enumerate}
\item
Si $F$ est lisse sur $S$, et dans $\cal QDM$, alors
$$\begin{array}{cccc}
\psi_{F} : & \bf K\mit_{*}(F) & \longrightarrow & \underline{\bf
G}^{rep}_{*}(F) \\
   & x & \mapsto & (\alpha_{F}^{rep})^{-1}.\phi_{F}(x)
\end{array}$$
\item
Pour tout morphisme propre de dimension cohomologique finie de $\cal
DM$, $f : F \longrightarrow F'$, on a
$$f_{*}\circ \psi_{F} = \psi_{F'}\circ f_{*}$$
dans l'un des cas suivants
\begin{itemize}
\item le morphisme $f$ est repr\'esentable
\item le champ $F$ est lisse sur $S$
\end{itemize}
\item
Si $f : F \longrightarrow F'$ est un morphisme repr\'esentable et \'etale
de champs de $\cal DM$, alors
$$\psi_{F}\circ f^{*} = f^{*} \circ \psi_{F'}$$
\item
Pour tout champ $F$ de $\cal DM$, tout $x \in \bf G\mit_{*}(F)$ et
tout $y \in \bf K\mit_{*}(F)$, on a
$$\psi_{F}(x.y)=\psi_{F}(x).\phi_{F}(y)$$
\item
Pour tout morphisme propre de dimension cohomologique finie, \\
$f~:~F~\longrightarrow~F'$ de champs de $\cal DM$, le diagramme suivant
commute
$$\xymatrix{
\bf G\mit_{0}(F) \ar[r]^{f_{*}} \ar[d]_{\psi_{F}} & \bf G\mit_{0}(F')
\ar[d]^{\psi_{F'}} \\
\underline{\bf G}^{rep}_{0}(F) \ar[r]_{f_{*}} & \underline{\bf
G}^{rep}_{0}(F')}$$
\end{enumerate}
\end{thm}

Par la suite, nous dirons "sch\'ema", pour "$S$-sch\'ema". Les termes
"projectifs", "fortement projectif", "lisse" ... feront r\'ef\'erence aux
notions relatives sur $S$.\\

\underline{\bf Preuve:} \rm Notons que la m\^eme preuve que \ref{th3.1}, montre
que la transformation naturelle
$$\psi_{F}:=(\alpha_{F}^{rep})^{-1}.\phi_{F}$$
commute avec les
images directes de morphismes fortement projectifs entre champs
lisses de $\cal QDM$.\\

\underline{\bf D\'efinition de $\psi_{F}$:} \\

Soit $F$ un champ de $\cal DM$. Nous dirons qu'un champ
simplicial $q~:~F_{\bullet}~\longrightarrow~F$
augment\'e et propre (~\ref{chsimp}~) sur $F$, est une enveloppe si~:
\begin{itemize}
\item
Pour tout $[m] \in \Delta^{op}$
$F_{m}$ est une gerbe
triviale sur un sch\'ema quasi-projectif et born\'ee par un faisceau
constant de groupes finis (donc des r\'eunions disjointes de $X_{n}\times BH_{n}$.
\item
Les morphismes de faces et d\'eg\'en\'erescences $F_{m}
\longrightarrow F_{n}$ sont propres repr\'esentables et de la forme
$f\times \rho : X_{n}\times BH_{n} \longrightarrow X_{m}\times BH_{m}$, pour $f : X_{n} \longrightarrow X_{m}$ et 
$\rho : H_{n} \longrightarrow H_{m}$. 
\item
Le morphisme induit (~\ref{chsimp}~)
$$q_{*} : \bf G\mit_{*}(F_{\bullet})_{\bf Q} \longrightarrow
\bf G\mit_{*}(F)_{\bf Q}$$
est un isomorphisme.
\end{itemize}

\begin{lem}\label{l3.7}
\begin{enumerate}
\item
Il existe toujours une enveloppe pour $F$.
\item
Deux enveloppes quelconques de  $F$ sont domin\'ees par une troisi\`eme.
\item
Pour tout $1$-morphisme de champs de $\cal DM$, $f : F \longrightarrow
F'$, il existe un diagramme commutatif de champs simpliciaux augment\'es
$$\xymatrix{
F_{\bullet} \ar[r]^{f_{\bullet}} \ar[d]_{b} & F_{\bullet}'
\ar[d]^{a} \\
F \ar[r]_{f} & F' }$$
avec $a$ et $b$ des enveloppes.
\end{enumerate}
\end{lem}

\underline{\bf Preuve:} \rm La preuve est la m\^eme que \cite[$3.1$,
$3.2$, $3.3$]{g3}, mais
en rempla\c{c}ant le mot "enveloppe" par "quasi-enveloppe de
Chow". $\Box$\\

Consid\'erons alors $F$ un champ de $\cal DM$, et $q : F_{\bullet}
\longrightarrow F$ une enveloppe.

\begin{lem}\label{l3.8}
Il existe un morphisme $\psi_{F_{\bullet}}$ dans la cat\'egorie homotopique des spectres simpliciaux, entre
$$[n] \mapsto \bf G\mit(F_{n})$$
et 
$$[n] \mapsto \underline{\bf G\mit}^{rep}(F_{n}),$$
telle que pour tout $[n]$, le morphisme induit $\bf G\mit(F_{n}) \longrightarrow 
\underline{\bf G}^{rep}(F_{n})$ soit le morphisme $\phi_{m}$.
\end{lem}

\underline{\bf Preuve:} \rm Pour chaque $[n] \in \Delta$, on a une \'equivalence de spectres
$$\underline{\bf G\mit}^{rep}(F_{n}) \longrightarrow \bf G\mit(MIF_{n})_{K}$$
o\`u $MIF_{n}$ est l'espace de modules du champ $I_{F_{n}}^{t}$. Ainsi, le spectre simplicial 
$[n] \mapsto \underline{\bf G\mit}^{rep}(F_{n})$ est \'equivalent \`a $[n] \mapsto
\bf G\mit(MIF_{n})_{K}$. 

D'un autre cot\'e, pour tout $[n] \in \Delta$, on dispose du morphisme naturel 
$$\bf G\mit(F_{n}) \longrightarrow \bf G\mit(X_{n})\otimes K(H_{n}),$$
o\`u $K(H_{n})$ est l'alg\`ebre des fonctions centrales sur $H_{n}$, nulles sur les \'el\'ements
d'ordre non-inversible sur $S$, \`a valeurs 
dans $K$. Il est d\'efini exactement de la m\^eme fa\c{c}on que le morphisme $\phi$.
Si $f\times \rho : X_{n}\times BH_{n} \longrightarrow X_{m}\times BH_{m}$ est un morphisme de 
transition, alors on dispose de
$$f_{*}\times Ind_{\rho} : \bf G\mit(X_{n})\otimes K(H_{n}) \longrightarrow \bf G\mit(X_{m})\otimes K(H_{m}).$$
Ici, $f_{*}$ est l'image directe pour la version strictement fonctorielle de $\bf G$, et 
$Ind_{\rho} : K(H_{n}) \longrightarrow K(H_{m})$ est \'egal par d\'efinition \`a $[H_{n}:H_{m}]$ fois le morphisme
d'induction des caract\`eres (remarquer qu'ici $H_{n} \longrightarrow H_{m}$ est injectif). 

Ce morphisme est de plus fonctoriel pour les images directes pour $F_{n} \longrightarrow F_{m}$, et 
donc d\'efinit un morphisme de spectres simpliciaux, entre $[n] \mapsto \bf G\mit(F_{n})$ et
$[n] \mapsto \bf G\mit(X_{n})\otimes K(H_{n})$. 

Enfin, il est clair que $\bf G\mit(X_{n})\otimes K(H_{n}) \simeq \bf G\mit(MIF_{n})_{K}$, et donc
$[n] \mapsto \bf G\mit(MIF_{n})_{K}$ est \'equivalent \`a $[n] \mapsto  \bf G\mit(X_{n})\otimes K(H_{n})$. 

En conclusion, on a un diagramme de spectres simpliciaux
$$\bf G\mit(F_{\bullet}) \longrightarrow \bf G\mit(MIF_{\bullet})_{K} \simeq \underline{\bf G}^{rep}(F_{\bullet}),$$ 
ce qui d\'efinit bien le morphisme cherch\'e dans la cat\'egorie homotopique. $\Box$\\

Le lemme pr\'ec\'edent implique que l'on peut d\'efinir un morphisme de spectres
(~\ref{chsimp}~)
$$\psi_{F_{\bullet}} : \bf G\mit_{*}(F_{\bullet}/F)_{\bf Q} \longrightarrow
\underline{\bf G}^{rep}_{*}(F_{\bullet}/F)_{K}$$
On d\'efinit alors $\psi_{F}$ par le diagramme commutatif suivant
dans $HoSp$
$$\xymatrix{
\bf G\mit_{*}(F_{\bullet}/F)_{\bf Q} \ar[d]_{q_{*}}
\ar[r]^-{\psi_{F_{\bullet}}} & \underline{\bf G}^{rep}_{*}(F_{\bullet}/F)_{K} \ar[d]^{Iq_{*}} \\
\bf G\mit_{*}(F)_{\bf Q}
\ar[r]^-{\psi_{F}} & \underline{\bf G}^{rep}_{*}(F)_{K} }$$
Cette d\'efinition poss\`ede un sens, car d'apr\`es \ref{c2.2}, on sait que
$q_{*}$ est un isomorphisme. 

Les propri\'et\'es des enveloppes rappel\'ees dans \ref{l3.7} et
une autre application du lemme pr\'ec\'edent montre que $\psi_{F}$
est ind\'ependant du choix de l'enveloppe.\\

Remarquons enfin, que si $q_{0} : F_{0} \longrightarrow F$ est une
quasi-enveloppe de Chow, avec $q_{0}$ un morphisme propre et
repr\'esentable, et $F_{0}$ une gerbe triviale, on peut construire
une enveloppe $q : F_{\bullet} \longrightarrow F$, avec \\
$q=q_{0}~:~F_{0}~\longrightarrow~F$. Alors, comme le champ simplicial constant
$F_{0}$ est une enveloppe pour $F_{0}$, on a
$$(Iq_{0})_{*}\psi_{F_{0}}=\psi_{F}(q_{0})_{*}$$

\underline{\bf Preuve de $(1)$:} \rm \\

Commen\c{c}ons par une petite digression sur les orbifolds lisses de
$\cal QDM$~: ce sont les champs $F$ de $\mathcal{QDM}$,lisses sur $S$, tel que
$\pi_{F} : I_{F} \longrightarrow F$ soit birationnel.

\begin{lem}\label{l3.9}
Soit $F$ un orbifold lisse et connexe de $\cal QDM$. Alors il existe un
sch\'ema lisse et quasi-projectif $X$, une action du $S$-sch\'ema en groupes $\bf
Gl\mit_{n}/S$ sur $X$, et une \'equivalence de champs
$$F \simeq [X/\bf Gl\mit_{n}/S]$$
\end{lem}

\underline{\bf Preuve:} \rm La preuve qui suit m'a \'et\'e communiqu\'ee
par A. Kresch. \\

Soit $q : T^{k} \longrightarrow F$ le $\bf
Gl\mit_{n}/S$-torseur correspondant au fibr\'e vectoriel des $k$-jets sur $F$
(~relatif \`a $S$~). Commen\c{c}ons par montrer que $T^{k}$ est repr\'esentable d\`es que $k$ est assez grand.

Comme ceci est local sur l'espace de modules de $F$, on peut supposer
que $F=[X/H]$, est le quotient d'un sch\'ema quasi-projectif lisse et
connexe par un groupe fini. Dire que $F$ est un orbifold est alors
\'equivalent \`a dire qu'aucun \'el\'ement de $H$ ne fixe le point g\'en\'erique
de $X$, ou encore que pour tout $h \in H$, le sous-sch\'ema des points
fixes de $h$ dans $X$ est un ferm\'e strict, $X^{h} \hookrightarrow X$.

Soit $x : Spec k(x) \longrightarrow X$ un point g\'eom\'etrique, et
consid\'erons $H_{x}~\longrightarrow~Gl(T^{k}_{X,x})$ l'action du sous-groupe
d'isotropie de $x$ sur l'espace des $k$-jets en $x$. Comme pour $h \in H_{x}$,
$T^{k}_{X^{h},x}\simeq (T^{k}_{X,x})^{h}$, cette action est fid\`ele pour $k$ assez grand. Ainsi, pour
tout point g\'eom\'etrique $x$ de $X$, le morphisme naturel
$$H_{x} \longrightarrow Gl(T^{k}_{X,x})\simeq \bf Gl\mit_{n}/k(x)$$
est injectif pour $k$ assez grand.

Soit $T^{k}_{X}$ le $\bf Gl\mit_{n}/S$-torseur des $k$-jets sur $X$
(~relativement \`a $S$~). Alors le champ $T^{k}$ est \'equivalent au champ
quotient $[T^{k}_{X}/H]$. Or, comme pour tout point g\'eom\'etrique $x$ de
$X$, le morphisme naturel de $H_{x}$ dans $\bf Gl\mit_{n}/k(x)$ est
injectif, $H$ op\`ere sans points fixes sur $T^{k}_{X}$. Ainsi, le champ
$T\simeq [T^{k}_{X}/H]$ est repr\'esentable par le $S$-sch\'ema
quotient $T^{k}_{X}/H$. \\

De plus, si $M$ est l'espace de modules de $F$, le
morphisme induit $T^{k} \longrightarrow M$ est affine. Comme $M$ est
quasi-projectif sur $S$, il en est de m\^eme de $T^{k}$.

On vient donc de d\'emontrer que $F \simeq [T^{k}/\bf Gl\mit_{n}/S]$, avec
$T^{k}$ un $S$-sch\'ema quasi-projectif et lisse. $\Box$\\

A l'aide du lemme pr\'ec\'edent, et de \cite{th2}, on sait que tout faisceau
coh\'erent sur un orbifold est quotient d'un fibr\'e vectoriel, et que $F$ admet
un faisceau inversible ample. En particulier, tout morphisme propre
et repr\'esentable $f : F' \longrightarrow F$ est fortement projectif. \\

Revenons au cas g\'en\'eral d'un champ lisse de $\cal QDM$.
On sait qu'il existe un morphisme
$$p : F \longrightarrow F'$$
o\`u $F'$ est un orbifold lisse de $\cal QDM$, et $p$ fait de $F$ une
gerbe de groupe fini sur $F'$ (~\cite{s2}~).
On choisit \`a l'aide de \ref{th1.4} une
quasi-enveloppe de Chow, repr\'esentable
$$r : F'_{0} \longrightarrow F'$$
avec $F'_{0}$ une gerbe triviale sur un sch\'ema quasi-projectif.
Posons \\
$q~:~F_{0}=F'_{0}\times_{F'}F~\longrightarrow~F$ le morphisme
induit. Alors $F_{0}$ est une gerbe sur $F'_{0}$, donc est encore une
gerbe.

Comme $r$ est fortement projectif (~par \ref{l3.9}~), $q$ l'est
aussi. En effectuant un changement de base fini de l'espace de modules
de $F_{0}$, ce qui ne change pas le caract\`ere fortement projectif de $q$,
on peut m\^eme supposer que c'est une gerbe triviale.\\

Montrons alors que $q_{*} : \bf G\mit_{*}(F_{0})_{\bf Q} \longrightarrow \bf
G\mit_{*}(F)_{\bf Q}$ est surjectif.

En effet, comme $q$ est fortement projectif, on peut utiliser le
th\'eor\`eme de Lefschetz-Riemann-Roch pour les
morphismes fortement projectifs, qui nous donne un diagramme
commutatif
$$\xymatrix{
\bf G\mit_{*}(F_{0})_{K} \ar[rr]^-{(\alpha_{F_{0}}^{rep})^{-1}.\phi_{F}}
\ar[d]_{q_{*}} & & \underline{\bf G}_{*}^{rep}(F_{0}) \ar[d]^{q_{*}} \\
\bf G\mit_{*}(F_{0})_{K}
\ar[rr]_{(\alpha_{F}^{rep})^{-1}.\phi_{F}} & & \underline{\bf
G}_{*}^{rep}(F) }$$
Comme les fl\`eches horizontales sont des isomorphismes, il suffit de
montrer que $q_{*} : \underline{\bf G}^{rep}(F_{0}) \longrightarrow
\underline{\bf G}^{rep}(F)$ est surjectif,
ce qui est une cons\'equence directe du lemme suivant.

\begin{lem}\label{l3.11}
Soit $f : F' \longrightarrow F$ un morphisme fini repr\'esentable et surjectif de
champs, avec $F$ lisse. Alors, le morphisme
$$f_{*} : \underline{\bf G}_{*}(F') \longrightarrow
\underline{\bf G}_{*}(F)$$
est surjectif.
\end{lem}

\underline{\bf Preuve:} \rm A l'aide de la dualit\'e de Poincar\'e \ref{p2.1}, on
peut utiliser la formule de projection pour $f$. Ceci implique que
$$f_{*}(f^{*}(x))=x.f_{*}(1)$$
pour $x \in \underline{\bf G}_{*}(F)$. Il suffit donc de montrer que
$f_{*}(1)$ est inversible dans $\underline{\bf G}_{0}(F)$, ce
qui, d'apr\`es le lemme \ref{l3.3}, est local sur $F_{et}$. On peut donc
supposer que $F$ et $F'$ sont des sch\'emas quasi-projectifs. Mais
$f_{*}(1)$ poss\`ede un rang \'egal au degr\'e de $f$, qui est non-nul, car
$f$ est dominant. Ainsi $f_{*}(1) \in \bf G\mit_{0}(F)_{\bf Q}$ a un
rang inversible, donc est inversible (~\cite{fl}~). $\Box$\\

Soit $x \in \bf G\mit_{*}(F)_{\bf Q}$. On vient de voir que l'on peut
\'ecrire $x=q_{*}(y)$, avec $y \in \bf G\mit_{*}(F_{0})_{\bf Q}$. On
applique alors une nouvelle fois la formule
pour le cas des morphismes fortement projectifs, \`a $q$ et $y$
$$(\alpha_{F}^{rep})^{-1}\phi_{F}(q_{*}(y)) = Iq_{*}\phi_{F_{0}}(y)$$
Ainsi, comme $Iq_{*}(\phi_{F_{0}}(y))=\psi_{F}(x)$, on a
$$\psi_{F}=(\alpha_{F}^{rep})^{-1}.\phi_{F}$$

\underline{\bf Preuve de $(2)$:} \rm \\

Soit $f : F \longrightarrow F'$ un morphisme propre de dimension
cohomologique finie. \\

\underline{Cas o\`u $f$ est repr\'esentable :} D'apr\`es \ref{l3.7}, on peut
trouver un diagramme commutatif de champs simpliciaux augment\'es
(~\ref{chsimp}~)
$$\xymatrix{
F_{\bullet} \ar[r]^{f_{\bullet}} \ar[d]_{b} & F_{\bullet}'
\ar[d]^{a} \\
F \ar[r]_{f} & F' }$$
avec $a$ et $b$ des enveloppes.
En prenant les images par les foncteurs $\bf G$ et $\underline{\bf
G}^{rep}$, ce diagramme induit deux diagrammes commutatifs dans $HoSp$
(~\ref{chsimp}~)
$$\xymatrix{
\bf G\mit(F_{\bullet}) \ar[r]^{(f_{\bullet})_{*}} \ar[d]_{b_{*}} &
\bf G\mit(F_{\bullet}') \ar[d]^{a_{*}} & &
\underline{\bf G}^{rep}(F_{\bullet}) \ar[r]^{(f_{\bullet})_{*}}
\ar[d]_{b_{*}} &
\underline{\bf G}^{rep}(F_{\bullet}') \ar[d]^{a_{*}} \\
\bf G\mit(F) \ar[r]_{f_{*}} & \bf G\mit(F') & &
\underline{\bf G}^{rep}(F) \ar[r]_{f_{*}} &
\underline{\bf G}^{rep}(F') }$$
Le morphisme $\phi$ d\'efinit alors un morphisme entre ces deux
diagrammes
$$\xymatrix{
 & \bf G\mit(F_{\bullet})_{\bf Q} \ar[dd]_(.3){\phi_{F_{\bullet}}}
 \ar[rr]^{(f_{\bullet})_{*}}
 \ar[dl]^{a_{*}} & & \bf G\mit(F_{\bullet}')_{\bf Q}
\ar[dd]^{\phi_{F_{\bullet}'}}
 \ar[dl]^{b'_{*}} \\
\bf G\mit(F)_{\bf Q} \ar[dd]_{\phi_{F}} \ar[rr]^(.7){f_{*}} & &
\bf G\mit(F')_{\bf Q} \ar[dd]^(.3){\phi_{F'}} \\
& \underline{\bf G}^{rep}(F_{\bullet}) \ar[rr]_(.3){(f_{\bullet})_{*}}
\ar[dl]^{a_{*}}
 & & \underline{\bf G}^{rep}(F'_{\bullet}) \ar[dl]^{b_{*}}\\
\underline{\bf G}^{rep}(F) \ar[rr]_{f_{*}} & &
\underline{\bf G}^{rep}(F') }$$
Ce qui, par d\'efinition, implique que
$$f_{*}\circ\psi_{F}=\psi_{F'}\circ f_{*}$$

\underline{Cas o\`u $F$ est lisse :} Soit $q' : F_{0}' \longrightarrow F'$
une quasi-enveloppe de Chow, avec $F_{0}'$ une gerbe triviale, et $q'$
un morphisme repr\'esentable fini (~\ref{th1.4}~). Notons $q :
F_{0}'\times_{F'}F \longrightarrow F$ la quasi-enveloppe induite,
\mbox{$r : F_{0} \longrightarrow F_{0}'\times_{F'}F$} une seconde
quasi-enveloppe de Chow par une gerbe triviale, et $r$ un morphisme
repr\'esentable fini. Consid\'erons le diagramme commutatif suivant
$$\xymatrix{
F_{0} \ar[r]^{g} \ar[d]_{p} & F_{0} \ar[d]^{q'} \\
F \ar[r]_{f} & F'}$$

Nous avons d\'ej\`a vu pr\'ec\'edemment que le morphisme induit
$$p_{*} : \bf G\mit_{*}(F_{0})_{\bf Q} \longrightarrow \bf
G\mit_{*}(F)_{\bf Q}$$
est surjectif.

Soit $x \in \bf G\mit_{*}(F)_{\bf Q}$, et $y \in \bf
G\mit_{*}(F_{0})_{\bf Q}$ tel que $p_{*}(y)=x$. Alors
$$\begin{array}{cl}
\psi_{F}(x) & = \psi_{F}(p_{*}(x)) \\
            & = p_{*}\psi_{F_{0}}(y)
\end{array}$$
par d\'efinition de $\psi_{F}$. Ainsi, si l'on savait que
$$g_{*}\psi_{F_{0}}=\psi_{F_{0}'}g_{*}$$
on aurait
$$\begin{array}{cl}
f_{*}\psi_{F}(x) & = f_{*}\psi_{F}(p_{*}(x)) \\
                 & = f_{*}p_{*}\psi_{F_{0}}(y) \\
                 & = q'_{*}g_{*}\psi_{F_{0}}(y) \\
                 & = q'_{*}\psi_{F_{0}'}(g_{*}(y)) \\
                 & = \psi_{F'}(q'_{*}g_{*}(y)) \\
                 & = \psi_{F'}(f_{*}p_{*}(y)) \\
                 & = \psi_{F'}(f_{*}(x))
\end{array}$$
Il nous reste donc \`a d\'emontrer le lemme suivant

\begin{lem}\label{l3.12}
Soit $f : F \longrightarrow F'$ un morphisme propre de dimension
cohomologique finie de gerbes
triviales de Deligne-Mumford, alors le diagramme suivant commute
$$\xymatrix{
\bf G\mit_{*}(F) \ar[r]^{f_{*}} \ar[d]_-{\phi_{F}} & \bf G\mit_{*}(F')
\ar[d]^-{\phi_{F'}} \\
\underline{\bf G}^{rep}_{*}(F) \ar[r]_{f_{*}} & \underline{\bf
G}^{rep}_{*}(F')}$$
\end{lem}

\underline{\bf Preuve :} \rm Factorisons $f$ par son "espace de
modules relatif"
$$f : \xymatrix{F \ar[r]^{k} & F'' \ar[r]^{h} & F'}$$
Le morphisme $k$ est, par d\'efinition, le morphisme universel
vers les champs qui sont repr\'esentable sur $F'$.

Remarquons qu'un tel champ existe. Comme l'existence de
ce champ est locale sur $F'_{et}$, il suffit de montrer qu'il existe
lorsque $F'$ est un sch\'ema. Mais dans ce cas, $F''$ est tout simplement
l'espace de modules de $F$.

On peut donc supposer que $f=k$. Dans ce cas les champs $F$ et $F'$
sont des gerbes sur $M$, l'espace de modules de $F$, et le morphisme induit
par $f$ est l'identit\'e sur $M$. De plus, le morphisme $f$ fait de $F$ une
gerbe \'etale sur $F'$, born\'ee par un groupe
fini $K$. Ce qui implique en particulier que $f$ est
\'etale (~\'eventuellement non repr\'esentable~).

Formons alors le carr\'e cart\'esien suivant
$$\xymatrix{F \ar[r]^{f} & F' \\
F'' \ar[u]^{u} \ar[r]_{u} & F \ar[u]_{f} }$$
Supposons que l'on ait d\'emontr\'e le lemme pour le morphisme $u$, alors,
par la formule de transfert (~\ref{p2.2}~), et la fonctorialit\'e de $\phi$
pour les images r\'eciproques, on peut \'ecrire
$$\begin{array}{cl}
f^{*}f_{*}\phi_{F} & = g_{*}g^{*}\phi_{F} \\
& = g_{*}\phi_{F''}g^{*} \\
& = \phi_{F}g_{*}g^{*} \\
& = \phi_{F}f^{*}f_{*} \\
& = f^{*}\phi_{F'}f_{*}
\end{array}$$
Or, on sait que $f_{*}f^{*}=\times\frac{1}{k}$, o\`u $k$ est l'ordre de
$K$. Ainsi, en multipliant par $k.f_{*}$, on trouve
$$f_{*}\phi_{F}=\phi_{F'}f_{*}$$
Il nous reste \`a d\'emontrer la formule pour $u$. Or, par construction,
$u$ poss\`ede une section $s : F \longrightarrow F''$. Ainsi, $F''$ est
une gerbe triviale sur $F$. On peut donc supposer que
$F''=F\times_{S} BK$, o\`u $BK=[S/K]$, et $u$ la projection
$$b : F\times_{S} BK \longrightarrow F$$
On utilise alors la "formule de Kunneth" pour $F''$.

\begin{lem}
Soit $a : F\times_{S} BK \longrightarrow BK$, et
$b : F\times_{S} BK \longrightarrow F$ les deux projections. Alors, le
morphisme
$$a^{*}.b^{*} : \bf K\mit_{0}(BK) \otimes_{\bf K\mit_{0}(S)} \bf
G\mit_{*}(F)\longrightarrow \bf G\mit_{*}(F\times BK)$$
est bijectif.
\end{lem}

\underline{\bf Preuve:} \rm Soit $Spec A \hookrightarrow Spec \bf Z$
l'image du morphisme canonique $S~\longrightarrow~Spec \bf Z$.
Comme $K$ est d'ordre premier aux
caract\'eristiques de $S$, le $Spec A$-sch\'ema en groupes qu'il d\'efinit est
r\'eductif. Notons alors $R(K)$ un syst\`eme de repr\'esentants des classes
d'isomorphie des objets simples de la cat\'egorie $\bf Coh\mit([Spec A/K])$.
Ce sont aussi les objets simples de la cat\'egorie des
$A[K]$-modules. Le
th\'eor\`eme de d\'evissage de Quillen \cite[$5.1$]{q}, implique alors qu'il
existe un
isomorphisme canonique
$$\bf K\mit_{0}(BK) \simeq \bigoplus_{R(K)}\bf K\mit_{0}(S)$$
De m\^eme, il existe un isomorphisme
$$\bf G\mit_{*}(F\times_{S} BK) \simeq \bigoplus_{R(K)}\bf G\mit_{*}(F)$$
et le lemme s'en suit. $\Box$\\

Soit $x \in \bf G\mit_{*}(F\times_{S} BK)_{\bf Q}$. On \'ecrit
$x=a^{*}(y).b^{*}(z)$, avec $y \in \bf K\mit_{0}(BK)$ et
$z \in \bf G\mit_{*}(F)$.
Notons $p : BK \longrightarrow S$, et
$q : F \longrightarrow S$ les deux projections.
Supposons que le lemme soit d\'emontr\'e pour $F=BK$,
$F'=S$, et pour les \'el\'ements
de $\bf G\mit_{0}(F)$, alors
$$\begin{array}{cl}
\phi_{F}(b_{*}(x)) & = \phi_{F}(b_{*}(a^{*}(y).b^{*}(z))) \\
& = \phi_{F}(b_{*}a^{*}(y).z) \\
& = \phi_{F}(q^{*}p_{*}(y).z) \\
& = q^{*}\phi_{Spec k}(p_{*}(y)).\phi_{F}(z) \\
& = q^{*}p_{*}(\phi_{BK}(y)).\phi_{F}(z) \\
& = b_{*}a^{*}(\phi_{BK}(y)).\phi_{F}(z) \\
& = b_{*}\phi_{F\times BK}(a^{*}(y)).\phi_{F}(z) \\
& = b_{*}(\phi_{F\times BK}(a^{*}(y)).b^{*}\phi_{F}(z) \\
& = b_{*}(\phi_{F\times BK}(a^{*}(y)).\phi_{F\times BK}(b^{*}(z)))\\
& = b_{*}(\phi_{F\times BK}(a^{*}(y).b^{*}(z))) \\
& = b_{*}(\phi_{F\times BK}(x))
\end{array}$$

On s'est donc ramen\'e au cas o\`u $F=BK$, et $f$ est le morphisme
structural
$$f : BK \longrightarrow S$$
En identifiant $\bf Coh\mit(F)$ avec la cat\'egorie des repr\'esentations
de $H$ dans $\bf Coh \mit(S)$, le morphisme
$$f_{*} : \bf G\mit_{0}(F)_{\bf Q} \longrightarrow \bf
G\mit_{0}(S)_{\bf Q}$$
associe \`a un $\cal O\mit_{S}$-module coh\'erent muni d'une action de $H$,
son sous-module des invariants.

De plus, $\underline{\bf G}^{rep}_{0}(F)$
est isomorphe \`a l'anneau des fonctions centrales de $H$ dans
$\bf G\mit_{0}(S)_{K}$.
Et par cet isomorphisme, le morphisme
$$\phi_{F} : \bf G\mit_{0}(F) \longrightarrow \bf G\mit_{0}(S)_{K}$$
associe \`a une repr\'esentation $V$ de $H$ dans $\bf Coh\mit(S)$ "son caract\`ere"
$$\begin{array}{cccc}
\chi : & K & \longrightarrow & \bf G\mit_{0}(S)_{K} \\
& h & \mapsto & \sum_{\zeta \in \mu_{\infty}(S)}\zeta.V^{h,(\zeta)}
\end{array}$$
o\`u $V^{h,(\zeta)}$ est le sous-module de $V$ sur lequel $h$ op\`ere par
multiplication par $\zeta$.

Enfin, le morphisme
$$f_{*} : \underline{\bf G}^{rep}(F) \longrightarrow
\bf G\mit_{0}(S)_{K}$$
est alors donn\'e par
$$f_{*}(\chi)=\frac{1}{k}.\sum_{h \in K}\chi(h)$$
pour toute fonction centrale
$$\chi : K \longrightarrow \bf G\mit_{0}(S)_{K}$$
La formule \`a d\'emontrer se traduit donc par
$$[V^{K}]=\frac{1}{k}.\sum_{h \in K}\zeta.[V^{h,(\zeta)}]$$
en tant qu'\'el\'ements de $\bf G\mit_{0}(S)_{K}$.
Ce qui est une formule bien connue.\\

\underline{\bf Preuve de $(3)$ :} \rm \\

Soit $f : F \longrightarrow F'$ un morphisme repr\'esentable \'etale,
$q' : F_{\bullet}' \longrightarrow F'$ une enveloppe,
$q : F_{\bullet}:=F'_{\bullet}\times_{F'}F \longrightarrow F$ l'enveloppe
induite, et
$f_{\bullet} : F_{\bullet} \longrightarrow F'_{\bullet}$
le morphisme induit.

Soit $x \in \bf G\mit_{*}(F')_{\bf Q}$, et  $x=q'_{*}(y)$, avec
$y \in \bf G\mit_{*}(F'_{\bullet}/F')_{\bf Q}$.
Alors, par la formule de transfert, on a
$$f^{*}(x)=q_{*}f^{*}_{\bullet}(y)$$
Ainsi,
$$\begin{array}{cl}
\psi_{F}(f^{*}x) & = q_{*}(\phi_{F_{\bullet}}(f_{\bullet}^{*}(y))) \\
& = q_{*}(f_{\bullet}^{*}\phi_{F_{\bullet}}(y)) \\
& = f^{*}(q'_{*}\phi_{F_{\bullet}}(y))\\
& = f^{*}\psi_{F'}(x)
\end{array}$$

\underline{\bf Preuve de $(4)$ :} \rm \\

Soit $x \in \bf G\mit_{*}(F)$ et $y \in \bf K\mit_{*}(F)$. Ecrivons
$$x=q_{*}(z)$$
avec $q : F_{\bullet} \longrightarrow F$ une enveloppe, et
$z \in \bf G\mit_{*}(F_{\bullet})_{\bf Q}$. Alors, comme
$q_{*}(z.q^{*}(y))=x.y$, on a
$$\begin{array}{cl}
\psi_{F}(x.y) & = q_{*}(\phi_{F_{\bullet}}(z.q^{*}(y))) \\
& = q_{*}(\phi_{F_{\bullet}}(z).\phi_{F_{\bullet}}(q^{*}(y))) \\
& = q_{*}(\phi_{F_{\bullet}}(z).q^{*}\phi_{F}(y)) \\
& = q_{*}(\phi_{F_{\bullet}}(z)).\phi_{F}(y) \\
& = \psi_{F}(x).\phi_{F}(y)
\end{array}$$

\underline{\bf Preuve de $(5)$ :} \rm C'est la m\^eme que le second cas
du point $(2)$, car on a le r\'esultat suivant.

\begin{lem}
Soit $q : F_{0} \longrightarrow F$ une quasi-enveloppe de Chow, avec
$q$ repr\'esentable et fini. Alors le morphisme
$$q_{*} : \bf G\mit_{0}(F_{0})_{\bf Q} \longrightarrow
\bf G\mit_{0}(F)_{\bf Q}$$
est surjectif.
\end{lem}

\underline{\bf Preuve :} \rm On proc\`ede par r\'ecurrence noeth\'erienne dans $F$.

Soit $i : U \hookrightarrow F$ un ouvert non-vide de $F$ lisse, et
$F' \hookrightarrow F$ son compl\'ementaire r\'eduit.
On dispose alors d'un morphisme de suites exactes longues
$$\xymatrix{
\bf G\mit_{0}(f^{-1}(F'))_{\bf Q} \ar[r] \ar[d]_{f_{*}} & \bf
G\mit_{0}(X)_{\bf Q} \ar[r] \ar[d]^{f_{*}} & \bf
G\mit_{0}(f^{-1}(U))_{\bf Q} \ar[r] \ar[d]_{f_{*}} & 0 \ar[d] \\
\bf G\mit_{0}(F')_{\bf Q} \ar[r] &
\bf G\mit_{0}(F)_{\bf Q} \ar[r] &
\bf G\mit_{0}(U)_{\bf Q} \ar[r] & 0 }$$
Le "lemme des cinq" implique donc qu'il nous suffit de consid\'erer le
cas de $U$. On suppose donc que $F$ est lisse. Mais alors on a d\'ej\`a vu que
$f_{*}$ est surjectif. $\Box$\\

$\Box$\\

Pour la suite on suppose que $S=Spec k$, avec $k$ un corps contenant
les racines de l'unit\'e.

\begin{thm}\label{th3.4}{(~Grothendieck-Riemann-Roch~)}
Pour chaque $F$ objet de $\cal QDM$, il existe un unique morphisme
$$\tau_{F} : \underline{\bf G}_{*}(F) \longrightarrow H_{\bullet}(F,*)$$
tel que~:
\begin{enumerate}
\item
Si $F$ est lisse dans $\cal QDM$, alors
$$\begin{array}{cccc}
\tau_{F} : & \underline{\bf G}_{*}(F) & \longrightarrow &
H_{\bullet}(F,*) \\
   & x & \mapsto & Td(T_{F}).Ch(x)
\end{array}$$
De m\^eme, si $X$ est un sch\'ema quasi-projectif, $\tau_{X}$ coincide
avec le morphisme d\'efini dans \cite[$4.1$]{g}.
\item
Pour tout morphisme propre de $\cal QDM$, $f : F \longrightarrow F'$, on a
$$f_{*}\circ \tau_{F} = \tau_{F'}\circ f_{*}$$
\item
Si $f : F \longrightarrow F'$ est un morphisme repr\'esentable et \'etale
de champs de $\cal QDM$, alors
$$\tau_{F}\circ f^{*}(x) = f^{*} \circ \tau_{F'}(x)$$
pour tout $x \in \underline{\bf G}_{0}(F')$.
\item
Pour tout champ $F$ de $\cal QDM$, tout $x \in \underline{\bf G}_{0}(F)$ et
tout $y \in \underline{\bf K}_{*}(F)$, on a
$$\tau_{F}(x.y)=\tau_{F}(x).Ch(y)$$
\item
La transformation naturelle $\tau_{F} : \underline{\bf G}_{0}(F)
\longrightarrow H_{\bullet}(F,*)$ se prolonge de fa\c{c}on
unique \`a une transformation naturelle de foncteurs covariants sur $(DM,pr.)$,
la sous-cat\'egorie des champs de Deligne-Mumford et morphismes propres.
\end{enumerate}
\end{thm}

\underline{\bf Preuve:} \rm La preuve suit le m\^eme principe que celle
de \ref{th3.3}. \\

\underline{\bf D\'efinition de $\tau_{F}$:} \\

Soit $F$ un champ de $\cal QDM$, et $p : F \longrightarrow M$ la
projection sur son espace de modules. On d\'efinit $\tau_{F}$ par le
carr\'e commutatif suivant
$$\xymatrix{
\underline{\bf G}_{*}(F) \ar[r]^-{\tau_{F}} \ar[d]_-{p_{*}} &
H_{\bullet}(F,*) \ar[d]^-{p_{*}}\\
\bf G\mit(M)_{\bf Q} \ar[r]_-{\tau_{M}} & H_{\bullet}(M,*)}$$
o\`u $\tau_{F}$ est le morphisme d\'efini dans \cite[$4.1$]{g}. Cette d\'efinition
poss\`ede un sens car, d'apr\`es \ref{c2.2}, le morphisme $p_{*}$ est
bijectif. \\

\underline{\bf Preuve de $(1)$ :} \\

Commen\c{c}ons par le cas o\`u $F$ est un orbifold. D'apr\`es \ref{l3.9}, on sait
que tout morphisme propre et repr\'esentable $F' \longrightarrow F$ est
en r\'ealit\'e fortement projectif.

Soit $f : X \longrightarrow F$ un morphisme propre et surjectif avec
$X$ un sch\'ema quasi-projectif et lisse (~\cite{jo}~). Alors, comme $F$ est
lisse, la formule de projection implique que le morphisme
$$f_{*} : \bf G\mit_{*}(X)_{\bf Q} \longrightarrow \underline{\bf
G}_{*}(F)$$
est surjectif. Notons $q : X \longrightarrow M$ le morphisme induit.

Si $x \in \underline{\bf G}_{*}(F)$, on peut \'ecrire
$x=f_{*}(y)$, avec $y \in \bf G\mit_{*}(X)_{\bf Q}$. Alors, par
d\'efinition de $\tau_{F}$, et par la formule de
Grothendieck-Riemann-Roch pour les sch\'emas (~\cite[$4.1$]{g}~)
$$\begin{array}{cl}
\tau_{F}(x) & = (p_{*})^{-1}\tau_{M}(q_{*}(y)) \\
& = (p_{*})^{-1}q_{*}(\tau_{X}(y)) \\
& = (p_{*})^{-1}q_{*}(Td(X).Ch(y)) \\
& = f_{*}(Td(X).Ch(y))
\end{array}$$
Or, une application de la formule de Grothendieck-Riemann-Roch
au morphisme $f$, qui est fortement
projectif, donne \\
$f_{*}(Td(X).Ch(y))=Td(T_{F}).Ch(x)$. \\

Dans le cas d'un champ lisse $F$ de $\cal QDM$, on peut trouver un
morphisme
$$q : F \longrightarrow F_{0}$$
qui fait de $F$ une gerbe sur $F_{0}$, et avec $F_{0}$ un orbifold
lisse. Remarquons alors que le morphisme naturel
$$Tq : T_{F} \longrightarrow q^{*}T_{F_{0}}$$
est un isomorphisme.

Le lecteur v\'erifiera que l'on peut utiliser le point $(2)$ sans tomber
dans un cercle vicieux.

Notons alors, $p : F \longrightarrow M$, et $p_{0} : F_{0}
\longrightarrow M$ les projections sur l'espace de modules. On a
\'evidemment $p_{0} \circ q = p$. Ainsi, par $(2)$ et le cas pr\'ec\'edent
$$\begin{array}{cl}
\tau_{F}(x) & = p_{*}^{-1}(\tau_{M}(p_{*}(x))) \\
            & =
            q_{*}^{-1}(p_{0})_{*}^{-1}(\tau_{M}((p_{0})_{*}q_{*}(x))) \\
            & = q_{*}^{-1}(\tau_{F_{0}}(q_{*}(x))) \\
            & = q_{*}^{-1}(Td(T_{F_{0}}).Ch(q_{*}(x)))
\end{array}$$

\begin{lem}\label{l3.13}
Soit $q : F \longrightarrow F_{0}$ un morphisme de champs de
Deligne-Mumford, lisses et connexes, qui fasse de $F$ une gerbe born\'ee
par un groupe fini $H$ sur $F_{0}$. Alors
$$q^{*}q_{*} : H^{\bullet}(F,*) \longrightarrow
H^{\bullet}(F,*)$$
est la multiplication par $\frac{1}{m}$, o\`u $m$ est l'ordre de $H$.
\end{lem}

\underline{\bf Preuve:} \rm Consid\'erons le
diagramme cart\'esien
$$\xymatrix{
F \ar[r]^{q} & F_{0} \\
F_{1} \ar[u]^{a} \ar[r]_{b} & F \ar[u]_{q} }$$
Comme $q$ est \'etale, la formule de transfert \ref{p2.2} implique que
$$q^{*}q_{*}(x)=b_{*}a^{*}(x)$$
pour tout $x \in H^{\bullet}(F,*)_{\bf Q}$.
Or $b$ est une gerbe triviale de groupe $H$. Ainsi, $F_{1}\simeq
F\times [Spec k/H]$. La formule de transfert nous ram\`ene alors au cas
o\`u $q : F=[Spec k/H] \longrightarrow F_{0}=Spec k$ est la projection
naturelle. Soit $e : F_{0} \longrightarrow F$ le $H$-torseur
universel, qui est une section de $q$.
Alors, $y=\frac{1}{m}.e^{*}(x)$ est un ant\'ec\'edent de $x$
par $q_{*}$. Ainsi
$$q^{*}q_{*}(x)=\frac{1}{m}.(e\circ q)^{*}(x)=\frac{1}{m}.x$$
$\Box$\\

Comme on sait que $q_{*}$ est un isomorphisme, le lemme ci-dessus
montre que
$$q^{*}=\frac{1}{m}.(q_{*})^{-1}$$
$$q_{*}=\frac{1}{m}.(q^{*})^{-1}$$
Donc
$$\begin{array}{cl}
\tau_{F}(x) & = q_{*}^{-1}(Td(T_{F_{0}}).Ch(q_{*}(x))) \\
            & = m.q^{*}(Td(T_{F_{0}}).Ch(q_{*}(x))) \\
            & = m.Td(q^{*}(T_{F_{0}})).q^{*}Ch(q_{*}(x)) \\
            & = m.Td(T_{F}).Ch(q^{*}q_{*}(x)) \\
            & = m.Td(T_{F}).Ch(\frac{1}{m}.x) \\
            & = Td(T_{F}).Ch(x)
\end{array}$$

Pour finir, le cas d'un sch\'ema $X$ est \'evident, car $M=X$ et la
projection naturelle est $p=Id : X \longrightarrow X$. \\

\underline{\bf Preuve de $(2)$ :} \\

Soit $f : F \longrightarrow F'$ un morphisme propre de champs de
$\cal QDM$. Notons $p : F \longrightarrow M$ et $p' : F'
\longrightarrow M'$ les projections sur les espaces de modules, ainsi
que $Mf : M \longrightarrow M'$ le morphisme induit. Alors, pour $x
\in \underline{\bf G}_{*}(F)$, on a
$$\begin{array}{cl}
\tau_{F'}(f_{*}(x)) & = (p'_{*})^{-1}\tau_{M'}(p'_{*}f_{*}(x)) \\
& = (p'_{*})^{-1}\tau_{M'}(Mf_{*}p_{*}(x)) \\
& = (p'_{*})^{-1}Mf_{*}\tau_{M}(p_{*}(x)) \\
& = f_{*}(p_{*})^{-1}\tau_{M}(p_{*}(x)) \\
& = f_{*}\tau_{F}(x)
\end{array}$$

\underline{\bf Preuve de $(3)$ :} \\

Soit $u : F \longrightarrow F'$ repr\'esentable et \'etale.

\begin{lem}\label{l3.14}
Soit $f : X \longrightarrow F$ un morphisme propre, tel que pour tout
sous-champ irr\'eductible $F' \hookrightarrow F$, il existe une
composante irr\'eductible $X'$ de $X$ au-dessus de $F'$,
telle que $f : X' \longrightarrow F'$ soit g\'en\'eriquement fini.

Alors le morphisme
$$f_{*} : \bf G\mit_{0}(X)_{\bf Q} \longrightarrow \underline{\bf
G}_{0}(F)$$
est surjectif.
\end{lem}

\underline{\bf Preuve :} \rm On proc\`ede par r\'ecurrence noeth\'erienne dans $F$.

Soit $i : U \hookrightarrow F$ un ouvert non-vide de $F$ lisse, et
$F' \hookrightarrow F$ son compl\'ementaire r\'eduit.
On dispose alors d'un morphisme de suite exactes longues
$$\xymatrix{
\bf G\mit_{0}(f^{-1}(F'))_{\bf Q} \ar[r] \ar[d]_{f_{*}} & \bf
G\mit_{0}(X)_{\bf Q} \ar[r] \ar[d]^{f_{*}} & \bf
G\mit_{0}(f^{-1}(U))_{\bf Q} \ar[r] \ar[d]_{f_{*}} & 0 \ar[d] \\
\underline{\bf G}_{0}(F') \ar[r] &
\underline{\bf G}_{0}(F) \ar[r] &
\underline{\bf G}_{0}(U) \ar[r] & 0 }$$
Le "lemme des cinq" implique donc qu'il nous suffit de consid\'erer le
cas de $U$. On suppose donc que $F$ est lisse. Mais alors $f : X
\longrightarrow F$ \'etant propre et surjectif, on a d\'ej\`a vu que
$f_{*}$ est surjectif. $\Box$\\

Soit $x \in \underline{\bf G}_{0}(F')$, et $f' : X'
\longrightarrow F'$ un morphisme comme dans le lemme, avec $y \in
\bf G\mit_{0}(X')_{\bf Q}$ et $f'_{*}(y)=x$.
Notons $f : X=X'\times_{F'}F \longrightarrow F$, et
$v : X \longrightarrow X'$ les morphismes induits. Alors, \`a l'aide de
$(2)$ et du cas des sch\'emas \cite{g}, on a
$$\begin{array}{cl}
u^{*}\tau_{F'}(x) & = u^{*}\tau_{F'}(f'_{*}(y)) \\
                  & = u^{*}f'_{*}\tau_{X'}(y) \\
                  & = f_{*}v^{*}\tau_{X'}(y) \\
                  & = f_{*}\tau_{X}(v^{*}(y)) \\
                  & = \tau_{F}(f_{*}v^{*}(y)) \\
                  & = \tau_{F}(u^{*}f'_{*}(y)) \\
                  & = \tau_{F}(u^{*}(x))
\end{array}$$

\underline{\bf Preuve de $(4)$ :} \\

Soit $x \in \underline{\bf G}_{0}(F)$ et $y \in \underline{\bf K}_{*}(F)$.
Choisissons un morphisme propre et surjectif
$$f : X \longrightarrow F$$
avec $X$ un sch\'ema quasi-projectif, tel que
$x=f_{*}(z)$, avec $z \in \bf G\mit_{0}(X)_{\bf Q}$.
Un tel morphisme existe d'apr\`es \cite[$2.6$]{vi2}, et le lemme \ref{l3.14}.
On a
alors $x.y=f_{*}(x.f^{*}(y))$, et donc, en utilisant le point $(2)$
et le cas des sch\'emas
$$\begin{array}{cl}
\tau_{F}(x.y)& = \tau_{F}(f_{*}(z.f^{*}(y))) \\
             & = f_{*}\tau_{X}(z.f^{*}(y)) \\
             & =f_{*}\tau_{X}(z).f^{*}Ch(y) \\
             & = f_{*}(\tau_{X}(z)).Ch(y) \\
             & = \tau_{F}(f_{*}(z)).Ch(y) \\
             & = \tau_{F}(x).Ch(y)
\end{array}$$

\underline{\bf Preuve de $(5)$ :} \\

Soit $F$ un champ de Deligne-Mumford, et $p : F \longrightarrow M$ la
projection sur son espace de modules. Alors,
on d\'efinit $\tau_{F}$ par le diagramme commutatif suivant
$$\xymatrix{
\underline{\bf G}_{0}(F) \ar[r]^-{\tau_{F}} \ar[d]_-{p_{*}} &
H_{\bullet}(F,*) \ar[d]^-{p_{*}}\\
\bf G\mit_{0}(M)_{\bf Q} \ar[r]_-{\tau_{M}} & H_{\bullet}(M,*)}$$
o\`u $\tau_{M}$ est d\'efini dans \cite[$7$]{g2}. Il est clair que cette
d\'efinition r\'epond aux conditions demand\'ees. $\Box$\\

Les th\'eor\`emes \ref{th3.3} et \ref{th3.4} peuvent se composer. On obtient de
cette fa\c{c}on le th\'eor\`eme de Grothendieck-Riemann-Roch sous sa forme
finale.

\begin{thm}\label{th3.5}{(~Grothendieck-Riemann-Roch~)}
Pour chaque $F$ objet de $\cal QDM$, il existe un unique morphisme
$$\tau^{rep}_{F} : \bf G\mit_{*}(F) \longrightarrow
H_{\bullet}^{rep}(F,*)$$
tel que~:
\begin{enumerate}
\item
Si $F$ est lisse dans $\cal QDM$, alors
$$\begin{array}{cccc}
\tau^{rep}_{F} : & \bf G\mit_{*}(F) & \longrightarrow &
H_{\bullet}^{rep}(F,*) \\
   & x & \mapsto & Td^{rep}(F).Ch^{rep}(x)
\end{array}$$
De m\^eme, si $X$ est un sch\'ema quasi-projectif, $\tau^{rep}_{X}$ coincide
avec le morphisme d\'efini dans \cite[$4.1$]{g}.
\item
Pour tout morphisme propre de $\cal QDM$, $f : F \longrightarrow F'$, on a
$$f_{*}\circ \tau^{rep}_{F} = \tau^{rep}_{F'}\circ f_{*}$$
dans l'un des deux cas suivants
\begin{itemize}
\item le morphisme $f$ est repr\'esentable
\item le champ $F$ est lisse
\end{itemize}
\item
Si $f : F \longrightarrow F'$ est un morphisme repr\'esentable et \'etale
de champs de $\cal QDM$, alors
$$\tau^{rep}_{F}\circ f^{*}(x) = f^{*} \circ \tau^{rep}_{F'}(x)$$
pour tout $x \in \bf G\mit_{0}(F')$.
\item
Pour tout champ $F$ de $\cal QDM$, tout $x \in \bf G\mit_{0}(F)$ et
tout $y \in \bf K\mit_{*}(F)$, on a
$$\tau^{rep}_{F}(x.y)=\tau^{rep}_{F}(x).Ch^{rep}(y)$$
\item
La transformation naturelle $\tau^{rep}_{F} : \bf G\mit_{0}(F)
\longrightarrow H_{\bullet}^{rep}(F,*)$ se prolonge de fa\c{c}on
unique \`a une transformation naturelle de foncteurs covariants sur
$(\cal DM\mit,pr.)$,
la sous-cat\'egorie des champs de Deligne-Mumford et morphismes propres.
\end{enumerate}
\end{thm}

\underline{\bf Preuve :} \rm On d\'efinit $\tau_{F}^{rep}$ par la
composition
$$\tau_{F}^{rep} : \xymatrix{
\bf G\mit_{*}(F) \ar[r]^{\psi_{F}} & \underline{\bf G}(I_{F}^{t})_{K}
\ar[r]^-{\tau_{I_{F}^{t}}} & H_{\bullet}(I_{F}^{t},*)=H^{rep}_{\bullet}(F,*)}$$
Les cinq assertions du th\'eor\`eme sont alors obtenues en composant
les cinq assertions des th\'eor\`emes \ref{th3.3} et \ref{th3.4}. $\Box$\\

En utilisant le th\'eor\`eme \ref{th2.3'}, ainsi que les lemmes sur
l'existence des quasi-enveloppes de Chow (~\ref{th1.4} $(2)$~), on peut
montrer le th\'eor\`eme
suivant, dont nous ferons l'\'economie de la d\'emonstration (~qui est
tout \`a fait analogue au cas des champs de Deligne-Mumford~).

Nous supposerons que $S$ contient les racines de l'unit\'es, et que nous
avons fix\'e un plongement $\mu_{\infty}(S)\hookrightarrow \bf C\mit^{*}$.
On posera alors $K:=\bf Q\mit(\mu_{\infty}(S))$.

Dans le th\'eor\`eme suivant, nous avons bien entendu pos\'e les d\'efinitions
suivantes~:
$$H^{p}_{rep}(F,q):=\pi_{dq-p}\bf H\mit((I_{F}^{t,f})_{li},\cal
H\mit^{q}\otimes K)$$
$$H_{p}^{rep}(F,q):=\pi_{dq-p}\bf H\mit((I_{F}^{t,f})_{li},\cal
H\mit'_{q}\otimes  K).$$

\begin{thm}\label{th3.5'}
Supposons que $S$ soit de caract\'eristique nulle.
Notons $\cal QGA$ la sous-cat\'egorie de $HoChAlg(S)$, des champs
lisses sur $S$, $\Delta$-\'equidimensionnels et pseudo-s\'epar\'es (\ref{ps}), et
$\cal QQGA$ celle de $\cal QGA$ form\'ees des champs dont l'espace de modules
est quasi-projectif. \\

Pour chaque $F$ objet de $\cal QGA$ notons
$$\begin{array}{cccc}
\tau^{rep}_{F} : & \bf G\mit_{*}(F) & \longrightarrow &
H_{\bullet}^{rep}(F,*) \\
   & x & \mapsto & Td^{rep}(F).Ch^{rep}(x)
\end{array}$$
Alors, on a les propri\'et\'es suivantes.
\begin{enumerate}
\item
Pour tout morphisme propre de $\cal QGA$, $f : F \longrightarrow F'$, on a
$$f_{*}\circ \tau^{rep}_{F} = \tau^{rep}_{F'}\circ f_{*}.$$
\item
Si $f : F \longrightarrow F'$ est un morphisme repr\'esentable et \'etale
de champs de $\cal QGA$, alors
$$\tau^{rep}_{F}\circ f^{*}(x) = f^{*} \circ \tau^{rep}_{F'}(x)$$
pour tout $x \in \bf G\mit_{0}(F')$.
\item
Pour tout champ $F$ de $\cal QGA$, tout $x \in \bf G\mit_{0}(F)$ et
tout $y \in \bf K\mit_{*}(F)$, on a
$$\tau^{rep}_{F}(x.y)=\tau^{rep}_{F}(x).Ch^{rep}(y)$$
\end{enumerate}
\end{thm}

\end{subsubsection}

\begin{subsubsection}{Exemples d'applications}
\hspace{5mm}
On notera $k$ un corps qui contient les racines de l'unit\'e. \\

Nous allons donner toute une s\'erie de corollaires du th\'eor\`eme
\ref{th3.5}. Nous nous int\'eresserons tout particuli\`erement aux
formules de type
Gauss-Bonnet. \\

Pour un champ alg\'ebrique $F$, nous noterons
$$A^{p}(F):=H^{p}(F,\cal K\mit_{p}\otimes \bf Q\mit)$$
$$A_{p}(F):=H^{p}(F,\cal R\mit^{p}\otimes \bf Q\mit)$$
$$A^{p}_{rep}(F):=H^{p}(I_{F}^{t},\cal K\mit_{p}\otimes K)$$
$$A_{p}^{rep}(F):=H^{p}(I_{F}^{t},\cal R\mit^{p}\otimes K)$$
o\`u $\cal R\mit^{p}$ est le complexe de Gersten de codimension
$p$. Notons, qu'aux graduations pr\`es, ces groupes sont des morceaux de
la cohomologie et de l'homologie \`a coefficients dans les repr\'esentations
de $F$.

\begin{cor}
Si $F$ est un champ de Deligne-Mumford, alors
les morphismes
$$\tau_{F} : \underline{\bf G}_{0}(F) \longrightarrow
\bigoplus_{p}A_{p}(F)$$
$$\tau_{F}^{rep} : \bf G\mit_{0}(F)_{K} \longrightarrow
\bigoplus_{p}A_{p}^{rep}(F)$$
sont des isomorphismes.

En particulier, si $F$ est lisse dans $\cal QDM$, alors
$$Ch : \underline{\bf G}_{0}(F) \longrightarrow
\bigoplus_{p}A^{p}(F)$$
$$Ch^{rep} : \bf G\mit_{0}(F)_{K} \longrightarrow
\bigoplus_{p}A^{p}_{rep}(F)$$
sont des isomorphismes.
\end{cor}

\underline{\bf Preuve:} \rm La seconde partie se d\'eduit de la premi\`ere
en remarquant que les classes de Todd $Td(T_{F})$ et $Td^{rep}(F)$
sont des \'el\'ements inversibles.

Comme $\tau_{F}^{rep}=\tau_{F}\circ \psi_{F}$, et que $\psi_{F}$ est
un isomorphisme, il suffit de d\'emontrer que $\tau_{F}$ est un
isomorphisme. Soit $p : F \longrightarrow M$ la projection sur
l'espace de modules. Le th\'eor\`eme de Riemann-Roch montre alors que l'on
a un carr\'e commutatif
$$\xymatrix{
\bf G\mit_{0}(F)_{\bf Q} \ar[r]^{\tau_{F}} \ar[d]_{p_{*}} &
A_{*}(F) \ar[d]^{p_{*}} \\
\bf G\mit_{0}(M)_{\bf Q} \ar[r]_{\tau_{M}} &
A_{*}(M)}$$
Or, on sait que $\tau_{M}$ est un isomorphisme, et donc $\tau_{F}$
aussi. $\Box$\\

\underline{Remarque:} En utilisant les complexes de Chow-Bloch
d\'efinis dans \cite{b}, on peut d\'emontrer des r\'esultats analogues pour
la $G$-th\'eorie sup\'erieure. \\

Si $F$ est un champ propre sur un $Spec k$, nous noterons
$$\int_{F}^{rep}:=p_{*} : H^{rep}_{\bullet}(F,*) \longrightarrow
H^{rep}_{\bullet}(Spec k,*)$$
o\`u $p$ est le morphisme structural.

De m\^eme, nous noterons
$$\int_{F}:=p_{*} : H_{\bullet}(F,*) \longrightarrow
H_{\bullet}(Spec k,*)$$

\begin{cor}\label{HRR}{(~Hirzebruch-Riemann-Roch~)}
Soit $F$ un champ alg\'ebrique de Deligne-Mumford mod\'er\'e et propre sur $k$.
Alors, pour tout faisceau coh\'erent $\cal F$ sur $F$, on a
$$\chi(F,\cal F\mit):=\sum_{i}(-1)^{i}Dim_{k}H^{i}(F,\cal F\mit)
=\int_{F}^{rep}\tau_{F}^{rep}(\cal F\mit)$$
En particulier si $F$ est lisse et dans $\cal QDM$, on a
$$\chi(F,\cal F\mit)=\int_{F}^{rep}Td^{rep}(F).Ch^{rep}(\cal F\mit)$$
\end{cor}

\underline{\bf Preuve : } \rm Il suffit d'appliquer \ref{th3.5} \`a la
projection $p : F \longrightarrow Spec k$, et remarquer que si un champ
de Deligne-Mumford est mod\'er\'e sur $Spec k$, alors $p$ est de
dimension cohomologique finie. $\Box$\\

Dans le cas des orbifolds de dimension $1$, cette formule peut
s'expliciter. On retrouve alors la formule de Riemann-Roch pour les
courbes orbifolds complexes d\'emontr\'ee dans \cite{ta}. Pour
simplifier, on supposera que $k$ est alg\'ebriquement clos.

Soit $F$ est un orbifold lisse propre de dimension $1$, et mod\'er\'e sur
$k$, et $M$ son espace de modules. Alors, $M$ est une courbe lisse et
projective sur $k$. Notons $x_{1}, \dots ,x_{m}$ les points ferm\'es de
$M$ au-dessus desquels $F$ n'est pas un sch\'ema, et
$A_{i}:=\cal O\mit_{M,x_{i}}^{h}$ l'hens\'elis\'e de l'anneau local de
$M$ en $x_{i}$. Alors, la restriction de $F$ au-dessus de
$Spec A_{i}$, est \'equivalent \`a un quotient de la forme
$[X_{i}/H_{i}]$, avec $X_{i}$ le spectre d'une $k$-alg\`ebre locale
hens\'elienne et r\'eguli\`ere de dimension $1$, et $H_{i}$ un groupe fini
op\'erant sur $X_{i}$, et ne fixant que le point ferm\'e $y_{i} \in X_{i}$.
Ainsi, $X_{i} \longrightarrow Spec A_{i}$ d\'efinit un rev\^etement fini,
mod\'er\'ement ramifi\'e en $x_{i}$, car $F$ est mod\'er\'e sur $k$. Il
correspond donc a un quotient de noyau ouvert
$$\widehat{\pi}_{1}^{t}(Spec A_{i},x_{i}) \longrightarrow H_{i}$$
Or, on sait que
$$\widehat{\pi}_{1}^{t}(Spec A_{i},x_{i})\simeq
\widehat{\bf Z\mit[\frac{1}{p}]}$$
Ainsi, on a forc\'ement $H_{i}\simeq \bf Z\mit/r_{i}$, o\`u $r_{i}$ est
l'ordre de ramification de $F$ en $x_{i}$ (~\ref{d1.11}~). Fixons nous des
sections $s_{i} : Spec k(x_{i}) \longrightarrow \widetilde{x_{i}}$,
ainsi que des g\'en\'erateurs $h_{i} \in H_{i}$.

Si $\pi_{i}$ est une uniformisante en $y_{i}$, l'action de
$H_{i}=<h_{i}>$, est donn\'ee par
$$\begin{array}{cccc}
h_{i} : & B_{i} & \longrightarrow & B_{i} \\
        & \pi_{i} & \mapsto & \zeta_{i}.\pi_{i}
\end{array}$$
o\`u $\zeta_{i}$ est une racine primitive $r_{i}$-\`eme de l'unit\'e dans
$k$.

Soit $T_{F,x_{i}}:=s_{i}^{*}i_{x_{i}}^{*}T_{F}$ la restriction
du fibr\'e tangent sur $Spec k(x_{i})$. Alors, $T_{F,x_{i}}$ est
un espace vectoriel de dimension $1$, sur lequel $h_{i}$ op\`ere par
multiplication par $\zeta_{i}$. Maintenant, si
$\cal L$ est un fibr\'e inversible sur $F$, sa restriction sur $Spec
k(x_{i})$, donne une repr\'esentations de $H_{i}$ sur un \\
$k$-espace
vectoriel de dimension $1$, $\cal
L\mit_{x_{i}}:=s_{i}^{*}i_{x_{i}}^{*}\cal L$. Notons $k_{i}$
l'entier compris entre $0$ et $r_{i}-1$, tel que $h_{i}$ op\`ere par
multiplication par $\zeta_{i}^{k_{i}}$ sur $\cal L\mit_{x_{i}}$.

\begin{df}{(~\cite{ta}~)}
L'entier $k_{i}$ est appel\'e la multiplicit\'e de $\cal L$ en $x_{i}$.
\end{df}

\begin{cor}{(~Riemann-Roch~)}{(~\cite{ta}~)}
Soit $F$ un orbifold lisse, propre, et mod\'er\'e sur $k$
alg\'ebriquement clos, de dimension $1$. Soit $M$ son espace de
modules, et $x_{i} \in M$ les points de ramifications de $F$, d'ordre
$r_{i}$.

Soit $\cal L$ un fibr\'e
inversible sur $F$, et $k_{i}$ la multiplicit\'e de $\cal L$ en $x_{i}$.
Alors on a
$$\chi(F,\cal L\mit)=\int_{F}C_{1}(\cal
L\mit)+1-g-\sum_{i}\frac{k_{i}}{r_{i}}$$
\end{cor}

\underline{\bf Preuve:} \rm On utilisera la th\'eorie de Gersten comme
th\'eorie cohomologique. \\

Le champ $I_{F}$ s'\'ecrit comme une
r\'eunion disjointe
$$I_{F}\simeq F \coprod_{i=1}^{m}BH_{i}$$
et donc
$$H_{rep}^{0}(F,0)\simeq H^{0}(F,0)
\bigoplus_{i=1}^{m}K(H_{i})$$
$$H_{rep}^{1}(F,1)\simeq H^{1}(F,1)$$
Ainsi, l'anneau gradu\'e $H^{*}_{rep}(F,*)$ est isomorphe \`a
$$H^{0}(F,0) \oplus H^{1}(F,1)
\bigoplus_{i=1}^{m}K(H_{i})$$
A travers cet isomorphisme, le caract\`ere de Chern de $\cal L$ est
donn\'e par la formule suivante
$$Ch^{rep}(\cal L\mit)=Ch(\cal L\mit)\bigoplus_{i=1}^{m}
\chi(\cal L\mit_{x_{i}})$$
o\`u $\chi(\cal L\mit_{x_{i}}) \in K(H_{i})$ est le
caract\`ere de la repr\'esentation de $H_{i}$ sur $\cal L\mit_{x_{i}}$.
Ainsi, avec les notations pr\'ec\'edentes, on a
$$\begin{array}{cccc}
\chi(\cal L\mit_{x_{i}}) : & H_{i} & \longrightarrow & K
\\
& h_{i}^{a} & \mapsto & \zeta_{i}^{a.k_{i}}
\end{array}$$
Quand \`a la classe de Todd de $F$, elle s'\'ecrit
$$Td^{rep}(F)=Td(T_{F})\bigoplus_{i=1}^{m}Td_{i}$$
o\`u
$$\begin{array}{cccc}
Td_{i} : & H_{i} & \longrightarrow & K
\\
& h_{i}^{a} & \mapsto & \frac{1}{1-\zeta_{i}^{a}}
\end{array}$$
Ainsi, la classe $\tau_{F}^{rep}(\cal L\mit)$ est
$$\tau_{F}(\cal L\mit)\bigoplus_{i=1}^{m}\chi(\cal
L\mit_{x_{i}}).Td_{i}$$
et donc, la formule \ref{HRR} implique que
$$\chi(F,\cal L\mit)=\int_{F}\tau_{F}(\cal L\mit)+
\sum_{i}\frac{1}{r_{i}}.\sum_{a=0}^{r_{i}-1}
\frac{\zeta_{i}^{a.k_{i}}}{1-\zeta_{i}^{a}}$$
Or, $\tau_{F}(\cal L\mit)=C_{1}(\cal L\mit)+\frac{1}{2}.C_{1}(T_{F})$.
La formule de Gauss-Bonnet \ref{gb1} (~ou une application de
\ref{HRR} au fibr\'e trivial~) donne
$$\int_{F}\tau_{F}(\cal L\mit)=\int_{F}C_{1}(\cal L\mit)+
1-g+\sum_{i}\frac{r_{i}-1}{2.r_{i}}$$
On conclut alors \`a l'aide de la formule
$$\sum_{a=0}^{r_{i}-1}\frac{\zeta_{i}^{a.k_{i}}}{1-\zeta_{i}^{a}}
=-\frac{r_{i}-1}{2}-k_{i}$$
ou encore
$$\sum_{a=0}^{r_{i}-1}\frac{\zeta_{i}^{a.k_{i}}+1}{1-\zeta_{i}^{a}}
=-k_{i}$$
$\Box$\\

\begin{cor}
Soit $F$ un champ de Deligne-Mumford propre et mod\'er\'e
sur $Spec k$, et $M$ son espace de modules. Alors, le genre
arithm\'etique de $M$ est donn\'e par la formule suivante
$$\chi^{a}(M):=\chi(M,\cal O\mit_{M})=\int_{F}^{rep}\tau_{F}^{rep}(1)$$
En particulier, si $F$ est lisse et dans $\cal QDM$, on a
$$\chi^{a}(M)=\int_{F}^{rep}Td^{rep}(F)$$
\end{cor}

\underline{\bf Preuve :} \rm Commen\c{c}ons par remarquer que, comme $F$
est mod\'er\'e, le foncteur $q_{*}$ (~o\`u $q$ est la projection canonique
$q : F \longrightarrow M$~) est exact.

En proc\'edant par r\'ecurrence sur la dimension du support des faisceaux
coh\'erents $\bf R\mit p_{*}^{i}(\cal F\mit)$, on peut se restreindre \`a
un sous-champ ouvert non-vide de $F$. Comme on peut aussi supposer que $F$ est
r\'eduit (~car $F_{red}~\hookrightarrow~F$ est acyclique pour les
faisceau coh\'erents~), on se ram\`ene au cas o\`u $F$ est une gerbe sur $M$.

Comme l'assertion est locale sur $M_{et}$, on peut supposer que
\mbox{$F=[X/H]$}, avec $X$ un sch\'ema et $H$ un groupe fini op\'erant
trivialement sur $X$.
En identifiant alors la cat\'egorie $\bf Coh\mit([X/H])$, \`a la cat\'egorie
des faisceaux coh\'erents sur $X$, munit d'une action de $H$, le foncteur
$$b_{*} : \bf Coh\mit([X/H]) \longrightarrow \bf Coh\mit(X/H)$$
associe \`a un faisceau coh\'erent $\cal F$ son sous-faisceau des
invariants par $H$,  $\cal F\mit^{H}$. Mais, comme $F$ est mod\'er\'e
sur $Spec k$, l'ordre de $H$ est premier avec la caract\'eristique de
$k$, et donc le foncteur $\cal F\mit \longrightarrow \cal F\mit^{H}$
est exact.\\

Comme $q_{*}$ est exact, on a
$$\chi(F,\cal O\mit_{F})=\chi(M,q_{*}\cal O\mit_{F})$$
mais $q_{*}(\cal O\mit_{F})=\cal O\mit_{M}$, et donc
$$\chi(F,\cal O\mit_{F})=\chi(M,\cal
O\mit_{M})=\int_{F}^{rep}\tau_{F}^{rep}(1)$$
$\Box$\\

\underline{Remarque :} La formule pr\'ec\'edente est une g\'en\'eralisation au cas des
champs singuliers, ainsi qu'\`a la caract\'eristique quelconque, de la formule
obtenue par \cite{ka}. Son int\'er\^et principal est qu'elle permet de
calculer le genre arithm\'etique d'un sch\'ema \'eventuellement singulier
$M$, en fonction de certaines classes de cohomologie sur un objet
lisse $F$.

Les formules de Gauss-Bonnet qui vont suivre poss\`edent le m\^eme int\'er\^et.

\begin{df}\label{d3.7}
\begin{itemize}
Soit $F$ un champ de Deligne-Mumford sur un corps $k$,
alg\'ebriquement clos et de caract\'eristique nulle.
\item
Nous d\'efinissons sa caract\'eristique d'Euler topologique par
$$\chi^{top}(F):=\sum_{p=0}^{p=Dim_{k}F}(-1)^{p}Dim_{k}H_{p}(F,p)$$
o\`u l'on utilise la cohomologie de De Rham.
\item
Notons $\{M_{i}\}_{i}$ une stratification de son espace de modules par
des sous-espaces localement ferm\'es, lisses et connexes,
tels que la fonction qui \`a un point $x \in M$ associe l'ordre de
ramification de $F$ en $x$ soit constante sur chacun des $M_{i}$,
\'egale \`a $m_{i}$.

La caract\'eristique d'Euler orbifold de $F$ est d\'efinie par
$$\chi^{orb}(F):=\sum_{i}\frac{\chi^{top}(M_{i})}{m_{i}}$$
\item
La caract\'eristique d'Euler des physiciens et d\'efinie par
$$\chi^{phy}(F):=\chi^{top}(I_{F})$$
\end{itemize}
\end{df}

Pour la suite $k$ est sera un corps alg\'ebriquement clos
de caract\'eristique nulle, et contenant
les racines de l'unit\'es.

Pour un \'el\'ement $x$ de $\underline{\bf K}_{0}(F)$, nous noterons
$C_{max}(x)$ sa classe de Chern maximale. Lorsque $F$ est connexe, et
$x$ de rang $r$, \mbox{$C_{max}(x)=C_{r}(x)$}. En g\'en\'eral, $C_{max}(x)$ est
la somme des $C_{r_{i}}(x_{i})$, o\`u $i$ parcourt l'ensemble des
composantes connexes de $F$, et $r_{i}$ est le rang de $x$ sur la
composante $i$.

\begin{cor}{(~Gauss-Bonnet $I$~)}\label{gb1} Soit $F$ un champ alg\'ebrique
de $\cal QDM$
propre et lisse sur $k$, et $M$ son espace de modules.
Alors
$$\chi^{orb}(F):=\int_{F}C_{max}(T_{F})$$
\end{cor}

\underline{\bf Preuve :} \rm On va appliquer le th\'eor\`eme \ref{th3.4}, au cas
de la projection $p : F \longrightarrow M$, et de l'\'el\'ement
$$x=can(\sum_{i}(-1)^{i}\Omega_{F}^{i}) \in \underline{\bf G}_{0}(F)$$
En utilisant le lemme \ref{l3.2}, on peut \'ecrire
$$\tau_{F}(x)=Td(T_{F}).Ch(can(\lambda_{-1}(T_{F})^{-1}))=C_{max}(T_{F})$$
On en d\'eduit alors
$$\int_{F}C_{max}(T_{F})=f_{*}(x)$$
o\`u $f : F \longrightarrow Spec k$. On conclut alors \`a l'aide de
\ref{c4.5} (~le lecteur v\'erifiera qu'il n'y pas de cercle vicieux~!~).
$\Box$\\

\begin{cor}{(~Gauss-Bonnet $II$~)}\label{gb2}
Soit $F$ un champ de Deligne-Mumford, lisse et propre sur $k$.
Alors, la caract\'eristique d'Euler topologique de $F$ est donn\'ee par
$$\chi^{top}(F)=\chi^{top}(M)=\int_{I_{F}}C_{max}(T_{I_{F}})$$
\end{cor}

\underline{\bf Preuve :} \rm On applique la formule
d'Hirzebruch-Riemann-Roch
$$\chi(F,x)=\int_{F}^{rep}Td^{rep}(F).Ch^{rep}(x)$$
avec $x=\lambda_{-1}(\Omega_{F}^{1})$.

Alors, par la suite spectrale d'hyper-cohomologie associ\'e au complexe
de De Rham sur $F$, on trouve
$$\begin{array}{cl}
\chi(F,x) & =\sum_{i}(-1)^{p+q}Dim_{k}H^{p}(F,\Omega_{F}^{q}) \\
& = \sum_{i}(-1)^{i}Dim_{k}H^{i}(F,\Omega_{F}^{\bullet}) \\
& = \chi^{top}(F)
\end{array}$$
Il nous reste \`a \'evaluer $Td^{rep}(F).Ch^{rep}(x)$.

Soit $\pi : I_{F} \longrightarrow F$ la projection naturelle, et
$$0 \longrightarrow \cal N\mit_{\pi}^{\vee} \longrightarrow
\pi^{*}\Omega^{1}_{F} \longrightarrow
\Omega_{I_{F}}^{1} \longrightarrow 0$$
la suite exacte courte des $1$-formes diff\'erentielles. Par le
morphisme $F$ d\'efini dans \ref{th2.3}
$$F : \bf K\mit_{0}(I_{F}) \longrightarrow \bf K\mit_{0}(I_{F})$$
on obtient
$$F(\pi^{*}\Omega_{F}^{1})=F(\Omega^{1}_{I_{F}})+F(\cal
N\mit^{\vee}_{\pi})$$
et donc
$$\lambda_{-1}(F(\pi^{*}\Omega_{F}^{1}))=
\lambda_{-1}(F(\Omega^{1}_{I_{F}})).\lambda_{-1}(F(\cal
N\mit^{\vee}_{\pi}))$$
Or, par d\'efinition
$$\lambda_{-1}(F(\cal N\mit^{\vee}_{\pi}))=\alpha_{F}$$
De plus, pour toute racine de l'unit\'e $\zeta\neq 1$, on a
$$F_{\zeta}(\Omega_{I_{F}}^{1})=0$$
et donc
$$F(\Omega_{I_{F}}^{1})=\Omega_{I_{F}}^{1}$$
Ainsi, on a
$$Ch^{rep}(x)=Ch(\Omega^{1}_{I_{F}}).Ch(\alpha_{F})$$
et donc, par le lemme \ref{l3.2}
$$Ch^{rep}(x).Td(T_{I_{F}})=Ch(\alpha_{F}).C_{max}(T_{I_{F}})$$
Ce qui implique que
$$\begin{array}{cl}
Td^{rep}(F).Ch^{rep}(x) & = Td(T_{I_{F}}).Ch(\alpha_{F})^{-1}.Ch^{rep}(x) \\
& = C_{max}(T_{I_{F}})
\end{array}$$
$\Box$\\

Ces deux derniers corollaires permettent de donner une relation de
r\'ecurrence sur les diff\'erentes caract\'eristiques d'Euler.

Pour cela, si $F$ est un champ, nous d\'efinirons par r\'ecurrence
$$I_{F}^{m}:=I_{I_{F}^{m-1}}$$
avec $I_{F}^{0}:=F$. On dispose alors des relations suivantes.

\begin{cor}
Si $F$ est un champ de Deligne-Mumford, lisse et propre sur $k$,
alors on a
$$\chi^{phy}(I_{F}^{m})=\chi^{top}(I_{F}^{m+1})=\chi^{orb}(I_{F}^{m+2})$$
pour tout $m \geq 0$.
\end{cor}

\underline{\bf Preuve :} \rm On applique les corollaires pr\'ec\'edents. $\Box$\\

Il est \`a noter que lorsque $F=[X/H]$ est un champ quotient par un
groupe fini, $I_{F}^{m}\simeq [X^{(m)}/H]$, o\`u
$$X^{(m)}=\{(x,h_{1},\dots ,h_{m}) \in X\times H^{m} \; / \;
h_{i}.h_{j}=h_{j}.h_{i} \; \forall \; i,j \; et \; h_{i}.x=x \;
\forall \; i\}$$
et l'action de $h \in H$ est donn\'ee par la formule
$$h.(x,h_{1},\dots ,h_{m}):=(h.x,h.h_{1}.h^{-1},\dots
,h.h_{m}.h^{-1})$$
On remarque alors que $\chi^{orb}(I_{F}^{m})$ est la caract\'eristique
d'Euler sup\'erieure $\chi_{m}(X,H)$ d\'efinie dans \cite{brfu}. Ainsi, on
pourrait d\'efinir une "cohomologie \`a coefficients dans les
repr\'esentations" sup\'erieure, pour un champ $F$ propre et lisse, par
$$H^{\bullet}_{rep,m}(F,*):=H^{\bullet}(I_{F}^{m},*)$$
ainsi que les caract\'eristiques d'Euler sup\'erieures
$$\chi_{m}(F):=\sum_{i}Dim_{k}H^{i}_{rep,m}(F,0)$$
Dans chacun des anneaux $H^{\bullet}_{rep,m}(F,*)$ on dispose d'une
classe d'Euler
$$Eu_{m}(F):=C_{max}(T_{I_{F}^{m}})\in H^{\bullet}(I_{F}^{m},*)$$
La relation du corollaire pr\'ec\'edent se traduit alors par
$$\chi_{m}(F)=\int_{F}^{m+1}Eu_{m+1}(F)$$
o\`u $\int_{F}^{m}=p_{*}^{(m)}$, avec
$p^{(m)} : I_{F}^{m} \longrightarrow Spec k$ le morphisme structural.

Dans \cite{brfu}, les nombres $\chi_{m}(X,H)$ sont \'etudi\'es \`a l'aide d'une s\'erie
g\'en\'eratrice
$$f(X,H;t):=\sum_{m}\chi_{m}(X,H).t^{m}$$
Suivant cette id\'ee, pour un champ $F$ propre et lisse,
on peut d\'efinir une cohomologie rassemblant toute
les cohomologies $H^{\bullet}_{rep,m}(F,*)$, en posant
$$H^{\bullet}_{\infty}(F,*):=\prod_{m}H^{\bullet}_{rep,m}(F,*)$$
Alors, la classe $\sum_{m}Eu_{m}(F)$ d\'efinit une classe "d'Euler
totale"
$$Eu_{\infty}(F):=\prod_{m}Eu_{m}(F)\in H^{\bullet}_{\infty}(F,*)$$
qui v\'erifie
$$\int_{F}^{\infty}Eu_{\infty}(F)=f(X,H;t)\in
H^{\bullet}_{\infty}(Spec k,*)\simeq k[[t]]$$

Il serait peut-\^etre int\'eressant d'\'etudier les propri\'et\'es de la
th\'eorie cohomologique $F \mapsto H^{\bullet}_{\infty}(F,*)$, en vue
d'obtenir des "formules de Riemann-Roch sup\'erieures". \\

Supposons maintenant que $k$ soit un sous-corps de $\bf C$. A tout
champ de Deligne-Mumford $F$ sur $k$ est alors associ\'e le champ
analytique \\
$F^{an}:=(F\times_{Spec k}Spec \bf C\mit)^{an}$
(~\ref{d5.5}~). Notons $M^{an}$ son espace de modules, et $M$ celui de $F$.
Nous noterons $F^{top}$ et $M^{top}$ les champs topologiques associ\'es
aux champs analytiques $F^{an}$ et $M^{an}$.

\begin{lem}
Le morphisme naturel de champs topologiques
$$p : F^{top} \longrightarrow M^{top}$$
induit un isomorphisme de $\bf Q$-alg\`ebres
$$p^{*} : H^{*}(M^{top};\bf Q\mit) \longrightarrow H^{*}(F^{top},\bf
Q\mit)$$
\end{lem}

\underline{\bf Preuve:} \rm Par Mayer-Vietoris, ceci est local sur
$M^{top}$. On peut donc supposer d'apr\`es \ref{p1.2}, que $F=[X/H]$, avec
$H$ un groupe fini op\'erant sur un espace topologique $X$. Il faut
alors montrer que le morphisme naturel
$$p^{*} : H^{*}(X/H^{top},\bf Q\mit) \longrightarrow
H^{*}(X^{top},\bf Q\mit)^{H}$$
est un isomorphisme. Mais ceci est un cas particulier du th\'eor\`eme de
changement de base propre. $\Box$\\

Dans le cas o\`u $F$ est lisse et propre sur $Spec k$,
ce lemme implique que le produit d'intersection
$$\phi : H^{n}(M^{top},\bf Q\mit)\times H^{n}(M^{top},\bf Q\mit)
\longrightarrow
H^{2n}(M^{top},\bf Q\mit)\simeq \bf Q\mit$$
o\`u $n=Dim_{k}F$, est non-d\'eg\'en\'er\'e.

\begin{df}
Soit $F$ un champ de Deligne-Mumford, propre et lisse sur $Spec k$,
avec $Dim_{k}F=2m$. La signature de $M$ (~ou de $F$~) est la signature
de la forme quadratique r\'eelle
$$\phi\otimes \bf R\mit : H^{2m}(M^{top},\bf R\mit) \times
H^{2m}(M^{top},\bf R\mit) \longrightarrow \bf R\mit$$
Elle est not\'ee $\sigma(M)$ (~ou bien $\sigma(F)$~).
\end{df}

Si $F$ est un champ lisse de $\cal QDM$, une application de la th\'eorie de
Hodge pour les champs alg\'ebriques complexes (~\ref{hodge}~) implique
imm\'ediatement que
$$\sigma(M)=\sigma(F)=\chi(F,\lambda_{1}(\Omega_{F}^{1}))$$
o\`u $\lambda_{1}(x):=\sum_{i}\lambda^{i}(x) \in \bf
K\mit_{0}(F)$, pour tout \'el\'ement $x \in \bf K\mit_{0}(F)$ de rang positif.
On peut alors appliquer la formule d'Hirzebruch-Riemann-Roch \`a
l'\'el\'ement $\lambda_{1}(\Omega_{F}^{1})$, pour en d\'eduire une formule de
la signature. Pour cela rappelons la d\'efinition du $L$-genre
$L(F)$ (~\cite[$8.2.2$]{hi}~), d'un champ lisse $F$.

Soit $T_{F}$ le fibr\'e tangent de $F$, et $a_{i}$ ses racines de Chern.
On a donc
$$\sum_{i}C_{i}(T_{F}).t^{i}=\prod_{i=1}^{i=n}(1+a_{i}t)$$
Alors, le "polyn\^ome"
$$\prod_{i=1}^{i=n}(\frac{a_{i}}{tanh(a_{i})})$$
est sym\'etrique en les $a_{i}$. Il s'exprime donc comme un polyn\^ome
$L(F)~\in~H^{\bullet}(F,*)$, en les classes de Chern $C_{i}(T_{F})$.

\begin{cor}{(~Formule de la signature~)}\label{sign}
Soit $F$ un champ de Deligne-Mumford lisse et propre sur $Spec k$, de
dimension paire. Notons
$$\beta_{F}:=can(F(\lambda_{1}(\cal N\mit_{\pi}^{\vee})))
\in \underline{\bf K}_{0}^{rep}(I_{F})$$
o\`u $\cal N\mit_{\pi}^{\vee}$ est le fibr\'e conormal du morphisme
$$\pi_{F} : I_{F} \longrightarrow F$$
Alors, on a
$$\sigma(M)=\sigma(F)=\int_{I_{F}}Ch(\alpha_{F}^{-1}.\beta_{F}).L(I_{F})$$
\end{cor}

\underline{\bf Preuve:} \rm C'est le m\^eme calcul que pour la formule
de Gauss-Bonnet \ref{gb2}. $\Box$\\

\underline{Remarque:} \rm Lorsque $F=[X/H]$ est un quotient d'un
sch\'ema $X$ par un groupe fini $H$ d'ordre $m$, la formule pr\'ec\'edente se r\'e\'ecrit
$$\sigma(X/H)=\frac{1}{m}.\sum_{h \in H}
\int_{X^{h}}Ch(\alpha_{h}^{-1}.\beta_{h}).L(X^{h})$$
o\`u $\alpha_{h}$ et $\beta_{h}$ sont les restrictions de $\alpha_{F}$
et $\beta_{F}$ sur $X^{h}$. On reconna\^it alors
la formule de la "$G$-signature" d'Atiyah-Singer. Ceci permet
d'interpr\'eter ce genre de formule
dans le cadre de la "cohomologie \`a coefficients dans le
repr\'esentations". Il pourrait-\^etre par ailleurs int\'eressant d'utiliser cette
cohomologie pour d\'emontrer un th\'eor\`eme d'indice \`a la Atiyah-Singer
pour des champs diff\'erentiels. Comme tout orbifold diff\'erentiel est un quotient
par un groupe de Lie compact, la formule d'indice
\'equivariant d'Atiyah-Singer, permet au moins de traiter ce cas
(~c'est ce qui est fait dans \cite{ka}~). Le passage aux
champs diff\'erentiels g\'en\'eraux ne devrait pas alors poser trop de
probl\`emes.

\end{subsubsection}

\end{subsection}

\begin{subsection}{Comparaison avec d'autres th\'eories cohomologiques}
\hspace{5mm}
Dans ce paragraphe nous allons comparer la cohomologie \`a coefficients
dans les caract\`eres (~ou dans les repr\'esentations~) avec deux autres
constructions : la cohomologie p\'eriodique de la cat\'egorie des
faisceaux coh\'erents (~\cite{k}~), et les "groupes de Chow" associ\'es
\`a la filtration par la codimension du support sur la $K$-th\'eorie des
faisceaux coh\'erents (~qui apparaissent d\'ej\`a dans \cite{g2}~). \\

\begin{subsubsection}{Groupes de Chow}
\hspace{5mm}
On supposera que $S=Spec k$ avec $k$ un corps contenant les racines de
l'unit\'e, et on utilisera la th\'eorie de Gersten comme th\'eorie cohomologique. \\

Soit $F$ un champ alg\'ebrique de Deligne-Mumford, et $\bf
Coh\mit^{p}(F)$ la cat\'egorie de faisceaux coh\'erents sur $F$ dont le
support est de codimension sup\'erieure ou \'egale \`a $p$. Si on note
$$\bf G\mit^{p}(F):=K(\bf Coh\mit^{p}(F))$$
le foncteur d'inclusion $\bf Coh\mit^{p+1}(F) \longrightarrow \bf
Coh\mit^{p}(F)$ induit un triangle dans $HoSp$
$$\xymatrix{
& \bf G\mit^{p+1}(F) \ar[rd] &  \\
\bigvee_{x \in |F|^{(p)}}\bf G\mit(\widetilde{x}) \ar[ru]^{+1} &
& \bf G\mit^{p}(F) \ar[ll]}$$
o\`u $|F|^{(p)}$ et l'ensemble des points de $F$ de codimension sup\'erieure
ou \'egale \`a $p$.
On d\'eduit de ces triangles des complexes de groupes ab\'eliens
$\cal R\mit^{p}(F)$, que l'on consid\`ere comme concentr\'es en degr\'e $[-p,0]$ \\

\hspace{-10mm}$\displaystyle{\cal R\mit^{p}(F) : 0 \rightarrow
\bigoplus_{x \in |F|^{(0)}}\bf G\mit_{p}(\widetilde{x}) \longrightarrow
\bigoplus_{x \in |F|^{(1)}}\bf G\mit_{p-1}(\widetilde{x})
\dots
\bigoplus_{x \in |F|^{(p-1)}}\bf G\mit_{1}(\widetilde{x})
\longrightarrow \bigoplus_{x \in |F|^{(p)}}\bf G\mit_{0}(\widetilde{x})
\rightarrow 0}$ \\

Remarquons que le terme de droite s'interpr\`ete comme le groupe des
cycles de $F$ \`a "coefficients dans les repr\'esentations des gerbes
r\'esiduelles".

\begin{df}
Les groupes de Chow \`a coefficients dans les repr\'esentations sont d\'efinis par
$$CH^{p}_{rep}(F,q):=H^{-q}(\cal R\mit^{p}(F))$$
Nous noterons
$$CH^{\bullet}_{rep}(F,*):=\bigoplus_{p,q}CH^{p}_{rep}(F,q)$$
\end{df}

Indiquons sans d\'emonstrations (~le lecteur pourra les d\'eduire des
propri\'et\'es g\'en\'erales du foncteur de $G$-th\'eorie~)
les principales propri\'et\'es de ces groupes.

\begin{itemize}
\item
Pour tout morphisme propre de dimension cohomologie finie
$$f : F \longrightarrow F'$$
entre deux champs de Deligne-Mumford, il existe une image directe
fonctorielle
$$f_{*} : CH^{\bullet}_{rep}(F,*) \longrightarrow
CH^{\bullet}_{rep}(F,*)$$
Si, de plus $F$ et $F'$ sont irr\'eductibles, alors on a
$$f_{*} : CH^{p}_{rep}(F,q) \longrightarrow
CH^{p+d}_{rep}(F,q)$$
o\`u $d=Dim_{k}F' - Dim_{k}F$.
\item
Pour tout morphisme repr\'esentable et plat
$$f : F \longrightarrow F'$$
entre deux champs de Deligne-Mumford, il existe une image r\'eciproque
fonctorielle
$$f^{*} : CH^{p}_{rep}(F',q) \longrightarrow
CH^{p}_{rep}(F,q)$$
\item
Si $j : F' \hookrightarrow F$ est une immersion ferm\'ee de
codimension $d$ de champs
de Deligne-Mumford, et $i : F-F' \hookrightarrow F$ l'immersion
compl\'ementaire, alors on a une suite exacte longue\\

\hspace{-10mm}
$\displaystyle{\xymatrix{
\dots \ar[r] &
CH^{p-d}_{rep}(F',q) \ar[r]^{j_{*}} & CH^{p}_{rep}(F,q)
\ar[r]^{i^{*}} & CH^{p}_{rep}(U,q) \ar[r] & CH^{p-d}_{rep}(F',q-1)
\ar[r] & \dots}}$
\item
Soit $p : F \longrightarrow M$ la projection d'un champ de
Deligne-Mumford sur son espace de modules. Alors, les complexes
$$p_{*}(\cal R\mit^{p})\otimes \bf Q$$
sont flasques sur $M_{et}$.
\end{itemize}

\begin{thm}\label{th3.6}
Il existe un isomorphisme compatible avec les images directes
$$\epsilon_{F} : CH^{\bullet}_{rep}(F,*)_{K} \longrightarrow
H_{\bullet}^{rep}(F,*)$$
\end{thm}

\underline{\bf Preuve:} \rm Comme les techniques utilis\'ees ci-dessous
sont en tout point similaires \`a celles utilis\'ees dans le paragraphe
pr\'ec\'edent, nous ne donnerons qu'une esquisse de preuve. \\

Remarquons tout d'abord que par
d\'efinition, on a un isomorphisme naturel
$$H_{p}(F,q)\simeq H^{p-q}((I_{F})_{et},\cal R\mit^{p}\otimes K)$$

On cherche donc \`a construire un isomorphisme fonctoriel pour les
images directes (~qui ne respectera pas la graduation~!~)
$$\epsilon_{F} : CH^{\bullet}_{rep}(F,*)_{K}\longrightarrow
H^{*}((I_{F})_{et},\cal R\mit^{\bullet}\otimes K)$$

Commen\c{c}ons par supposer que
$F$ est une gerbe sur un espace alg\'ebrique $X$. Comme $\pi_{F} :
I_{F} \longrightarrow F$ est plat, il existe
$$\pi_{F}^{*} : CH^{p}_{rep}(F,q) \longrightarrow
CH^{p}_{rep}(I_{F},q)$$
Consid\'erons alors pour toute racine de l'unit\'e
$\zeta \in \mu_{\infty}(k)$,
le foncteur d\'efini dans la preuve de \ref{th2.3}
$$F_{\zeta} : \bf Coh\mit(I_{F}) \longrightarrow \bf Coh\mit(I_{F})$$
Comme ce foncteur pr\'eserve la filtration par la codimension du
support, on obtient un morphisme de complexes de pr\'efaisceaux sur
$(I_{F})_{et}$
$$F_{\zeta} : \cal R\mit^{p} \longrightarrow
\cal R\mit^{p}\otimes K$$
On pose alors
$$\epsilon_{F}:=\sum_{\zeta}can\circ F_{\zeta}\circ \pi_{F}^{*} :
\cal R\mit^{p}(F) \longrightarrow \cal R\mit^{p}(I_{F})\otimes K$$
C'est un morphisme de complexes.
De plus, si $f : F \longrightarrow F'$ est un morphisme propre et
repr\'esentable de gerbes, alors on a
$$\epsilon_{F'}\circ f_{*}=If_{*}\circ \epsilon_{F}
: \cal R\mit^{p}(F) \longrightarrow \cal R\mit^{p}(I_{F})\otimes K$$
en tant que morphisme de complexes. \\

Dans le cas d'un champ de Deligne-Mumford g\'en\'eral, choisissons une
hyper-quasi-enveloppe de Chow
$$q : F_{\bullet} \longrightarrow F$$
tel que $F_{m}$ soit une gerbe triviale.

En appliquant les foncteurs covariants $\cal R\mit^{*}_{\bf Q}$ au champ
simplicial $F_{\bullet}$, on obtient un complexe simplicial, dont le
complexe simple associ\'e sera not\'e
$$\cal R\mit^{*}(F_{\bullet})_{\bf Q}$$
Ce complexe est muni d'un morphisme de complexe
$$q_{*} : \cal R\mit^{*}(F_{\bullet}/F)_{\bf Q} \longrightarrow
\cal R\mit^{*}(F)_{\bf Q}$$
On montre de fa\c{c}on analogue \`a \ref{th2.2}, que le morphisme $q_{*}$ est
un quasi-isomorphisme.

On applique alors le morphisme $\epsilon$ construit ci-dessus \`a
chaque "\'etage" de $F_{\bullet}$, ce qui induit un morphisme de
complexes
$$\epsilon_{F_{\bullet}} : \cal R\mit^{*}(F_{\bullet}/F)_{K} \longrightarrow
\cal R\mit^{*}(I_{F_{\bullet}}/I_{F})_{K}$$
On a ainsi construit un diagramme de complexes
$$\xymatrix{
\cal R\mit^{*}(F_{\bullet}/F)_{K} \ar[r]^-{\epsilon_{F_{\bullet}}}
\ar[d]_{q_{*}} &
\cal R\mit^{*}(I_{F_{\bullet}}/I_{F})_{K}
\ar[d]^{Iq_{*}} \\
\cal R\mit^{*}(F)_{K} & \cal R\mit^{*}(F)_{K}}$$
Comme $q_{*}$ est un quasi-isomorphisme, on peut d\'efinir
un morphisme
$$\epsilon_{F}:=(Iq_{*})\epsilon_{F_{\bullet}}(q_{*})^{-1} :
CH^{\bullet}_{rep}(F,*)_{K}
\longrightarrow CH^{\bullet}_{rep}(I_{F},*)_{K}$$
On le compose alors avec le morphisme canonique
$$can : CH^{\bullet}_{rep}(I_{F},*)_{K}\longrightarrow
H^{*}((I_{F})_{et},\cal R\mit^{\bullet}\otimes K)$$
pour obtenir le morphisme cherch\'e
$$\epsilon_{F} : CH^{\bullet}_{rep}(F,*)_{K}
\longrightarrow H^{rep}_{\bullet}(F,*)$$
Par construction, ce morphisme est compatible avec les images
directes, et les changements de base
\'etales de l'espace de modules. Pour montrer que c'est un
isomorphisme, on peut donc raisonner par r\'ecurrence sur la dimension
de $F$ en appliquant les suites exactes longues, et donc se ramener
au cas o\`u $F$ est une gerbe triviale sur un sch\'ema. Comme le morphisme
$\epsilon_{F}$ est clairement un isomorphisme dans ce cas, on en
d\'eduit le th\'eor\`eme. $\Box$\\

\underline{Remarque:} L'isomorphisme \ref{th3.6} ne pr\'eserve pas la
graduation.

\end{subsubsection}

\begin{subsubsection}{Cohomologie P\'eriodique}
\hspace{5mm} Rappelons que pour une cat\'egorie exacte $\cal E$, il
existe une construction d'un spectre de cohomologie p\'eriodique de
$\cal E$ (~\cite{k}~). On se propose dans ce paragraphe de
comparer la cohomologie p\'eriodique des fibr\'es vectoriels et des
faisceaux coh\'erents sur un champ, et sa cohomologie de De Rham \`a
coefficients dans les caract\`eres et dans les repr\'esentations.
Pour cela on supposera que $S=Spec k$, avec $k$ un corps de caract\'eristique
nulle contenant les racines de l'unit\'e, et que la th\'eorie cohomologique
utilis\'ee
est la cohomologie de De Rham. \\

Soit $\cal E$ une cat\'egorie $k$-lin\'eaire et exacte. Dans \cite{k}, Keller lui
associe un complexe double $CP(\cal E\mit)$ de $k$-espaces vectoriels. Le
spectre associ\'e au complexe simple $sCP(\cal E\mit)$ par la
construction de Dold-Puppe, sera not\'e $\cal
CP\mit(E)$. On l'appellera le spectre de cohomologie p\'eriodique de $\cal
E\mit$. Il sera consid\'er\'e comme un objet de $HoSp$.
On dispose alors des propri\'et\'es suivantes : \\

\begin{itemize}

\item $\cal E\mit \mapsto \cal CP\mit(E)$ est un foncteur de la cat\'egorie
des cat\'egories exactes $k$-lin\'eaires (~et foncteurs exacts
$k$-lin\'eaires~), vers la cat\'egorie des spectres.

\item si $R$ est une $k$-alg\`ebre, il existe un isomorphisme naturel
$$CP(R) \stackrel{\simeq}{\longrightarrow} CP(\cal E\mit)$$
o\`u $\cal E$ est la cat\'egorie des $R$-modules projectifs de type fini
et $CP(R)$ est le complexe d'homologie p\'eriodique de la $k$-alg\`ebre
$R$.

\item si $A$ est une cat\'egorie ab\'elienne et $B$ une
sous-cat\'egorie de Serre, alors les foncteurs naturels
$$B \longrightarrow A \longrightarrow A/B$$
induisent un triangle dans $HoSp$
$$\xymatrix{
& \cal CP\mit(B) \ar[rd] & \\
\cal CP\mit(A/B) \ar[ru]^{-1} & & \cal CP\mit(A) \ar[ll]}$$

\item si $F : \cal E\mit \longrightarrow \cal E'\mit$ est un foncteur exact qui
induit une \'equivalence sur les cat\'egories d\'eriv\'ees, alors le
morphisme induit
$$F : \cal CP\mit (\cal E\mit) \longrightarrow \cal CP\mit (\cal E\mit ')$$
est un isomorphisme

\item il existe un morphisme fonctoriel en $\cal E$
$$Ch : K(\cal E\mit) \longrightarrow \cal CP\mit(\cal E\mit)$$

\item si on note $(Li/k)_{li}$ le site des espaces alg\'ebriques lisses
muni de la topologie lisse, alors
il existe un isomorphisme dans $HoSp((Li/k)_{li})$
$$\cal CP\mit \simeq \cal H\mit$$
\end{itemize}

Bien entendu, la liste n'est pas exhaustive. En particulier la troisi\`eme
propri\'et\'e que l'on cite n'est qu'un cas particulier d'une
propri\'et\'e de localisation bien plus g\'en\'erale.\\

\begin{df}
Pour tout champ alg\'ebrique $F$, on d\'efinit sa cohomologie p\'eriodique par
$$HP_{*}(F):=\pi_{*}\cal CP\mit(\bf Vect\mit(F))$$
\end{df}

Par construction, $F \mapsto HP_{*}(F)$ est clairement un foncteur
contravariant.

\begin{prop}
Soit $F$ un champ lisse et bien ramifi\'e.
\begin{enumerate}
\item
Il existe une transformation naturelle de foncteurs contravariants
$$Ch^{per} : \bf K\mit_{*} \longrightarrow HP_{*}$$
\item
Si $F$ est un champ alg\'ebrique lisse, il existe un morphisme
$$\epsilon_{F} : HP_{*}(F) \longrightarrow H^{\bullet}_{\chi}(F,*)$$
telle que
\begin{itemize}
\item $\epsilon$ est compatible avec les images r\'eciproques de
morphismes entre champs lisses
\item on a
$$\epsilon_{F} \circ Ch^{per} = Ch^{\chi}\circ \epsilon_{F}$$
\end{itemize}
\end{enumerate}
\end{prop}

\underline{\bf Preuve:} \rm Nous ne donnerons qu'une esquisse de preuve. \\

L'existence du caract\`ere de Chern provient directement des propri\'et\'es
rappel\'ees ci-dessus. \\

De fa\c{c}on tout \`a fait analogue \`a \ref{th2.4}, on construit un morphisme
$$\rho : \cal CP\mit(D_{F}) \longrightarrow
\bf H\mit((D_{F})_{li},\cal CP\mit\otimes \cal A\mit_{F})$$
que l'on compose avec l'\'equivalence rappel\'ee ci-dessus
$$\bf H\mit((D_{F})_{li},\cal CP\mit\otimes \cal A\mit_{F})
\longrightarrow \bf H\mit((D_{F})_{li},\cal H\mit\otimes \cal
A\mit_{F})$$
et l'image r\'eciproque $d_{F}^{*} : \cal CP\mit(F) \longrightarrow
\cal CP\mit(D_{F})$ pour obtenir le morphisme cherch\'e
$$HP_{*}(F) \longrightarrow H^{\bullet}_{\chi}(F,*)$$
Par la propri\'et\'e de naturalit\'e de cette construction et du caract\`ere
de Chern construit dans \cite{k}, on a clairement
$$\epsilon_{F} \circ Ch^{per} = Ch^{\chi}\circ \epsilon_{F}$$
$\Box$\\

Dans le cas o\`u $F$ est un champ de Deligne-Mumford, on peut aller un
peu plus loin.

\begin{prop}
Soit $p : F \longrightarrow M$ la projection d'un champ de
Deligne-Mumford lisse, sur son espace de modules. Consid\'erons le
pr\'efaisceau en spectres sur $M_{et}$
$$\begin{array}{cccc}
p_{*}\cal CP\mit : & M_{et} & \longrightarrow & Sp \\
& U & \mapsto & \cal CP\mit (\bf Vect\mit(p^{-1}U))
\end{array}$$
Posons alors
$$HP_{*}'(F):=\bf H\mit^{-*}(M_{et},p_{*}\cal CP\mit)$$
Nous noterons encore $Ch^{per}$ le morphisme compos\'e
$$\xymatrix{
\bf K\mit_{*}(F) \ar[r]^-{Ch^{per}} & HP_{*}(F) \ar[r]^-{can}&
HP_{*}'(F)}$$
Alors, il existe un isomorphisme
$$\epsilon_{F}' : HP_{*}'(F) \longrightarrow H^{\bullet}_{rep}(F,*)$$
tel que
\begin{itemize}
\item $\epsilon'$ est compatible avec les images r\'eciproques de
morphismes entre champs lisses
\item on a
$$\epsilon_{F} \circ Ch^{per} = Ch^{rep}\circ \epsilon_{F}$$
\end{itemize}
\end{prop}

\underline{\bf Preuve:} \rm L'existence du caract\`ere de Chern, et du
morphisme $\epsilon_{F}'$ se d\'emontre comme dans la proposition
pr\'ec\'edente. Il ne reste qu'\`a d\'emontrer que $\epsilon_{F}'$ est un
isomorphisme.

Mais, par d\'efinition de $HP_{*}'(F)$, ceci est local sur $M_{et}$. On
peut donc supposer que $F=[X/H]$ est un quotient d'un sch\'ema affine
lisse $X$ par un groupe fini. Le fait que $\epsilon_{F}'$ est un
isomorphisme dans ce cas est alors d\'emontr\'e dans \cite[$A6.1$, $A6.9$]{ft}.
$\Box$\\

\underline{Remarque:} Je ne sais pas si le morphisme naturel
$$can : HP_{*}(F) \longrightarrow HP_{*}'(F)$$
est un isomorphisme, mais cela me semblerait moral. On aurait alors
un isomorphisme pour un champ de Deligne-Mumford lisse
$$HP_{*}(F) \simeq H^{\bullet}_{rep}(F,*)$$
Cet isomorphisme serait alors un analogue alg\'ebrique de la description de la
cohomologie p\'eriodique d'un orbifold (~\cite[$6.12$]{cm}~).

\end{subsubsection}

\end{subsection}

\begin{subsection}{Conclusion sur les th\'eor\`emes de Riemann-Roch}
\hspace{5mm}
A la lumi\`ere des th\'eor\`emes \ref{c3.3}, \ref{th3.5} et \ref{th3.5'}, on peut
conclure que le probl\`eme
initial consistant \`a d\'emontrer des th\'eor\`emes de Riemann-Roch pour les
champs alg\'ebriques est r\'esolu dans le cas des champs de
Deligne-Mumford, et partiellement r\'esolu dans le cas g\'en\'eral. Il me
semble cependant que la d\'efinition de la cohomologie \`a coefficients
dans les repr\'esentations est une bonne notion pour traiter le cas des
champs $\Delta$-affines. Je pense que l'on pourra dire que la
situation est satisfaisante lorsque l'on aura d\'emontr\'e que le
th\'eor\`eme de Grothendieck-Riemann-Roch reste vrai (~au moins au niveau
du $\bf G\mit_{0}$~) pour un morphisme propre de champs alg\'ebriques
$\Delta$-affines et lisses sur un corps. Pour arriver \`a un tel
r\'esultat, il manque essentiellement deux choses. Premi\`erement, il
faudrait montrer que les champs $\Delta$affines et lisses sur un corps
sont bien ramifi\'es. Le second point important est celui de
l'existence de "lemmes de Chow" pour des morphismes propres de champs
alg\'ebriques, qui nous fait d\'efaut en toute g\'en\'eralit\'e, et nous oblige
donc \`a ne traiter que le cas des champs qui sont localement des
quotients g\'eom\'etriques uniformes affines (~pour lesquels on dispose
des quasi-enveloppe de Chow~). Un autre point de vu consisterait \`a
montrer que l'on peut "localiser en bas" le th\'eor\`eme de Riemann-Roch, et donc
se ramener syst\'ematiquement \`a ce dernier cas.

\end{subsection}

\end{section}

\end{part}

\newpage
\begin{part}*{\Huge Seconde Partie :\\
\vspace{25mm}$\cal D$-modules et th\'eor\`eme "GAGA"\large}
\newpage
\begin{section}{Chapitre $4$ : Th\'eor\`emes de Grothendieck-Riemann-Roch
pour les $\cal D$-modules et formules d'indices}
\hspace{5mm}
Comme il est expliqu\'e dans \cite{l}, une importante application de la
formule de Riemann-Roch est le calcul d'indices de $\cal D$-modules
coh\'erents. Dans ce chapitre nous nous inspirons de ces constructions
pour comparer les spectres de $K$-th\'eorie des faisceaux coh\'erents aux
spectres de \\
$K$-th\'eorie des $\cal D$-modules coh\'erents, et nous en
d\'eduisons des formules de Riemann-Roch pour les $\cal D$-modules
coh\'erents dans le cadre des champs alg\'ebriques. \\

Le point essentiellement technique de ce chapitre est la d\'efinition
des images directes en $K$-th\'eorie. En effet, dans la litt\'erature les
images directes de $\cal D$-modules ne sont d\'efinies qu'au niveau des
cat\'egories d\'eriv\'ees. Or, ceci n'est pas suffisant pour induire un
morphisme au niveau des spectres de $K$-th\'eorie. Pour contourner
cette difficult\'e nous avons choisi de travailler avec les m\'ethodes de
\cite{th}. Cependant, par souci de l\'eg\`eret\'e nous n'avons pas donn\'e
toutes les d\'emonstrations en d\'etails. Ce choix est justifi\'e par le
fait que nous nous int\'eressons plus ici aux applications du th\'eor\`eme
de Riemann-Roch de la section pr\'ec\'edente, qu'\`a l'\'etude d\'etaill\'ee des
$\cal D$-modules sur les champs alg\'ebriques. \\

Dans un premier temps, le th\'eor\`eme de
Grothendieck-Riemann-Roch au niveau du $K_{0}$ nous permet d'obtenir
des formules calculant des caract\'eristiques d'Euler pond\'er\'ees,
analogues \`a ce qui est d\'emontr\'e dans \cite{mac}. En particulier, on
compl\`ete ainsi la preuve de la formule de Gauss-Bonnet
\ref{gb1}. Par la suite, la formule de Riemann-Roch appliqu\'ee \`a certains
\'el\'ements de la $K$-th\'eorie des $\cal D$-modules holonomes, permet
d'obtenir des formules calculant des caract\'eristiques d'Euler pond\'er\'ees
par des fonctions constructibles \`a valeurs dans la $K$-th\'eorie
sup\'erieure du corps de base.

Notons enfin, que la description de la $K$-th\'eorie des $\cal
D$-modules holonomes r\'eguliers nous fait esp\'erer qu'il serait
possible d'\'etendre la formule de Riemann-Roch de S. Bloch et H. Esnault
\cite{be} au cas des $\cal D$-modules holonomes r\'eguliers. Ce th\'eor\`eme
serait alors \`a valeurs dans les fonctions constructibles \`a valeurs
dans la cohomologie de Chern-Simons d\'efinie dans \cite{be}.\\

Pour tout ce chapitre, on notera $k$ un corps
de caract\'eristique nulle et contenant les racines de l'unit\'e. Un
champ sera un champ sur $(Esp/Spec k)_{et}$.

Nous noterons $(Li/k)_{et}$ (~resp. $(Li/k,li)_{et}$~) le site
des espaces alg\'ebriques lisses sur $k$ (~resp. des espaces alg\'ebriques
lisses et morphismes lisses~), muni de la topologie \'etale.

Nous noterons aussi $\cal LDM$ (~resp.
$\cal QLDM$~) la cat\'egorie homotopiques des champs de Deligne-Mumford
lisses et s\'epar\'es sur $k$ (~resp. lisses et s\'epar\'es sur $k$ et dont
l'espace de modules est
quasi-projectif~).\\

Par la suite on se fixera une th\'eorie cohomologique avec images directes
(~\ref{s4}~). Les groupes de cohomologie et d'homologie associ\'es seront not\'es
$H^{\bullet}(-,*)$ et $H_{\bullet}(-,*)$. On utilisera aussi la
cohomologie \`a coefficients dans les repr\'esentations d\'efinie dans
\ref{s3}. Elle sera not\'ee comme pr\'ec\'edemment $H^{\bullet}_{rep}(-,*)$
et $H_{\bullet}^{rep}(-,*)$.

\begin{subsection}{Les champs de $\cal D$-modules}
\hspace{5mm}
Pour chaque objet $X \in Ob(Li/k)$, on dispose du
faisceau des op\'erateurs diff\'erentiels $\cal D\mit_{X}$ (~\cite{bo}~).
La cat\'egorie des $\cal D\mit_{X}$-modules  (~resp. des $\cal
D\mit_{X}$-modules $\cal D\mit_{X}$-coh\'erents, resp.
des $\cal D\mit_{X}$-modules holonomes, resp. des
$\cal D\mit_{X}$-modules holonomes r\'eguliers, resp. des
$\cal D\mit_{X}$-modules $\cal O\mit_{X}$-coh\'erents, resp.
des $\cal D\mit_{X}$-modules $\cal O\mit_{X}$-coh\'erents r\'eguliers~)
sera not\'ee $\cal D\mit-\bf Mod\mit(X)$
(~resp. $\cal D\mit-\bf Coh\mit(X)$, resp. $\cal D\mit_{h}-\bf Coh\mit(X)$,
resp.
$\cal D\mit_{hr}-\bf Coh\mit(X)$, resp. $\nabla(X)$,
resp. $\nabla_{r}(X)$~). Par abus de langage on dira
"$\cal D\mit_{X}$-module coh\'erent"
pour "$\cal D\mit_{X}$-modules $\cal D\mit_{X}$-coh\'erent". On appelera
aussi $\nabla(X)$ la cat\'egorie des connexions sur $X$.

Si $f : X \longrightarrow Y$ est un morphisme dans $(Li/k)$, on
dispose d'images r\'eciproques (~\cite{bo}~)
$$f^{!} : \cal D\mit-\bf Mod\mit(Y) \longrightarrow \cal D\mit-\bf
Mod\mit(X)$$
qui pr\'eservent les propri\'et\'es d'\^etre
holonome (~resp. holonome r\'egulier, resp. $\cal O$-coh\'erent,
resp. $\cal O$-coh\'erent r\'egulier~). Ceci
permet de d\'efinir la cat\'egorie cofibr\'ee $\cal D\mit-\bf Mod\mit$, des $\cal
D$-modules (~resp. $\cal D\mit_{h}-\bf Coh\mit$, des \\
$\cal D$-modules holonomes, resp.
$\cal D\mit_{hr}-\bf Coh\mit$, des $\cal D$-modules holonomes r\'eguliers,
resp. $\nabla$,
des $\cal D$-modules $\cal O$-coh\'erents, resp.
des $\nabla_{r}$, des $\cal D$-modules \\
$\cal O$-coh\'erents r\'eguliers~).
Ce sont des cat\'egories cofibr\'ees sur $(Li/k)_{et}$.

De m\^eme, lorsque $f$ est lisse, le foncteur $f^{!}$ pr\'eserve la
propri\'et\'e d'\^etre $\cal D$-coh\'erent. Ainsi, on peut d\'efinir
$\cal D\mit-\bf Coh\mit$, la cat\'egorie cofibr\'ee des \\
$\cal D$-modules coh\'erents sur $(Li/k,li)_{et}$.

\begin{prop}\label{p4.1}
Les cat\'egories cofibr\'ees $\cal D\mit_{h}-\bf Coh\mit$,
$\cal D\mit_{hr}-\bf Coh\mit$ et $\nabla$ sont des champs sur $(Li/k)_{et}$.

La cat\'egorie cofibr\'ee $\cal D\mit-\bf Coh\mit$ est un champ sur
$(Li/k,li)_{et}$.

Se sont de plus des champs en cat\'egories exactes sur $(Li/k,li)_{et}$.
\end{prop}

\underline{\bf Preuve:} \rm Cette proposition provient directement
du fait qu'un $\cal D$-module sur un espace alg\'ebrique lisse $X$, est la donn\'ee
d'une connexion int\'egrable sur un $\cal O$-module quasi-coh\'erent. Ce qui est
\'equivalent \`a la donn\'ee d'un faisceau de $\cal O\mit_{X}$-modules
quasi-coh\'erent $\cal M$, et d'un scindage de la suite exacte d'Atiyah
$$0 \longrightarrow \Omega^{1}_{X}\otimes_{\cal O\mit_{X}}\cal M\mit
\longrightarrow \cal P\mit(\cal M\mit) \longrightarrow
\cal M\mit \longrightarrow 0$$
Il est clair que ces donn\'ees poss\`edent la
propri\'et\'e de descente pour des morphismes \'etales et surjectifs. $\Box$\\

\begin{df}\label{d4.1}
Soit $F$ un champ de $\cal LDM$, $F_{et}$ son site
\'etale, et \\
$\cal D\mit-\bf Mod\mit_{F}$, $\cal D\mit-\bf Coh\mit_{F}$, $\cal
D\mit_{h}-\bf Coh\mit_{F}$, $\cal D\mit_{hr}-\bf Coh\mit_{F}$,
$\nabla_{r,F}$ et $\nabla_{F}$, les champs restreints
sur $F_{et}$. Ce sont des champs en cat\'egories exactes sur $F_{et}$.

\begin{itemize}
\item
On d\'efinit
$$\cal D\mit-\bf Mod\mit(F):=\int_{F_{et}}\cal D\mit-\bf Mod\mit_{F}$$
$$\cal D\mit-\bf Coh\mit(F):=\int_{F_{et}}\cal D\mit-\bf Coh\mit_{F}$$
$$\cal D\mit_{h}-\bf Coh\mit(F):=\int_{F_{et}}\cal D\mit_{h}-\bf
Coh\mit_{F}$$
$$\cal D\mit_{hr}-\bf Coh\mit(F):=\int_{F_{et}}\cal D\mit_{hr}-\bf
Coh\mit_{F}$$
$$\nabla(F):=\int_{F_{et}}\nabla_{F}$$
$$\nabla_{r}(F):=\int_{F_{et}}\nabla_{r,F}$$
Ce sont respectivement les cat\'egories exactes des $\cal D$-modules
sur $F$, des $\cal D$-modules coh\'erents sur $F$, des $\cal D$-modules
holonomes sur $F$, des $\cal D$-modules holonomes r\'eguliers sur $F$,
et des connexions sur $F$.
\item
Les spectres de $K$-th\'eorie associ\'es (~\ref{d1.4}~) seront
not\'es respectivement
$$\bf K\mit^{\cal D}(F):=K(\cal D\mit-\bf Coh\mit(F))$$
$$\bf K\mit^{h}(F):=K(\cal D\mit_{h}-\bf Coh\mit(F))$$
$$\bf K\mit^{h,r}(F):=K(\cal D\mit_{hr}-\bf Coh\mit(F))$$
$$\bf K\mit^{\nabla}(F):=K(\nabla(F))$$
$$\bf K\mit^{\nabla,r}(F):=K(\nabla_{r}(F))$$
Les spectres de $K$-cohomologie associ\'es (~\ref{d1.4}~) seront not\'es
respectivement
$$\underline{\bf K}^{\cal D}(F):=\bf H\mit(F_{li},\underline{K}^{\cal D}_{\bf Q})$$
$$\underline{\bf K}^{h}(F):=\bf H\mit(F_{li},\underline{K}^{h}_{\bf Q})$$
$$\underline{\bf K}^{h,r}(F):=\bf H\mit(F_{li},\underline{K}^{h,r}_{\bf Q})$$
$$\underline{\bf K}^{\nabla}(F):=\bf H\mit(F_{li},\underline{K}^{\nabla}_{\bf Q})$$
$$\underline{\bf K}^{\nabla,r}(F):=\bf H\mit(F_{li},\underline{K}^{\nabla,r}_{\bf Q})$$
\end{itemize}
\end{df}

\underline{Remarque:} \rm Sur $F_{et}$ on dispose du faisceaux des
op\'erateurs diff\'erentiels $\cal D\mit_{F}$.
Il y a alors une \'equivalence naturelle entre la
cat\'egorie des $\cal D$-modules sur $F$ au sens pr\'ec\'edent, et celle des
$\cal D\mit_{F}$-modules sur $F_{et}$. De m\^eme, il existe des
d\'efinitions \'evidentes de $\cal D\mit_{F}$-modules coh\'erents,
holonomes, holonomes r\'eguliers ... Ces notions se correspondent
\'evidemment par l'\'equivalence de cat\'egories \'evoqu\'ee ci-dessus.

Avant d'aller plus loin, rappelons quelques propri\'et\'es concernant les
$\cal D\mit_{F}$-modules coh\'erents.

\begin{prop}
Soit $F$ un champ de $\cal LDM$.
\begin{enumerate}
\item
Pour tout $\cal D\mit_{F}$-module coh\'erent $\cal M$, il existe un
sous-$\cal O\mit_{F}$-module coh\'erent $\cal F\mit \hookrightarrow \cal M$,
tel que
$$\cal M\mit=\cal D\mit_{F}.\cal F$$
\item
Tout $\cal D\mit_{F}$-module coh\'erent admet une bonne filtration
(~\cite{l}~).
\end{enumerate}
\end{prop}

\underline{\bf Preuve:} \rm $(1)$ Comme $\cal M$ est $\cal
O\mit_{F}$-quasi-coh\'erent, il est limite inductive de ses sous-$\cal
O\mit_{F}$-modules coh\'erents. Donc $\cal M$ est limite inductive de ses
sous-modules de la forme $\cal D\mit_{F}.\cal F$, o\`u $\cal F$ est un
sous-$\cal O\mit_{F}$-module coh\'erent. Comme $\cal M$ est coh\'erent, on
en d\'eduit qu'il existe un sous-$\cal O\mit_{F}$-module coh\'erent $\cal
F\mit_{0}$, tel que $\cal M\mit=\cal D\mit_{F}.\cal F\mit_{0}$. \\

$(2)$ Soit $\cal M\mit=\cal D\mit_{F}.\cal F$, avec $\cal F$ un
sous-$\cal O\mit_{F}$-module coh\'erent. Alors, si $\cal D\mit_{F}^{i}$
est le sous-$\cal O\mit_{F}$-module de $\cal D\mit_{F}$,
des op\'erateurs diff\'erentiels d'ordre inf\'erieur \`a $i$, la filtration
$$\cal M\mit^{i}:=\cal D\mit_{F}^{i}.\cal F\mit \hookrightarrow
\cal M$$
est une bonne filtration. $\Box$ \\

Il est clair que les correspondances
$$F \mapsto \bf K\mit^{h}(F) \qquad
F \mapsto \underline{\bf K}^{h}(F) \qquad
F \mapsto \bf K\mit^{h,r}(F) \qquad
F \mapsto \underline{\bf K}^{h,r}(F)$$
$$F \mapsto \bf K\mit^{\nabla}(F) \qquad
F \mapsto \underline{\bf K}^{\nabla}(F)\qquad
F \mapsto \bf K\mit^{\nabla,r}(F)\qquad
F \mapsto \underline{\bf K}^{\nabla,r}(F)$$
d\'efinissent des foncteurs contravariants de la cat\'egorie
$\cal LDM$ vers $HoSp$.

De m\^eme les correspondances
$$F \mapsto \bf K\mit^{\cal D}(F)$$
$$F \mapsto \underline{\bf K}^{\cal D}(F)$$
d\'efinissent des foncteurs contravariants de la sous-cat\'egorie
de $\cal LDM$ des champs et morphisme lisses, vers $HoSp$.

\begin{subsubsection}{Images directes de $\cal D$-modules}
\hspace{5mm}

\begin{prop}\label{p4.3}
\begin{enumerate}
\item
La correspondance
$$F \mapsto \bf K\mit^{\cal D}(F)$$
d\'efinit un foncteur covariant de la sous-cat\'egorie
$(\cal LDM\mit,pr.)$ de
$HoChAlg(k)$, des champs de Deligne-Mumford lisses et morphismes
propres, vers $HoSp$.
\item
Les correspondances
$$F \mapsto \bf K\mit^{h}(F)$$
$$F \mapsto \bf K\mit^{h,r}(F)$$
d\'efinissent des foncteurs covariants de
$\cal LDM$, vers $HoSp$.
\end{enumerate}
\end{prop}

\underline{\bf Preuve:} \rm Nous allons utiliser les constructions de
Thomason (~\cite{th}~), afin de d\'efinir un analogue des images directes
de $\cal D$-modules d\'efinies dans \cite[$5$]{bo}, suffisamment
fonctorielles pour
passer aux spectres de $K$-th\'eorie.\\

Soit $f : F \longrightarrow F'$ un $1$-morphisme de champs
de Deligne-Mumford lisses. Sur le faisceau d'anneaux $f^{*}\cal
D\mit_{F'}$, il existe une structure naturelle de $\cal
D\mit_{F}$-module \`a gauche. Elle est d\'efinie en consid\'erant le
morphisme naturel sur les faisceaux tangents
$T_{F} \longrightarrow f^{*}T_{F'}$, qui induit un morphisme sur les
faisceaux des op\'erateurs diff\'erentiels $\cal D\mit_{F}
\longrightarrow f^{*}\cal D\mit_{F'}$. En tant que $(\cal
D\mit_{F},f^{-1}\cal D\mit_{F'})$-bi-module,
$f^{*}\cal D\mit_{F'}$ sera not\'e
$$\cal D\mit_{F\rightarrow F'}$$
Le $\cal D\mit_{F}$-module \`a droite dual (~\cite[$3.3$]{bo}~)
de $\cal D\mit_{F\rightarrow F'}$ sera alors not\'e
$$\cal D\mit_{f}:=\cal D\mit_{F\rightarrow F'}\otimes_{\cal O
\mit_{F}}\omega_{F}$$
C'est un $(f^{-1}\cal D\mit_{F'},\cal D\mit_{F})$-bi-module. De fa\c{c}on
explicite on a
$$\cal D\mit_{f}=f^{-1}(\cal D\mit_{F'}\otimes_{\cal
O\mit_{F'}}\omega_{F}^{\vee})\otimes_{f^{-1}\cal
O\mit_{F'}}\omega_{F'}$$

Introduisons les notations suivantes.
\begin{itemize}
\item
Soit $B(f)$ la cat\'egorie bi-compliciale de Waldhausen des complexes de
$f^{-1}\cal D\mit_{F'}$-modules injectifs, cohomologiquement born\'es.
\item
Soit $A(f)$ la cat\'egorie bi-compliciale de Waldhausen des triplets
$(E_{\bullet},F_{\bullet},u)$, o\`u $E_{\bullet}$ est un complexe de
$\cal D\mit_{F}$-modules plats sur $\cal D\mit_{F}$, et \`a cohomologie born\'ee
et coh\'erente, $F_{\bullet}$ est un objet de $B(f)$, et
$$u : \cal D\mit_{f}\otimes_{\cal D\mit_{F}} E_{\bullet}
\longrightarrow F_{\bullet}$$
est un quasi-isomorphisme.
\item
Soit $C(f)$ la cat\'egorie des complexes de $\cal D\mit_{F'}$-modules \`a
cohomologie
born\'ee.
\end{itemize}
Notons alors que le foncteur naturel
$$(E_{\bullet},F_{\bullet},u) \mapsto E_{\bullet}$$
de $A(f)$ dans la cat\'egorie
des complexes de $\cal D\mit_{F}$-modules \`a cohomologie born\'ee et coh\'erente,
induit une \'equivalence au niveau des cat\'egories d\'eriv\'ees
(~appliquer \cite[$1.9.7$]{th}~), et donc un isomorphisme
naturel dans $HoSp$ (~\cite[$1.9.8$]{th}~)
$$K(A(f)) \simeq \bf K\mit^{\cal D}(F)$$
De plus, le foncteurs
$$\begin{array}{cccc}
\cal D\mit_{f}\otimes_{\cal O\mit_{F}}- : &
A(f) & \longrightarrow & B(f) \\
& (E_{\bullet},F_{\bullet},u) & \mapsto & F_{\bullet}
\end{array}$$
$$\begin{array}{cccc}
f_{*} : &
B(f) & \longrightarrow & C(f) \\
& F_{\bullet} & \mapsto & f_{*}(F_{\bullet})
\end{array}$$
sont "exacts" au sens de \cite{th}. Par composition on obtient un foncteur exact
$$f_{+} : A(f) \longrightarrow C(f)$$
qui est une r\'ealisation de l'image directe d\'efinie dans \cite[$5$]{bo}, au
niveau des
cat\'egories de complexes. Nous noterons
$$\bf R\mit f_{+} : D^{b}(\cal D\mit_{F}) \longrightarrow D^{b}(\cal
D\mit_{F'})$$
le morphisme induit sur les cat\'egories d\'eriv\'ees des complexes
de \\
$\cal D$-modules cohomologiquement born\'es.

\underline{Preuve de $(1)$:}

\begin{lem}\label{l4.2}
Si $f$ est un morphisme propre, alors
le foncteur $\bf R\mit f_{+}$ pr\'eserve la propri\'et\'e d'\^etre "\`a cohomologie
coh\'erente".
\end{lem}

\underline{\bf Preuve:} \rm La preuve donn\'ee
dans \cite[$2.2.1.1$]{al} se traduit mot pour mot au cas des champs
alg\'ebriques.
Elle sera r\'e\'ecrite au cours de la d\'emonstration de \ref{th4.1}. Le lecteur
v\'erifiera qu'il n'y a pas de cercle vicieux. $\Box$\\

Le lemme pr\'ec\'edent permet donc de dire que le foncteur $\bf R\mit f_{+}$ se
factorise par
$$\bf R\mit f_{+} : A(f) \longrightarrow C'(f)$$
o\`u $C'(f)$ est la cat\'egorie des complexes de $\cal D\mit_{F'}$-modules \`a
cohomologie coh\'erente et born\'ee. Il induit donc un morphisme dans $HoSp$
$$f_{+} : K(A(f)) \longrightarrow \bf K\mit^{\cal D}(F')$$
Et en composant avec l'isomorphisme canonique $\bf K\mit^{\cal D}(F)\simeq
K(A(f))$, on obtient le morphisme cherch\'e dans $HoSp$
$$f_{+} : \bf K\mit^{\cal D}(F) \longrightarrow \bf K\mit^{\cal D}(F')$$
Il reste \`a montrer que cela d\'efinit bien un foncteur, $F \mapsto
\bf K\mit^{\cal D}(F)$.

\begin{lem}
\begin{enumerate}
\item
Soit $i : F \longrightarrow F'$ un morphisme repr\'esentable fini et non-ramifi\'e,
$p : F' \longrightarrow F''$ un morphisme lisse, et $f=p\circ i$. Alors, on
a une
\'egalit\'e dans $HoSp$
$$f_{+}=p_{+}\circ i_{+} : \bf K\mit^{\cal D}(F)
\longrightarrow \bf K\mit^{\cal D}(F')$$
\item
Soit $i : F' \longrightarrow F''$ un morphisme repr\'esentable fini et
non-ramifi\'e,
$p : F \longrightarrow F'$ un morphisme lisse, et $f=i\circ p$. Alors, on a une
\'egalit\'e dans $HoSp$
$$f_{+}=i_{+}\circ p_{+} : \bf K\mit^{\cal D}(F)
\longrightarrow \bf K\mit^{\cal D}(F')$$
\end{enumerate}
\end{lem}

\underline{\bf Preuve:} \rm On ne d\'emontrera que le premier point, le
second se traitant de la m\^eme fa\c{c}on. \\

On va utiliser une construction
diff\'erente de $i_{+}$ et $p_{+}$.\\

Commen\c{c}ons par le cas de $i : F \longrightarrow F'$, un morphisme
repr\'esentable,
fini et non-ramifi\'e.

On sait alors que $\cal D\mit_{i}$ est un $\cal D\mit_{F}$-module
localement libre (~\cite[$7.7$]{bo}~). Ainsi, comme $i_{*}$ est exact, on
peut d\'efinir
$i_{+}$ au niveau des complexes de $\cal D\mit_{F}$-modules \`a cohomologie
born\'ee, par la formule
$$\begin{array}{cccc}
i_{+} : & A'(i) & \longrightarrow & C(p) \\
& E_{\bullet} & \mapsto &
i_{*}(\cal D\mit_{i}\otimes_{\cal D\mit_{F}}E_{\bullet})
\end{array}$$
o\`u $A'(i)$ est la cat\'egorie des complexes de $\cal D\mit_{F}$-modules
\`a cohomologie born\'ee et coh\'erente. Il est alors clair, que le morphisme induit
sur les spectres de $K$-th\'eorie
$$i_{+} : \bf K\mit^{\cal D}(F)\simeq K(A'(i)) \longrightarrow
K(C(p))\simeq \bf
K\mit^{\cal D}(F')$$
est \'egal au morphisme d\'efini pr\'ec\'edemment. \\

Passons au cas d'un morphisme lisse (~non-n\'ecessairement
repr\'esentable~)
$$p : F' \longrightarrow F''$$
Notons $\phi : \cal D\mit_{p,\bullet} \longrightarrow \cal D\mit_{p}$ la
r\'esolution $\cal D\mit_{F'}$-localement libre de De Rham
(~\cite[$5.3$, $(ii)$]{bo}~).
Notons alors $A'(p)$ la cat\'egorie bi-compliciale de Waldhausen des triplets
$$(E_{\bullet},F_{\bullet},u)$$
o\`u $E_{\bullet}$ est un complexe de $\cal D\mit_{F'}$-module \`a cohomologie
coh\'erente et born\'ee, $F_{\bullet}$ un complexe de $p^{-1}\cal
D\mit_{F''}$-modules acycliques, et
$$u : \cal D\mit_{p,\bullet}\otimes_{\cal D\mit_{F'}}E_{\bullet}
\longrightarrow F_{\bullet}$$
un quasi-isomorphisme.
On pose alors
$$\begin{array}{cccc}
p_{+} : & A'(p) & \longrightarrow & C(f) \\
& (E_{\bullet},F_{\bullet},u) & \mapsto &
p_{*}(F_{\bullet})
\end{array}$$
Il est alors clair, que le morphisme induit sur les spectres
de $K$-th\'eorie
$$p_{+} : \bf K\mit^{\cal D}(F')\simeq K(A'(p)) \longrightarrow
K(C(f))\simeq \bf
K\mit^{\cal D}(F'')$$
coincide avec celui d\'efini pr\'ec\'edemment. \\

Soit $A(f,i,p)$ la cat\'egorie bi-compliciale de Waldhausen dont les objets
sont les
$$(E_{\bullet},F_{\bullet},u,G_{\bullet},v,w)$$
avec
\begin{itemize}
\item $(E_{\bullet},F_{\bullet},u)$ un objet de $A(f)$
\item $v : \cal D\mit_{p,\bullet}\otimes_{\cal D\mit_{F'}}i_{*}(\cal
D\mit_{i}\otimes_{\cal D\mit_{F}}E_{\bullet}) \longrightarrow G_{\bullet}$
est un quasi-isomorphisme, et $G_{\bullet}$ un complexe de
$p^{-1}(\cal D\mit_{F''})$-modules acycliques
\item $w : i_{*}(F_{\bullet}) \longrightarrow G_{\bullet}$
est un quasi-isomorphisme de complexes de \\
$\cal D\mit_{F''}$-modules
\end{itemize}
Consid\'erons le diagramme "exact" de cat\'egories bi-compliciales de
Waldhausen
$$\xymatrix{
A(f) \ar[rr]^{f_{+}} \ar[dd]_{a} &  & C(f) \\
 & A(f,i,p) \ar[ul]_{b} \ar[dl]_{c} \ar[dr]^{d} & \\
A'(i) \ar[r]_{i_{+}} & C(p) & A'(p) \ar[uu]_-{p_{+}}\ar[l]_{e} }$$
Les foncteurs $a$, $b$, $c$, $d$ et $e$ \'etant les foncteurs canoniques
$$\begin{array}{cccc}
a : & A(f) & \longrightarrow & A'(i) \\
 & (E_{\bullet},F_{\bullet},u) & \mapsto & E_{\bullet}
\end{array}$$
$$\begin{array}{cccc}
b : & A(f,i,p) & \longrightarrow & A(f) \\
 & (E_{\bullet},F_{\bullet},u,G_{\bullet},v,w) & \mapsto &
 (E_{\bullet},F_{\bullet},u)
\end{array}$$
$$\begin{array}{cccc}
c : & A(f,i,p) & \longrightarrow & A'(i) \\
 & (E_{\bullet},F_{\bullet},u,G_{\bullet},v,w) & \mapsto &
 E_{\bullet}
\end{array}$$
$$\begin{array}{cccc}
d : & A(f,i,p) & \longrightarrow & A'(p) \\
 & (E_{\bullet},F_{\bullet},u,G_{\bullet},v,w) & \mapsto &
 (\cal D\mit_{p,\bullet}\otimes_{\cal D\mit_{F'}}i_{*}(\cal
D\mit_{i}\otimes_{\cal D\mit_{F}}E_{\bullet},G_{\bullet},v))
\end{array}$$
$$\begin{array}{cccc}
e : & A'(p) & \longrightarrow & C(p) \\
 & (E_{\bullet},F_{\bullet},u) & \mapsto & E_{\bullet}
\end{array}$$
Il est alors \'evident que $a\circ b=c$, et que
$i_{+}\circ c=e\circ d$.

De plus
$$(E_{\bullet},F_{\bullet},u,G_{\bullet},v,w) \mapsto
p_{*}w$$
induit un quasi-isomorphisme entre les foncteurs
$f_{+}\circ b$ et $p_{+}\circ d$.
Enfin, \`a l'aide de \cite[$1.9.7$, $1.9.8$]{th}, on montre que les
morphismes $e$ et $a$ induisent des
\'equivalences sur les cat\'egories d\'eriv\'ees, et donc des isomorphismes
dans $HoSp$ sur les spectres de $K$-th\'eorie.
Ainsi, lorsque que l'on consid\`ere
le diagramme induit au niveau des spectres de $K$-th\'eorie, on obtient
un diagramme commutatif dans $HoSp$
$$\xymatrix{
\bf K\mit^{\cal D}(F)\simeq K(A(f)) \ar[rr]^-{f_{+}} \ar[dd]_-{a} &  &
K(C(f))\simeq \bf K\mit^{\cal D}(F'') \\
 & K(A(f,i,p)) \ar[ul]_-{b} \ar[dl]_-{c} \ar[dr]^-{d} & \\
\bf K\mit^{\cal D}(F)\simeq K(A'(i)) \ar[r]^-{i_{+}} &
\bf K\mit^{\cal D}(F')\simeq K(C(p))\simeq \bf K\mit^{\cal D}(F') & K(A'(p))
\simeq \bf K\mit^{\cal D}(F')\ar[l]_-{e} \ar[uu]_-{p_{+}} }$$
Pour terminer la preuve, il nous reste \`a montrer que le morphisme
$$b : K(A(f,i,p)) \longrightarrow K(A(f))$$
est un isomorphisme. Pour cela il nous suffit de montrer que
le foncteur $b : A(f,i,p) \longrightarrow A(f)$ induit une \'equivalence
au niveau des cat\'egories d\'eriv\'ees
$$\bf R\mit b : D(A(f,i,p)) \longrightarrow D(A(f))$$
ce qui provient d'une nouvelle application des m\'ethodes de
\cite[$1.9.7$, $1.9.8$]{th}.
Montrons \`a titre d'exemple que $\bf R\mit  b$ est essentiellement
surjectif.

Soit $(E_{\bullet},F_{\bullet},u)$ un objet de $D(A(f))$. Choisissons un
r\'esolution injective dans la cat\'egorie des $p^{-1}\cal
D\mit_{F''}$-modules
$$v : \cal D\mit_{p,\bullet}\otimes_{\cal D\mit_{F'}}i_{*}(\cal
D\mit_{i}\otimes_{\cal D\mit_{F}}E_{\bullet}) \longrightarrow G_{\bullet}$$
Comme $\cal D\mit_{p,\bullet}$ est un complexe de $\cal
D\mit_{F'}$-modules localement libres, le morphisme naturel
$$can : \cal D\mit_{p,\bullet}\otimes_{\cal D\mit_{F'}}i_{*}(\cal
D\mit_{i}\otimes_{\cal D\mit_{F}}E_{\bullet}) \longrightarrow
i_{*}(i^{-1}\cal D\mit_{p,\bullet}\otimes_{i^{-1}\cal D\mit_{F'}}\cal
D\mit_{i}\otimes_{\cal D\mit_{F}}E_{\bullet})$$
est un isomorphisme dans $C(p)$.

De plus, comme $\cal D\mit_{i}$ est plat sur $\cal D\mit_{F}$, le
quasi-isomorphisme
$$\phi : \cal D\mit_{p,\bullet} \longrightarrow \cal D\mit_{p}$$
induit un quasi-isomorphisme
$$\phi : i_{*}(i^{-1}\cal D\mit_{p,\bullet}\otimes_{i^{-1}\cal
D\mit_{F'}}i_{*}(\cal
D\mit_{i}\otimes_{\cal D\mit_{F}}E_{\bullet})) \longrightarrow
i_{*}(i^{-1}\cal D\mit_{p}\otimes_{i^{-1}\cal D\mit_{F'}}\cal
D\mit_{i}\otimes_{\cal D\mit_{F}}E_{\bullet})$$
Or, on sait que $\cal D\mit_{f}$ est isomorphe \`a
$i^{-1}\cal D\mit_{p}\otimes_{i^{-1}\cal D\mit_{F'}}\cal
D\mit_{i}$. On a donc un quasi-isomorphisme
naturel
$$\phi \circ can^{-1} : \cal D\mit_{p,\bullet}\otimes_{\cal
D\mit_{F'}}i_{*}(\cal
D\mit_{i}\otimes_{\cal D\mit_{F}}E_{\bullet}) \longrightarrow
i_{*}(\cal D\mit_{f}\otimes_{\cal D\mit_{F}}E_{\bullet})$$
Enfin, comme $i_{*}$ est exact,
le morphisme
$$i_{*}(u) : i_{*}(\cal D\mit_{f}\otimes_{\cal D\mit_{F}}E_{\bullet})
\longrightarrow i_{*}(F_{\bullet})$$
est un quasi-isomorphisme. Par composition, ceci induit un
quasi-isomorphisme
$$\cal D\mit_{p,\bullet}\otimes_{\cal D\mit_{F'}}i_{*}(\cal
D\mit_{i}\otimes_{\cal D\mit_{F}}E_{\bullet}) \longrightarrow
i_{*}(F_{\bullet})$$
Donc, comme $v$ est une r\'esolution injective du membre de gauche, on
en d\'eduit qu'il existe un quasi-isomorphisme
$$w : i_{*}(F_{\bullet}) \longrightarrow G_{\bullet}$$
Ainsi, l'objet $(E_{\bullet},F_{\bullet},u,G_{\bullet},v,w)$ est un
"ant\'ec\'edent" de $(E_{\bullet},F_{\bullet},u)$ par
$\bf R\mit b$. $\Box$\\

\begin{lem}
\begin{enumerate}
\item Soit $i : F \longrightarrow F'$ et $j : F' \longrightarrow F''$
deux morphismes repr\'esentables, fini et non-ramifi\'es. Alors
$$(j\circ i)_{+}=j_{+}\circ i_{+} :
\bf K\mit^{\cal D} \longrightarrow \bf K\mit^{\cal D}(F'')$$
\item
Soit $p : F \longrightarrow F'$ et $q : F' \longrightarrow F''$
deux morphismes propres et lisses. Alors
$$(q\circ p)_{+}=q_{+}\circ p_{+} :
\bf K\mit^{\cal D} \longrightarrow \bf K\mit^{\cal D}(F'')$$
\end{enumerate}
\end{lem}

\underline{\bf Preuve:} \rm \\

$(1)$ Comme $\cal D\mit_{i}$ (~resp. $\cal D\mit_{j}$~) sont
localement libre sur $\cal D\mit_{F}$ (~resp. $\cal D\mit_{F'}$~), on
a, pour tout complexe de $\cal D\mit_{F}$-modules \`a cohomologie born\'ee
et coh\'erente $E_{\bullet}$, des isomorphismes fonctoriels en $E_{\bullet}$
$$\begin{array}{cl}
j_{+}\circ i_{+}(E_{\bullet}) & \simeq j_{*}(\cal
D\mit_{j}\otimes_{\cal D\mit_{F'}}i_{+}(E_{\bullet})) \\
 & \simeq j_{*}(\cal
D\mit_{j}\otimes_{\cal D\mit_{F'}}i_{*}(\cal D\mit_{i}\otimes_{\cal
D\mit_{F}}E_{\bullet}) )\\
 & \simeq j_{*}(i_{*}(i^{-1}\cal
D\mit_{j}\otimes_{i^{-1}\cal D\mit_{F'}}\cal D\mit_{i}\otimes_{\cal
D\mit_{F}}E_{\bullet})) \\
 & \simeq (j\circ i)_{*}(\cal D\mit_{j\circ i}\otimes_{\cal
D\mit_{F}}E_{\bullet}) \\
 & \simeq (j\circ i)_{+}(E_{\bullet})
\end{array}$$
Ce qui montre $(1)$. \\

$(2)$ Soit $\cal D\mit_{p,\bullet} \longrightarrow \cal D\mit_{p}$
(~resp. $\cal D\mit_{q,\bullet} \longrightarrow \cal D\mit_{q}$~) la
r\'esolution de De Rham (~\cite[$5.3$ $(ii)$]{bo}~). Alors, comme
$\cal D\mit_{q}$ est plat sur $q^{-1}\cal D\mit_{F'}$, le morphisme
naturel
$$\phi : p^{-1}\cal D\mit_{q,\bullet}\otimes_{q^{-1}\cal D\mit_{F'}}\cal
D\mit_{p,\bullet} \longrightarrow
p^{-1}\cal D\mit_{q}\otimes_{q^{-1}\cal D\mit_{F'}}\cal
D\mit_{p} \simeq \cal D\mit_{q\circ p}$$
est un quasi-isomorphisme.

Notons $A''(p)$ (~resp. $A''(q)$~) la cat\'egorie bi-compliciale de
Waldhausen des
complexes de $\cal D\mit_{F}$-modules (~resp. $\cal
D\mit_{F'}$-modules~) acycliques, et \`a cohomologie coh\'erente et born\'ee.

Remarquons que si $E$ est un $\cal D\mit_{F}$-module acyclique,
$\cal D\mit_{p,\bullet}\otimes_{\cal D\mit_{F}}E$ est un complexe
de $p^{-1}\cal D\mit_{F'}$-modules, localement isomorphes \`a une somme
directe de $E$. C'est donc un complexe de $p^{-1}\cal
D\mit_{F'}$-modules acycliques. Ainsi, on peut d\'efinir des foncteurs
exacts
$$\begin{array}{cccc}
p_{+} : & A''(p) & \longrightarrow & A''(q) \\
        & E_{\bullet} & \mapsto & p_{*}(\cal D\mit_{p,\bullet}
        \otimes_{\cal D\mit_{F}}E)
\end{array}$$
$$\begin{array}{cccc}
q_{+} : & A''(q) & \longrightarrow & C(q\circ p) \\
        & E_{\bullet} & \mapsto & q_{*}(\cal D\mit_{q,\bullet}
        \otimes_{\cal D\mit_{F'}}E)
\end{array}$$
$$\begin{array}{cccc}
(q\circ p)_{+} : & A''(p) & \longrightarrow & C(q\circ p) \\
        & E_{\bullet} & \mapsto & (q\circ p)_{*}(
p^{-1}\cal D\mit_{q,\bullet}\otimes_{p^{-1}\cal D\mit_{F'}}
\cal D\mit_{p,\bullet}\otimes_{\cal D\mit_{F}}E)
\end{array}$$
Il est imm\'ediat que ces foncteurs induisent les morphismes d\'efinis plus haut
$$p_{+} : \bf K\mit^{\cal D}(F) \longrightarrow \bf K\mit^{\cal
D}(F')$$
$$q_{+} : \bf K\mit^{\cal D}(F') \longrightarrow \bf K\mit^{\cal
D}(F'')$$
$$(q\circ p)_{+} : \bf K\mit^{\cal D}(F) \longrightarrow \bf K\mit^{\cal
D}(F'')$$
Or, comme $\cal D\mit_{q,\bullet}$ est localement libre sur
$\cal D\mit_{F'}$, on a des isomorphismes canoniques
$$\begin{array}{cl}
q_{+}\circ p_{+}(E_{\bullet}) & \simeq
q_{*}(\cal D\mit_{q,\bullet}\otimes_{\cal D\mit_{F'}}p_{*}(\cal
D\mit_{p,\bullet}\otimes_{\cal D\mit_{F}}E_{\bullet})) \\
& \simeq
q_{*}(p_{*}(p^{-1}\cal D\mit_{q,\bullet}\otimes_{p^{-1}\cal D\mit_{F'}}\cal
D\mit_{p,\bullet}\otimes_{\cal D\mit_{F}}E_{\bullet})) \\
& \simeq
(q\circ p)_{*}(p^{-1}\cal D\mit_{q,\bullet}\otimes_{p^{-1}\cal D\mit_{F'}}\cal
D\mit_{p,\bullet}\otimes_{\cal D\mit_{F}}E_{\bullet}) \\
& \simeq
(q\circ p)_{+}(E_{\bullet})
\end{array}$$
$\Box$ \\

On termine alors la preuve du point $(1)$ en utilisant un lemme formel
suivant.

\begin{lem}
Soit $H : Ob(\cal DM\mit) \longrightarrow Ob(HoSp)$, une application.
On suppose que pour tout morphisme propre
$f : F \longrightarrow F'$ dans $\cal LDM$, il existe un morphisme dans $HoSp$
$$f_{+} : H(F) \longrightarrow H(F')$$
qui v\'erifie les propri\'et\'es suivantes
\begin{itemize}
\item si $f=p\circ i$, avec $i$ un morphisme repr\'esentable fini et
non-ramifi\'e de $\cal LDM$, et $p$ un morphisme lisse, alors
$$f_{+}=p_{+}\circ i_{+}$$
\item si $f=i\circ p$, avec $i$ un morphisme repr\'esentable fini et
non-ramifi\'e de $\cal LDM$, et $p$ un morphisme lisse, alors
$$f_{+}=i_{+}\circ p_{+}$$
\item si $i : F \longrightarrow F'$ et $j : F' \longrightarrow
F''$ sont deux morphismes repr\'esentables finis, et non-ramifi\'es de
$\cal LDM$, alors
$$(j\circ i)_{+} = j_{+}\circ i_{+}$$
\item si $p : F \longrightarrow F'$ et $q : F' \longrightarrow
F''$ sont deux morphismes propres et lisses de
$\cal LDM$, alors
$$(q\circ p)_{+} = q_{+}\circ p_{+}$$
\end{itemize}
Alors, pour tout morphisme propre $f : F \longrightarrow F'$ et $g :
F' \longrightarrow F''$ de $\cal LDM$, on a
$$(g\circ f)_{+} = g_{+}\circ f_{+}$$
\end{lem}

\underline{\bf Preuve:} \rm On consid\`ere le diagramme commutatif
suivant de $\cal LDM$
$$\xymatrix{
F \ar[r]^-{f} \ar[d]_-{i} & F' \ar[d]^-{j} \ar[r]^-{g} & F'' \\
F\times F' \ar[d]_-{k} \ar[ru]_-{p} \ar[r]^-{h} & F'\times F'' \ar[ru]_-{q} & \\
F\times F'\times F' \times F'' \ar[ru]_-{r}& & }$$
o\`u $i$, $k$ et $j$ sont les graphes de $f$, $j\circ p$ et $g$, et
$p$, $q$, $r$ les projections. On a alors
$$\begin{array}{cl}
g_{+} \circ f_{+} & = q_{+} \circ j_{+} \circ p_{+} \circ i_{+} \\
& = q_{+} \circ h_{+} \circ i_{+} \\
& = q_{+} \circ r_{+} \circ k_{+} \circ i_{+} \\
& = (q\circ r)_{+} \circ (k\circ i)_{+} \\
& = (g \circ f)_{+}
\end{array}$$
$\Box$\\

\underline{Preuve de $(2)$:} Il faut montrer que
$\bf R\mit f_{+}$ pr\'eserve les propri\'et\'es d'\^etre \`a cohomologie
holonome et \`a cohomologie holonome r\'egulier. La fonctorialit\'e se
d\'emontre alors exactement comme dans le point $(1)$. \\

\begin{lem}
Le foncteur $\bf R\mit f_{+}$ preserve la propri\'et\'e d'\^etre "\`a
cohomologie holonome (~resp. "\`a cohomologie holonome r\'egulier"~).
\end{lem}

\underline{\bf Preuve:} \rm Comme c'est une assertion locale sur
$F'_{et}$, le cas o\`u $f$ est repr\'esentable se traite comme
celui des sch\'emas (~\cite[$10.1$, $12.2$]{bo}~). De plus, on peut toujours
supposer que
$F'$ est un sch\'ema affine $Y$.\\

Soit $F_{1} \hookrightarrow F$ le support r\'eduit d'un $\cal
D\mit_{F}$-module holonome (~resp. holonome r\'egulier~) $\cal M$. Soit
$U \hookrightarrow F_{1}$ un sous-champ ouvert dense, lisse sur $k$,
qui est une gerbe sur son espace de modules $V$, que l'on peut
supposer lisse et affine. Notons
$$i : U \hookrightarrow F$$
l'immersion canonique. On consid\`ere alors le triangle dans la
cat\'egorie d\'eriv\'ee des $\cal D\mit_{F}$-modules \`a cohomologie born\'ee et
holonome (~resp. holonome r\'egulier~)
$$\xymatrix{
& \cal N\mit \ar[rd]^{a} & \\
i_{+}\circ i^{!}\cal M\mit \ar[ur]^{+1} & & \ar[ll] \cal M\mit }$$
o\`u $\cal N$ est le cone du morphisme naturel
$$a : \cal M\mit \longrightarrow \bf R\mit i_{+}\circ i^{!} \cal M$$
Alors, comme $a$ et un isomorphisme sur $U$, $\cal N$ poss\`ede un
support de dimension strictement plus petite que celui de $\cal M$. En
raisonnant par r\'ecurrence, il nous suffit donc de d\'emontrer le lemme
pour des $\cal D\mit_{F}$-modules de la forme
$i_{+}\circ \cal M$, o\`u $i : U \hookrightarrow F$ est une immersion
localement ferm\'ee, $U$ une gerbe sur un sch\'ema lisse et affine $V$,
et $\cal M$ un $\cal D\mit_{U}$-module holonome (~resp. holonome
r\'egulier~). Ainsi, on peut supposer que $F$ est une gerbe sur un
sch\'ema affine et lisse $X$. \\

Le morphisme $f : F \longrightarrow Y$ se factorise alors par
$$f : \xymatrix{F \ar[r]^{p} & X \ar[r]^{g} & Y}$$
o\`u $p$ est la projection sur l'espace de modules. Le cas de
$\bf R\mit g_{+}$ est d\'ej\`a connu, car $g$ est repr\'esentable. Il nous reste
donc \`a
traiter le cas de  $\bf R\mit p_{+}$. En localisant sur $X_{et}$, on
se ram\`ene au cas o\`u $F=[X/H]$, avec $H$ un groupe fini op\'erant
trivialement sur $X$. Alors, la donn\'ee d'un $\cal D\mit_{F}$-module est
\'equivalente \`a la donn\'ee d'un $\cal D\mit_{X}$-module $H$-\'equivariant.
De plus le foncteur $\bf R\mit p_{+}\simeq p_{+}$ est exact, et
associe \`a un $\cal D\mit_{X}$-module \'equivariant $\cal M$, le sous $\cal
D\mit_{X}$-module des invariants par $H$, $\cal M\mit^{H}$.
Comme ce sous-modules est un facteur direct de $\cal M$, il est
holonome (~resp. holonome r\'egulier~) si $\cal M$ l'est. $\Box$\\

\end{subsubsection}

\begin{subsubsection}{Comparaison entre $\bf G$ et $\bf K\mit^{\cal D}$}
\hspace{5mm}
Le but de ce paragraphe est de d\'emontrer le th\'eor\`eme suivant.

\begin{thm}\label{th4.1}
Pour tout champ $F$ dans $\cal LDM$, il existe un isomorphisme dans $HoSp$
$$\delta_{F}^{-1} : \bf G\mit(F) \longrightarrow \bf K\mit^{\cal D}(F)$$
qui commute avec les images directes de morphismes propres.
\end{thm}

\underline{\bf Preuve:} \rm On d\'efinit un foncteur exact
$$\begin{array}{cccc}
\delta_{F}^{-1} : & \bf Coh\mit(F) & \longrightarrow & \bf Coh\mit^{\cal D} \\
       & \cal E\mit     & \mapsto & \cal D\mit_{F}\otimes_{\cal
       O\mit_{F}}(\cal E\mit\otimes_{\cal O\mit_{F}}\omega_{F}^{\vee})
\end{array}$$
o\`u $\omega_{F}^{\vee}$ est le dual du faisceau inversible canonique
$\omega_{F}:=Det\Omega_{F}^{1}$.

Ce foncteur induit un morphisme sur les spectres de $K$-th\'eorie
$$\delta_{F}^{-1} : \bf G\mit(F) \longrightarrow \bf K\mit^{\cal D}(F)$$

\begin{lem}
Le morphisme $\delta_{F}^{-1}$ est un isomorphisme dans $HoSp$
\end{lem}

\underline{\bf Preuve:} \rm Comme tout $\cal D$-module
admet une bonne filtration, le d\'evissage de \cite[$6.7$]{q}, implique
que le foncteur $\cal E\mit \mapsto \cal D\mit_{F}\otimes_{\cal
O\mit_{F}}\cal E$ induit un isomorphisme entre les spectres
$\bf G\mit(F)$ et $\bf K\mit^{\cal D}(F)$. De plus, comme
$\omega_{F}^{\vee}$ est inversible, le foncteur
$\cal E\mit \mapsto \cal E\mit\otimes_{\cal
O\mit_{F}}\omega_{F}^{\vee}$ induit aussi une \'equivalence sur le
spectre $\bf G\mit(F)$. Ainsi, par composition $\delta_{F}^{-1}$ est
un isomorphisme dans $HoSp$. $\Box$\\

Pour la covariance de $\delta_{F}^{-1}$, on factorise tout morphisme
propre $f~:~F~\longrightarrow~F'$ par son graphe
$$f : \xymatrix{F \ar[r]^{i} & F\times F' \ar[r]^{p} & F'}$$
On se ram\`ene donc \`a deux cas, celui o\`u $f=i$ est repr\'esentable fini et
non-ramifi\'e, et celui o\`u $f=p$ est un morphisme lisse. \\

Commen\c{c}ons par le cas o\`u $f$ est fini et non-ramifi\'e. Comme $i_{+}$
est alors exact, on peut r\'ealiser l'image directe en $K$-th\'eorie par
la foncteur
$$\begin{array}{cccc}
i_{+} : & \bf Coh\mit^{\cal D}(F) & \longrightarrow & \bf
Coh\mit^{\cal D}(F') \\
& \cal M\mit & \mapsto & i_{*}(\cal D\mit_{i}\otimes_{\cal
D\mit_{F}}\cal M \mit)
\end{array}$$
Comme $\cal D\mit_{i}$ est localement libre sur $\cal D\mit_{F}$, on a
des isomorphismes canoniques
$$\begin{array}{cl}
i_{+}\circ \delta_{F}^{-1}(\cal E\mit) & \simeq i_{*}(\cal
D\mit_{i}\otimes_{\cal
D\mit_{F}}\rho(\cal E \mit)) \\
& \simeq i_{*}(\cal D\mit_{i}\otimes_{\cal D\mit_{F}}\cal
D\mit_{F}\otimes_{\cal O\mit_{F}}\cal E \mit\otimes_{\cal
O\mit_{F}}\omega_{F}^{\vee}) \\
& \simeq i_{*}(\cal D\mit_{i}\otimes_{\cal O\mit_{F}}\cal E \mit\otimes_{\cal
O\mit_{F}}\omega_{F}^{\vee}) \\
& \simeq i_{*}(i^{-1}(\cal D\mit_{F'}\otimes_{\cal
O\mit_{F'}}\omega_{F'}^{\vee})\otimes_{f^{-1}\cal
O\mit_{F'}}\omega_{F}\otimes_{\cal O\mit_{F}}
\cal E \mit\otimes_{\cal O\mit_{F}}\omega_{F}) \\
& \simeq i_{*}(i^{-1}(\cal D\mit_{F'}\otimes_{\cal
O\mit_{F'}}\omega_{F'}^{\vee})\otimes_{f^{-1}\cal
O\mit_{F'}}\cal E \mit) \\
& \simeq \cal D\mit_{F'}\otimes_{\cal
O\mit_{F'}}\omega_{F'}^{\vee}\otimes_{\cal O\mit_{F'}}i_{*}(\cal
E\mit) \\
& \simeq \delta_{F}^{-1}(i_{*}(\cal E\mit))
\end{array}$$

Le cas o\`u $f$ est un morphisme lisse se traite d'une mani\`ere
analogue, en utilisant la r\'esolution de De Rham de $\cal D\mit_{f}$. $\Box$\\

\begin{df}\label{d4.2}
On pose
$$\delta_{F}:=(\delta_{F}^{-1})^{-1} : \bf K\mit^{\cal D}(F)
\longrightarrow \bf G\mit(F)$$
\end{df}

En localisant cette construction pour la topologie \'etale, on obtient une
comparaison entre $\underline{\bf G}$ et $\underline{\bf K}^{\cal D}$.

\begin{thm}\label{th4.2}
Pour tout champ $F$ dans $\cal LDM$, il existe un isomorphisme dans $HoSp$
$$\delta_{F}^{-1} : \underline{\bf G}(F) \longrightarrow
\underline{\bf K}^{\cal D}(F)$$
qui commute avec les images directes de morphismes propres
repr\'esentables.

De plus, on a un diagramme commutatif dans $HoSp$
$$\xymatrix{
\bf G\mit(F) \ar[r]^-{\delta_{F}^{-1}} \ar[d]_{can} &
\bf K\mit^{\cal D}(F) \ar[d]^{can} \\
\underline{\bf G}(F) \ar[r]^-{\delta_{F}^{-1}} &
\underline{\bf K}^{\cal D}(F) }$$
\end{thm}

\underline{\bf Preuve:} \rm C'est la m\^eme que celle de \ref{th4.2}.
$\Box$\\

Le point fondamental que l'on utilisera par la suite, est le calcul
explicite de $\delta_{F}$ au niveau du $K_{0}$.

\begin{prop}\label{p4.4}
Soit $\cal M$ un $\cal D\mit_{F}$-module coh\'erent sur un champ
$F$ de $\cal LDM$. Soit $\cal M\mit^{i} \hookrightarrow \cal
M\mit^{i+1} \hookrightarrow \cal M$ une bonne filtration de $\cal M$,
et $Gr\cal M$ son module gradu\'e, vu comme faisceau coh\'erent sur le
champ cotangent $T^{*}F$. Alors, on a une \'egalit\'e dans $\bf
G\mit_{0}(F)$
$$\delta_{F}([\cal M\mit])=e^{*}[Gr\cal M\mit]$$
o\`u $e : F \hookrightarrow T^{*}F$ est la section nulle.

On a aussi
$$\delta_{F}([\cal M\mit])=DR(\cal M\mit)$$
o\`u $DR(\cal M\mit)$ est la classe du complexe de De Rham de
$\cal M$ dans $\bf G\mit_{0}(F)$ (~au sens de \cite{a}~).
\end{prop}

\underline{\bf Preuve:} \rm C'est la m\^eme que \cite[$1.2.1$]{al}. $\Box$ \\

\end{subsubsection}

\end{subsection}

\begin{subsection}{Etude des spectres de $K$-th\'eorie
des $\cal D$-modules holonomes}

\begin{subsubsection}{Le lemme de Kashiwara}
\hspace{5mm} Dans cette section nous allons rappeler Le lemme de
Kashiwara et en d\'eduire de tr\`es fortes propri\'et\'es pour les spectres
de $K$-th\'eorie des $\cal D$-modules holonomes.

\begin{thm}{(~"Lemme de Kashiwara"~)}\label{kas}
Soit $j : F' \hookrightarrow F$ une immersion ferm\'ee de
champs de $\cal LDM$, et $i : U \hookrightarrow F$ l'immersion
ouverte compl\'ementaire. Alors on a
$$j^{!}\circ j_{+}=Id : \bf K\mit^{h}(F)
\longrightarrow \bf K\mit^{h}(F)$$
et
$$j_{+}\circ j^{!}+i_{+}\circ i^{!}=Id : \bf K\mit^{h}(F)
\longrightarrow \bf K\mit^{h}(F)$$
dans $HoSp$.

Il en est de m\^eme pour les foncteurs $\underline{\bf K}^{h}$,
$\bf K\mit^{h,r}$ et $\underline{\bf K}^{h,r}$.
\end{thm}

\underline{\bf Preuve:} \rm Donnons la d\'emonstration pour le cas du
foncteur $\bf K\mit^{h}$. Le cas de $\bf K\mit^{h,r}$ se traite de la
m\^eme fa\c{c}on. Par localisation on obtiendra alors le r\'esultat pour les
foncteurs $\underline{\bf K}^{h}$ et $\underline{\bf
K}^{h,r}$. \\

Soit $\cal HOL\mit(F \; on \; F')$ la cat\'egorie des $\cal D\mit_{F}$-modules
holonome dont le support est contenu dans $F'$. Les m\^emes calculs que
ceux faits dans \cite[$7.11$]{bo} montrent que les foncteurs
$j_{+}$ et $j^{!}$ d\'efinissent des \'equivalences, "inverses l'une de
l'autres", entre les cat\'egories
$\cal HOL \mit(F')$ et $\cal HOL\mit(F \; on \; F')$ (~c'est le "lemme
de Kashiwara"~).

Soit $\bf K\mit^{h}(F \; on \; F')$ le spectre de $K$-th\'eorie de
la cat\'egorie $\cal HOL\mit(F \; on \; F')$. On notera
$$c : \bf K\mit^{h}(F \; on \; F') \longrightarrow \bf K\mit(F)$$
le morphisme induit par l'inclusion naturelle.
Notons aussi
$$j_{+}^{F'} : \bf K\mit_{*}^{h}(F') \longrightarrow \bf
K\mit^{h}_{*}(F \; on \; F')$$
$$j_{F'}^{!} : \bf K\mit_{*}^{h}(F \; on \; F') \longrightarrow
\bf K\mit^{h}_{*}(F')$$
les morphisme induits par les foncteurs exacts $j_{+}$ et $j^{!}$. Ils
v\'erifient clairement $j_{F'}^{!}\circ j^{F'}_{+}=Id$ et
$j^{F'}_{+} \circ j^{!}_{F'}=Id$. De plus, on a
$j_{+}=c \circ j_{+}^{F'}$ et $j^{!}\circ c=j^{!}_{F'}$.

On peut alors \'ecrire
$$\begin{array}{cl}
j^{!}\circ j_{+}& = j^{!}\circ c \circ j_{+}^{F'} \\
& = j^{!}_{F'} \circ j_{+}^{F'} \\
& = Id
\end{array}$$

Le th\'eor\`eme de localisation de \cite[$5.5$]{q}, implique que l'on a un
triangle dans $HoSp$
$$\xymatrix{
& \bf K\mit^{h}(F \; on \; F') \ar[rd]^{c} & \\
\bf K\mit^{h}(U) \ar[ru]^{-1} & & \ar[ll]_{i^{!}} \bf K\mit(F)}$$
De plus, comme $i^{!}$ poss\`ede une section $i_{+} : \bf K\mit^{h}(U)
\longrightarrow \bf K\mit^{h}(F)$, ce triangle se scinde, et donne
lieu \`a une suite exacte courte scind\'ee
$$\xymatrix{
0 \ar[r] & \bf K\mit_{*}^{h}(F \; on \; F') \ar[r]^{c}
& \bf K\mit^{h}_{*}(F) \ar[r]^{i^{!}} & \bf K\mit^{h}_{*}(U)
\ar[r] & 0}$$
Et en utilisant l'\'equivalence entre $\cal HOL\mit(F \; on \; F')$ et
$\cal HOL\mit(F')$, on trouve une suite exacte scind\'ee
$$\xymatrix{
0 \ar[r] & \bf K\mit_{*}^{h}(F') \ar[r]^{j_{+}}
& \bf K\mit^{h}_{*}(F) \ar[r]^{i^{!}} & \bf K\mit^{h}_{*}(U)
\ar[r] & 0}$$
Ecrivons alors
$$Id=\alpha + i_{+}i^{!}$$
Ainsi, $i^{!}=i^{!}\circ \alpha + i^{!}$, et donc
$i^{!}\circ \alpha=0$. Ceci montre que l'on peut \'ecrire
$\alpha=j_{+} \circ \beta$, avec
$\beta : \bf K\mit_{*}^{h}(F') \longrightarrow \bf K\mit_{*}^{h}(F)$.
Le fait que $j^{!}\circ j_{+}=Id$ implique alors que
$j^{!}=\beta + j^{!}\circ i_{+}\circ i^{!}$.

Il ne nous reste qu'\`a montrer que $j^{!}\circ i_{+}=0$. \\

Soit $A$ la cat\'egorie bi-compliciale Waldhausen des triplets
$(E_{\bullet},F_{\bullet},u)$, o\`u $E_{\bullet}$ est un complexe
de $\cal D\mit_{U}$-modules acycliques \`a cohomologie
holonome et born\'ee, et
$$u : i_{*}E_{\bullet}
\longrightarrow F_{\bullet}$$
est un quasi-isomorphisme de complexes de $\cal D\mit_{F}$-modules
avec $F_{\bullet}$ un complexe de $\cal D\mit_{F}$-modules plats.
Soit $B$ la cat\'egorie des complexes de $\cal D\mit_{F'}$-modules
\`a cohomologie holonome et born\'ee. On veut montrer que l'unique
transformation naturelle de foncteurs exacts
$$\begin{array}{cccc}
0 \Rightarrow i_{+}\circ j^{!} : & A & \longrightarrow & B \\
& (E_{\bullet},F_{\bullet},u)  & \mapsto &
j^{*}(F_{\bullet})
\end{array}$$
est un quasi-isomorphisme. Comme ceci est local sur $F_{et}$, on peut
se ramener au cas o\`u $F$ est un sch\'ema affine. Le r\'esultat provient
alors de \cite[$VI$ $8.5$]{bo}. $\Box$\\

\begin{cor}\label{c4.1}
Soit
$$\xymatrix{
F'_{1} \ar[r]^{j'} \ar[d]_{f_{1}} & F' \ar[d]^{f} \\
F_{1} \ar[r]_{j} & F}$$
un diagramme cart\'esien de champ de $\cal LDM$, avec $j$ une immersion
localement ferm\'ee. Alors on a
$$j^{!}\circ f_{+} = (f_{1})_{+} \circ (j')^{!} :
\bf K\mit^{h}(F) \longrightarrow \bf K\mit^{h}(F)$$

Il en est de m\^eme pour les foncteurs  $\underline{\bf K}^{h}$,
$\bf K\mit^{h,r}$ et $\underline{\bf K}^{h,r}$.
\end{cor}

\underline{\bf Preuve:} \rm Comme pour le th\'eor\`eme, nous ne donnerons
la preuve que dans le cas du foncteur $\bf K\mit^{h}_{*}$.\\

De la cas o\`u $j$ est une immersion
ouverte, le corollaire provient du fait que les images directes sont
compatibles avec les changements de bases \'etales. \\

Supposons maintenant que $j$ est une immersion ferm\'ee. Notons \\
$i~:~F~-~F_{1}~\hookrightarrow~F$ et $i' : F'-F'_{1} \hookrightarrow
F'$, les immersions ouvertes compl\'ementaires, et
$g : F'-F'_{1} \longrightarrow F-F_{1}$ le morphisme induit.

Soit $x \in \bf K\mit_{*}^{h}(F')$. A l'aide du th\'eor\`eme \ref{kas}
on peut \'ecrire
$$x=j'_{+}(j')^{!}(x) + i'_{+}(i')^{+}(x)$$
On a donc
$$\begin{array}{cl}
j^{!}f_{+}(x) & = j^{!}f_{+}j'_{+}(j')^{!}(x) +
j^{!}f_{+}i'_{+}(i')^{!}(x) \\
& = j^{!}j_{+}(f_{1})_{+}(j')^{!}(x) +
j^{!}i_{+}g_{+}i'_{+}(x) \\
& = (f_{1})_{+}(j')^{!}(x)
\end{array}$$
car $j^{!}i_{+}=0$. \\

Dans le cas g\'en\'eral, on factorise $j$ en une immersion ferm\'ee suivie
d'une immersion ouverte, et on utilise les deux cas pr\'ec\'edents. $\Box$\\

\begin{cor}\label{c4.2}
Soit $F$ un champ de $\cal LDM$, et $j_{i} : F_{i} \hookrightarrow F$
des sous-champs localement ferm\'es de $\cal LDM$ formant une
stratification de $F$. Alors
$$\sum_{i} : (j_{i})_{+}\circ j_{i}^{!} = Id : \bf K\mit^{h}(F)
\longrightarrow \bf K\mit^{h}(F)$$

Il en est de m\^eme pour les foncteurs $\underline{\bf K}^{h}$,
$\bf K\mit^{h,r}$ et $\underline{\bf K}^{h,r}$.
\end{cor}

\underline{\bf Preuve:} \rm Ce corollaire se d\'eduit imm\'ediatement du
th\'eor\`eme pr\'ec\'edent par une r\'ecurrence sur le nombre de sous-champs
$F_{i}$. $\Box$\\

\end{subsubsection}

\begin{subsubsection}{Fonctions constructibles \`a valeurs dans
un foncteur et th\'eor\`emes de d\'evissage}

Commen\c{c}ons par introduire la notion de fonctions
constructibles \`a valeurs dans un foncteur.

Rappelons pour cela que pour tout point $x \in |F|$, on dispose d'un morphisme
canonique
$$i_{x} : \widetilde{x} \longrightarrow F$$
Si $x$ est un point de $F$, et $H : \cal LDM\mit \longrightarrow Ab$
un foncteur contravariant, nous noterons
$$H(\widetilde{x}):=colim_{U}H(U)$$
o\`u $U$ parcourt les sous-champs ouverts non-vides et lisses de
$\overline{\left\{x\right\}}$. On dispose alors d'un morphisme de
restriction
$$i_{x}^{*} : H(F) \longrightarrow H(\widetilde{x})$$

\begin{df}
Soit $H : \cal LDM\mit \longrightarrow Ab$ un foncteur contravariant.
On d\'efinit le groupe des fonctions constructibles sur $F$ \`a coefficients
dans $H$ par
\begin{multline*}
F_{ct}H(F):=\{ (\sigma_{x}) \in \prod_{x \in |F|}H(\widetilde{x})
\; / \; \forall x \in |F|, \\
\exists U \; ouvert \; dense \; de \;
\overline{\{x\}}, \; et \; \sigma_{U} \in H(U) \; / \;
\forall y \in |U| \; \sigma_{y}=i_{y}^{*}(\sigma_{U}) \}	
\end{multline*}
\end{df}

Par exemple, si $H$ est un foncteur constant d\'efini par un groupe
ab\'elien $A$, $F_{ct}H(F)$ coincide avec le groupe des fonctions
constructibles de $M$ \`a valeurs dans $A$, o\`u $M$ est l'espace de
modules de $F$.

Notons qu'il existe un morphisme naturel
$$\begin{array}{ccc}
H(F) & \longrightarrow & F_{ct}H(F) \\
\sigma & \mapsto & ((i_{x}^{*}(\sigma))_{x \in |F|})
\end{array}$$

\begin{prop}\label{p4.5}
Pour tout foncteur contravariant
$$H : \cal LDM\mit \longrightarrow Ab$$
on peut munir $F \mapsto F_{ct}H(F)$ d'une unique structure de foncteur
contravariant compatible avec le morphisme naturel
$$H(F) \longrightarrow F_{ct}H(F)$$
\end{prop}

\underline{\bf Preuve:} \rm Soit $f : F \longrightarrow F'$ un
morphisme de $\cal LDM$, et $(\sigma_{x})$ un \'el\'ement de
$F_{ct}H(F')$. Soit $y \in |F|$ un point de $F$ au-dessus d'un point
$x$ de $F'$, et $f_{y} : \widetilde{y} \longrightarrow \widetilde{x}$
le morphisme induit. On pose alors
$$f^{*}(\sigma_{x})_{y}:=f_{y}^{*}(\sigma_{x})$$
Cette d\'efinition v\'erifie clairement les hypoth\`eses demand\'ees. $\Box$\\

\begin{thm}\label{th4.3}
Il existe des isomorphismes compatibles avec les images
r\'eciproques sur $\cal LDM$
$$\bf K\mit_{*}^{h}(F) \simeq F_{ct}\bf K\mit_{*}^{\nabla}(F)$$
$$\bf K\mit_{*}^{h,r}(F) \simeq F_{ct}\bf K\mit_{*}^{\nabla,r}(F)$$
$$\underline{\bf K}_{*}^{h}(F) \simeq
F_{ct}\underline{\bf K}_{*}^{\nabla}(F)$$
$$\underline{\bf K}_{*}^{h,r}(F) \simeq F_{ct}\underline{\bf
K}_{*}^{\nabla,r}(F)$$
\end{thm}

\underline{\bf Preuve:} \rm Comme pour le th\'eor\`eme \ref{kas}, on
ne d\'emontrera ce r\'esultat que pour le foncteur
$\bf K\mit_{*}^{h}$. \\

Pour un champ $F$ (~suppos\'e connexe~) dans $\cal LDM$, on consid\`ere
$$\gamma_{F}:=\prod_{x \in |F|} j_{x}^{!} : \bf K\mit_{*}^{h}(F)
\longrightarrow \prod_{x \in |F|}\bf K\mit_{*}^{h}(\widetilde{x})$$
Comme la cat\'egorie des $\cal D\mit_{\widetilde{x}}$-modules est \'equivalente
\`a la cat\'egorie des connexions sur $\widetilde{x}$, la continuit\'e du
foncteur de $K$-th\'eorie implique que l'inclusion naturelle induit
un isomorphisme
$$\bf K\mit_{*}^{\nabla}(\widetilde{x})\simeq \bf K\mit_{*}^{h}(\widetilde{x})$$
Le morphisme obtenu
$$\bf K\mit_{*}^{h}(F) \longrightarrow \prod_{x \in |F|}\bf
K\mit_{*}^{\nabla}(\widetilde{x})$$
est contravariant par d\'efinition.

Remarquons alors que l'image d'un \'el\'ement $\sigma \in \bf K\mit_{*}^{h}(F)$ par
$\gamma_{F}$ est un \'el\'ement constructible dans
$\prod_{x \in |F|}\bf K\mit_{*}^{\nabla}(\widetilde{x})$. Ceci est \'equivalent
au fait que pour tout \'el\'ement $\sigma \in \bf K\mit_{*}^{h}(F)$, il existe un
sous-champ ouvert $i : U \hookrightarrow F$ tel que $i^{!}\sigma$ est dans
l'image de $\bf K\mit^{\nabla}_{*}(U)$.
Soit $\xi$ le point g\'en\'erique de $F$. Comme le morphisme naturel
$$\bf K\mit^{\nabla}_{*}(\widetilde{\xi}) \longrightarrow
\bf K\mit^{h}_{*}(\widetilde{\xi})$$
est un isomorphisme, pour tout $\sigma \in \bf K\mit_{*}^{h}(F)$, il existe un
sous-champ ouvert lisse $U \hookrightarrow F$ non vide, un \'el\'ement $\sigma' \in
\bf K\mit^{\nabla}(U)$, tel que
$i_{\xi}^{!}(\sigma)=i_{\xi}^{!}(\sigma')$. Ce qui implique
 qu'il existe un ouvert lisse $j : V \hookrightarrow U
\hookrightarrow F$, non vide, avec
$j^{!}(\sigma)=j^{!}(\sigma')$. Ainsi, $j^{!}(\sigma)$ est dans
l'image de $\bf K\mit^{\nabla}_{*}(V)$. \\

Montrons alors que $\gamma_{F}$ induit un isomorphisme de $\bf
K\mit_{*}^{h}(F)$ sur $F_{ct}\bf K\mit_{*}^{\nabla}(F)$.

Soit $\sigma \in \bf K\mit_{*}^{h}(F)$, et supposons que
$\gamma_{F}(\sigma)=0$. En particulier si $\xi$ est le point
g\'en\'erique de $F$, $i_{\xi}^{!}(\sigma)=0$, et
donc il existe un sous-champ ouvert non-vide lisse $j_{0} : F_{0}
\hookrightarrow
F$ tel que $j_{0}^{!}(\sigma)=0$. Soit $x$ un des points g\'en\'eriques
du compl\'ementaire r\'eduit $F'=(F-U)_{red}$. Alors $i_{x}^{!}(\sigma)=0$.
On peut donc trouver un ouvert non-vide lisse
de $F_{1}$, $j_{1} : F_{1} \hookrightarrow
F_{0} \hookrightarrow F$, tel que $j_{1}^{!}(\sigma)=0$. Ainsi, par
r\'ecurrence noeth\'erienne on trouve une stratification de $F$ par des
sous-champs localement ferm\'es et lisses, \\
$j_{i} : F_{i} \longrightarrow F$, tel que
$$j_{i}^{!} (\sigma)=0 \; \forall \; i$$
Le corollaire \ref{c4.2} implique alors que
$$\sigma=\sum_{i}(j_{i})_{+}\circ j_{i}^{!}(\sigma)=0$$
Ceci montre que $\gamma$ est injective.

Soit $(\sigma_{x}) \in F_{ct}\bf K\mit^{\nabla}_{*}(F)$. Par
d\'efinition, il existe une stratification de $F$ par des sous-champs
localement ferm\'es et lisses, $j_{i} : F_{i} \hookrightarrow F$, et des
\'el\'ements $\sigma_{i} \in \bf K\mit^{\nabla}(F_{i})$, tels que
$$i_{x}^{!}(\sigma_{i})=\sigma_{x} \quad \forall \; x \in |F_{i}|$$
Posons
$$\sigma=\sum_{i}(j_{i})_{+}(\sigma_{i})$$
Alors, le corollaire \ref{c4.1}, implique que, pour $x \in  |F_{i}|$
$$i_{x}^{!}(\sigma)=i_{x}^{!}(\sigma_{i})=\sigma_{x}$$
Ce qui montre que $\gamma_{F}(\sigma)=(\sigma_{x})$. $\Box$\\

Par la suite, si $f : F \longrightarrow F'$ est un morphisme
de $\cal LDM$, nous noterons
$$f_{+} : F_{ct}\bf K\mit_{*}^{\nabla}(F) \longrightarrow
F_{ct}\bf K\mit_{*}^{\nabla}(F')$$
le morphisme d\'efini par le carr\'e commutatif suivant
$$\xymatrix{
\bf K\mit^{h}_{*}(F) \ar[r]^{f_{+}} \ar[d]_{\gamma_{F}} &
\bf K\mit^{h}_{*}(F') \ar[d]^{\gamma_{F'}} \\
F_{ct}\bf K\mit^{\nabla}_{*}(F) \ar[r]_{f_{+}} & F_{ct}\bf
K\mit^{\nabla}_{*}(F')}$$
On d\'efinit de la m\^eme mani\`ere des images directes sur les foncteurs
$F_{ct}\bf K\mit_{*}^{\nabla,r}$, $F_{ct}\underline{\bf
K}_{*}^{\nabla}$ et $F_{ct}\underline{\bf K}_{*}^{\nabla}$.
Remarquons que pour les deux derniers, ces images directes ne sont
d\'efinies que pour des morphismes repr\'esentables. \\

Par la suite nous identifierons $\bf K\mit^{h}_{*}$ avec
$F_{ct}\bf K\mit^{\nabla}_{*}$ et
$\bf K\mit^{h,r}_{*}$ avec $F_{ct}\bf K\mit^{\nabla,r}_{*}$
par les isomorphismes \ref{th4.3}.

\begin{prop}\label{p4.6}
Les foncteurs covariants $\underline{\bf K}^{h}_{*}(-)$ et
$\underline{\bf K}^{h,r}_{*}(-)$, de la cat\'egorie
$\cal LDM$ et morphismes repr\'esentables, s'\'etendent de fa\c{c}on unique
en des foncteurs covariants sur $\cal LDM$.

On a alors, pour tout morphisme de $\cal LDM$,
$f : F \longrightarrow F'$ un diagramme commutatif
$$\xymatrix{
\underline{\bf K}^{h}_{*}(F) \ar[r]^-{\gamma_{F}} \ar[d]_{f_{+}}
& \underline{\bf G}_{*}(F) \ar[d]^{f_{*}} \\
\underline{\bf K}^{h}_{*}(F') \ar[r]^-{\gamma_{F'}}
& \underline{\bf G}_{*}(F')}$$
\end{prop}

\underline{\bf Preuve:} \rm Nous ne d\'emontrerons cette proposition
que pour le foncteur $\underline{\bf K}^{h}_{*}(-)$, le
second cas se traitant de fa\c{c}on similaire.\\

Commen\c{c}ons par l'unicit\'e. \\

Soit $f : F \longrightarrow F'$ un morphisme de $\cal LDM$, et
$\sigma~\in~\underline{\bf K}^{h}_{*}(F)$. On peut choisir
une stratification $j : F_{i} \hookrightarrow F$ de $F$, telle que
chaque $F_{i}$ soit une gerbe triviale sur un sch\'ema lisse $X_{i}$, et de
plus qu'il existe $\sigma_{i} \in \underline{\bf
K}^{h}_{*}(F_{i})$, avec $\sigma=\sum_{i}(j_{i})_{+}(\sigma_{i})$.
Or, si on note $p_{i} : X_{i} \longrightarrow F_{i}$
une section, le morphisme induit
$$(p_{i})_{+} : \underline{\bf K}^{h}_{*}(X_{i}) \longrightarrow
\underline{\bf K}^{h}_{*}(F_{i})$$
est surjectif. On \'ecrit alors $\sigma_{i}=(p_{i})_{+}(\alpha_{i})$.
Mais comme les morphismes $g_{i} : X_{i} \longrightarrow F$, et
$h_{i} : X_{i} \longrightarrow F'$ sont repr\'esentables, l'image
directe de $\sigma$, si elle existe, doit v\'erifier
$$f_{+}(\sigma)=\sum_{i}f_{+}(g_{i})_{+}(\alpha_{i})=
\sum_{i}(h_{i})_{+}(\alpha_{i})$$
Elle est donc d\'etermin\'ee uniquement par $(h_{i})_{+}$ et
$(g_{i})_{+}$. \\

Notons toujours $f : F \longrightarrow F'$ un morphisme de $\cal
LDM$, et \\
$\sigma\in \underline{\bf K}^{h}_{*}(F)$. \\

Si $F$ et $F'$ sont des gerbes lisses sur des espaces alg\'ebriques
lisses $X$ et $X'$, les morphismes naturels induits par les
projections $p : F \longrightarrow X$ et $p' : F' \longrightarrow X'$
$$p^{*} : \underline{\bf K}^{h}_{*}(X) \longrightarrow
\underline{\bf K}^{h}_{*}(F)$$
$$(p')^{*} : \underline{\bf K}^{h}(X') \longrightarrow
\underline{\bf K}^{h}(F')$$
sont des isomorphismes. Notons $m$ et $m'$ les ordres de
ramifications de $F$ et $F'$ sur $X$ et $X'$. On d\'efinit alors
l'image directe par le diagramme commutatif suivant
$$\xymatrix{
\underline{\bf K}^{h}_{*}(F) \ar[r]^{f_{+}} &
\underline{\bf K}^{h}_{*}(F') \\
\underline{\bf K}^{h}_{*}(X) \ar[u]^-{\frac{1}{m}.p^{*}}
\ar[r]_{Mf_{+}}&
\underline{\bf K}^{h}_{*}(X') \ar[u]_{\frac{1}{m'}.(p')^{*}} }$$
o\`u $Mf : X \longrightarrow X'$ est le morphisme induit sur les
espaces de modules.
Dans le cas o\`u $f$ est quelconque, on trouve des stratifications
$j'_{i} : F'_{i} \hookrightarrow F'$ et $j_{k} : F_{k} \hookrightarrow F$
par des gerbes lisses et compatible avec $f$ (~i.e. une strate est
envoy\'ee dans une strate~). Pour chaque $k$, on note $f(k)$ un entier
tel que $f(F_{k}) \hookrightarrow F'_{f(k)}$, et
$$f_{k} : F_{k} \longrightarrow F_{f(k)}$$
le morphisme induit. On pose alors
$$f_{+}(\sigma):=\sum_{k}(j'_{f(k)})_{+}(f_{k})_{+}j_{k}^{!}(\sigma)$$
Comme l'ensemble des stratifications est filtrant, on v\'erifie que
cette d\'efinition ne d\'epend pas des choix des $F'_{i}$, $F_{k}$,
$f(k)$, et d\'efinit sur $\underline{\bf K}^{h}_{*}(-)$ une
structure de foncteur covariant, compatible avec celle d\'efinie
pr\'ec\'edemment pour les morphismes repr\'esentables. $\Box$\\

\end{subsubsection}

\end{subsection}

\begin{subsection}{Les th\'eor\`emes de Grothendieck-Riemann-Roch pour
les $\cal D$-modules}

\begin{subsubsection}{Le th\'eor\`eme de Grothendieck-Riemann-Roch en\\
$\bf K\mit^{\cal D}$-th\'eorie}

\begin{thm}\label{th4.4}
Pour chaque $F$ objet de $\cal LDM$, il existe un unique morphisme
$$\tau^{rep,\cal D}_{F} : \bf K\mit^{\cal D}_{*}(F) \longrightarrow
H_{\bullet}^{rep}(F,*)_{\bf Q}$$
tel que
\begin{enumerate}
\item
si $F$ est dans $\cal QLDM$, alors
$$\begin{array}{cccc}
\tau^{rep,\cal D}_{F} : & \bf K\mit^{\cal D}_{*}(F) & \longrightarrow &
H_{\bullet}^{rep}(F,*) \\
   & x & \mapsto & Td^{rep}(F).Ch^{rep}(\delta_{F}(x))
\end{array}$$
De m\^eme, si $X$ est un sch\'ema quasi-projectif, $\tau^{rep,\cal D}_{X}$ coincide
avec le morphisme d\'efini dans \cite[$3.4$]{al}.
\item
pour tout morphisme propre de $\cal LDM$, $f : F \longrightarrow F'$, on a
$$f_{*}\circ \tau^{rep,\cal D}_{F} = \tau^{rep,\cal D}_{F'}\circ f_{+}$$
\item
si $f : F \longrightarrow F'$ est un morphisme repr\'esentable et \'etale
de champs de $\cal LDM$, alors
$$\tau^{rep,\cal D}_{F}\circ f^{!}(x) = f^{*} \circ \tau^{rep,\cal D}_{F'}(x)$$
pour tout $x \in \bf K\mit^{\cal D}_{0}(F')$.
\end{enumerate}
\end{thm}

\underline{\bf Preuve:} \rm On pose
$$\tau_{F}^{rep,\cal D}:=\tau_{F}^{rep}\circ \delta_{F}$$
et on compose le th\'eor\`eme pr\'ec\'edent avec le
th\'eor\`eme de Grothendieck-Riemann-Roch \ref{th3.5}. $\Box$\\

Remarquons que l'on peut composer ce th\'eor\`eme avec le morphisme
canonique
$$\bf K\mit^{h}_{*}(F) \longrightarrow \bf K\mit^{\cal D}_{*}(F)$$
On notera la transformation de Riemann-Roch obtenue par
$$\tau_{F}^{rep,h} : \bf K\mit^{h}_{*}(F) \longrightarrow
H^{\bullet}_{rep}(F,*)$$

\end{subsubsection}

\begin{subsubsection}{Le th\'eor\`eme de Grothendieck-Riemann-Roch en\\
$\bf K\mit^{h}$-cohomologie}

\begin{thm}\label{th4.5}
Pour chaque $F$ objet de $\cal LDM$, il existe un unique morphisme
$$\tau^{h}_{F} : \underline{\bf K}^{h}_{*}(F) \longrightarrow
H_{\bullet}(F,*)_{\bf Q}$$
tel que
\begin{enumerate}
\item
si $F$ est dans $\cal QLDM$, alors
$$\begin{array}{cccc}
\tau^{h}_{F} : & \underline{\bf K}^{h}_{*}(F) & \longrightarrow &
H_{\bullet}(F,*)_{\bf Q} \\
   & x & \mapsto & Td(F).Ch(\delta_{F}(x))
\end{array}$$
si $X$ est un espace alg\'ebrique lisse, alors $\tau^{h}_{X}$ coincide
avec le morphisme d\'efini dans \cite[$3.4$]{al}.
\item
pour tout morphisme propre de $\cal LDM$, $f : F \longrightarrow F'$, on a
$$f_{*}\circ \tau^{h}_{F} = \tau^{h}_{F'}\circ f_{+}$$
\item
si $f : F \longrightarrow F'$ est un morphisme repr\'esentable et \'etale
de champs de $\cal LDM$, alors
$$\tau^{h}_{F}\circ f^{!}(x) = f^{*} \circ \tau^{h}_{F'}(x)$$
pour tout $x \in \underline{\bf K}^{h}_{0}(F')$.
\end{enumerate}
\end{thm}

\underline{\bf Preuve:} \rm On pose $\tau_{F}^{h}:=\tau_{F}\circ
\delta_{F}$, o\`u on a encore not\'e $\delta_{F}$ le morphisme compos\'e
$$\xymatrix{\underline{\bf K}^{h}(F) \ar[r]^{nat} & \underline{\bf
K}^{\cal D}(F) \ar[r]^-{\delta_{F}} & \underline{\bf K}(F)}$$
Le seul point non trivial est le point $(2)$, dans le cas o\`u $f$ est non
repr\'esentable.

Soit $f : F \longrightarrow F'$ un morphisme propre, et
$\sigma \in \underline{\bf K}^{h}_{*}(F)$. A l'aide de
\cite{dm} et de la r\'esolution des singularit\'e, on peut trouver un
sch\'ema lisse $X$, muni d'un morphisme propre, surjectif et
g\'en\'eriquement fini\\
$p  : X \longrightarrow F$. En raisonnant par
r\'ecurrence noeth\'erienne, on peut supposer que $p$ est tel que, pour tout
sous-champ ferm\'e irr\'eductible $F_{1}~\hookrightarrow~F$, il existe une
composante connexe de $X$ au-dessus de $F_{1}$, sur laquelle $p$ est
g\'en\'eriquement fini.

Sous ces hypoth\`ese, le morphisme
$$p_{+} : \underline{\bf K}^{h}_{*}(X) \longrightarrow
\underline{\bf K}^{h}_{*}(F)$$
est surjectif. En effet, soit $j_{i} : F_{i} \hookrightarrow F$ une
stratification par des sous-champs lisses, et $k_{i} : X_{i} \hookrightarrow
X$ des sous-espace localement ferm\'es, tel que
$$p : X_{i} \longrightarrow F_{i}$$
soit \'etale et fini. Une utilisation de \ref{c4.2} nous ram\`ene donc au cas
o\`u $p$ est un morphisme \'etale et fini. Mais alors, on a
$$p_{+}p^{!}(\sigma)=m.\sigma$$
o\`u $m$ est le degr\'e de $X$ sur $F$. Ceci montre que $p_{+}$ est
surjectif. \\

Soit alors $\sigma' \in \underline{\bf K}^{h}_{*}(X)$, avec
$\sigma=p_{+}(\sigma')$. Alors, comme le point $(2)$ est vrai pour les
morphismes repr\'esentables, il est vrai pour $p$ et $f\circ p$. On a
donc
$$\begin{array}{cl}
\tau_{F'}^{h}(f_{+}(\sigma)) & = \tau_{F'}^{h}(f_{+}p_{+}(\sigma')) \\
& = f_{+}p_{+}\tau_{X}^{h}(\sigma') \\
& = f_{+}\tau_{F}^{h}(p_{+}(\sigma')) \\
& = f_{+}\tau_{F}^{h}(\sigma)
\end{array}$$
$\Box$ \\

\end{subsubsection}

\end{subsection}

% \begin{subsection}{Classes Caract\'eristiques Secondaires
% des $\cal D$-modules Holonomes R\'eguliers}
%
% \begin{subsubsection}{Classes Caract\'eristiques Secondaires}
%
% \end{subsubsection}
%
% \begin{subsubsection}{Formule d'Indice pour les $\cal D$-modules
% Holonomes R\'eguliers}
%
% \end{subsubsection}
%
% \end{subsection}

\begin{subsection}{Exemples d'application}
\hspace{5mm}

Si $p : F \longrightarrow F$ est le morphisme structural d'un champ
propre $F$ de $\cal LDM$, et $\cal M$ un $\cal D\mit_{F}$-module sur
$F$, nous noterons
$$Ind(F,\cal M\mit):=p_{+}(\cal M\mit) \in \bf K\mit^{\cal D}_{0}(Spec
k)=\bf Z$$

\underline{Remarque:} Par d\'efinition des images directes de $\cal
D$-modules, $Ind(F,\cal M\mit)$ est aussi la caract\'eristique d'Euler
du complexe de De Rham de $\cal M$ sur $F$. \\

\begin{cor}\label{c4.3}
Soit $F$ un champ propre de $\cal QLDM$, et $\cal M\mit$ un \\
$\cal D\mit_{F}$-module coh\'erent sur $F$. Notons
$[Gr\cal M\mit] \in \bf G\mit_{0}(T^{*}F)$ la classe du gradu\'e
associ\'e. Alors
$$Ind(F,\cal M\mit)=\int_{T^{*}F}^{rep}Ch^{rep}([Gr\cal
M\mit]).Ie_{*}(Td^{rep}(F))$$
o\`u $e : F \hookrightarrow T^{*}F$ est la section nulle du fibr\'e
cotangent \`a $F$.
\end{cor}

\underline{\bf Preuve:} \rm On utilise la proposition \ref{p4.4}.
$\Box$\\

Par exemple, dans le cas o\`u $\cal M\mit=\cal O\mit_{F}$ est la
connexion triviale, on trouve
$$Ind(F,\cal
O\mit_{F})=\int_{T^{*}I_{F}}(Ie_{*}(1))^{2}=\int_{I_{F}}C_{max}(T_{I_{F}})$$
On retrouve ainsi la formule de Gauss-Bonnet \ref{gb2}.

Remarquons aussi, que le corollaire pr\'ec\'edent appliqu\'e \`a une
connexion $(V,\nabla)$ sur un fibr\'e vectoriel de rang $r$ sur $F$, donne
$$Ind(F,(V,\nabla)=r.\chi^{top}(F)$$
et m\^eme, plus pr\'ecisemment
$$\tau_{F}^{h}(V,\nabla)=r.C_{max}(T_{I_{F}}) \in
H^{\bullet}_{rep}(F,*)$$

On remarquera que, contrairement \`a ce qui se passe pour un sch\'ema
(~\cite[$6.6.4$]{l}~), la formule pr\'ec\'edente ne se simplifie pas dans le cas o\`u
$\cal M$ est holonome. En effet, le fait que $Gr\cal M$ poss\`ede un
support de dimension moiti\'e dans $T^{*}F$, n'implique pas \`a priori
que $\pi_{F}^{*}(Gr\cal M\mit) \in \bf G\mit_{0}(T^{*}I_{F})$ poss\`ede
encore cette propri\'et\'e. \\

Notons que l'on a toujours un morphisme canonique
$$p^{!} : \bf K\mit_{*}(Spec k) \longrightarrow \bf K\mit_{*}^{\nabla}(F)$$
induit par le morphisme structural $p : F \longrightarrow Spec k$.
Ce morphisme induit donc un morphisme sur les
fonctions constructibles
$$F_{ct}\bf K\mit_{*}(k)(F) \longrightarrow F_{ct}\bf
K\mit_{*}^{\nabla}(F)\simeq \bf K\mit_{*}^{h}(F)$$
En composant avec la transformation de Riemann-Roch, on obtient
des "morphismes d'Euler-MacPherson"
$$E_{rep,M}^{i} : F_{ct}\bf K\mit_{i}(k)(F) \longrightarrow
H^{\bullet}_{rep}(F,*)$$
De la m\^eme fa\c{c}on on peut aussi d\'efinir
$$F_{ct}\bf K\mit_{*}(k)(F) \longrightarrow F_{ct}\underline{\bf
K}_{*}^{\nabla}(F)\simeq \underline{\bf K}_{*}^{h}(F)$$
et donc
$$E_{M}^{i} : F_{ct}\bf K\mit_{i}(k)(F) \longrightarrow
H^{\bullet}(F,*)_{\bf Q}$$
Par exemple, si $i=0$, on trouve deux morphisme
$$E_{rep,M}^{0} : F_{ct}\bf Z\mit(F) \longrightarrow H^{\bullet}_{rep}(F,*)$$
$$E_{M}^{0} : F_{ct}\bf Z\mit(F) \longrightarrow
H^{\bullet}(F,*)_{\bf Q}$$
Pour $i=1$, on a
$$E_{rep,M}^{1} : F_{ct}k^{*}(F) \longrightarrow H^{\bullet}_{rep}(F,*)$$
$$E_{M}^{1} : F_{ct}k^{*}(F) \longrightarrow
H^{\bullet}(F,*)_{\bf Q}$$

\begin{df}\label{d4.3}
Soit $F$ un champ dans $\cal LDM$.
Soit $\sigma \in F_{ct}\bf Z\mit(F)$.
Les caract\'eristiques d'Euler topologique et orbifold de $F$ pond\'er\'ees
par $\sigma$, sont d\'efinies par
$$\chi^{orb}(F,\sigma):=\sum_{i}\sigma_{i}.\chi^{orb}(F_{i})$$
$$\chi^{top}(F,\sigma):=\sum_{i}\sigma_{i}.\chi^{top}(F_{i})$$
o\`u les sous-champs $F_{i} \hookrightarrow F$ forment
une stratification de $F$, telle que $\sigma$ soit constante \'egale
\`a $\sigma_{i}$ sur $|F_{i}|$.
\end{df}

\begin{cor}\label{c4.4}
Soit $F$ un champ de $\cal LDM$.
Pour tout $\sigma \in F_{ct}\bf Z\mit(F)$, on a
$$\chi^{top}(F,\sigma)=\int_{F}^{rep}E_{rep,M}^{0}(\sigma)$$
$$\chi^{orb}(F,\sigma)=\int_{F}E_{M}^{0}(\sigma)$$
\end{cor}

\underline{\bf Preuve:} \rm Fixons nous une stratification $j_{i} : F_{i}
\hookrightarrow F$, telle que $\sigma$ soit constant \'egale \`a
$\sigma_{i} \in \bf Z$ sur $|F_{i}|$.
Notons $\bf 1\mit_{i}$ la classe de la connexion triviale sur $F_{i}$.
Alors, par d\'efinition de l'isomorphisme $\gamma$ (~\ref{th4.3}~), on a
$$E_{rep,M}^{0}(\sigma)=\sum_{i}\sigma_{i}.\tau_{F}^{rep,h}((j_{i})_{+}(\bf
1\mit_{i})) \in
H^{\bullet}(F,*)$$
Ainsi, le th\'eor\`eme de Grothendieck-Riemann-Roch \ref{th4.4}, appliqu\'e
\`a  $\sum_{i}\sigma_{i}.(j_{i})_{+}(\bf 1\mit_{i})$, donne
$$\sum_{i}\sigma_{i}.Ind(F_{i},\bf
1\mit_{i})=\int_{F}^{rep}E_{rep,M}^{0}(\sigma)$$
Ce qui est la formule que l'on voulait d\'emontrer.

Le cas de la caract\'eristique d'Euler orbifold se traite de la m\^eme
fa\c{c}on. $\Box$\\

Ce dernier corollaire explique le choix de la terminologie de "classe
d'Euler-MacPherson" pour l'application $E_{M}^{0}$. Notons qu'elle
n'est pas \'egale \`a l'application d\'efinie dans \cite{mac}. En effet,
si $F=X$ est un sch\'ema lisse et propre, on a
$E_{M}^{0}(\bf 1\mit_{X})=C_{max}(T_{X})=Eu(X)$, alors que
$C_{*}(\bf 1\mit_{X})$ est la classe de Chern total du fibr\'e tangent
$T_{X}$. De plus, nous n'avons pas su d\'emontrer que $E_{M}^{0}$
commute avec les images directes. \\

\begin{cor}\label{c4.5}
Soit $F$ un champ propre de $\cal LDM$, et
$x=can(\lambda_{-1}(\Omega^{1}_{F}))~\in~\underline{\bf K}_{0}(F)$.
Alors
$$p_{*}(x)=\chi^{orb}(F)$$
o\`u $p : F \longrightarrow Spec k$ est le morphisme structural.
\end{cor}

\underline{\bf Preuve:} \rm On applique le corollaire \ref{c4.4} \`a la
fonction constructible constante $\bf 1$. On a donc
$$\chi^{orb}(F,\bf 1\mit)=\chi^{orb}(F)=\int_{F}E_{M}(\bf 1\mit)$$
Or, $E_{M}(\bf 1\mit)=\tau_{F}^{h}(\cal O\mit_{F})$, o\`u
$\cal O\mit_{F}$ est muni de la connexion triviale. Les \'enonc\'es
\ref{th4.1} et \ref{p4.4} impliquent alors que
$$\tau_{F}^{h}(\cal O\mit_{F})=\tau_{F}(\lambda_{-1}(\Omega_{F}^{1}))$$
Ainsi, d'apr\`es \ref{th3.4}
$$\chi^{orb}(F)=\int_{F}\tau_{F}(\lambda_{-1}(\Omega_{F}^{1}))=p_{*}
(\lambda_{-1}(\Omega_{F}^{1}))$$
$\Box$ \\

On peut aussi donner une interpr\'etation de
$\int_{F}^{rep}E_{M}^{1}(\sigma)$, et de $\int_{F}E_{M}^{1}(\sigma)$,
pour $\sigma \in F_{ct}k^{*}(F)$, en terme de "determinants d'Euler".

\begin{df}
Soit $F$ un champ de $\cal LDM$, et $\sigma \in F_{ct}k^{*}(F)$. Notons
$j_{i} : F_{i} \hookrightarrow F$, une stratification de $F$, telle
que $\sigma$ soit constante \'egale \`a $\sigma_{i} \in k^{*}$ sur
$|F_{i}|$. Le determinant d'Euler de $F$ pond\'er\'e par $\sigma$ est
$$Det(F,\sigma):=\prod_{i}\sigma_{i}^{\chi^{top}(F_{i})} \in k^{*}$$
Le determinant d'Euler orbifold de $F$ pond\'er\'e par $\sigma$ est
$$Det^{orb}(F,\sigma):=\prod_{i}\sigma_{i}^{\chi^{orb}(F_{i})} \in
k^{*}\otimes_{\bf Z}\bf Q$$
\end{df}

Pour le corollaire suivant, nous supposerons que la th\'eorie
cohomologique utilis\'e est la th\'eorie de Gersten. Nous la noterons
alors avec les indices usuels des groupes de Chow
$$A^{p}(F,q):=H^{p}(F_{et},\underline{K}_{p+q})$$
$$A_{rep}^{p}(F,q):=H^{p}((I_{F})_{et},\underline{K}_{p+q})$$

\begin{cor}\label{c4.6}
Soit $F$ un champ de $\cal LDM$ et $\sigma \in F_{ct}\bf K\mit_{p}(k)(F)$.
Notons $j_{i} : F_{i} \hookrightarrow F$ une stratification de $F$,
telle que $\sigma$ soit constante \'egale \`a $\sigma_{i} \in \bf
K\mit_{p}(k)$ sur $F_{i}$.
Alors
$$\sum_{i}\chi^{top}(F_{i}).\sigma_{i}=\int_{F}^{rep}E_{rep,M}^{p}(\sigma)$$
$$\sum_{i}\chi^{orb}(F_{i}).\sigma_{i}=\int_{F}^{rep}E_{M}^{p}(\sigma)$$
En particulier
$$Det(F,\sigma)=\int_{F}^{rep}E_{rep,M}^{1}(\sigma)$$
$$Det^{orb}(F,\sigma)=\int_{F}E_{M}^{1}(\sigma)$$
\end{cor}

\underline{\bf Preuve:} \rm Commen\c{c}ons par remarquer que la
transformation de Riemann-Roch
$$\tau_{Spec k}^{h} : \bf K\mit^{h}_{p}(Spec k) \simeq
A^{0}(Spec k,p)=\bf K\mit_{p}(Spec k)$$
est un isomorphisme. \\

Comme dans la preuve du corollaire
\ref{c4.4}, on a
$$E_{M}^{p}(\sigma)=\sum_{i}\tau_{F}^{rep,h}((j_{i})_{+}(p_{i}^{!}\sigma_{i}))$$
o\`u $p_{i} : F_{i} \longrightarrow Spec k$ est le morphisme structural.
Ainsi, par le th\'eor\`eme de Grothendieck-Riemann-Roch \ref{th4.5}, et la
formule de projection pour les morphismes structuraux, on trouve
$$\begin{array}{cl}
\int_{F}^{rep}E_{M}^{p}(\sigma)&
=\sum_{i}(p_{i})_{+}p_{i}^{!}(\sigma_{i}) \\
& =\sum_{i}(p_{i})_{+}(\cal O\mit_{F_{i}}).\sigma_{i} \\
& =\sum_{i}Ind(F_{i},\cal O\mit_{F_{i}}).\sigma_{i} \\
& =\sum_{i}\chi^{top}(F_{i}).\sigma_{i} \in \bf K\mit_{p}(k)
\end{array}$$
$\Box$\\

Pour terminer cette partie, notons qu'il serait int\'eressant de
poss\'eder des formules analogues aux formules pr\'ec\'edentes, mais dans un
cas relatif. On disposerait ainsi de "classes d'Euler-MacPherson pour
la $K$-th\'eorie sup\'erieure", bien que la signification de telles
formules pour les $\bf K\mit_{p}$ avec $p>1$ me semble un peu
myst\'erieuses.

\end{subsection}

\end{section}

\newpage
\begin{section}{Chapitre $5$ : Champs alg\'ebriques et champs
analytiques}
\hspace{5mm}
Ce dernier chapitre est ind\'ependant des pr\'ec\'edents. Nous nous
int\'eresserons au probl\`eme de comparaison entre champs analytiques et
champs alg\'ebriques complexes. Notre but est d'essayer de g\'en\'eraliser
les r\'esultats d'alg\'ebrisation d'Artin (~\cite[$7.3$]{a}~), au cas des champs
analytiques.

Pour cela nous commencerons par d\'emontrer les th\'eor\`emes GAGA pour des
champs alg\'ebriques complexes de Deligne-Mumford. Ils seront utilis\'es
par la suite pour d\'emontrer un cas particulier du th\'eor\`eme d'Artin. \\

Dans cette section nous appellerons sch\'ema un sch\'ema localement de type
fini et s\'epar\'e sur $Spec \bf C$. Nous noterons $(Sch/\bf C\mit)_{et}$ le
site des
sch\'emas, muni de la
topologie \'etale. Le petit site \'etale d'un sch\'ema $X$ sera not\'e
$X_{et}$. Un champ alg\'ebrique sera un champ de Deligne-Mumford,
localement de type fini et s\'epar\'e sur $\bf C\mit$ (~\ref{d1.8}~). Si
$X$ est un sch\'ema, $X^{an}$ sera l'espace analytique associ\'e.\\

Un espace analytique, est un espace analytique complexe, localement
de type fini et s\'epar\'e sur $\bf C$. L'espace topologique sous jacents \`a
un espace analytique $X$ sera not\'e $X^{top}$.

\begin{subsection}{Analytification des champs alg\'ebriques}

\begin{subsubsection}{Champs analytiques}
\hspace*{5mm}
Rappelons pour commencer les d\'efinitions suivantes.

\begin{df}\label{d5.1}
Un morphisme d'espaces analytiques
$$f : X \longrightarrow Y$$
est
\begin{itemize}
\item \'etale, si pour chaque point $x \in X$, le morphisme d'anneaux
locaux
$$f_{x}^{*} : \cal O\mit_{Y,f(x)} \longrightarrow \cal O\mit_{X,x}$$
est un isomorphisme.
\item non ramifi\'e, si localement sur $X^{top}$, c'est une immersion ferm\'ee.
\item fini, s'il est propre et \`a fibres finies.
\end{itemize}
\end{df}

Avec ces d\'efinitions, le site analytique est d\'efini par la cat\'egorie
des espaces analytiques, dont les morphismes couvrant sont les
morphismes \'etales et surjectifs. Il sera not\'e $(An/\bf C\mit)_{et}$.
Sur ce site on dispose alors de la notion de champs (~\cite{lm}~),
ainsi que de morphismes de champs repr\'esentables poss\'edant certaines
propri\'et\'es locales, comme \'etales, finis, non ramifi\'es, lisses ...
(~\ref{d1.5}~).
Un champ \'equivalent \`a un espace analytique, sera encore appel\'e un
espace analytique. \\

\begin{df}\label{d5.2}
Un champ analytique (~de Deligne-Mumford et s\'epar\'e~)
est un champ $F$ sur le site $(An/\bf C\mit)_{et}$
tel que
\begin{enumerate}
\item le morphisme diagonal
$$\Delta : F \longrightarrow F\times F$$
est repr\'esentable et fini.
\item il existe un espace analytique $X$
et un morphisme \'etale et surjectif
$$X \longrightarrow F$$
\end{enumerate}
\end{df}

Nous noterons $Ch^{an}(\bf C\mit)$
(~resp. $ChAn(\bf C\mit)$~) la $2$-cat\'egorie
des champs (~resp. champs analytiques~) sur $(An/\bf C\mit)_{et}$.
Les cat\'egories homotopiques associ\'ees, seront not\'ees respectivement
$HoCh^{an}(\bf C\mit)$ et $HoChAn(\bf C\mit)$.

\underline{Remarque:} Comme dans le cas des sch\'emas (~\ref{d1.8}~),
le morphisme diagonal d'un champ analytique est automatiquement non ramifi\'e.\\

L'exemple standard de champ analytique est le champ quotient par un
groupe fini. Pour un action d'un groupe fini $H$ sur un espace
analytique $X$, nous noterons $[X/H]$ le champ quotient. Son groupoide
des sections au-dessus d'un espace $Y$, est le groupoide des
couples $(P,f)$, o\`u $P \longrightarrow Y$ est un $H$-torseur, et
$f : P \longrightarrow X$ est un morphisme $H$-\'equivariant.\\

\begin{df}\label{d5.3}
Le champ des ramifications d'un champ analytique $F$, est d\'efini
par
$$I_{F}:=F\times_{F\times F}F$$
Si la projection naturelle
$$\pi_{F} : I_{F} \longrightarrow F$$
est \'etale, on dit que $F$ est une gerbe.
\end{df}

Les champs analytiques poss\`edent de nombreuses propri\'et\'es analogues \`a
celles des champs alg\'ebriques. Nous ne retiendrons que les suivantes.

\begin{prop}\label{p5.1}
Soit $F$ un champ analytique.
\begin{enumerate}
\item Il existe un espace analytique $M$, et un morphisme propre
$$p : F \longrightarrow M$$
qui est universel vers les espaces analytiques, et qui induit une
bijection
$$\pi_{0}F(Spec \bf C\mit^{an}) \simeq M(Spec \bf C\mit^{an})$$
On dira que $M$ est un espace de modules pour $F$.
\item Soit $M$ l'espace de modules de $F$. Alors, il existe un recouvrement
\'etale $U \longrightarrow M$, un
espace analytique $X$, un faisceau en groupes finis $H$ constant sur
$X$, et une op\'eration de $H$ sur $X$, tel que
le champ
$$F_{U}:=F\times_{M}U$$
soit \'equivalent au champ classifiant $[X/H]$.
\item Si $F$ est r\'eduit, il existe un sous-champ analytique
ferm\'e $F_{1} \hookrightarrow F$, tel que $F-F_{1}$ soit une gerbe.
\end{enumerate}
\end{prop}

Un champ analytique sera propre si son espace de
modules est un espace analytique
compact.\\

Notons pour finir que, comme dans le cadre alg\'ebrique (~\ref{d2.1}~), il y
a une
notion de faisceaux analytiques coh\'erents sur un champ analytique $F$. Ce
sont les sections globales cart\'esiennes sur $F$, du champ en cat\'egories
ab\'eliennes $\bf Coh\mit$ sur $(An/\bf C\mit)_{et}$, dont les sections
au-dessus d'un espace analytique $X$, sont les faisceaux
de $\cal O\mit_{X}^{an}$-modules coh\'erents.

\end{subsubsection}

\begin{subsubsection}{Analytification}

Soit $f : C \longrightarrow C'$ un foncteur entre deux sites. On
suppose que $f$ est continu, dans le sens o\`u
$$f(U) \longrightarrow F(X)$$
est couvrant dans $C'$, si $U \longrightarrow X$ est couvrant dans $C$.

Soit $p : \cal C\mit \longrightarrow C$ un champ. On dispose alors
d'un champ "image r\'eciproque" $f^{*}(\cal C\mit)$ (~\cite[$3.2$]{gi}~).

\begin{lem}\label{l5.1}
Soit
$$\alpha : (Esp/\bf C\mit)_{et} \longrightarrow (An/\bf C\mit)_{et}$$
le foncteur continu qui \`a un sch\'ema $X$ associe son espace analytique
$X^{an}$. Alors l'image par $\alpha^{*}$ d'un
champ alg\'ebrique est un champ analytique.
\end{lem}

\underline{\bf Preuve:} \rm En effet, si $F$ est un champ alg\'ebrique
\'equivalent \`a un pr\'efaisceau simplicial repr\'esent\'e par un espace
alg\'ebrique en groupoides $X_{\bullet}$, alors le champ $F^{an}$ est
\'equivalent au pr\'efaisceau repr\'esent\'e par l'objet en groupoides
$X_{\bullet}^{an}$, qui, par d\'efinition, est un champ analytique. $\Box$\\

Comme le foncteur $\alpha^{*}$ pr\'eserve les \'equivalences de champs, il
passe aux cat\'egories homotopiques.

\begin{df}\label{d5.5}
On notera
$$\begin{array}{ccc}
HoChAlg(\bf C\mit) & \longrightarrow & HoChAn(\bf C\mit) \\
F & \mapsto & F^{an}:=\alpha^{*}(F)
\end{array}$$
le foncteur d'analytification.
\end{df}

Comme on sait que les champs de Deligne-Mumford sont localement des
quotients par des groupes finis (~\ref{p5.1}, \ref{p1.2}~), il est naturel de
commencer par comparer les quotients alg\'ebriques, et les quotients
analytiques.

\begin{prop}\label{p5.3}
Soit $X$ un sch\'ema et $H$ un groupe fini op\'erant sur $X$. Alors il
existe une \'equivalence faible naturelle
$$[X/H]^{an} \simeq [X^{an}/H]$$
\end{prop}

\underline{\bf Preuve:} \rm Soit $X \longrightarrow [X/H]$ le
$H$-torseur universel. Comme le foncteur d'analytification transforme
$H$-torseurs en $H$-torseurs (~\cite{sga1}~), le morphisme induit
$$X^{an} \longrightarrow [X/H]^{an}$$
est un $H$-torseur. Il d\'efini donc une \'equivalence faible
$$[X/H]^{an} \longrightarrow [X^{an}/H]$$
$\Box$\\

\begin{cor}\label{c5.1}
Soit $F$ un champ alg\'ebrique et $p : F \longrightarrow M$ la
projection sur son espace de
modules (~\cite{km}~). Alors le morphisme induit
$$p : F^{an} \longrightarrow M^{an}$$
fait de $M^{an}$ l'espace de modules de $F^{an}$.
\end{cor}

\underline{\bf Preuve:} \rm Soit $X$ l'espace de modules de $F^{an}$.
Par la propri\'et\'e universelle des espaces de modules, il existe un
morphisme naturel
$$X \longrightarrow M^{an}$$
Pour montrer que c'est un isomorphisme, on peut effectuer un
changement de base de $M^{an}$ par un morphisme $u^{an} : (U \longrightarrow
M)^{an}$, avec $u$ \'etale et surjectif. Par \ref{p1.2}, on peut donc
supposer que $F=[X/H]$, avec $X$ un sch\'ema affine, et $H$ un groupe
fini op\'erant sur $X$. Mais alors, comme $F^{an}\simeq [X^{an}/H]$,
l'assertion \`a d\'emontrer est que le morphisme naturel
$$X^{an}/H \longrightarrow (X/H)^{an}$$
est un isomorphisme. Ce qui est vrai. $\Box$\\

Soit $\alpha_{*}\bf Coh\mit$ l'image r\'eciproque du champ des faisceaux
analytiques coh\'erents par le foncteur d'analytification. C'est le
champ sur $(Esp/\bf C\mit)_{et}$, dont la cat\'egorie des sections au-dessus d'un
espace alg\'ebrique $X$, est la cat\'egorie $\bf Coh\mit(X^{an})$ des
faisceaux analytiques coh\'erents sur $X^{an}$. Alors, on dispose d'un
foncteur naturel d'analytification (~\cite[$XII$]{sga1}~)
$$\begin{array}{ccc}
\bf Coh\mit(X) & \longrightarrow & \bf Coh\mit(X^{an})\\
\cal F\mit & \mapsto & \cal F\mit^{an}
\end{array}$$
Ce foncteur d\'efini un morphisme de champs sur $(Esp/\bf
C\mit)_{et}$
$$\bf Coh\mit \longrightarrow \alpha_{*}\bf Coh\mit$$
ou bien, par adjonction, un morphisme de champs sur $(An/\bf
C\mit)_{et}$
$$-^{an} : \alpha^{*}\bf Coh \mit \longrightarrow \bf Coh\mit$$
Par composition, ce foncteur d\'efini un foncteur sur les sections
cart\'esiennes
\hspace{-20mm}
$$\begin{array}{ccccc}
Hom_{Cart}(F,\bf Coh\mit) & \longrightarrow &
Hom_{Cart}(F^{an},\alpha^{*}\bf Coh\mit) & \longrightarrow &
Hom_{Cart}(F^{an},\bf Coh\mit)\\
\cal F\mit & \mapsto & \alpha^{*}\cal F\mit & \mapsto &
\cal F\mit^{an}
\end{array}$$
Et donc un foncteur
$$\begin{array}{ccc}
\bf Coh\mit(F):=Hom_{Cart}(F,\bf Coh\mit) & \longrightarrow &
\bf Coh\mit(F^{an}):=Hom_{Cart}(F^{an},\bf Coh\mit) \\
\cal F\mit & \mapsto & \cal F\mit^{an}
\end{array}$$

\begin{df}\label{d5.6}
Soit $F$ un champ alg\'ebrique, le foncteur d'analytification ci-dessus
est not\'e
$$\begin{array}{ccc}
\bf Coh\mit(F) & \longrightarrow & \bf Coh\mit(F^{an})\\
\cal F\mit & \mapsto & \cal F\mit^{an}
\end{array}$$
\end{df}

Notons, que $\cal F\mit \mapsto \cal F\mit^{an}$ est un foncteur
exact. En effet, comme ceci est local sur $F_{et}$, on se ram\`ene au
cas des sch\'emas \cite[$XII$ $4.4$]{sga1}.

\end{subsubsection}

\end{subsection}

\begin{subsection}{Th\'eor\`emes GAGA}
\hspace*{5mm}
Le th\'eor\`eme principal est le suivant.

\begin{thm}{(~"GAGA"~)}\label{GAGA}
Soit $F$ un champ alg\'ebrique propre. Alors le foncteur
$$\begin{array}{ccc}
\bf Coh\mit(F) & \longrightarrow & \bf Coh\mit(F^{an}) \\
\cal F\mit & \mapsto & \cal F\mit^{an}
\end{array}$$
est une \'equivalence de cat\'egories.
\end{thm}

\underline{\bf Preuve:} \rm

\begin{lem}\label{l5.2}
Soit $f : F \longrightarrow F'$ un morphisme propre et repr\'esentable de
champs alg\'ebriques, et $\cal F\mit$ un faisceau coh\'erent sur $F$.
Alors le morphisme naturel
$$R^{i}f^{an}_{*}(\cal F\mit^{an}) \longrightarrow R^{i}f_{*}(\cal
F\mit)^{an}$$
est un isomorphisme.
\end{lem}

\underline{\bf Preuve:} \rm Comme ceci est local sur $(F')^{an}_{et}$,
le lemme est une cons\'equence directe de \cite[$XII$ $4.2$]{sga1}. $\Box$\\

\begin{lem}\label{l5.3}
Soit $p : F \longrightarrow M$ le projection d'un champ alg\'ebrique sur
son espace de modules, et $\cal F$ un faisceau coh\'erent sur $F$. Alors le
morphisme naturel
$$p^{an}_{*}(\cal F\mit^{an}) \longrightarrow p_{*}(\cal F\mit)^{an}$$
est un isomorphisme.
\end{lem}

\underline{\bf Preuve:} \rm Soit $U \longrightarrow M$ un recouvrement
\'etale alg\'ebrique tel que $F_{U}$ soit un champ quotient sur $U$. En
localisant sur $U$, on peut donc supposer que $F=[X/H]$, avec
$H$ un groupe fini op\'erant sur un sch\'ema $X$. Le faisceau $\cal F\mit$
est alors donn\'e par un faisceau coh\'erent $\cal F\mit_{X}$ sur $X$, muni d'une
action de $H$.

Soit $q : X \longrightarrow X/H$ la projection. D'apr\`es le
lemme pr\'ec\'edent on a alors
$$q_{*}^{an}(\cal F\mit_{X}^{an}) \simeq q_{*}(\cal F\mit_{X})^{an}$$
De plus cet isomorphisme est un isomorphisme \'equivariant de
$H$-faisceaux analytiques coh\'erents sur $X^{an}/H$. En prenant les
invariants sous $H$ on obtient
$$p_{*}^{an}(\cal F\mit_{X}^{an}) \simeq
q_{*}^{an}(\cal F\mit_{X}^{an})^{H} \simeq
(q_{*}(\cal F\mit_{X})^{an})^{H} \simeq
(q_{*}(\cal F\mit_{X})^{H})^{an} \simeq p_{*}(\cal F\mit_{X})^{an}$$
$\Box$\\

Revenons \`a la preuve du th\'eor\`eme. Elle suit exactement le m\^eme sch\'ema
que la preuve donn\'ee dans \cite[$XII$ $4.4$]{sga1}.\\

\underline{$(1)$ Le foncteur est pleinement fid\`ele:}\\

Soit $\cal F$ un faisceau coh\'erent sur $F$. Alors, comme
$p_{*}$ est un foncteur exact, on a un isomorphisme canonique
$$H^{i}(F,\cal F\mit) \simeq H^{i}(M,p_{*}\cal F\mit)$$
De plus, $M$ est un espace alg\'ebrique propre, et $p_{*}\cal F\mit$ un
faisceau coh\'erent sur $M$, donc le th\'eor\`eme GAGA pour les espaces
alg\'ebriques (~\cite[$XII$ $4.4$]{sga1}~) nous dit que le morphisme naturel
$$H^{i}(M,p_{*}\cal F\mit) \longrightarrow
H^{i}(M^{an},(p_{*}\cal F\mit)^{an})$$
est un isomorphisme. De plus le lemme \ref{l5.3} implique que
$$H^{i}(M^{an},(p_{*}\cal F\mit)^{an}) \simeq
H^{i}(M^{an},p^{an}_{*}(\cal F\mit^{an}))$$
Or
$$H^{i}(M^{an},p^{an}_{*}(\cal F\mit^{an})) \simeq H^{i}(F^{an},\cal
F\mit^{an})$$
On obtient donc l'isomorphisme cherch\'e
$$H^{i}(F,\cal F\mit) \simeq H^{i}(F^{an},\cal F\mit^{an})$$
Si $\cal F$ et $\cal G$ sont deux faisceaux coh\'erents sur $F$, on
dispose du faisceau coh\'erent $\underline{Hom}_{\cal O\mit_{F}}(\cal
F\mit,\cal G\mit)$ sur $F$. Nous pouvons donc lui appliquer l'isomorphisme
pr\'ec\'edent avec $i=0$
$$Hom_{\cal O\mit_{F}}(\cal F\mit, \cal G\mit) \simeq
H^{0}(F,\underline{Hom}_{\cal O\mit_{F}}(\cal F\mit,\cal G\mit)) \simeq
H^{0}(F^{an},\underline{Hom}_{\cal O\mit_{F}}(\cal
F\mit,\cal G\mit)^{an})$$
Or
$$H^{0}(F^{an},\underline{Hom}_{\cal O\mit_{F}}(\cal
F\mit,\cal G\mit)^{an}) \simeq
Hom_{\cal O\mit_{F^{an}}}(\cal F\mit^{an},\cal G\mit^{an})$$
Ceci ach\`eve la preuve de l'assertion $(1)$.\\

\underline{$(2)$ Le foncteur est essentiellement surjectif:}\\

Commen\c{c}ons par le cas o\`u $F=X\times BH$ est une gerbe triviale de
groupe fini $H$ et d'espace de modules $X$ un sch\'ema propre. Alors un
faisceau coh\'erent $\cal F$ sur $F^{an}$ est donn\'e par un faisceau coh\'erent
$\cal F\mit_{X}$ sur $X^{an}$ muni d'une action du groupe $H$. D'apr\`es le
th\'eor\`eme GAGA pour $X$, on sait qu'il existe un faisceau coh\'erent
$\cal M\mit_{X}$ sur $X$ tel que $\cal F\mit_{X}\simeq \cal M\mit_{X}^{an}$.
L'action de $H$ est alors donn\'ee par une repr\'esentation
$$H \longrightarrow Aut_{\cal O\mit_{X^{an}}}(\cal M\mit_{X}^{an})$$
Or, une seconde application du th\'eor\`eme GAGA implique que
$$Aut_{\cal O\mit_{X^{an}}}(\cal M\mit_{X}^{an})\simeq Aut_{\cal
O\mit_{X}}(\cal M\mit_{X})$$
On munit ainsi $\cal M\mit_{X}$ d'une action de $H$, ce qui d\'efinit le
faisceau coh\'erent $\cal M$ sur $F$ tel que
$\cal M\mit^{an}\simeq \cal F\mit$.\\

Passons au cas g\'en\'eral. Soit $\cal F$ un faisceau analytique coh\'erent
sur $F^{an}$. On raisonne par r\'ecurrence sur la dimension $d$ du
support de $\cal F$.\\

Soit $\cal A$ l'id\'eal annulateur de $\cal F$ dans $\cal O\mit_{F^{an}}$.
En se restreignant au ferm\'e d\'efini par $\cal A$, on peut supposer que
le support de $\cal F$ est $F$.\\

Soit $\cal I$ l'id\'eal de $F_{red}$ dans $F$, et $k$ un entier tel que
$\cal I\mit^{k}=0$. On dispose de la filtration suivante
$$0 \hookrightarrow \cal I\mit^{k-1}.\cal F\mit \hookrightarrow \dots
\cal I\mit.\cal F\mit \hookrightarrow \cal F\mit$$
dont les quotients successifs sont des images directes de modules
coh\'erents sur $F_{red}$ par l'immersion canonique $F_{red}
\hookrightarrow F$. Remarquons alors qu'une extension analytique de
faisceaux coh\'erents alg\'ebrisables est alg\'ebrisable. En effet, le fait
que le foncteur d'analytification $\cal G\mit \mapsto \cal G\mit^{an}$
soit exact et pleinement fid\`ele, et la formule
$$\underline{Hom}_{\cal O\mit_{F}}(\cal F\mit, \cal G\mit)^{an}\simeq
\underline{Hom}_{\cal O\mit_{F^{an}}}(\cal F\mit^{an},\cal
G\mit^{an})$$
montre que
$$Ext^{1}_{\cal O\mit_{F^{an}}}(\cal F\mit^{an},\cal
G\mit^{an}) \simeq Ext^{1}_{\cal O\mit_{F}}(\cal F\mit,\cal
G\mit)$$
On peut donc se restreindre
au cas o\`u $F$ est r\'eduit. \\

D'apr\`es \cite[Thm. $4.12$]{dm}, on peut trouver un
sch\'ema projectif $X$ normal, et un morphisme propre et g\'en\'eriquement \'etale
$$X \longrightarrow F$$
Soit $F_{X}=F\times_{M}X$ le champ induit sur $X$, et $F_{0}$ la
normalisation de $F_{X}$. Alors, d'apr\`es \ref{l1.1}, on sait que
$F_{0}$ est une gerbe triviale sur $X$. Notons $q : F_{0}
\longrightarrow F$ la projection. Alors $q^{*}\cal F$ est un faisceau
analytique coh\'erent sur $F_{0}^{an}$, qui est alg\'ebrisable d'apr\`es la
premi\`ere partie. Ecrivons
$$\cal M\mit^{an} \simeq q^{*}\cal F$$
avec $\cal M$ un faisceau coh\'erent sur $F_{0}$. On consid\`ere le
faisceau coh\'erent sur $F^{an}$
$q_{*}(\cal M\mit^{an})$. Il est alg\'ebrisable d'apr\`es le lemme
\ref{l5.2}. De plus, par la formule de la projection, on a
$$q_{*}q^{*}\cal F \simeq \cal F\mit\otimes_{\cal O\mit_{F^{an}}}
q_{*}(\cal O\mit_{F_{0}^{an}})$$
Comme $q$ est g\'en\'eriquement un changement de base par un rev\^etement
\'etale de l'espace de modules, on sait que
$q_{*}\cal O\mit_{F_{0}}$ est g\'en\'eriquement isomorphe \`a $\cal
O\mit_{F}^{m}$. Cet isomorphisme g\'en\'erique donne lieu \`a un diagramme
de faisceaux coh\'erents sur $F$
$$\xymatrix{
\cal N\mit \ar[r]^{u} \ar[d]_{v} & q_{*}\cal O\mit_{F_{0}} \\
\cal O\mit_{F}^{m} }$$
o\`u $u$ et $v$ sont des isomorphismes g\'en\'eriques. Consid\'erons le
diagramme obtenu sur $F^{an}$ en tensorisant par $\cal F$, et en
compl\'etant avec les noyaux et conoyaux
$$\xymatrix{
              &           & 0 \ar[d] &           &             & \\
              &           & \cal K\mit_{1} \ar[d] &            & \\
 0 \ar[r] & \cal K\mit_{2} \ar[r] & \cal N\mit \otimes \cal F\mit
 \ar[d]_{v} \ar[r]^{u} & q_{*}(\cal O\mit_{F_{0}^{an}})\otimes \cal F\mit
\ar[r]
 & \cal C\mit_{2} \ar[r] & 0 \\
              &           & \cal F\mit^{m} \ar[d] &           & \\
              &           & \cal C\mit_{1} \ar[d] & & \\
              &           & 0 & & }$$
Comme $u$ et $v$ sont des isomorphismes sur un ouvert Zariski dense, on conclut
que $\cal K\mit_{i}$ et $\cal C\mit_{i}$ ont un support de dimension
strictement plus petit que $d$, et sont donc alg\'ebrisables
par r\'ecurrence. De plus, on a
vu que
$$q_{*}(\cal O\mit_{F_{0}^{an}})\otimes \cal F\mit \simeq (q_{*}\cal
M\mit)^{an}$$
Ainsi, on en d\'eduit que $\cal N\mit \otimes \cal F\mit$ est
alg\'ebrisable comme extension de faisceaux coh\'erents alg\'ebrisables. La
colonne verticale nous dit alors que $\cal F\mit^{m}$ est
alg\'ebrisable. On termine en remarquant par exemple, que $\cal F$ est le
noyau du
morphisme $\cal F\mit^{m} \longrightarrow \cal F\mit^{m}$, qui
d\'eplace les facteurs d'un cran vers la droite. $\Box$\\

On d\'eduit de ce th\'eor\`eme, les corollaires habituels.

\begin{cor}\label{c5.2}
Soit $F$ un champ alg\'ebrique propre. Alors le foncteur
$$G \mapsto G^{an}$$
induit une bijection entre les sous-champs alg\'ebriques ferm\'es de $F$, et
les sous-champs analytiques ferm\'es de $F^{an}$
\end{cor}

\underline{\bf Preuve:} \rm Le foncteur d'analytification induit une
bijection entre les faisceaux coh\'erents d'id\'eaux de $\cal O\mit_{F}$, et les
faisceaux coh\'erents d'id\'eaux de $\cal O\mit_{F^{an}}$. $\Box$\\

\begin{cor}\label{c5.3}
Soit $F$ un champ alg\'ebrique propre. Alors le foncteur
$$G \mapsto G^{an}$$
induit une \'equivalence de la cat\'egorie homotopique des
champs alg\'ebriques finis et repr\'esentables sur $F$, et celle des
champs analytiques finis et repr\'esentables sur $F^{an}$.
\end{cor}

\underline{\bf Preuve:} \rm En effet, le foncteur en question induit
une \'equivalence entre la cat\'egorie des faisceaux en $\cal O\mit_{F}$-alg\`ebres
coh\'erentes et celle des faisceaux en $\cal O\mit_{F^{an}}$-alg\`ebres
coh\'erentes. $\Box$\\

\begin{cor}\label{c5.4}
Soit $F$ et $G$ deux champs alg\'ebriques propres. Alors le foncteur
d'analytification induit un foncteur pleinement fid\`ele de la cat\'egorie
homotopique des champs alg\'ebriques propres, vers celle des champs
analytiques.
\end{cor}

\underline{\bf Preuve:} \rm Il faut remarquer qu'un morphisme
$$f : F^{an} \longrightarrow G^{an}$$
est d\'etermin\'e \`a homotopie pr\`es par son graphe
$$\gamma_{f} : F^{an} \longrightarrow F^{an}\times G^{an}$$
qui est un morphisme fini, et appliquer le corollaire pr\'ec\'edent.
$\Box$\\

\end{subsection}

\begin{subsection}{Alg\'ebrisation des champs analytiques}
\hspace*{5mm}
Dans cette derni\`ere section nous analysons le cas des champs
analytiques dont les espaces de modules sont des espaces alg\'ebriques.\\

\begin{prop}\label{p5.4}
Soit $X$ un espace alg\'ebrique, et $H$ un groupe fini. Alors le
foncteur
$$F \mapsto F^{an}$$
induit une \'equivalence entre le \\
$2$-groupoide des gerbes alg\'ebriques
born\'ees par le  groupe $H$ sur $X$, et celle des gerbes analytiques
born\'ees par le groupe $H$ sur $X^{an}$.
\end{prop}

\underline{\bf Preuve:} \rm Remarquons d'abord qu'il est clair que le
foncteur d'analytification transforme gerbes en gerbes. \\

Soit $BH$ l'ensemble simplicial classifiant de $H$, et
$Aut(BH)$ l'ensemble simplicial des auto-\'equivalences de $BH$. Ce
dernier poss\`ede un ensemble simplicial classifiant $BAut(BH)$, qui est
un ensemble simplicial connexe et $2$-tronqu\'e. Ainsi, \`a travers
l'\'equivalence d\'emontr\'ee dans \cite{tan}, nous le verrons comme un
$2$-groupoide. Ces groupes
d'homotopie sont donn\'es par
$$\pi_{1}(BAut(BH)) \simeq Out(H)$$
$$\pi_{2}(BAut(BH)) \simeq Z(G)$$
o\`u $Out(H)=Aut(H)/Int(H)$ est le groupe des automorphismes ext\'erieurs
de $H$, et $Z(H)$ le centre de $H$. Nous noterons $\cal G$ le
$2$-champ associ\'e au pr\'efaisceau simplicial constant sur $(Sch/\bf
C\mit)_{et}$, de fibre $BAut(BH)$. Alors $\cal G\mit^{an}$ est canoniquement
\'equivalent au champ sur $(An/\bf C\mit)_{et}$, associ\'e au pr\'efaisceau
simplicial
constant $BAut(BH)$.

Enfin nous savons d'apr\`es \cite{s2}, qu'il existe une \'equivalence
entre le \\
$2$-groupoide des gerbes alg\'ebriques de groupes $H$ sur $X$
(~resp. analytiques de groupes $H$ sur $X^{an}$~) et
$\cal G\mit(X)$ (~resp. $\cal G\mit^{an}(X^{an})$~). Pour d\'emontrer
la proposition, il nous suffit donc de montrer que le morphisme
d'analytification induit une \'equivalence faible
$$\cal G\mit(X) \simeq \cal G\mit^{an}(X^{an})$$
Mais ceci provient du lemme g\'en\'eral suivant.

\begin{lem}\label{l5.4}
Soit $F$ un pr\'efaisceau simplicial $n$-tronqu\'e sur $(Sch/\bf C\mit)_{et}$,
tel que pour chaque
section $s : X \longrightarrow F$, les faisceaux d'homotopie
$\pi_{m}(F,s)$ soient localement constants \`a fibres finies sur
$X_{et}$. Alors pour chaque espace alg\'ebrique $X$,
le morphisme naturel
$$F(X) \longrightarrow F^{an}(X^{an})$$
est une \'equivalence faible.
\end{lem}

\underline{\bf Preuve:} \rm Il suffit de regarder, pour chaque
section $s \in F(X)$, le morphisme induit sur les suites spectrales
$$E_{2}^{p,q}=H^{q}(X_{et},\pi_{p}(F,s)) \longrightarrow
(E')_{2}^{p,q}=H^{q}(X^{an}_{top},\pi_{p}(F^{an},s))$$
qui convergent vers $\pi_{*}(F(X),s)$ et $\pi_{*}(F^{an}(X^{an}),s)$.
Le th\'eor\`eme de comparaison entre cohomologie \'etale et cohomologie
transcendante pour des groupes finis (~\cite{sga1}~), permet de conclure que le
morphisme est un isomorphisme sur les termes $E_{2}$. $\Box$\\

\begin{df}
Un $1$-morphisme propre et repr\'esentable de champs analytiques
$$f : F_{0} \longrightarrow F$$
est une modification, s'il existe un sous-champ ferm\'e $F'
\hookrightarrow F$, tel que $f$ induise une \'equivalence
$$f : F_{0}-f^{-1}(F') \longrightarrow F-F'$$
\end{df}

Le r\'esultat suivant est un analogue du r\'esultat de Moisezon
\cite[$7.16$]{a}.

\begin{thm}\label{th5.1}
Soit $F$ un champ analytique r\'eduit tel que son espace de modules soit
alg\'ebrique et propre. Alors il existe une modification
$$f : F_{0} \longrightarrow F$$
avec $F_{0}$ \'equivalent \`a un champ alg\'ebrique.
\end{thm}

\underline{\bf Preuve:} \rm Comme $M$ est alg\'ebrisable, nous \'ecrirons
encore $M$ pour l'espace alg\'ebrique correspondant.\\

En consid\'erant la r\'eunion disjointe des
composantes irr\'eductibles de $F$, on peut supposer que $F$ est
irr\'eductible, et donc int\`egre.\\

Soit $M$ l'espace de modules de $F$, et
$S \hookrightarrow M$ un sous-espace ferm\'e tel que $F$ soit une gerbe sur
$M-S$.

Commen\c{c}ons par fixer quelques d\'efinitions. Si $H$ est un groupe
fini, un rev\^etement (~resp. rev\^etement analytique~) d'un sch\'ema $X$
de groupe $H$ non-ramifi\'e en dehors de $S$, est la donn\'ee d'une
action de $H$ sur un sch\'ema $Y$ (~resp. sur un espace analytique
$Y$~) et d'un morphisme \'equivariant
$Y \longrightarrow M$ (~resp. $Y \longrightarrow M^{an}$~) \'etale sur
$M-S$, et qui fasse de $M$ le quotient de $Y$ par $H$ (~resp. \'etale
sur $(M-S)^{an}$, et qui fasse de $M^{an}$ le quotient de $Y$ par
$H$~).

Si $x \in M^{an}$, nous parlerons de germe analytique de rev\^etements
de groupe $H$ de $M$ en $x$ et non-ramifi\'e en dehors de $S$, pour
d\'esigner une classe de tels rev\^etements, o\`u deux rev\^etements sont
\'equivalents s'ils sont isomorphes sur un voisinage analytique
de $x$. Pour simplifier le vocabulaire nous dirons simplement
"rev\^etement de groupe $H$" pour rev\^etement de groupe $H$ et
non-ramifi\'e en dehors de $S$. De m\^eme, nous parlerons de germes \'etales
de rev\^etements de groupe $H$ de $M$ en un point ferm\'e $x \in M$.

Le foncteur d'analytification transforme un germe \'etale de rev\^etements de
groupe $H$, en un germe analytique de rev\^etements de groupe $H$. Nous
dirons que l'image d'un germe \'etale est un germe analytique
alg\'ebrisable.\\

\underline{Etape $(1)$:} \\

Commen\c{c}ons par supposer que $F$ et $M$ sont normaux,
et que $M$ v\'erifie l'hypoth\`ese suivante.

\begin{hyp}
Pour tout point $x \in M^{an}$, tout germe analytique de rev\^etement
de groupe $H$ est alg\'ebrisable.
\end{hyp}

Soit $x \in M$, et $x \in U \hookrightarrow M$ un voisinage
analytique de $x$ tel que $F\times_{M}U=:F_{U}$ soit \'equivalent \`a un champ
quotient par un groupe fini
$$F_{U} \simeq [V/H]$$
On peut donc, d'apr\`es l'hypoth\`ese, supposer qu'il existe un
sch\'ema $X$ et un morphisme \'etale au voisinage de $x$
$$X \longrightarrow M$$
une action de $H$ sur un sch\'ema $Y$, et un morphisme
$$p : Y \longrightarrow X$$
qui fasse de $X$ le quotient de $Y$ par $H$, et tel que dans un voisinage
analytique de $x$, le morphisme $p$ soit isomorphe \`a la projection canonique
$$q : V \longrightarrow U$$
Compactifions le morphisme $Y \longrightarrow M$
$$\xymatrix{
Y \ar[r]^{j} \ar[d] & Z \ar[ld]^{p} \\
M & }$$
o\`u $j$ est une immersion ouverte d'image Zariski dense, et $p$ propre.
Consid\'erons le changement de base sur $Z$,
$F_{Z}:=F\times_{M}Z \longrightarrow Z^{an}$, ainsi que $F_{Z}^{\circ}
\longrightarrow
(F_{Z})_{red}$ sa normalisation r\'eduite.
Par construction, le morphisme
$$F_{Z}^{\circ} \longrightarrow Z^{an}$$
poss\`ede au voisinage analytique d'un point $z \in Z^{an}$
au-dessus de $x$, une section. Or, comme $Z^{an}$ est normal en $x$,
l'analogue de \ref{l1.1} pour
des champs analytiques implique que
$F_{Z}^{\circ}$ est une gerbe au voisinage analytique de $z$. Or, comme
$Z$ est un sch\'ema propre, il existe un ouvert Zariski dense de $Z$,
$W$ contenant $x$,
tel que $F_{Z}^{\circ}$ soit une gerbe sur $W^{an}$. Mais par la proposition
\ref{p5.4}, cette gerbe est alg\'ebrisable. Il existe donc un sch\'ema
$\widetilde{W}$, un morphisme \'etale et surjectif $\widetilde{W} \longrightarrow
W$,
et un diagramme $1$-commutatif
$$\xymatrix{
\widetilde{W}^{an} \ar[d] \ar[r] & W^{an} \\
F_{Z}^{\circ} \ar[ru] & }$$
Le $1$-morphisme compos\'e
$$\widetilde{W}^{an} \longrightarrow F_{Z}^{\circ} \longrightarrow
F_{Z}\longrightarrow F$$
est par construction un prolongement de la section canonique
$$V \longrightarrow F_{U} \longrightarrow F$$
Comme celle-ci est \'etale au voisinage analytique de $x$, le $1$-morphisme
$$\widetilde{W}^{an} \longrightarrow F$$
est \'etale dans un voisinage Zariski d'un point $w \in \widetilde{W}$
au-dessus de $x$. Ainsi, quitte \`a restreindre $W$, on montre
qu'il existe un $1$-morphisme \'etale
$$X(x)^{an} \longrightarrow F$$
o\`u $X(x)$ est un sch\'ema, et dont l'image dans $M$ contient le point ferm\'e $x$.

Lorsque $x$ parcourt $M^{an}$ tout entier, on construit ainsi un sch\'ema
$X=\coprod_{x}X(x)$ et un morphisme \'etale et surjectif
$$X^{an} \longrightarrow F$$
Or, par quasi-compacit\'e de la topologie de Zariski sur $M$, on peut
supposer que $X$ est un sch\'ema de type fini sur $\bf C$. Enfin, comme
$F$ est normal, $X$ est aussi normal.

Il nous reste alors \`a d\'emontrer que le morphisme naturel
$$r : Y_{1}:=X^{an}\times_{F}X^{an} \longrightarrow X^{an}\times_{M} X^{an}$$
est alg\'ebrisable. En effet, dans ce cas, on pourra \'ecrire
$Y_{1}=X_{1}^{an}$, et le champ $F$ sera alors \'equivalent au champ associ\'e
\`a l'analytifi\'e du groupoide alg\'ebrique
$$s,b : X_{1} \longrightarrow X$$

Comme $F$ est un champ analytique propre, il est s\'epar\'e. Ainsi, $r$
est un morphisme fini. De plus, son image
est $X^{an}\times_{M}X^{an}$,
et comme $F$ est une gerbe sur $M-S$,
$r$ est \'etale en dehors de $S^{an}\times_{M}X^{an} \cup
X^{an}\times_{M}S$. Ainsi, $r$ est un morphisme fini et \'etale
en dehors d'un sous-espace alg\'ebrisable.
Comme $Y_{1}$ est normal, ceci implique que
$r$ est alg\'ebrisable (~\cite[$XII$]{sga1}~). \\

\underline{Etape $(2)$:} \\

Soit $F$ comme dans l'\'enonc\'e du th\'eor\`eme. On consid\`ere un
\'eclatement
$$p : X \longrightarrow M$$
tel que $X$ soit lisse, et que $p^{-1}(S)$ soit un diviseur \`a
croisements normaux
dans $X$. Un tel morphisme existe d'apr\`es la r\'esolution des
singularit\'es. Soit $F_{0}$ la
normalisation de $F\times_{M}X$. Alors, le $1$-morphisme canonique
$$f : F_{0} \longrightarrow F$$
est une modification. De plus,
lemme suivant implique que le champ $F_{0}$ v\'erifie
les hypoth\`eses de la premi\`ere \'etape. Il est donc alg\'ebrisable.

\begin{lem}\label{l5.5}
Soit $X$ un sch\'ema lisse est $S \hookrightarrow X$ un sous-sch\'ema
ferm\'e, tel que
$S$ soit un diviseur \`a croisements normaux dans $X$.
Alors tout germe analytique de rev\^etement de
groupe $H$ non-ramifi\'e en dehors de $S$ sur $X$ est alg\'ebrisable.
\end{lem}

\underline{\bf Preuve:} \rm Soit $x \in X^{an}$. Alors il existe un
voisinage analytique $U$ de $x$, tel que la paire $(U,U\cap S)$ soit
hom\'eomorphe \`a $(\bf C\mit^{m},H)$, o\`u $H$ est une r\'eunion finie
d'hyperplans complexes. Ainsi, on a
$$\pi_{1}(U-U\cap S)\simeq \bf Z\mit^{r}$$
o\`u $r$ est le nombre de composantes irr\'eductibles de $S\cap U$.

De plus, il existe un voisinage \'etale de $x$, $f : V \longrightarrow X$,
et une immersion ouverte $V \hookrightarrow \bf A\mit^{m}$, tel que
la paire $(V,f^{-1}(S))$ soit isomorphe \`a $(V,V\cap H)$, o\`u
$H$ est une r\'eunion finie d'hyperplans de $\bf A\mit^{m}$. \\

Soit $q : V \longrightarrow U$ un germe analytique de rev\^etement de groupe
$H$. Quitte \`a restreindre $U$ on peut supposer qu'il est comme
ci-dessus. Comme le rev\^etement est non-ramifi\'e en-dehors de
$S$, il donne lieu \`a un homomorphisme
$$\pi_{1}(U-U\cap S) \longrightarrow H$$
Comme on a un isomorphisme canonique
$$\pi_{1}(U-U\cap S) \simeq \pi_{1}(\bf C\mit^{m}-H) \simeq \bf
Z\mit^{r}$$
on a un homomorphisme induit
$$\pi_{1}((\bf A\mit^{m}-H)^{an}) \longrightarrow H$$
Ce $H$-torseur est alg\'ebrisable d'apr\`es \cite{sga1}. On obtient donc un \\
$H$-torseur alg\'ebrique
$$Y \longrightarrow \bf A\mit^{m}-H$$
Si $V$ est  un voisinage \'etale de $x$ comme ci-dessus, alors par
restriction on a un $H$-torseur
$$Y_{V} \longrightarrow V$$
Soit $Z$ la normalisation de $V$ dans $Y_{V}$. Alors
$$p : Z \longrightarrow V$$
d\'efinit un germe \'etale de rev\^etement de groupe $H$ en $x$, qui, par
construction, est une alg\'ebrisation du germe analytique $q : V
\longrightarrow U$. $\Box$\\
$\Box$\\

On tire de ce th\'eor\`eme le corollaire suivant sur "les champs de
Moisezon".

\begin{cor}\label{c5.5}
Si $F$ est un champ analytique r\'eduit tel que son corps des
fonctions m\'eromorphes est de dimension de transcendance \'egal \`a sa
dimension, alors il existe une modification
$$f : F_{0} \longrightarrow F$$
o\`u $F_{0}$ est alg\'ebrisable.
\end{cor}

\underline{\bf Preuve:} \rm En effet, les fonctions m\'eromorphes sur
$F$, sont les fonctions m\'eromorphes sur son espace de modules $M$.
L'espace analytique $M$ est donc alg\'ebrisable d'apr\`es \cite[$7.3$]{a}.
On termine alors par le th\'eor\`eme pr\'ec\'edent.
$\Box$\\

\begin{cor}\label{c5.6}
Tout champ analytique propre et lisse de dimension $1$ est alg\'ebrique.
\end{cor}

\underline{\bf Preuve:} \rm C'est en r\'ealit\'e un corollaire de la
preuve du th\'eor\`eme. En effet, comme $F$ est lisse, son espace de
modules est une surface de Riemann compacte, et v\'erifie donc les
hypoth\`eses de la premi\`ere partie de la preuve. $\Box$\\

A la vue de ces derniers r\'esultats, on peut
l\'egitimement poser la question suivante. \\

\underline{\bf Question:} \rm Tout champ analytique propre $F$,
dont l'espace de modules est un espace alg\'ebrique, est alg\'ebrisable. \\

Cette conjecture implique \`a titre d'exemple le fait (~certainement bien
connu~)
suivant.

\begin{cor}\label{c5.7}
Supposons que la conjecture pr\'ec\'edente soit vraie.

Soit $X$ un espace analytique lisse muni d'une action propre d'un groupe
de Lie affine $H$, tel que le quotient $X/H$ soit alg\'ebrisable
et propre.
Soit $Z$ le $\bf Gl\mit_{m}$-torseur tangent \`a $X$, muni de son
action  naturelle de $H$. Alors, le quotient $Y=Z/H$ est repr\'esentable
par un espace alg\'ebrique lisse, et l'action naturelle de $\bf Gl\mit_{m}$ sur
$Y$ est alg\'ebrisable.
\end{cor}

\underline{\bf Preuve:} \rm Commen\c{c}ons par supposer que $F=[X/H]$ est
un orbifold. La conjecture implique que $F$ est un champ alg\'ebrisable.
Ainsi, le lemme \ref{l3.9} implique que le torseur tangent \`a $F$ est
repr\'esentable par un espace alg\'ebrique lisse $Y$. Comme
$Y^{an}\simeq [Z/H]$, ceci montre que $Y^{an}$ est un quotient pour
$Z$ par $H$. De plus, comme $Y$ est le torseur tangent \`a $F$,
l'action de $\bf Gl\mit_{m}$ sur $Y^{an}$ est alors l'analytification
de celle de $\bf Gl\mit_{m}$ sur $Y$. \\

Dans le cas g\'en\'eral, on peut \'ecrire $q : F\longrightarrow F'$, o\`u
$F'$ un orbifold alg\'ebrisable. On termine alors en remarquant que
$$T_{F}=T_{F'}\times_{F'}F$$
$\Box$\\

\underline{Remarques:}
\begin{itemize}
\item
Lorsque $X/H$ est une vari\'et\'e projective, P. Essydieux m'a fait
remarquer qu'il est facile de d\'emontrer le corollaire \ref{c5.7} "\`a la main".
\item
La moralit\'e de ce corollaire est qu'\`a toute construction d'un champ alg\'ebrique
lisse et propre de Deligne-Mumford par des m\'ethodes analytiques, on associe
naturellement
une construction alg\'ebrique de ce champ. En effet, en gardant les
m\^emes notations que dans le corollaire, on a
$$[X/H]\simeq [Y/Gl_{m}]^{an}$$
\end{itemize}

Nous terminerons par une esquisse de d\'emonstration de cette conjecture.
Nous pr\'ef\'erons cependant la laisser sous forme de conjecture, car, comme nous
n'avons pas v\'erifi\'e tous les d\'etails, il se pourrait qu'une
difficult\'e technique nous ait \'echapper. \\

\underline{Etape $(1)$ :} On peut supposer que $F$ est r\'eduit. \\

Supposons que $F_{red}$ est alg\'ebrisable, et montrons alors que $F$
aussi.

Soit $\cal I$ l'id\'eal nilpotent de $F_{red}$ dans $F$. Comme $\cal
I\mit^{k}=0$, un raisonnement par r\'ecurrence sur $k$ montre qu'il
suffit de d\'emontrer le lemme suivant.

\begin{lem}
Soit $F$ un champ analytique propre, et $F_{0}$ un sous-champ d\'efini
par un id\'eal de carr\'e nul dans $\cal O\mit_{F}$. Si $F_{0}$ est
alg\'ebrisable, $F$ aussi.
\end{lem}

\underline{\bf "Preuve:"\mit} \rm En s'inspirant de \cite{bel}, on esp\`ere
pouvoir d\'efinir un complexe
cotangent analytique de $F_{0}$.
On le notera $L_{F_{0}^{an}}$. C'est un
objet de la cat\'egorie d\'eriv\'ee des complexes de pr\'efaisceaux de
\mbox{$\cal O\mit_{F_{0}}^{an}$-modules} \`a cohomologie coh\'erente, tel que
$H^{1}(F_{0}^{an},L_{F_{0}^{an}})$ est en bijection avec les classes
d'isomorphie des d\'eformations infinit\'esimales de $F_{0}^{an}$. En
particulier, le champ $F$ est d\'etermin\'e \`a \'equivalence pr\`es
par une classe dans ce groupe.

D'autre part les d\'eformations infinit\'esimales du champ alg\'ebrique
$F_{0}$ sont classifi\'ees par le groupe $H^{1}(F_{et},L_{F_{0}})$, o\`u
$L_{F_{0}}$ est le complexe cotangent relatif au morphisme structural
$F \longrightarrow Spec \bf C$ (~\cite[$9$]{lm}~). \\

Comme l'analytifi\'e de $L_{F_{0}}^{an}$ est quasi-isomorphe
\`a $L_{F_{0}^{an}}$, le th\'eor\`eme GAGA implique que toute d\'eformation
infinit\'esimale du champ $F_{0}^{an}$ est alg\'ebrisable. Ainsi,
$F$ est alg\'ebrisable. $\Box$\\

\underline{Etape $(2)$ :} Soit $f : F \longrightarrow F'$ un
morphisme propre et birationnel de champs analytiques. Alors, $F$ est
alg\'ebrisable si et seulement si $F$ est alg\'ebrisable. \\

Pour d\'emontrer cette partie il faudrait r\'e\'ecrire l'article d'Artin
\cite{a} dans le cadre plus g\'en\'eral des champs analytiques de
Deligne-Mumford. Cela implique en particulier qu'il faut d\'emontrer des
th\'eor\`emes GAGF dans ce cadre. Cependant, la preuve du th\'eor\`eme GAGA
donn\'ee pr\'ec\'edemment donne bon espoir que de tels \'enonc\'es soient
v\'erifi\'es.

La preuve du corollaire \cite[$7.15$]{a} se g\'en\'eralise alors au cas des champs
analytiques, au d\'etail pr\`es que nous n'avons d\'emontr\'e \ref{th5.1} que pour
des champs r\'eduits. Cependant l'\'etape $(1)$ implique que cela suffit.
$\Box$\\

\end{subsection}

\end{section}

\end{part}
\newpage
\begin{section}{Appendice}

\begin{subsection}{Spectres}\label{spectre}

Le formalisme des spectres permet de ne prendre en compte que
"l'information stable" de la th\'eorie de l'homotopie. Pour
l'utilisation que nous en ferons, leur int\'er\^et r\'eside dans le fait que
les limites et colimites homotopiques de spectres sont des foncteurs "exacts".
Ceci nous permet en particulier de faire commuter ces limites
homotopiques avec la formation des fibres et des cofibres. Cet
argument est utilis\'e implicitement tout au long de ce travail. \\

Un spectre $E$ est la donn\'ee d'une famille d'ensembles simpliciaux
fibrants point\'es $(E_{[n]},x_{n})$, avec $n \in \bf Z$, et d'une famille de
morphismes
$$e_{n} : (E_{[n]},x_{n}) \longrightarrow
(\Omega_{x_{n+1}}(E_{[n+1]}),x_{n+1}).$$

Un morphisme de spectres $f : E \longrightarrow E'$ est la donn\'ee
d'une famille de morphismes point\'es
$$f_{n} : (E_{[n]},x_{n}) \longrightarrow (E_{[n]}',x_{n}')$$
tel que $e_{n}' \circ f_{n} = \Omega(f_{n+1}) \circ e_{n}$.

On peut d\'efinir pour
tout $i\geq 0$ un morphisme de groupes
$$\pi_{i}(E_{[n]},x_{n})\longrightarrow \pi_{i}(\Omega_{x_{n+1}}(
E_{[n]}),x_{n+1}) \simeq\pi_{i+1}(E_{[n+1]},x_{n+1}).$$
On notera la limite de ce syst\`eme inductif par
$$\pi_{i}(E,x):=colim_{m\geq n}(\pi_{i+m-n}(E_{[m]},x_{m})).$$
Un morphisme de spectres $f : E \longrightarrow E'$, est appel\'e une
\'equivalence faible, si pour tout $i \in \bf Z$, le morphisme
induit
$$f_{*} : \pi_{i}(E,x) \longrightarrow \pi_{i}(E',f(x))$$
est un isomorphisme.

Un th\'eor\`eme important ( \cite[$2.53$]{j} ) est qu'il existe sur la
cat\'egorie des spectres
une structure de cat\'egorie de mod\`eles ferm\'ee compatible avec la
d\'efinition pr\'ec\'edente d'\'equivalence faible. \\

Si $X$ est un ensemble simplicial, et $E$ un spectre, on peut d\'efinir
le spectre des morphismes $Hom_{sp}(X,E)$ par
$$Hom_{sp}(X,E)_{[n]}:=Hom_{s}(X,E_{[n]}),$$
les morphismes de transition \'etant donn\'es par
$$\xymatrix{
Hom_{s}(X,E_{[n]}) \ar[r]^{e_{n}} & Hom_{s}(X,\Omega_{x_{n+1}}(E_{[n+1]}))
\ar[r]^{nat} & \Omega_{x_{n+1}}(Hom_{s}(X,E_{[n+1]})) }$$

De plus, cette cat\'egorie poss\`ede des $Hom$ internes, not\'es $Hom_{sp}$.
Tout comme la cat\'egorie des ensembles simpliciaux, elle poss\`ede aussi
des limites et colimites homotopiques. En particulier, si $f : E
\longrightarrow E'$ est un morphisme de spectres, on dispose
de la fibre et cofibre homotopique de $f$
$$Fib(f):= holim (
E \longrightarrow E' \longleftarrow \bullet )$$
$$Cof(f):= hocolim (
\bullet \longleftarrow  E \longrightarrow E' ).$$
La principale propri\'et\'e, qui distingue la cat\'egorie des spectres de
celle des ensembles simpliciaux, est que si l'on note
$$u : E' \longrightarrow Cof(f)$$
le morphisme naturel, alors il existe une \'equivalence faible
fonctorielle
$$E \simeq Fib(u).$$
De m\^eme, si on note
$$v : Fib(f) \longrightarrow E$$
le morphisme naturel, il existe une \'equivalence faible naturelle
$$Cof(v) \simeq E'.$$
De cette fa\c{c}on, on voit que la formation des fibres et cofibres
homotopiques commute avec les limites et colimites
homotopiques. C'est une propri\'et\'e tout \`a fait remarquable, et propre \`a
l'homotopie stable, que nous utiliserons tout au long de ce travail.\\

Pour finir, la construction de Dold-Puppe permet d'associer \`a tout
complexe de groupes ab\'eliens $C : \dots C_{n} \longrightarrow C_{n+1}
\dots$, un spectre point\'e $(K(C),0)$. Cette construction est fonctorielle, et
est telle que
$$\pi_{i}(K(C),0)\simeq H^{-i}(C).$$
Notons que "$n$-\`eme \'etage" $K(C)_{[n]}$, est par d\'efinition l'ensemble
simplicial obtenu par la construction de Dolp-Puppe appliqu\'ee au
complexe $\tau_{\leq 0}C[n]$.

\end{subsection}

\begin{subsection}{Descente}

Dans cet appendice on d\'emontre le th\'eor\`eme de descente.
La preuve est d\'ej\`a dans la preuve de \cite[Thm. $3-10$]{j2}, mais
le r\'esultat n'\'etant pas explicit\'e sous cette forme nous avons tenu
\`a en donner une d\'emonstration compl\`ete.\\

Pour tout l'appendice, $C$ est un site. On utilisera les notations de
la section $1$, ainsi que la proposition suivante :

\begin{prop}{\cite[$I-1$ Cor. $1$]{q2}}\label{a1}
Soit $F$ un pr\'efaisceau simplicial fibrant, et $f : A \rightarrow B$
une \'equivalence faible. Alors le morphisme induit
$$f^{*} : Hom_{s}(B,F) \rightarrow Hom_{s}(A,F)$$
est une \'equivalence faible.
\end{prop}

Si $X\in C$, alors on dispose d'un foncteur image r\'eciproque
$$j^{*} : SPr(C) \rightarrow SPr(C/X)$$
Ce foncteur poss\`ede un adjoint \`a gauche
$$j_{!} : SPr(C/X) \rightarrow SPr(C)$$
qui est l'extension par le pr\'efaisceau vide. Il est d\'efini
par :\\
\hspace*{5mm} pour $F\in SPr(C/X)$ et $U\in C$, alors
$$(j_{!}F(U)):=\coprod_{Hom_{C}(U,X)}F(U\rightarrow X)$$
Il est clair que $j_{!}$ pr\'eserve les cofibrations ainsi que
les \'equivalences faibles.

\begin{lem}\label{a2}
Soit $C$ et $C'$ deux sites et un foncteur
$$\bf a\mit  : SPr(C') \rightarrow SPr(C)$$
poss\`edant un adjoint \`a gauche
$$\bf b\mit : SPr(C) \rightarrow SPr(C')$$
qui pr\'eserve les cofibrations et les \'equivalences faibles, et tel que
$\bf b\mit(*)=*$. Alors le foncteur $\bf a\mit$ transforme objets fibrants
en objets
fibrants.
\end{lem}

\underline{\bf Preuve:} \rm Soit $F$ un pr\'efaisceau simplicial fibrant sur
$C'$, et $i : A \hookrightarrow B$ une cofibration triviale de $SPr(C)$. Il
faut montrer que le morphisme induit
$$i^{*} : Hom(B,a(F)) \rightarrow Hom(A,a(F))$$
est surjectif. Mais par adjonction, on dispose d'un carr\'e commutatif
$$\begin{array}{ccc}
Hom(B,a(F)) & \rightarrow & Hom(A,a(F)) \\
 \wr \downarrow & & \downarrow \wr\\
Hom(b(B),F) & \rightarrow & Hom(b(B),F)
\end{array}$$
Or par hypoth\`ese $b(i) : b(A) \rightarrow b(B)$ est une cofibration
triviale de $SPr(C')$, et donc le morphisme du bas est surjectif. Ce qui
implique que celui du haut aussi. $\Box$\\

On vient de voir que si $F$ est fibrant sur $C$, alors le
pr\'efaisceau $j^{*}F$ (~que l'on notera $F_{X}$  par la suite~) est fibrant
sur $C/X$. De plus, comme $j^{*}$ pr\'eserve les cofibrations triviales,
on obtient une \'equivalence faible canonique
$$H(C/X,F_{X}) \stackrel{\simeq}{\rightarrow} F^{\circ}(X)$$
pour chaque r\'esolution injective $F \hookrightarrow F^{\circ}$. De cette
fa\c{c}on on identifiera toujours les espaces $F^{\circ}(X)$ et
$H(C/X,F_{X})$, que l'on notera $H(X,F)$.\\

\begin{df}\label{a3}
Un objet simplicial $X_{\bullet}$ de $C$ est un foncteur
$$X_{\bullet} : \Delta^{op} \rightarrow C$$
o\`u $\Delta$ est la cat\'egorie simpliciale standard.
On notera $X_{m}$ pour l'objet $X_{\bullet}([m])$.\\
\hspace*{5mm} Si $X_{\bullet}$ est un objet simplicial de $C$, le site
induit sur
$X_{\bullet}$ est le site suivant :
\begin{itemize}
\item les objets sont les morphismes de $C$
$$U \rightarrow X_{m}$$
pour $m$ un entier positif.
\item un morphisme de $f : U \rightarrow X_{m}$ vers $g : V \rightarrow
X_{n}$ est la donn\'ee d'un morphisme $a : [n] \rightarrow [m]$ dans
$\Delta$ et d'un diagramme commutatif dans $C$
$$\begin{array}{rcl}
U & \rightarrow & V \\
f \downarrow & & \downarrow g \\
X_{m} & \stackrel{X_{\bullet}(a)}{\rightarrow} & X_{n}
\end{array}$$
\item un morphisme est couvrant si le morphisme induit
$$U \rightarrow V$$
est couvrant dans $C$.
\end{itemize}
Ce site est not\'e $C/X_{\bullet}$
\end{df}

Remarquons que l'on a un foncteur de restriction
$$j^{*} : SPr(C) \rightarrow SPr(C/X_{\bullet})$$
A travers ce foncteur, tout pr\'efaisceau simplicial $F$ sera aussi
consid\'er\'e comme pr\'efaisceau sur $X_{\bullet}$. Ainsi
$H(X_{\bullet},F)$ d\'esignera $H(C/X_{\bullet},j^{*}F)$.\\

Soit $X$ un objet de $C$ et $U \rightarrow X$ un morphisme couvrant. Le nerf
du recouvrement $U/X$ est l'objet simplicial de $C$ d\'efini par
$$\begin{array}{ccc}
\Delta & \rightarrow & C \\
 \
[m] & \mapsto &
U^{(m)}=\underbrace{U_{\stackrel{\times}{X}}U \dots
_{\stackrel{\times}{X}}U}_{m+1 \; fois}
\end{array}$$
et les morphismes $U^{(m)} \rightarrow U^{(n)}$ sont induits par les
projections et les diagonales. On le note $\cal N\mit(U/X)$.\\

Si $F$ est un pr\'efaisceau simplicial et $U \rightarrow X$ un recouvrement
de $C$, on obtient un $\Delta$-diagramme (~espace cosimplicial~)
d'ensembles simpliciaux
$$\begin{array}{ccc}
\Delta & \rightarrow & SEns \\
 \ [m] & \mapsto & F(U^{(m)})
\end{array}$$

Rappelons que l'espace de cohomologie de $\Check{C}ech$ du recouvrement
$U/X$ \`a coefficients dans le pr\'efaisceau simplicial $F$ est
$$\Check{H}(U/X,F):=Holim_{[m] \in \Delta}F(U^{(m)})$$

\begin{df}\label{a4}
Un espace cosimplicial $Z$ est un foncteur
$$\begin{array}{cccc}
Z : & \Delta & \rightarrow & SEns \\
              & [m]   & \mapsto & Z([m])
\end{array}$$
La cat\'egorie des espaces cosimpliciaux est not\'ee $CSEns$
\end{df}

\underline{\bf Exemples}\rm
\begin{itemize}
\item $*$ est l'espace cosimplicial constant
\item l'espace cosimplicial $\Delta/-$ est d\'efini par
$$[m] \mapsto (\Delta/-)([m])=B(\Delta/[m])$$
o\`u $B(I)$ est l'espace classifiant de la cat\'egorie $I$,
et $I/i$ est la cat\'egorie des morphismes de $I$ de but $i$
\end{itemize}

Un espace cosimplicial $Z$ peut \^etre vu comme un
pr\'efaisceau simplicial sur $\Delta$ (~site trivial~). Ainsi, si
$Y$ et $Z$ sont deux espaces cosimpliciaux, on d\'efinit l'espace des
morphismes de $Y$ vers $Z$ par
$$Hom_{cs}(Y,Z):=Hom_{s}(Y,Z)$$
o\`u $Y$ et $Z$ sont consid\'er\'es comme pr\'efaisceaux sur $\Delta$.\\

Avec ces notations, la limite homotopique de $Z$ est donn\'ee par
(~\cite[Ch. $XI$ $3-2$]{bk}~)
$$Holim_{\Delta}Z=Hom_{cs}(\Delta/-,Z)$$

\begin{thm}\label{a5}
Soit $F$ un pr\'efaisceau simplicial sur $C$, et $X_{\bullet}$ un objet
simplicial de $C$. Alors il existe une \'equivalence faible fonctorielle en
$F$
$$H(X_{\bullet},F) \simeq Holim_{[m]\in \Delta}H(X_{m},F)$$
\end{thm}

\begin{lem}\label{a6}
Si $F$ est un objet fibrant de $SPr(C)$, alors le pr\'efaisceau induit $F$
sur le site $C/X_{\bullet}$ est flasque.
\end{lem}

\underline{\bf Preuve:} \rm Soit $j^{*} : SPr(C) \rightarrow
SPr(C/X_{\bullet})$ le morphisme de restriction. On consid\`ere $j^{*}F
\rightarrow F^{\circ}$ une r\'esolution injective sur $C/X_{\bullet}$.
Alors, d'apr\`es le lemme \ref{a2}, le morphisme induit sur $C/X_{m}$
$$F \rightarrow F^{\circ}$$
est une cofibration triviale d'objets fibrants. C'est donc une \'equivalence
d'homotopie. Ainsi, pour tout objet $U$ de $C/X_{m}$, le morphisme induit
$$F(U) \rightarrow F^{\circ}(U)$$
est une \'equivalence faible. Comme ceci est vrai pour tout $U$ et tout $m$,
on en d\'eduit que pour tout objet $U$ de $C/X_{\bullet}$
$$F(U) \rightarrow F^{\circ}(U)$$
est une \'equivalence faible. $\Box$\\

\begin{lem}\label{a7}
Le foncteur
$$\begin{array}{ccc}
SPr(C/X_{\bullet}) & \rightarrow & CSEns \\
F & \mapsto & F(X_{\bullet})
\end{array}$$
admet un adjoint \`a gauche not\'e $Z \mapsto \widetilde{Z}$
\end{lem}

\underline{\bf Preuve:} \rm Soit $Z$ un espace cosimplicial. On d\'efinit
$$\begin{array}{cccc}
\widetilde{Z} : & C/X_{\bullet} & \rightarrow & SEns \\
                & (~U \rightarrow X_{m}~) & \mapsto & Z([m])
\end{array}$$
$\Box$ \\

\underline{\bf Preuve du th\'eor\`eme:} \rm Soit $F$ un pr\'efaisceau
simplicial sur $C$. En rempla\c{c}ant $F$ par $F^{\circ}$ on peut supposer
que $F$ est fibrant sur $C$.
Soit $F~\hookrightarrow~F^{\circ}$ une r\'esolution injective sur
$C/X_{\bullet}$. On sait que \mbox{$H(X_{\bullet},F)=Hom_{s}(*,F^{\circ})$}.
Notons $\bf1\mit=\widetilde{(\Delta/-)}$. \\

\begin{lem}\label{a8}
Le morphisme canonique $\bf 1\mit \rightarrow *$ est une \'equivalence faible
dans $SPr(C/X_{\bullet})$.
\end{lem}

\underline{\bf Preuve:} \rm Il suffit de voir que pour chaque $m$ l'espace
$\bf 1\mit(X_{m})$ est contractile. Or par d\'efinition, on a
$$\bf 1\mit(X_{m})=B(\Delta/[m])$$
Mais le classifiant d'une cat\'egorie qui poss\`ede un objet final est
contractile. $\Box$ \\

Comme $F^{\circ}$ est fibrant, la proposition \ref{a1} montre que le morphisme
naturel
$$Hom_{s}(*,F^{\circ}) \stackrel{\simeq}{\rightarrow} Hom_{s}(\bf
1\mit,F^{\circ})$$
est une \'equivalence faible.
Mais, par adjonction, il existe une \'equivalence faible
fonctorielle
$$Hom_{s}(\bf 1\mit,F^{\circ}) \stackrel{\simeq}{\rightarrow}
Holim_{\Delta}F^{\circ}(X_{\bullet})$$
On obtient ainsi une \'equivalence naturelle
$$H(X_{\bullet},F)
\stackrel{\simeq}{\rightarrow}Holim_{\Delta}F^{\circ}(X_{\bullet})$$
De plus, par le lemme \ref{a6}, le morphisme naturel
$$F \rightarrow F^{\circ}$$
est une \'equivalence faible \em objet par objet\rm . Comme
les limites homotopiques pr\'eservent les \'equivalences faibles, le
morphisme induit
$$Holim_{\Delta}F(X_{\bullet}) \rightarrow
Holim_{\Delta}F^{\circ}(X_{\bullet})$$
est une \'equivalence faible. Enfin, comme $F$ est fibrant, et par la
remarque suivant \ref{a2}, on obtient un diagramme fonctoriel en $F$
$$H(X_{\bullet},F) \stackrel{\simeq}{\rightarrow}
Holim_{\Delta}F^{\circ}(X_{\bullet}) \stackrel{\simeq}{\leftarrow}
Holim_{\Delta}H(X_{m},F)$$
$\Box$\\

\begin{prop}\label{a9}
Si $F$ est un pr\'efaisceau simplicial, alors le morphisme naturel
$$H(X,F) \rightarrow H(\cal N\mit(U/X),F)$$
est une \'equivalence faible
\end{prop}

Pour cela on a besoin d'un lemme que nous ne d\'emontrerons pas (~voir
\cite[Cor. $2-7$]{j2}~).

\begin{lem}\label{a10}
Si $F$ est un faisceau simplicial sur un site $C$, alors il existe une
r\'esolution injective $F \hookrightarrow F^{\circ}$, o\`u $F^{\circ}$ est
un faisceau d'ensembles simpliciaux.
\end{lem}

\underline{\bf Preuve de la proposition:} \rm Remarquons que, si l'on note
$SS(X)$ la cat\'egorie des faisceaux simpliciaux sur le site $C/X$, alors
on dispose d'une \'equivalence de cat\'egories
$$b : SS(X) \rightarrow SS(\cal N\mit(U/X))$$
De plus, si $F \in SS(X)$, qui est aussi fibrant
comme pr\'efaisceau simplicial, $b(F)$ est alors fibrant en tant qu'objet
de $SPr(\cal N\mit(U/X))$. En effet, notons $a$ le foncteur
de faisceautisation.\\
\hspace*{5mm} Soit $A \hookrightarrow B$ une cofibration triviale de $SPr(\cal
N\mit(U/X))$, et un diagramme commutatif sur $\cal N\mit(U/X)$
$$\begin{array}{ccc}
A & \rightarrow & b(F) \\
\downarrow & & \downarrow \\
B & \rightarrow & *
\end{array}$$
Comme $F$ est un faisceau on obtient un carr\'e commutatif
$$\begin{array}{ccc}
Hom_{\cal N\mit(U/X)}(B,b(F)) & \rightarrow & Hom_{\cal N\mit(U/X)}(A,b(F))
\\\| & & \| \\
Hom_{\cal N\mit(U/X)}(a(B),b(F)) & \rightarrow & Hom_{\cal
N\mit(U/X)}(a(A),b(F))
\end{array}$$
Puis, par l'\'equivalence de cat\'egories $SS(X) \rightarrow SS(\cal
N\mit(U/X))$, un diagramme commutatif
$$\begin{array}{ccc}
Hom_{\cal N\mit(U/X)}(B,b(F)) & \rightarrow & Hom_{\cal N\mit(U/X)}(A,b(F))
\\
\| & & \| \\
Hom_{X}(a(B),F) & \rightarrow & Hom_{X}(a(A),F) \\
\| & & \| \\
Hom_{X}(B,F) & \rightarrow & Hom_{X}(A,F)
\end{array}$$
et le morphisme du bas est surjectif par hypoth\`ese.\\

Soit $F$ un pr\'efaisceau simplicial sur $C$. Quitte \`a remplacer $F$ par
son faisceau associ\'e, on peut supposer que $F$ est un faisceau. En effet
le morphisme naturel $F \rightarrow a(F)$ est une \'equivalence faible, donc
$F$ et $a(F)$ ont la m\^eme cohomologie.\\

Soit $F \hookrightarrow F^{\circ}$ une r\'esolution injective dans
$SPr(C)$, avec
$F^{\circ}$ un faisceau. Alors
$$b(F) \rightarrow b(F^{\circ})$$
est encore une r\'esolution injective. On en d\'eduit donc
$$H(\cal N\mit(U/X),F)=Hom_{s}(*,b(F^{\circ}))=Hom_{s}(*,F^{\circ})
=H(X,F)$$
$\Box$\\

\begin{cor}\label{desc}
Si $U \rightarrow X$ est un recouvrement de $C$ et $F$ un pr\'efaisceau
simplicial, alors il existe une \'equivalence faible fonctorielle
en $F$, $X$ et $U$
$$H(X,F) \stackrel{\simeq}{\rightarrow} \Check{H}(U/X,H(-,F))$$
\end{cor}

\underline{\bf Preuve:} \rm On applique la proposition \ref{a9} et le
th\'eor\`eme \ref{a5}. $\Box$\\

Le th\'eor\`eme \ref{a5} sera appliqu\'e de la fa\c{c}on suivante.

\begin{cor}\label{a11}
Soit $X_{\bullet}$ un objet simplicial de $C$, et $F_{X}$ le pr\'efaisceau
simplicial qu'il repr\'esente. Alors, pour tout pr\'efaisceau simplicial
$F$, on a un isomorphisme naturel dans $HoSp$
$$\bf R\mit Hom_{s}(F_{X},F)\simeq \bf H\mit(X_{\bullet},j^{*}F)$$
o\`u $j : SPr(C/X_{\bullet}) \longrightarrow SPr(C)$ est le morphisme de
restriction.
\end{cor}

\underline{\bf Preuve:} \rm Quitte \`a remplacer $F$ par une r\'esolution
injective, on peut supposer que $F$ est fibrant.

Alors, on a par d\'efinition
$$\bf R\mit Hom_{s}(F_{X},F)\simeq Hom_{s}(F_{X},F)\simeq
Tot(F(X_{\bullet})$$
o\`u $Tot(F(X_{\bullet}))$ est l'espace total de l'espace cosimplicial
$[n] \mapsto F(X_{n})$ (~\cite[]{bk}~). Mais, il est d\'emontr\'e dans
\cite{j2}, que le morphisme canonique
$$Tot(F(X_{\bullet})) \longrightarrow holim_{\Delta}F(X_{\bullet})$$
est un isomorphisme. On conclut alors par le th\'eor\`eme \ref{a5}. $\Box$\\

% Le cas particulier qui nous interesse est celui o\`u $X_{\bullet}$ est
% l'objet simplicial associ\'e \`a un pr\'efaisceau simplicial. Rappelons que
% $\widetilde{X_\bullet}$ est le pr\'efaisceau simplicial repr\'esent\'e par
% $X_{\bullet}$.
%
% \begin{prop}\label{a}
% Il existe un foncteur
% $$\begin{array}{ccc}
% Spr(C) & \longrightarrow & SC \\
% F & \mapsto & EF_{\bullet}
% \end{array}$$
% et une transformation naturelle
% $$\widetilde{EF_{\bullet}} \longrightarrow F$$
% telle que, pour chaque pr\'efaisceau simplicial $F$, et chaque objet $U
% \in C$, le morphisme naturel
% $$\widetilde{EF_{\bullet}}(U) \longrightarrow F(U)$$
% est une \'equivalence faible.
% \end{prop}
%
% \underline{\bf Preuve:} \rm Soit $F \in Spr(C)$. Commen\c{c}ons par
% d\'efinir un objet bi-simplicial $bF$ de $C$, de la fa\c{c}on suivante.
% $$(bF)_{i,j}=\coprod_{\phi_{\bullet} : X_{0} \rightarrow X_{i}}
% X_{0}\times F(X_{i})_{j}$$
% o\`u $\phi_{\bullet}$ parcours les diagrammes de $i+1$-fl\`eches
% composables
% $$X_{0} \stackrel{\phi_{1}}{\longrightarrow} X_{1} \dots
% \stackrel{\phi_{i}}{\longrightarrow} X_{i}$$
% On d\'efinit alors $F_{\bullet}$ par l'objet simplicial diagonal de $bF$
% $$F_{\bullet}:=diag(bF)$$
% Le fait que le morphisme naturel
% $$\widetilde{F_{\bullet}} \longrightarrow F$$
% soit une \'equivalence faible objet par objet, est d\'emontr\'e dans
% \cite[]{s}. $\Box$\\
%
% \underline{Remarque:} Par analogie avec l'alg\`ebre homologique,
% l'objet $EF_{\bullet}$ doit \^etre pens\'e comme "la r\'esolution projective"
% canonique de $F$.

En stabilisant ces r\'esultats, on obtient des r\'esultats analogues sur
les spectres.

\begin{cor}\label{a12}
\begin{enumerate}
\item Soit $K$ un pr\'efaisceau en spectres fibrant dans $Sp(C)$. Alors
$F$ est flasque.
\item Soit $K$ un pr\'efaisceau en spectres, et $X_{\bullet}$ un objet
simplicial de $C$. Alors, il existe un isomorphisme naturel dans
$HoSp$
$$\bf H\mit(X_{\bullet},j^{*}K) \simeq holim_{[m]\in \Delta}\bf H\mit
(X_{m},K)$$
o\`u $j : SPr(C/X_{\bullet}) \longrightarrow SPr(C)$ est le morphisme de
restriction.
\end{enumerate}
\end{cor}

\end{subsection}

\begin{subsection}{Strictification}\label{strict}
\hspace{5mm}
Les objets simpliciaux que nous avons manipul\'e au cours de cette th\`ese
ne sont pas
de vrais objets simpliciaux, mais des "pseudo-objets simpliciaux".
Dans cette section nous justifions les d\'efinitions des spectres de
$K$-th\'eorie de tels objets. \\

Soit $I$ une petite cat\'egorie, et $F : I \longrightarrow Cat$ un
pseudo-foncteur contravariant (~\cite[$3.1$]{ma}~). On lui associe, \`a l'aide de
la construction de Grothendieck (~\cite[$3.4$]{ma}~), une strictification
$$SF : I \longrightarrow Cat$$
qui est un foncteur, muni de pseudo-transformations naturelles
$$SF \longrightarrow F$$
$$F \longrightarrow SF$$
adjoint l'une de l'autre. \\

Si pour chaque objet $i$, $F(i)$ est une cat\'egorie exacte (~ou
bi-complicial de Waldhausen~), et que pour chaque morphisme
$i \rightarrow j$, $F(i) \rightarrow F(j)$ est un foncteur exact, alors
les $SF(i)$ sont encore des cat\'egories exactes (~ou
bi-compliciales de Waldhausen~), et les $SF(i) \rightarrow SF(j)$
des foncteurs exacts. On dira alors que $F$ est un $I$-pseudo-diagramme de
cat\'egories exactes.

\begin{df}
Soit $F : I \longrightarrow Cat$ un $I$-pseudo-diagramme de cat\'egories
exactes (~ou bi-compliciale de Waldhausen~). Notons
$K~:~CatEx~\longrightarrow~Sp$ le foncteur de $K$-th\'eorie.
On d\'efinit alors
$$holim_{I}K(F(i)):=holim_{I}K(SF(i))$$
$$hocolim_{I}K(F(i)):=hocolim_{I}K(SF(i))$$
en tant qu'objet de $HoSp$.
\end{df}

Il est \`a noter que si $F$ est d\'ej\`a un foncteur, la transformation
naturelle $SF \longrightarrow F$ induit une \'equivalence sur les
spectres de $K$-th\'eorie (~car elle poss\`ede un adjoint~)
objet par objet. Ainsi, la d\'efinition pr\'ec\'edente coïncide avec la
d\'efinition classique (~\cite[$X$ $3$]{bk}~) \`a isomorphisme canonique pr\`es dans
$HoSp$. \\

Supposons maintenant que l'on dispose d'une pseudo-transformation
naturelle de pseudo-foncteurs sur $I$ (~\cite[$3.2$]{ma} )
$$f : F \longrightarrow F'$$
Elle induit alors, une transformation naturelle de foncteurs
$$Sf : SF \longrightarrow SF'$$
et donc des morphismes dans $HoSp$
$$Sf : holim_{I}K(F(i)) \longrightarrow holim_{I}K(F'(i))$$
$$Sf : hocolim_{I}K(F(i)) \longrightarrow hocolim_{I}K(F'(i))$$
Par le proc\'ed\'e pr\'ec\'edent, tous les pseudo-diagrammes de cat\'egories
exactes que nous rencontrerons seront strictifi\'es pour pouvoir en prendre
le spectre de $K$-th\'eorie.

\end{subsection}

\begin{subsection}{Extension des coefficients}\label{extcoeff}
\hspace{5mm}
Lorsque $K$ est un pr\'efaisceau en spectres sur un site $C$, et
$\cal A$ un pr\'efaisceau de $\bf Q$-espaces vectoriels, nous souhaitons donner un
sens \`a l'expression $K\otimes \cal A$. % Ce doit \^etre un objet de
% $HoSp(C)$, qui est un $\cal A$-module, tel que si $K$ est associ\'e par
% la construction de Dold-Puppe \`a un complexe de pr\'efaisceaux ab\'eliens
% $C_{\bullet}$ sur $C$, $K\otimes \cal A$ est \'equivalent au spectre associ\'e
% au complexe $C_{\bullet}\otimes \cal A$. \\

Il est alors naturel de poser la d\'efinition suivante.

\begin{df}
Soit $C$ un site, muni d'un pr\'efaisceau de $\bf Q$-epaces vecoriels $\cal A$.
Notons $\cal H\mit_{\cal A}=K(\cal A\mit,0)$ le pr\'efaisceau en spectres d'Eilenberg-Maclane
associ\'e. C'est un pr\'efaisceau en spectres tel que ses pr\'efaisceaux
d'homotopie v\'erifient
$$\begin{array}{ll}
\pi_{k}(\cal H\mit_{\cal A})=& 0 \; pour \; k\neq 0 \\
\pi_{0}(\cal H\mit_{\cal A})\simeq H_{0}(\cal H\mit_{\cal A},\bf
Z\mit)& \simeq \cal A 
\end{array}$$

Pour tout objet $K \in Sp(C)$, nous noterons
$$K\otimes\cal A\mit:=K\wedge \cal H\mit_{\cal A}$$
en tant qu'objet de $HoSp(C)$.
\end{df}

Remarquons que comme $\cal A$ est un pr\'efaisceau en $\bf Q$-espaces vectoriels, $K(\cal A\mit,0)$ est aussi 
un pr\'efaisceau en spectres de Moore pour $\cal A$. 

\begin{prop}\label{a13}
\begin{enumerate}
\item
La correspondance $K \mapsto K\otimes \cal A$ est un
foncteur de $HoSp(C)$ dans $HoSp(C)$.
\item
Soit $K$ un objet en anneaux de $HoSp(C)$.
Si $\cal A$ est un pr\'efaisceau de $\bf Q$-alg\`ebres, alors $K\otimes \cal A$ est un
objet en anneaux dans $HoSp(C)$.
De plus, il existe un morphisme naturel dans $HoSp(C)$
$$K \longrightarrow K\otimes \cal A$$
\item
Si $\cal A$ est le pr\'efaisceau en $\bf Q$-espaces vectoriels libres
engendr\'e par un pr\'efaisceau d'ensembles $X$, alors il
existe un isomorphisme canonique
$$K\otimes \cal A\mit \simeq \bigvee_{X}K_{\bf Q}$$
\item
Si $\cal A$ est un pr\'efaisceau constant en $\bf Q$-alg\`ebres et acyclique, 
alors pour chaque objet $X~\in~Ob(C)$, il
existe un isomorphisme naturel
$$H^{*}(X,K\otimes \cal A\mit)\simeq H^{*}(X,K)\otimes H^{0}(\cal A\mit(X))$$
\end{enumerate}
\end{prop}

\underline{\bf Preuve:} \rm Les propri\'et\'es $(1)$ et $(2)$ sont
imm\'ediates par d\'efinition, et $(3)$ et $(4)$ sont une application de la
formule de Kunneth.
$\Box$\\

\underline{Remarque:} \rm Si $\cal A$ est un pr\'efaisceau constant de
fibre $A$, et $a \in A$, le morphisme
$$\begin{array}{cccc}
a : & A & \longrightarrow & A \\
 & b& \mapsto & a.b
\end{array}$$
induit un morphisme de pr\'efaisceaux de $\bf Q$-espaces vectoriels
$$a : \cal A\mit \longrightarrow \cal A$$
Ainsi, par fonctorialit\'e, on a un morphisme de spectres ab\'eliens
$$a : K\otimes \cal A\mit \longrightarrow K\otimes \cal A$$
Ceci justifie en particulier la notation
$$\frac{1}{m} : K\otimes \bf Q\mit \longrightarrow K\otimes \bf Q$$
ou encore si $K$ est un spectre en anneau, et
$\phi : K \longrightarrow K$ un morphisme dans $HoSp(C)$
$$\frac{1}{m}.\phi : K\otimes \bf Q\mit \longrightarrow K\otimes \bf Q$$
qui est un morphisme dans $HoSp(C)$.

\begin{cor}
Soit $F$ un champ de Deligne-Mumford, et $\Lambda$ le faisceau
de $\bf Q$-alg\`ebres d\'efini dans \ref{katrep}. Alors, si $S$ contient les racines de
l'unit\'e, il existe un isomorphisme naturel
$$H^{-*}(I_{F}^{t},\underline{K}\otimes \Lambda) \simeq
\underline{\bf K}_{*}(I_{F}^{t})\otimes \Lambda(S)$$
\end{cor}

\underline{\bf Preuve:} \rm D'apr\`es \ref{a13} $(5)$, il suffit de montrer que
$\Lambda$ est acyclique sur $(Esp/S)_{et}$.

Comme la cohomologie commute avec les limites inductives filtrantes,
il suffit de montrer que $\Lambda_{m}$ est acyclique pour
la topologie \'etale. Comme c'est un pr\'efaisceau constant de $\bf
Q$-espaces vectoriels, il est acyclique pour la topologie \'etale.

\end{subsection}

\begin{subsection}{Th\'eorie de Hodge pour les champs alg\'ebriques
complexes}\label{hodge}
\hspace{5mm}
Pour d\'emontrer la formule de la signature \ref{sign}, nous avons
utilis\'e que pour un champ alg\'ebrique complexe lisse $F$, dont
l'espace de modules est projectif, on a
$$\sigma(F)=\sum_{p,q}(-1)^{p}h^{p,q}$$
o\`u $h^{p,q}:=Dim_{\bf C}H^{p}(F,\Omega^{q}_{F})$. Pour le cas o\`u $F$
est un sch\'ema, cette formule se d\'emontre \`a l'aide de la th\'eorie de
Hodge (~\cite[$15.8.2$]{hi}~). Nous allons montrer dans cette partie,
que cette th\'eorie est applicable aux champs alg\'ebriques, et donc, que
l'on peut d\'emontrer la formule pr\'ec\'edente par un raisonnement tout \`a
fait analogue \`a celui fait dans \cite[$15.8.2$]{hi}. Nous utiliserons
pour cela un argument de Deligne (~\cite{de}~).\\

Fixons nous un champ alg\'ebrique complexe et lisse $F$, tel que son
espace de modules $M$ soit projectif. Nous noterons $p : F
\longrightarrow M$ la projection. Soit $h_{M} \in H^{2}(M^{an},\bf
Z\mit)$ une section hyperplane de $M$. Comme $p$ est un morphisme
fini, la classe $h:=p^{*}(h_{M})\in H^{2}(F^{an},\bf Z\mit)$ est
positive, et d\'efinit donc une m\'etrique K$\ddot{a}$hl\'erienne sur $F^{an}$. On
associe \`a cette m\'etrique les notions usuelles de formes harmoniques,
et des op\'erateurs $L$, $\Lambda$, $h$, $d$, $d'$, et $d''$. Les identit\'es
de K$\ddot{a}$hler impliquent que les formes harmoniques sur $F^{an}$, sont
les formes $\alpha$ v\'erifiant $d'\alpha=d''\alpha=0$. \\

Soit $q : X \longrightarrow F$ un morphisme propre et g\'en\'eriquement
fini, avec $X$ un sch\'ema lisse et projectif. En utilisant les
arguments de transferts utilis\'es dans \cite[$4.3$]{de}, on montre que les
morphismes
$$q^{*} : H^{p}(F,\Omega^{q}) \longrightarrow H^{p}(X,\Omega^{q})$$
sont injectifs. On suit alors l'argument utilis\'e dans \cite[$5.3$]{de},
pour aboutir \`a la proposition suivante.

\begin{prop}{ (~Deligne \cite[$5.3$, $5.4$]{de}~)}
\begin{itemize}
\item
La suite spectrale de Hodge vers de Rham d\'eg\'en\`ere.
\item
La filtration induite sur $H^{*}(F^{an},\bf C\mit)$ v\'erifie
$$F^{p}(H^{n}(F^{an},\bf C\mit))=\bigoplus_{i\geq p}
F^{p}(H^{p+q}(F^{an},\bf C\mit))\cap \overline{F^{q}(H^{p+q}(F^{an},\bf
C\mit))}$$
En particulier, on a un isomorphisme naturel
$$H^{n}(F^{an},\bf C\mit)\simeq \bigoplus_{p+q=n}H^{p}(F,\Omega^{q})$$
\item
Toute classe de cohomologie dans $H^{*}(F^{an},\bf C\mit)$ peut se
repr\'esenter par un forme $\alpha$ avec $d'\alpha=d''\alpha=0$. En
particulier l'inclusion des formes harmoniques dans les formes
exactes induit un isomorphisme
$$\cal H\mit^{n}(F^{an}) \simeq H^{n}(F^{an},\bf C\mit)$$
\end{itemize}
\end{prop}

En corollaire de cette proposition, on voit que tout le formalisme de
la th\'eorie de Hodge est valable pour le champ $F$. En particulier, la
preuve de \cite[$15.8.2$]{hi} se g\'en\'eralise \`a notre situation, et
montre la formule cherch\'ee
$$\sigma(F)=\sum_{p,q}(-1)^{p}h^{p,q}$$

Par exemple, dans le cas des champs alg\'ebrique de dimension $2$, on trouve la
formule d'indice de Hodge suivante.

\begin{cor}
Soit $F$ un champ alg\'ebrique complexe lisse, tel que son espace de
modules soit une surface projective. Alors, la forme d'intersection
$$H^{1}(F,\Omega^{1}_{F})\cap H^{2}(F^{an},\bf R\mit)\times
H^{1}(F,\Omega^{1}_{F})\cap H^{2}(F^{an},\bf R\mit)
\longrightarrow \bf R$$
poss\`ede un unique vecteur propre positif.
\end{cor}

\underline{Remarque:} On doit aussi pouvoir d\'emontrer plus
g\'en\'eralement un th\'eor\`eme de Hodge pour des champs diff\'erentiels
riemanniens compacts.

\end{subsection}

\end{section}

\end{document}